\pdfoutput=1

\documentclass[oneside,reqno,11pt]{amsart}

\usepackage{preamble}

\myexternaldocument[I:]{corank1_ASW_I}
\myexternaldocument[II:]{corank1_ASW_II}
\myexternaldocument[IV:]{corank1_ASW_IV}

\title{Co-rank \texorpdfstring{$1$}{1} Arithmetic Siegel--Weil III: \\ Geometric local-to-global}
\author{Ryan C. Chen}
\date{February 8, 2025}
\address{Department of Mathematics, Massachusetts Institute of Technology, 182 Memorial Drive, Cambridge, MA 02139, USA}
\email{rcchen.math@gmail.com}

\begin{document}
    
    \begin{abstract}
This is the third in a sequence of four papers, where we prove the arithmetic Siegel--Weil formula in co-rank $1$ for Kudla--Rapoport special cycles on exotic smooth integral models of unitary Shimura varieties of arbitrarily large even arithmetic dimension. 
Our arithmetic Siegel--Weil formula implies that degrees of Kudla--Rapoport arithmetic special $1$-cycles are encoded in the first derivatives of unitary Eisenstein series Fourier coefficients.

In this paper, we make the reduction from global arithmetic intersection numbers for special cycles to the local geometric quantities in our companion papers. 
Building on our previous companion papers, we also propose a construction for arithmetic special cycle classes associated to possibly singular matrices of arbitrary co-rank. 
\end{abstract}

    \maketitle
    
    \tableofcontents
    
    \clearpage


        \section{Introduction}
        \label{sec:part_III:intro}
            This paper is a continuation of our companion papers \cite{corank1_ASW_I.pdf,corank1_ASW_II.pdf}, and will be followed by \cite{corank1_ASW_IV.pdf}.
We refer the reader to the introductions of \cite{corank1_ASW_I.pdf,corank1_ASW_IV.pdf} for further motivation, overview, and strategy for our four-part series of papers. The discussion in \cref{ssec:part_III:intro:arith_Siegel--Weil} is an abridged version of loc. cit..

As explained in \cref{ssec:part_III:intro:global_decomp,ssec:part_III:intro:local-to-global_Archimedean,ssec:part_III:intro:local-to-global_non-Arch},
the main objective of this paper is to explain how to reduce global intersection numbers to quantities appearing on the geometric sides of our local main theorems from \cite{corank1_ASW_I.pdf,corank1_ASW_II.pdf}.
We also discuss a construction of arithmetic cycle classes; see \cref{ssec:part_III:intro:arith_cycle_classes}.
In \cref{ssec:part_III:intro:outline}, we outline the structure of this paper and its relation with our companion papers \cite{corank1_ASW_I.pdf, corank1_ASW_II.pdf, corank1_ASW_IV.pdf}.

            \subsection{Arithmetic Siegel--Weil}
            \label{ssec:part_III:intro:arith_Siegel--Weil}
            For the introduction, fix an imaginary quadratic field $F / \Q$ with ring of integers $\mc{O}_F$ and odd discriminant $\Delta$. Let $V$ be a non-degenerate $F / \Q$ Hermitian space of signature $(n - 1, 1)$ with pairing $(-,-)$. Set $G = U(V)$.
For the introduction, we assume $V$ contains a full-rank self-dual\footnote{Our convention is that \emph{self-duality} for Hermitian lattices is always with respect to the trace pairing, here $\mrm{tr}_{F / \Q}(-,-)$.} $\mc{O}_F$-lattice,
and assume that the prime $p = 2$ is split in $\mc{O}_F$.

Since the work of Kudla--Rapoport \cite{KR14} (also Rapoport--Smithling--Zhang \cite{RSZ21}), it has been customary to define special cycles $\mc{Z}(T) \ra \mc{M}$ over (stacky) integral models $\mc{M} \ra \Spec \mc{O}_F$ for Shimura varieties associated to $G' \coloneqq \Res_{F / \Q} \G_m \times G$. 
In this paper, we are mainly interested in the case when $\mc{M} \ra \Spec \mc{O}_F$ is the ``exotic smooth'' Rapoport--Smithling--Zhang (RSZ) integral model of odd relative dimension $n - 1$ \cite[{\S 6}]{RSZ21} (empty if $n \equiv 0 \pmod{4}$). The stack $\mc{M}$ admits a moduli description: it parameterizes tuples $(A_0, \iota_0, \lambda_0, A, \iota, \lambda)$ where $A_0$ and $A$ are abelian schemes (dimensions $1$ and $n$ respectively) with $\mc{O}_F$-actions $\iota_0$ and $\iota$, and with compatible quasi-polarizations $\lambda_0$ and $\lambda$. The datum $(A_0, \iota_0, \lambda_0, A, \iota, \lambda)$ satisfies a few additional conditions, which we suppress in the introduction (see \cref{ssec:ab_var:integral_models}). We are able to prove versions of our main global results for more general $\mc{M}$ (including odd arithmetic dimension $n$) at the price of discarding finitely many primes (particularly ramified primes for odd $n$); see \crefext{IV:remark:arithmetic_Siegel-Weil:main_results:other_levels}.

The moduli stack $\mc{M}$ carries a natural family of Hermitian $\mc{O}_F$-lattices 
    \begin{equation}
    \mc{L}at \ra \mc{M} \quad \quad \mc{L}at \coloneqq \underline{\Hom}_{\mc{O}_F}(A_0, A),
    \end{equation}
e.g. as explained in \cref{ssec:ab_var:special_cycles,appendix:quasi-compactness_lemma:terminology}.
Given\footnote{The notation $\mrm{Herm}_m$ denotes a scheme over $\Spec \Z$, e.g. $\mrm{Herm}_m(\Q)$ denotes $m \times m$ Hermitian matrices with entries in $F$.} any $T \in \mrm{Herm}_m(\Q)$, the associated \emph{Kudla--Rapoport} special cycle $\mc{Z}(T) \ra \mc{M}$ is defined as the substack
    \begin{equation}\label{equation:part_III:intro:arith_Siegel-Weil:KR_cycle}
    \mc{Z}(T) \coloneqq \{\underline{x} \in \mc{L}at^m : (\underline{x}, \underline{x}) = T\} \subseteq \mc{L}at^m
    \end{equation}
consisting of $m$-tuples with Gram matrix $T$. More precisely, see \cref{ssec:ab_var:integral_models}. The morphism $\mc{Z}(T) \ra \Spec \mc{O}_F$ is smooth of relative dimension $n - 1 - \rank(T)$ in the generic fiber over $\Spec F$. If $T$ is not positive semi-definite, then $\mc{Z}(T)$ is empty.

We propose a new candidate definition of arithmetic cycle classes
    \begin{equation}\label{equation:intro:arith_Siegel--Weil:cycle_class_def}
    [\widehat{\mc{Z}}(T)] \coloneqq [\widehat{\mc{Z}}(T)_{\ms{H}}] + \sum_{p \text{ prime}} [{}^{\mathbb{L}}\mc{Z}(T)_{\ms{V},p}] \in \arithCh^m(\mc{M})_{\Q}.
    \end{equation}
Here, $[\widehat{\mc{Z}}(T)_{\ms{H}}] \in \arithCh^m(\mc{M})_{\Q}$ is intended to describe ``horizontal'' contributions,
involving both the flat part\footnote{Given an algebraic stack $\mc{X}$ over a Dedekind domain $R$, its \emph{flat part} or \emph{horizontal part} $\mc{X}_{\ms{H}}$ is the largest closed substack $\mc{X}_{\ms{H}} \subseteq \mc{X}$ which is flat over $\Spec R$. The stack $\mc{X}_{\ms{H}}$ is also the scheme-theoretic image of the generic fiber of $\mc{X}$. Given a formal algebraic stack $\mc{X}$ over $\Spf R$ for a complete discrete valuation ring $R$, its \emph{flat part} or \emph{horizontal part} $\mc{X}_{\ms{H}}$ is the largest closed substack $\mc{X}_{\ms{H}} \subseteq \mc{X}$ which is flat over $\Spf R$ (in the sense discussed in \cref{footnote:flatness_formal_alg_stacks}).} $\mc{Z}(T)_{\ms{H}}$ of $\mc{Z}(T)$
and a Green current $g_{T,y}$ on the Archimedean fiber of $\mc{M}$ (where $g_{T,y}$ depends on an auxiliary parameter $y \in \opn{Herm}_m(\R)_{>0}$, as is typical in Kudla's program).
The class ${}^{\mathbb{L}}\mc{Z}(T)_{\ms{V},p} \in \mrm{gr}^m_{\mc{M}} K'_0(\mc{Z}(T)_{\F_p})_{\Q}$ is intended to describe ``vertical'' contributions from positive characteristic.
This construction is discussed further in \cref{ssec:part_III:intro:arith_cycle_classes}.
The actual construction appears in \cref{sec:arith_cycle_classes} of the paper.

Due to non-properness of $\mc{M} \ra \Spec \mc{O}_F$ in general, one should likely modify $[\widehat{\mc{Z}}(T)]$ on a suitable compactification of $\mc{M}$. 
If $\mc{Z}(T) \ra \Spec \mc{O}_F$ is proper, however, we consider certain ``arithmetic degrees without boundary contributions'' (a real number) when $m \leq n$, given by
    \begin{align}\label{equation:part_III:intro:results:if_proper:Archimedean}
        \widehat{\deg}([\widehat{\mc{Z}}(T)] \cdot \widehat{c}_1(\widehat{\mc{E}}^{\vee})^{n-m}) & \coloneqq \left ( \int_{\mc{M}_{\C}} g_{T,y} \wedge c_1(\widehat{\mc{E}}^{\vee}_{\C})^{n-m} \right ) 
        \\
        & \hphantom{\coloneqq} + \widehat{\deg}((\widehat{\mc{E}}^{\vee})^{n - \rank(T)}|_{\mc{Z}(T)_{\ms{H}}}) \label{equation:part_III:intro:results:if_proper:horizontal}
        \\
        & \hphantom{\coloneqq} + \sum_{p \text{ prime}} \deg_{\F_p}({}^{\mathbb{L}}\mc{Z}(T)_{\ms{V},p} \cdot (\mc{E}^{\vee})^{n-m}) \log p \label{equation:part_III:intro:results:if_proper:vertical}
    \end{align}
conditional on convergence of the integral, for a certain metrized tautological bundle $\widehat{\mc{E}}$ on $\mc{M}$ (\cref{ssec:ab_var:integral_models}); we do check convergence of the integral in the settings of our arithmetic Siegel--Weil results. Here we set $\mc{M}_{\C} \coloneqq \mc{M} \times_{\Spec \mc{O}_F} \Spec \C$ for either embedding $F \ra \C$. The middle term is mixed characteristic in nature: for $\rank (T) = n - 1$, it is (essentially) a weighted sum of (stable) Faltings heights of abelian varieties \crefext{IV:remark:arithmetic_Siegel-Weil:main_results:Faltings_height}.
For proper $\mc{Z}(T) \ra \mc{O}_F$,
the total quantity from the formulas
\cref{equation:part_III:intro:results:if_proper:Archimedean,equation:part_III:intro:results:if_proper:horizontal,equation:part_III:intro:results:if_proper:vertical}
should coincide with the arithmetic degree (without boundary contributions) of a version of $[\widehat{\mc{Z}}(T)]$ on any reasonable compactification of $\mc{M}$.

We consider the (normalized) $U(m,m)$ Siegel Eisenstein series
    \begin{equation}\label{equation:part_III:intro:results:Eisenstein}
    E^*(z,s)^{\circ}_n \coloneqq \Lambda_m(s)^{\circ}_n \sum_{\begin{psmallmatrix} a & b \\ c & d \end{psmallmatrix} \in P_1(\Z) \backslash SU(m,m)(\Z)} \frac{\det(y)^{s - s_0}}{\det(c z + d)^n |\det (c z + d)|^{2(s - s_0)}}.
    \end{equation}
Here $\Lambda_m(s)^{\circ}_n$ is a certain normalizing factor \crefext{IV:equation:global_normalized_Fourier:global_normalization:normalizing_factor}, consisting of various $L$-functions and $\Gamma$ functions, etc..
We wrote $P_1 \coloneqq P \cap SU(m,m)$ for $P \subseteq U(m,m)$ denoting the Siegel parabolic ($m \times m$ block upper triangular), the element $z = x + i y$ lies in Hermitian upper-half space (i.e. $x \in \mrm{Herm}_m(\R)$ and $y \in \mrm{Herm}_m(\R)_{>0}$, meaning $y$ is positive definite), and $s_0 = (n - m)/2$. Given $T \in \mrm{Herm}_m(\Q)$, the Eisenstein series $E^*(z,s)^{\circ}_n$ has \emph{$T$-th Fourier coefficient}
    \begin{equation}\label{equation:part_III:intro:results:Fourier_coeff}
    E^*_T(y,s)^{\circ}_n \coloneqq 2^{m(m-1)/2} |\Delta|^{-m(m-1)/4} \int_{\mrm{Herm}_m(\Z) \backslash \mrm{Herm}_m(\R)} E^*(z, s)^{\circ}_n e^{-2 \pi i \mrm{tr}(Tz)} ~dx
    \end{equation}
for $z = x + i y$ in Hermitian upper-half space, where this integral is taken with respect to the Euclidean measure on $\mrm{Herm}_m(\R)$.

In our four-part sequence of papers, our main result is an \emph{arithmetic Siegel--Weil formula} in co-rank $1$, i.e. we show that the formula
    \begin{equation}\label{equation:part_III:intro:results:conjecture}
    \frac{h_F}{w_F} \frac{d}{ds} \bigg|_{s = s_0} \frac{2 \Lambda_n(s-s_0)^{\circ}_n}{\kappa \Lambda_m(s)^{\circ}_n} E^*_T(y,s)^{\circ}_n \overset{?}{=} \widehat{\mrm{vol}}_{\widehat{\mc{E}}^{\vee}}([\widehat{\mc{Z}}(T)])
    \end{equation}
holds when $m = n$ and $T$ is singular of co-rank $1$. Here, we write $h_F$ (resp. $w_F$) for the class number of (resp. number of roots of unity in) $\mc{O}_F$, and $\kappa = 1$ (resp. $\kappa = 2$) if $m < n$ (resp. $m = n$). We also prove (the closely related version of) \cref{equation:part_III:intro:results:conjecture} when $T$ is nonsingular of rank $m = n - 1$; in that case, the special value at the non-central point $s = s_0$ simultaneously has geometric meaning (``geometric Siegel--Weil'').
These proofs are completed in our companion paper \cite{corank1_ASW_IV.pdf},
and we refer to the introduction of \cite{corank1_ASW_IV.pdf} for further discussion.

            \subsection{Global decomposition}
            \label{ssec:part_III:intro:global_decomp}
                Our proof of arithmetic Siegel--Weil in co-rank $1$ is local in nature. In our previous companion papers \cite{corank1_ASW_I.pdf,corank1_ASW_II.pdf}, we formulated and proved the key ``local arithmetic Siegel--Weil'' formulas at all places, which are local analogues of (the co-rank $1$ case of) \crefext{equation:part_III:intro:results:conjecture}. At non-Archimedean places (resp. the Archimedean place), the geometric side of our local theorem takes place on appropriate Rapoport--Zink spaces (resp. Hermitian symmetric domains).

The main objective of this paper is to explain how to patch these local geometric quantities to the (global) arithmetic intersection number considered in the formulas \cref{equation:part_III:intro:results:if_proper:Archimedean,equation:part_III:intro:results:if_proper:horizontal,equation:part_III:intro:results:if_proper:vertical},
using uniformization (non-Archimedean and Archimedean).
The main difference between this ``local-to-global'' argument in our work compared with previous works
(such as \cite{KR14,LZ22unitary}) is that we have a new mixed characteristic contribution
from \eqref{equation:part_III:intro:results:if_proper:horizontal}
(e.g. Faltings heights of abelian varieties, coming from components of $\mc{Z}(T)$ flat over $\Spec \mc{O}_F$)
which do not admit an obvious \emph{canonical} local decomposition,
especially if one seeks a local interpretation in terms of e.g. Rapoport--Zink spaces.
Dealing with this issue and giving such a local interpretation is the main new contribution of this paper.

To state our ``local-to-global'' results, we set up some additional notation.
Consider $T \in \opn{Herm}_m(\Q)$.
In \cref{ssec:part_III:intro:global_decomp,ssec:part_III:intro:local-to-global_Archimedean,ssec:part_III:intro:local-to-global_non-Arch}, we assume $\rank(T) \geq n - 1$ and either $m = n - 1$ or $m = n$,
or assume that $T$ is nonsingular and not positive definite.
Using the expressions in
\cref{equation:part_III:intro:results:if_proper:Archimedean,equation:part_III:intro:results:if_proper:horizontal,equation:part_III:intro:results:if_proper:vertical},
we set
    \begin{align}
        \label{equation:intro:local-to-global:Archimedean}
        \mrm{Int}_{\infty, \mrm{global}}(T, y) & \coloneqq \int_{\mc{M}^{\mrm{an}}_{\C}} g_{T,y} \wedge c_1(\widehat{\mc{E}}_{\C}^{\vee})^{n - m}
        \\
        \label{equation:intro:local-to-global:horizontal}
        \sum_{p} \mrm{Int}_{\ms{H},p,\mrm{global}}(T) & \coloneqq \widehat{\deg}((\widehat{\mc{E}}^{\vee})^{n - \rank(T)}|_{\mc{Z}(T)_{\ms{H}}}) - (\deg_{\Z} \mc{Z}(T)_{\ms{H}}) \cdot h_{\widehat{\mc{E}}^{\vee}}^{\mrm{CM}}
        \\
        \mrm{Int}_{\ms{V},p,\mrm{global}}(T) & \coloneqq \deg_{\F_p} ({}^{\mbb{L}} \mc{Z}(T)_{\ms{V},p} \cdot (\mc{E}^{\vee})^{n - m}) \log p
    \end{align}
for primes $p$ and the height constant
    \begin{equation}
    h^{\mrm{CM}}_{\widehat{\mc{E}}^{\vee}} \coloneqq \frac{L'(1,\eta)}{L(1, \eta)} + \frac{\Gamma'(1)}{\Gamma(1)} + \frac{1}{2} \log |\Delta| - \log(2 \pi) - \frac{1}{2} \log(4 \pi e^{\gamma})
    \end{equation}
from \eqref{equation:arith_cycle_classes:Hodge_bundles:taut_height_constants},
where $\eta$ is the quadratic character associated to $F / \Q$ and $\gamma$ is the Euler--Mascheroni constant.
This height constant $h^{\mrm{CM}}_{\widehat{\mc{E}}^{\vee}}$
essentially comes from Faltings height for elliptic curves with CM by $\mc{O}_F$.
We have
    \begin{equation}
        \mrm{Int}_{\infty, \mrm{global}}(T, y) \in \R 
        \quad \quad
        \mrm{Int}_{\ms{H},p,\mrm{global}}(T) \in \Q \cdot \log p
        \quad \quad
        \mrm{Int}_{\ms{V},p,\mrm{global}}(T) \in \Q \cdot \log p
    \end{equation}
for all primes $p$. For fixed $T$ and $y$, all but finitely many of these quantities equal zero.

In this paper, we explain how to calculate these quantities locally via uniformization,
on Hermitian symmetric domains and Rapoport--Zink spaces.
The trickiest piece for us is the horizontal contribution $\mrm{Int}_{\ms{H},p,\mrm{global}}(T)$ (``mixed characteristic'').
The subtraction in \eqref{equation:intro:local-to-global:horizontal} occurs for roughly the following reason:
while Faltings height does not admit a \emph{canonical} local decomposition,\footnote{At least without picking a section
of the appropriate metrized line bundle. Perhaps the analytic side should not ``know'' which section was picked!}
the \emph{difference} of Faltings heights
for isogenous abelian varieties \emph{does} admit a canonical local decomposition over primes $p$,
as a $\Q$-linear combination of $\log p$ for primes $p$.
Indeed, given isogenous abelian varieties $A_1, A_2$ over a number field $E$, with everywhere semi-abelian reduction,
the difference of Faltings heights is given by
\begin{align}\label{equation:intro:local-to-global:change_Faltings_height}
    h_{\mrm{Fal}}(A_2) - h_{\mrm{Fal}}(A_1) & = \sum_{p \mid \deg \phi} a_p \log p = \frac{1}{[E : \Q]} \sum_p \sum_{\breve{w} \mid p} [\breve{E}_{\breve{w}} : \breve{\Q}_p] \breve{\delta}_{\mrm{Fal}}(\phi_{\breve{w}})
    \\
    \breve{\delta}_{\mrm{Fal}}(\phi_{\breve{w}}) & \coloneqq \frac{1}{2} \log(\deg \phi_{\breve{w}}) - \frac{1}{[\breve{E}_{\breve{w}} : \breve{\Q}_p]} \mrm{length}_{\mc{O}_{\breve{E}_{\breve{w}}}}(e^* \Omega^1_{\ker \phi/ \mc{O}_{\breve{E}_{\breve{w}}}}) \log p
    \label{equation:intro:local-to-global:change_Faltings_height_local}
\end{align}
for some $a_p \in \Q$ and any choice of isogeny $\phi \colon A_1 \ra A_2$
(extend $\phi$ over $\mc{O}_E$). Here, the symbol $e$ means the identity section,
the inner sum in \cref{equation:intro:local-to-global:change_Faltings_height}
runs over all prime ideals $\breve{w}$ of $E \otimes_{\Q} \breve{\Q}_p$
(with associated residue field $\breve{E}_{\breve{w}}$),
and $\phi_{\breve{w}}$ means the associated isogeny of $p$-divisible groups
$A_1[p^{\infty}] \ra A_2[p^{\infty}]$ base-changed to $\Spf \mc{O}_{\breve{E}_{\breve{w}}}$.

Since Faltings height is intrinsic, the coefficients $a_p$ do not depend on the choice of $\phi$,
by linear independence of $\log p$ for different $p$.
This is also the reason why a difference of Faltings heights appears in e.g. \crefext{I:corollary:intro:results:n=2}. 
We then argue that these numbers $a_p \in \Q$ can be calculated in a purely local way for our situation,
in terms of local special cycles on Rapoport--Zink spaces.
The argument we give is somewhat delicate,
as we wish to avoid writing down explicit (global) isogenies $\phi \colon A_1 \ra A_2$
in order to obtain a more ``purely local'' formulation.
The quantities $\breve{\delta}_{\mrm{Fal}}(\phi_{\breve{w}})$ are indeed amenable to
(local) calculation on Rapoport--Zink spaces.
Unlike the global coefficients $a_p \in \Q$, however,
the $\breve{\delta}_{\mrm{Fal}}(\phi_{\breve{w}})$
may in general depend on the choice of $\phi_{\breve{w}}$, i.e. the local isogeny of $p$-divisible groups.
This issue is a key difficulty for us, which needs to be resolved in this paper.
We show that some well-chosen local isogenies $\phi_{\breve{w}}$ 
lift to global isogenies of $p$-divisible groups (\cref{ssec:qcan_heights:minimal_isogenies})
and can be used to calculate contributions to (differences of) Faltings height.

There is also the issue that the tautological bundle $\widehat{\mc{E}}^{\vee}$
is not the same as the metrized Hodge bundle (but is known to behave similarly,
as first observed by Gross \cite{Gross78} and studied further in \cite{BHKRY20II}),
so the (more natural) version with ``tautological height'' needs additional argument.
The ``tautological height'' and Faltings height are treated in parallel in \cref{part:local_change_heights}.

While previous work for special points on Shimura curves \cite{KRY04} also studied the change in Faltings heights along isogenies,
our insistence on a purely local formulation is an important difference for our method.
We only observe the limiting phenomena in \crefext{I:figure:intro:strategy:local_limit_method} on a local level;
this is what allows us to prove a theorem on Shimura varieties of arbitrarily large dimension.

            \subsection{Local-to-global: Archimedean}
            \label{ssec:part_III:intro:local-to-global_Archimedean}

We first consider Archimedean intersection numbers, which are more standard;
this is still ``purely vertical'' from the perspective of this paper.

Given any $m$-tuple $\underline{x} \in V^m$,
with image $\underline{x}_{\infty} \in (V \otimes \R)^m$,
recall that there is an Archimedean
\emph{local special cycle} $\mc{D}(\underline{x}_{\infty}) \subseteq \mbb{D}$
\crefext{II:ssec:Hermitian_domain:local_cycles}
on a Hermitian symmetric domain $\mc{D}$ arising in complex uniformization for $\mc{M}$
(\cref{sec:Arch_uniformization}). Here $\mc{D}$ parameterizes maximal negative definite subspaces
of $V$, and $\mc{D}(\underline{x}_{\infty})$ consists of those subspaces contained in the orthogonal complement of
$\opn{span}(\underline{x}_{\infty})$.
If $\underline{x}_f$ denotes the image of $\underline{x}$ in $(V \otimes_{\Q} \A_f)^m$,
there is an associated discrete set
    \begin{equation}
        \mc{D}(\underline{x}_f) \coloneqq
        \{ g \in U(V)(\A_f) / K_f :
        g^{-1} x_i \in L \otimes_{\Z} \hat{\Z}^p \text{ for all $x_i \in \underline{x}_f$}\}
        \subseteq U(V)(\A_f) / K_f
    \end{equation}
(see \cref{ssec:Arch_uniformization:away_infty_special_cycle})
where $K_f$ is the stabilizer of some full rank self-dual lattice $L \subseteq V$.
We call $\mc{D}(\underline{x}_f)$ an ``away-from-$\infty$'' special cycle,
because it appears as an ``away-from-$\infty$'' component of complex uniformization
of global special cycles, from an ad\`elic point of view.

Consider any $y \in \mrm{Herm}_m(\R)_{>0}$ and $\underline{x} \in V^m$ with Gram matrix $(\underline{x}, \underline{x}) = T$.
For $T$ as in \cref{ssec:part_III:intro:global_decomp},
we gave currents $[\xi(\underline{x},y)]$ on the Hermitian symmetric domain $\mc{D}$
\crefext{II:ssec:part_II:Hermitian_domain:corank1_modification}, using Kudla Green functions and star products
(this needs some additional care when $T$ is singular; see loc. cit.).
These currents $[\xi(\underline{x},y)]$ will descend to currents $g_{T,y}$ on the (analytification of the) complex Shimura variety $\mc{M}_{K'_f,\C}^{\mrm{an}}$.
These $g_{T,y}$ are (by definition) the currents appearing in our global intersection number \cref{equation:intro:local-to-global:Archimedean}.
We have a local Archimedean intersection number
    \begin{equation}\label{equation:part_II:Arch_uniformization:Archimedean}
    \mrm{Int}_{\infty}(T,y) \coloneqq \int_{\mc{D}} [\xi(\underline{x}, y)] \wedge c_1(\widehat{\mc{E}}^{\vee})^{n - m}.
    \end{equation}
Via complex uniformization, this is related to the global Archimedean intersection number by the following
``local-to-global'' formula
(\cref{ssec:Arch_uniformization:Archimedean}).
\begin{propositionLetter}
    \label{proposition:intro:local-to-global_Archimedean}
With notation as above, we have
    \begin{equation}
    \mrm{Int}_{\infty, \mrm{global}}(T, y) = \mrm{Int}_{\infty}(T,y) \frac{h_F}{w_F} \cdot \deg \Biggl [ U(V)(\Q) \backslash \coprod_{\substack{\underline{x} \in V^m \\ (\underline{x}, \underline{x}) = T}} \mc{D}(\underline{x}_f) \Biggr ]
    \end{equation}
where $\deg [ - ]$ denotes groupoid cardinality.
\end{propositionLetter}

Note that the actual body of this paper proves something slightly stronger,
e.g. not requiring that $V$ contain a self-dual lattice (\cref{ssec:Arch_uniformization:Archimedean}).

Our main Archimedean local theorem from \cite{corank1_ASW_II.pdf} related $\opn{Int}_{\infty}(T,y)$
with the derivative of an Archimedean local Whittaker function.
The term $\deg [ \cdots ]$ essentially arises from a ``uniformization degree''
and essentially unfolds as a product of local theta integrals, which may be related with special values of
(normalized) local Whittaker functions (away from $\infty$). This is explained (with precise constants)
in our companion paper \crefext{IV:ssec:local_Siegel-Weil:uniformization_degree}.
These local Whittaker functions appear in ``product formula local decompositions'' for Eisenstein series Fourier coefficients,
and $\opn{Int}_{\infty, \mrm{global}}(T,y)$ then corresponds to a derivative at the Archimedean place,
for an application of the Leibniz derivative rule to a product of local Whittaker functions.
This is explained in \cite{corank1_ASW_IV.pdf}.

            \subsection{Local-to-global: Non-Archimedean}
            \label{ssec:part_III:intro:local-to-global_non-Arch}
                We next consider ``non-Archimedean'' intersection numbers, which will involve both ``vertical'' positive characteristic contributions (e.g. from components $\mc{Z}(T)$ in positive characteristic) and ``horizontal'' parts of mixed characteristic contributions (from components of $\mc{Z}(T)$ flat over $\mc{O}_F$, giving Faltings heights of abelian varieties).
Set $F_p \coloneqq F \otimes_{\Q} \Q_p$.
If $p$ is nonsplit in $F$, we let $\breve{F}_p$ be the completion of the maximal unramified extension of $F_p$.
If $p$ is split in $F$, we let $\breve{F}_p \coloneqq \breve{\Q}_p$ and pick a morphism $F_p \ra \breve{F}_p$.

There is a space of \emph{global special quasi-homomorphisms} $\mbf{W}$ depending on $p$
(\cref{ssec:non-Arch_uniformization:framing})
which is an $F / \Q$ Hermitian space of dimension $n$ (resp. $n - 1$) if $F_p / \Q_p$ is nonsplit (resp. split).
Given any $m$-tuple $\underline{\mbf{x}} \in \mbf{W}^m$ with image $\underline{\mbf{x}}_p \in \mbf{W}_p^m$
(where $\mbf{W}_p \coloneqq \mbf{W} \otimes_{\Q} \Q_p$),
recall that there is a \emph{local special cycle}
$\mc{Z}(\underline{\mbf{x}}_p) \subseteq \mc{N}$ on a Rapoport--Zink space $\mc{N}$
\crefext{I:ssec:moduli_pDiv:special_cycles}.
Here $\mc{N}$ is a moduli space of $p$-divisible groups appearing in (non-Archimedean)
Rapoport--Zink uniformization of $\mc{M}$ and its special cycles.\footnote{Strictly speaking, when $F_p / \Q_p$ is split,
we only uniformize the special cycles we need; these special cycles need not be contained in supersingular loci.}
There is an inclusion of Hermitian spaces $\mbf{W} \otimes_{\Q} \A_f^p \hookrightarrow V \otimes_{\Q} \A_f^p$,
which is an isomorphism if $F / \Q_p$ is nonsplit.
If $\underline{\mbf{x}}^p$ denotes the image of $\underline{\mbf{x}}$ in $(V \otimes_{\Q} \A_f^p)^m$,
there is an associated discrete set
    \begin{equation}
        \mc{Z}(\underline{\mbf{x}}^p) \coloneqq
        \{ g \in U(V)(\A_f^p) / K_f^p :
        g^{-1} x_i \in L \otimes_{\Z} \hat{\Z}^p \text{ for all $x_i \in \underline{\mbf{x}}^p$}\}
        \subseteq U(V)(\A^p) / K^p
    \end{equation}
(see \cref{ssec:Arch_uniformization:away_infty_special_cycle})
where $K_f^p$ is the stabilizer of some full rank self-dual lattice $L \subseteq V$.
As above, we call $\mc{Z}(\underline{\mbf{x}}^p)$ an ``away-from-$p$'' special cycle,
because it appears as an ``away-from-$p$'' component of Rapoport--Zink uniformization
of global special cycles, from an ad\`elic point of view.
We also ``complete'' $\mbf{W}$ to a certain $n$-dimensional Hermitian space
$\mbf{V} = \mbf{W} \oplus \mbf{W}^{\perp}$ (where $\mbf{W}^{\perp} = 0$ if $F_p / \Q_p$ is nonsplit);
see \cref{ssec:non-Arch_uniformization:framing}.
Set $I_1 \coloneqq U(\mbf{W}) \times U(\mbf{W}^{\perp})$.

Consider any $\underline{\mbf{x}} \in \mbf{W}^{m}$ with Gram matrix $(\underline{\mbf{x}}, \underline{\mbf{x}}) = T$.
As studied in \cite{corank1_ASW_I.pdf},
we then have local ``vertical'' and local ``horizontal'' intersection numbers at $p$ given by
    \begin{align}
    \opn{Int}_{\ms{V},p}(T) & \coloneqq
    2 [\breve{F}_p : \breve{\Q}_p]^{-1}
    \deg_{\overline{\F}_p}( {}^{\mbb{L}} \mc{Z}(\underline{\mbf{x}}_p)_{\ms{V}} \cdot (\mc{E}^{\vee})^{n - m} ) \log p
    \\
    \opn{Int}_{\ms{H},p}(T) & \coloneqq
    2 \left ( \sum_{\mc{Z} \hookrightarrow \mc{Z}(\underline{\mbf{x}}_p)_{\ms{H}}} \deg(\mc{Z}) \delta_{\mrm{tau}}(s_{\mc{Z}}) \right ) \log p
    \end{align}
respectively. The quantities $\opn{Int}_{\ms{V},p}(T)$ and $\opn{Int}_{\ms{H},p}(T)$
appeared previously on the geometric side of our non-Archimedean local main theorem \crefext{I:theorem:intro:non_Arch_local_theorems},
but we quickly recall some of the objects involved.
Here ${}^{\mbb{L}} \mc{Z}(\underline{\mbf{x}}_p)_{\ms{V}}$ is a certain ``derived vertical local special cycle class''
supported on the special fiber of the local special cycle $\mc{Z}(\underline{\mbf{x}}_p) \subseteq \mc{N}$,
and $\mc{E}^{\vee}$ is some (local) dual tautological bundle on $\mc{N}$.
The sum runs over ``components'' of the flat part
$\mc{Z}(\underline{\mbf{x}}_p)_{\ms{H}} \subseteq \mc{Z}(\underline{\mbf{x}}_p)$
of the local special cycle.
The symbol $\deg(\mc{Z})$ denotes the degree of the adic finite flat morphism $\mc{Z} \ra \Spf \mc{O}_{\breve{F}_p}$,
and $\delta_{\mrm{tau}}(s_{\mc{Z}})$ is a certain ``local change of tautological height''.
More precisely, the universal $p$-divisible group over each $\mc{Z}$ involves a quasi-canonical lifting $\mf{X}_s$ of level
$s \coloneqq s_{\mc{Z}} \in \Z_{\geq 0}$ in the sense of Gross \cite{Gross86},
and $\delta_{\mrm{tau}}(s_{\mc{Z}}) \log p$ is defined to be $- \breve{\delta}_{\mrm{Fal}}(\psi_s)$
with $\breve{\delta}_{\mrm{Fal}}(-)$ as in \eqref{equation:intro:local-to-global:change_Faltings_height_local}
and $\psi_s \colon \mf{X}_0 \ra \mf{X}_s$ any isogeny of minimal degree from a canonical lifting $p$-divisible group
$\mf{X}_0$ to a level $s$ quasi-canonical lifting $\mf{X}_s$.

The main result of this paper is the following non-Archimedean ``local-to-global'' formula
(\cref{ssec:non-Arch_uniformization:vertical,ssec:non-Arch_uniformization:horizontal}).
The notation below $K_{1,\mbf{L}^{\perp}_p} \subseteq U(\mbf{W}_p^{\perp})$ means the unique maximal compact subgroup,
where $\mbf{W}_p^{\perp} \coloneqq \mbf{W}^{\perp} \otimes_{\Q} \Q_p$ (this term vanishes when $F_p / \Q_p$ is nonsplit).

\begin{theoremLetter} We have
    \begin{equation}
    \mrm{Int}_{\ms{H},p,\mrm{global}}(T) = \mrm{Int}_{\ms{H},p}(T) \frac{h_F}{w_F} \cdot \deg \Biggl [ I_1(\Q) \backslash \Biggl ( \coprod_{\substack{\underline{x} \in \mbf{W}^m \\ (\underline{x}, \underline{x}) = T}} U(\mbf{W}^{\perp}_p) / K_{1,\mbf{L}^{\perp}_{p}} \times \mc{Z}(\underline{\mbf{x}}^p) \Biggr ) \Biggr ],
    \end{equation}
and the same holds when $\ms{H}$ (``horizontal'') is replaced by $\ms{V}$ (``vertical'').
\end{theoremLetter}

We again emphasize that the main difficulty for us is the horizontal case,
due to the issues explained in \cref{ssec:part_III:intro:global_decomp}.
Additional arguments are needed here, besides just Rapoport--Zink uniformization; this is the subject of
\cref{part:local_change_heights}.

We also note that the actual body of this paper (\cref{ssec:non-Arch_uniformization:vertical,ssec:non-Arch_uniformization:horizontal})
proves a slightly stronger statement,
e.g. we only need to require that $V$ contains a lattice which is self-dual at $p$ (as opposed to self-dual everywhere,
as assumed for simplicity in this introduction).

Our main non-Archimedean local theorem from \cite{corank1_ASW_I.pdf}
related $\opn{Int}_{\ms{H},p}(T) + \opn{Int}_{\ms{V},p}(T)$
with the derivative of a local Whittaker function at $p$.
As in the Archimedean case, the term $\deg [ \cdots ]$ essentially arises
as a ``uniformization degree'' and may be expressed in terms of local theta integrals
and special values of local Whittaker functions away from $p$.
The total intersection number $\mrm{Int}_{\ms{H},p,\mrm{global}}(T)+\mrm{Int}_{\ms{V},p,\mrm{global}}(T)$ will
correspond to an application of the Leibniz derivative rule to product of local Whittaker functions on the analytic side,
for a derivative at $p$. Again, this is explained (with precise constants) in our companion paper \cite{corank1_ASW_IV.pdf}.


            \subsection{Arithmetic cycle classes}
            \label{ssec:part_III:intro:arith_cycle_classes}
                As mentioned above, construction of arithmetic cycle classes
    \begin{equation}
    [\widehat{\mc{Z}}(T)] \coloneqq [\widehat{\mc{Z}}(T)_{\ms{H}}] + \sum_{p \text{ prime}} [{}^{\mathbb{L}}\mc{Z}(T)_{\ms{V},p}] \in \arithCh^m(\mc{M})_{\Q}
    \end{equation}
from \cref{equation:intro:arith_Siegel--Weil:cycle_class_def} needs additional work.
We quickly sketch our contributions in this direction; the relevant material appears is \cref{sec:arith_cycle_classes} below.

For general $T$, there is no precise definition of $[\widehat{\mc{Z}}(T)]$ which has been proposed in the prior literature \cite[{Remark 4.4.2}]{Li22IHES}.
Our candidate definition may need modification on a compactification, but we expect it to apply in already-compact situations
(e.g. the Rapoport--Smithling--Zhang \cite{RSZ21} setup for CM extensions of totally real fields $\neq \Q$).
In general, it may also be necessary to modify the Green currents differently than in \cite[{Definition 4.7}]{GS19}; see discussion below.

We suggest treating the ``Archimedean'' and ``horizontal'' part of $\mc{Z}(T)$ by Garcia--Sankaran's idea (see below),
and suggest treating the ``vertical'' part of $\mc{Z}(T)$ by some derived/$K$-theoretic methods.
It seems unknown whether these two approaches can be combined, e.g. there seems to be no
``derived arithmetic intersection theory'' at the moment.

As we explain in \cref{ssec:arith_cycle_classes:horizontal} of this paper, the class $[\widehat{\mc{Z}}(T)_{\ms{H}}]$ may be constructed using currents $g_{T,y}$ (associated to $T$ and allowed to vary with a parameter $y \in \mrm{Herm}_m(\R)_{>0}$) satisfying a \emph{modified current equation}, i.e. that
    \begin{equation}\label{equation:part_III:global_theorem_statement:arithmetic_Siegel-Weil:modified_current_equation}
    -\frac{1}{2 \pi i} \partial \overline{\partial} g_{T,y} + \delta_{\mc{Z}(T)_{\C}} \wedge [c_1(\widehat{\mc{E}}^{\vee}_{\C})^{m - \rank(T)}]
    \end{equation}
is represented by a smooth $(m,m)$-form. For such currents, we copy the idea of \cite[{\S 5.4}]{GS19} and apply it to the flat part $\mc{Z}(T)_{\ms{H}}$ of $\mc{Z}(T)$.

Currents satisfying \cref{equation:part_III:global_theorem_statement:arithmetic_Siegel-Weil:modified_current_equation} were constructed by Garcia and Sankaran \cite{GS19}, using the Mathai--Quillen theory of superconnections \cite[{(4.38)}]{GS19}. For their arithmetic Siegel--Weil results, however, they need a non ``linearly invariant'' modification of their current \cite[{Definition 4.7}]{GS19} (see discussion below).

We choose to instead use the star-product approach of Kudla \cite{Kudla97a} (as formulated by Liu for unitary groups \cite{Liu11}) to define the currents $g_{T,y}$ for our arithmetic Siegel--Weil results. Traditionally, the star product approach was used for nonsingular $T$ (or at least block diagonal $T$, with diagonal entries $0$ or nonsingular). In \crefext{II:ssec:part_II:Hermitian_domain:corank1_modification}, we gave a (new) \emph{linearly invariant} modification in the case of singular $T \in \mrm{Herm}_n(\Q)$ with rank $n - 1$, which will appear in our arithmetic Siegel--Weil result for singular $T$. The discussion in \crefext{II:ssec:part_II:Hermitian_domain:corank1_modification} focused on the local version (on the Hermitian symmetric domain); we descend it to the complex Shimura variety in \cref{ssec:Arch_uniformization:Archimedean} of this paper.

The ``vertical'' part of $[\widehat{\mc{Z}}(T)]$ at a prime $p$ will be given by a
class ${}^{\mbb{L}} \mc{Z}(T)_{\ms{V},p} \in \opn{gr}^m_{\mc{M}} K'_0( \mc{Z}(T)_{\F_p})_{\Q}$
defined by derived tensor products
(\cref{ssec:arith_cycle_classes:vertical}).
One might hope to reduce to the case when $T = \opn{diag}(0,T^{\flat})$ for a nonsingular $T^{\flat}$,
but there is a slight complication due to possible non-diagonalizability of $T$
(because of a class number phenomenon). We get around this with a certain ``$p$-local linear invariance''
which is enough for our purpose.

As part of the expected automorphic behavior of $[\widehat{\mc{Z}}(T)]$, it is expected that these classes should satisfy a certain ``linear invariance'' property for the action\footnote{For any $\gamma \in \GL_m(\mc{O}_F)$ and any Hermitian matrix $T \in \mrm{Herm}_m(\Q)$, we say e.g. that $T$ and ${}^t \overline{\gamma} T \gamma$ are \emph{$\GL_m(\mc{O}_F)$-equivalent}, and that they lie in the same \emph{$\GL_m(\mc{O}_F)$-equivalence class}.} of $\GL_m(\mc{O}_F)$ on Hermitian matrices $T$. We verify this for the classes we define: for any $g_{T,y}$ satisfying 
    \begin{equation}\label{intro:part_III:results:Archimedean_linear_invariance}
    g_{T,y} = g_{{}^t \overline{\gamma} T \gamma, \gamma^{-1} y {}^t \overline{\gamma}^{-1}},
    \end{equation}
we show
    \begin{equation}
    [\widehat{\mc{Z}}(T)] = [\widehat{\mc{Z}}({}^t \overline{\gamma} T \gamma)]
    \end{equation}
where $[\widehat{\mc{Z}}(T)]$ is formed with respect to $y$ and $[\widehat{\mc{Z}}({}^t \overline{\gamma} T \gamma)]$ is formed with respect to $\gamma^{-1} y {}^t \overline{\gamma}^{-1}$.
In fact, we prove refined results: we show that the vertical part at each prime $p$ is linearly invariant on the level of Grothendieck groups \cref{equation:arith_cycle_classes:vertical:linear_invariance}, and we will see that the horizontal part is linear invariant on its own (\cref{ssec:arith_cycle_classes:horizontal}).
The currents $g_{T,y}$ appearing in our main arithmetic Siegel--Weil results satisfy the linear invariance property in \eqref{intro:part_III:results:Archimedean_linear_invariance}; see \cref{ssec:Arch_uniformization:Archimedean} (also \crefext{II:ssec:part_II:Hermitian_domain:corank1_modification}). Note that the Garcia--Sankaran Green currents in \cite[{(4.38)}]{GS19} also satisfy the same linear invariance property (but the modified currents in \cite[{Definition 4.7}]{GS19} do not).

            \subsection{Outline}
            \label{ssec:part_III:intro:outline}
                We briefly summarize the remaining content in this paper, and discuss the relation with our companion papers \cite{corank1_ASW_I.pdf,corank1_ASW_II.pdf,corank1_ASW_IV.pdf}.
Further explanations may be found at the beginning of some sections.

In \cref{part:part_III:global_cycles} ``Global special cycles'', we set up the global moduli stacks (RSZ) and special cycles (KR) appearing in our main global theorems.
In the present paper, we also allow a more general setup at the price of possibly discarding finitely many primes (particularly the ramified primes in the case of odd arithmetic dimension). Our main theorems will also apply in these situations \crefext{IV:remark:arithmetic_Siegel-Weil:main_results:other_levels}.
In Section \ref{sec:arith_cycle_classes}, we define associated arithmetic special cycle classes and discuss arithmetic degrees.

In \cref{part:local_change_heights} ``Local change of heights'', we give a key argument in our reduction process from global heights in mixed characteristic to local quantities. We formulate these local quantities to be computable in terms of local special cycles and certain ``local change of heights'' along ``minimal isogenies'' (sketched in \cref{ssec:part_III:intro:global_decomp}). The ``local change of height'' constants appearing in our non-Archimedean ``local arithmetic Siegel--Weil'' main theorems in \cite{corank1_ASW_I.pdf} emerge from this global-to-local reduction process.

In \cref{part:part_III:uniformization} ``Uniformization'' we discuss complex and Rapoport--Zink uniformization of special cycles in our setup, and finish the reduction process from global heights/intersection numbers to the geometric quantities appearing in the ``local arithmetic Siegel--Weil'' theorems of our companion papers \cite{corank1_ASW_I.pdf,corank1_ASW_II.pdf}. Strictly speaking, the Rapoport--Zink uniformization we need at split places does not seem covered by the literature (not supersingular locus). We treat inert/ramified/split in parallel. For the most part, we disallow $p = 2$ only in the ramified case. Using the analogous local construction from \crefext{II:ssec:part_II:Hermitian_domain:corank1_modification}, we explain a modified Green current for singular $T$ (of rank $n - 1$ and size $n \times n$) in Section \ref{ssec:Arch_uniformization:Archimedean}.

The analytic side of the local-to-global reduction process will appear in our final companion paper \cite{corank1_ASW_IV.pdf}, where we complete the proof of our (global) co-rank $1$ arithmetic Siegel--Weil formula using essentially all preceding results.

Appendix \ref{appendix:K0} explains the setup we use for $K_0$ groups of Deligne--Mumford stacks. This setup is roughly as in \cite[{Appendix A}]{YZ17}, but they work over a field (we need schemes over $\mc{O}_F$, or at least over discrete valuation rings). Appendix \ref{appendix:quasi-compactness_lemma} contains some notation on abelian schemes, and records a proof for quasi-compactness of special cycles (which does not seem explicitly available in the literature). Appendix \ref{appendix:pDiv_prelim} concerns $p$-divisible groups, where we fix some notation and record some (presumably standard) facts.

Our algebro-geometric conventions follow the Stacks project \cite{stacks-project} unless stated otherwise.

            \subsection{Acknowledgements}
            \label{part_III:acknowledgements}
                I thank my advisor Wei Zhang for suggesting this topic, for his dedicated support and constant enthusiasm, for insightful discussions throughout the entire course of this project, and for helpful comments on earlier drafts. I thank Tony Feng, Qiao He, Benjamin Howard, Ishan Levy, Chao Li, Keerthi Madapusi, Andreas Mihatsch, Siddarth Sankaran, Ananth Shankar, Yousheng Shi, Tonghai Yang, Shou-Wu Zhang, and Zhiyu Zhang for helpful comments or discussions.

This work was partly supported by the National Science Foundation Graduate Research Fellowship under Grant Nos. DGE-1745302 and DGE-2141064.
Parts of this work were completed at the Mathematical Sciences Research Institute (MSRI), now becoming the Simons Laufer Mathematical Sciences Institute (SLMath), and the Hausdorff Institute for Mathematics. I thank these institutes for their support and hospitality. The former is supported by the National Science Foundation (Grant No. DMS-1928930), and the latter is funded by the Deutsche Forschungsgemeinschaft (DFG, German Research Foundation) under Germany's Excellence Strategy – EXC-2047/1 – 390685813.

    \clearpage


\part{Global special cycles}
\label{part:part_III:global_cycles}

    \section{Moduli stacks of abelian varieties} 
    \label{sec:ab_var}
        We discuss Kudla--Rapoport (KR) global special cycles on Rapoport--Smithling--Zhang (RSZ) smooth integral models of unitary Shimura varieties (which may be stacks). 
Fix an imaginary quadratic field extension $F / \Q$ with ring of integers $\mc{O}_F$ and write $a \mapsto a^{\s}$
for the nontrivial automorphism $\s$ of $F$. We write $\Delta \in \Z_{<0}$ and $\sqrt{\Delta} \in \mc{O}_F$ (pick a square root) for (generators of the) discriminant and different, respectively.
        
        \subsection{Integral models} 
        \label{ssec:ab_var:integral_models}
            \begin{definition}\label{definition:ab_var:integral_models:Hermitian_abelian_scheme}
Let $S$ be a scheme over $\Spec \mc{O}_F$. By a \emph{Hermitian abelian scheme} over $S$, we mean a tuple $(A, \iota, \lambda)$ where
    \begin{align*}
    & A & & \text{is an abelian scheme over $S$ of constant relative dimension $n$} 
    \\
    & \iota \colon \mc{O}_F \ra \End(A) & & \parbox[t]{0.75\textwidth}{is a ring homomorphism
    } 
    \\
    & \lambda \colon A \ra A^{\vee} & & \parbox[t]{0.75\textwidth}{is a quasi-polarization satisfying:
    \begin{itemize}[label={}]
    \item (Action compatibility) The Rosati involution $\dagger$ on $\End^{0}(A)$ satisfies $\iota(a)^{\dagger} = \iota(a^{\sigma})$ for all $a \in \mc{O}_F$.
    \end{itemize}
    }
\end{align*}
\end{definition}
An \emph{isomorphism} of Hermitian abelian schemes is an isomorphism of abelian schemes which respects the $\mc{O}_F$-actions and polarizations (exactly). For fixed $n \geq 1$, the \emph{moduli stack of Hermitian abelian schemes} $\ms{M}$ is the stack\footnote{By a \emph{stack in groupoids} over some base scheme $S$, we always mean a (not necessarily algebraic) stack in groupoids as in \cite[\href{https://stacks.math.columbia.edu/tag/02ZI}{Definition 02ZI}]{stacks-project} over the fppf site $(\Sch/S)_{fppf}$.} 
in groupoids over $\Spec \mc{O}_F$ with
    \begin{equation}
    \ms{M}(S) \coloneqq \{ \text{groupoid of relative $n$-dimensional Hermitian abelian schemes over $S$} \}
    \end{equation}
for $\mc{O}_F$-schemes $S$.

For an integer $r$ with $0 \leq r \leq n$, we next consider
    \begin{itemize}[label = {}]
    \item (Kottwitz $(n - r, r)$ signature condition) For all $a \in \mc{O}_F$, the characteristic polynomial of $\iota(a)$ acting on $\Lie A$ is $(x - a)^{n - r}(x - a^{\s})^r \in \mc{O}_S[x]$
    \end{itemize}
for pairs $(A, \iota)$, where $A \ra S$ is a relative $n$-dimensional abelian scheme with $\mc{O}_F$-action $\iota$, and $S$ is an $\mc{O}_F$-scheme. Here we view $\mc{O}_S$ as an $\mc{O}_F$-algebra via the structure map $S \ra \Spec \mc{O}_F$. This defines a substack\footnote{A \emph{substack} will always mean a strictly full substack.}
    \begin{equation}
    \ms{M}(n - r, r) \subseteq \ms{M}
    \end{equation}
consisting of Hermitian abelian schemes of signature $(n - r, r)$. The inclusion $\ms{M}(n - r, r) \ra \ms{M}$ is representable by schemes (in the sense of \cite[\href{https://stacks.math.columbia.edu/tag/04ST}{Section 04ST}]{stacks-project}) and is a closed immersion. There is an isomorphism\footnote{As in the Stacks project (e.g. \cite[\href{https://stacks.math.columbia.edu/tag/04XA}{Section 04XA}]{stacks-project}), we often abuse terminology and say ``isomorphism'' of stacks instead of ``equivalence''.} $\ms{M}(n - r, r) \ra \ms{M}(r, n - r)$ given by $(A, \iota, \lambda) \mapsto (A, \iota \circ \s, \lambda)$.

For any integer $d \geq 1$, there is a substack $\ms{M}^{(d)} \subseteq \ms{M}$ consisting of Hermitian abelian schemes $(A, \iota, \lambda)$ where $\lambda$ is polarization of constant degree $\deg \lambda \coloneqq \deg \ker \lambda = d$. If $\mc{A}_{n,d}$ (over $\Spec \mc{O}_F$) denotes the moduli stack of (relative) $n$-dimensional abelian schemes equipped with a polarization of degree $d$, the forgetful map $\ms{M}^{(d)} \ra \mc{A}_{n,d}$ is representable by schemes, finite, and unramified (e.g. via Lemma \ref{lemma:quasi-compactness_lemma}).
Hence $\ms{M}^{(d)}$ is a Noetherian Deligne--Mumford stack which is separated and finite type over $\Spec \mc{O}_F$ (because this is true of $\mc{A}_{n,d}$ as proved with level structure in the classical \cite[{\S 7.2 Theorem 7.9}]{MFK94}; one can deduce the stacky version upon inverting primes dividing the level, taking stack quotients, and patching over $\Spec \mc{O}_F$).

We set
    \begin{equation}
    \ms{M}(n - r, r)^{(d)} \coloneqq \ms{M}(n - r, r) \cap \ms{M}^{(d)}
    \end{equation}
where the right-hand side is an intersection of substacks of $\ms{M}$. There is an open and closed disjoint union decomposition\footnote{Here, the notation $\ms{M}^{(d)}[1/\Delta]$ means $\ms{M}^{(d)} \times_{\Spec \mc{O}_F} \Spec \mc{O}_F[1/\Delta]$. We often use such shorthand, along with subscripts for base change, e.g. $\ms{M}^{(d)}_S \coloneqq \ms{M}^{(d)} \times_{\Spec \mc{O}_F} S$ over an understood base.}
    \begin{equation}
    \ms{M}^{(d)}[1/\Delta] = \coprod_{(n-r,r)} \ms{M}(n-r,r)^{(d)}[1/\Delta]
    \end{equation}
over $\Spec \mc{O}_F[1/\Delta]$, where the disjoint union runs over all possible signatures $(n-r,r)$. 

The structure morphism $\ms{M}(n - r, r)^{(d)}[1/(d \Delta)] \ra \Spec \mc{O}_F[1/\Delta]$ is smooth of relative dimension $(n - r) r$ (e.g. by Remark \ref{remark:smoothness_moduli_stack} below; recall that being smooth of some relative dimension may be checked fppf locally on the target for morphisms of algebraic stacks).
We set $\ms{M}_0 \coloneqq \ms{M}(1,0)^{(1)}$. The structure morphism $\ms{M}_0 \ra \Spec \mc{O}_F$ is proper, quasi-finite,\footnote{Following the Stacks project \cite[\href{https://stacks.math.columbia.edu/tag/0CHU}{Definition 0CHU}]{stacks-project}, we require that finite morphisms of algebraic stacks are by definition (relatively) representable by schemes.
The morphism $\ms{M}_0 \ra \Spec \mc{O}_F$ is not finite in this sense, because $\ms{M}_0$ is not a scheme. Nevertheless, we continue to use terminology like ``representable by schemes and finite'' for morphisms of stacks which are not necessarily algebraic.} and \'etale by \cite[{Proposition 3.1.2}]{Howard12} or \cite[{Proposition 2.1.2}]{Howard15}.

Given any non-degenerate Hermitian $\mc{O}_F$-lattice $L$ of rank $n$ and signature $(n - r, r)$, we define an associated substack
    \begin{equation}
    \mc{M} \subseteq \ms{M}_0 \times_{\Spec \mc{O}_F} \ms{M}(n - r, r)
    \end{equation}
as follows (cf. \cite[{Proposition 2.12}]{KR14}, there in a principally polarized situation). 
Write $(-,-)$ for the pairing on $L$.
Let $b_L$ be the smallest positive integer such that $b_L \cdot (-,-)$ is $\mc{O}_F$-valued. Let $L'$ be the Hermitian $\mc{O}_F$ lattice which is the $\mc{O}_F$-module $L$ but with Hermitian pairing $b_L \cdot (-,-)$. Form the dual lattice $L^{\prime \vee}$ of $L'$ with respect to the Hermitian pairing, and set $d'_L \coloneqq |L^{\prime \vee} / L'|$.

If $L$ is self-dual of signature $(n - 1, 1)$ and $2 \nmid \Delta$, we set $d_L \coloneqq 1$ (the \emph{exotic smooth} setup for even $n$, see Section \ref{ssec:ab_var:exotic_smooth}). Otherwise, let $d_L \in \Z_{>0}$ be the product of ramified primes and the primes $p$ for which $L \otimes_{\Z} \Z_p$ is not self-dual.

\begin{definition}\label{definition:stack_ass_to_lattice}
Let $\mc{M} \subseteq \ms{M}_0 \times_{\Spec \mc{O}_F} \ms{M}(n - r, r)[1/(d_L \Delta)]$ be the substack
    \begin{equation}
    \mc{M}(S) \coloneqq \left \{ (A_0, \iota_0, \lambda_0, A, \iota, \lambda) : \begin{array}{l} \Hom_{\mc{O}_F \otimes \hat{\Z}^p} (T^p(A_{0,\overline{s}}), T^p(A_{\overline{s}})) \cong L \otimes_{\Z} \hat{\Z}^p \\ \text{for every geometric point $\overline{s}$ of $S$, with $p = \opn{char}(\overline{s})$,} \\ \text{and $b_L \cdot \lambda$ is a polarization of degree $d'_L$} \end{array} \right \}
    \end{equation}
for schemes $S$ over $\Spec \mc{O}_F[1/(d_L \Delta)]$, where
    \begin{equation}
    (A_0, \iota_0, \lambda_0) \in \ms{M}_0(S) \quad \quad (A, \iota, \lambda) \in \ms{M}(n - r, r)(S).
    \end{equation}
\end{definition}

\begin{warning}\label{warning:exotic_smooth_notation}
Whenever $L$ satisfies the even rank \emph{exotic smooth} setup (Section \ref{ssec:ab_var:exotic_smooth}), we will extend $\mc{M}$ to a smooth Deligne--Mumford stack which surjects onto $\Spec \mc{O}_F$ (see loc. cit.). In that case, we will override the notation here: after Section \ref{ssec:ab_var:exotic_smooth}, the notation $\mc{M}$ will always denote the exotic smooth moduli stack for such $L$. The restriction of the exotic smooth moduli stack over $\Spec \mc{O}_F[1/\Delta]$ will recover the stack in Definition \ref{definition:stack_ass_to_lattice}.
\end{warning}

In the definition of $\mc{M}$, the notation $\Hom_{\mc{O}_F \otimes \hat{\Z}^p} (T^p(A_{0,\overline{s}}), T^p(A_{\overline{s}})) \cong L \otimes_{\Z} \hat{\Z}^p$ asserts the existence of isomorphisms of Hermitian lattices, and the elements of $\Hom_{\mc{O}_F \otimes \hat{\Z}^p} (T^p(A_{0,\overline{s}}), T^p(A_{\overline{s}}))$ are not required to respect Hermitian pairings. As usual, $T^p(-)$ is the away-from-$p$ ad\`elic Tate module (if $p = 0$, this is over the full finite ad\`eles) and $\hat{\Z}^p = \prod_{\ell \neq p} \Z_{\ell}$. Note that $\mc{M}$ depends only on the \emph{ad\`elic isomorphism class}\footnote{We say that non-degenerate Hermitian $\mc{O}_F$ lattices $L$ and $L'$ are \emph{ad\`elically isomorphic} (or are in the same \emph{ad\`elic isomorphism class}) if there exist isomorphisms of $\mc{O}_F \otimes_{\Z} \Z_p$-Hermitian lattices $L \otimes_{\Z} \Z_p \cong L' \otimes_{\Z} \Z_p$ for every prime $p$, as well as isomorphisms of $\mc{O}_F \otimes_{\Z} \R$-Hermitian spaces $L \otimes_{\Z} \R \cong L' \otimes_{\Z} \R$ (classical terminology: \emph{genus}).} of $L$. The stack $\mc{M}$ (also the extension in Section \ref{ssec:ab_var:exotic_smooth}) and its special cycles will be the global moduli stacks of main interest in this work. We generally suppress $L$ from notation, but sometimes write $\mc{M}^L$ instead of $\mc{M}$ to emphasize $L$ dependence.

We claim that $\mc{M}$ is a Noetherian Deligne--Mumford stack which is separated and smooth of relative dimension $(n - r) r$ over $\Spec \mc{O}_F[1/(d_L \Delta)]$.
Indeed, there is an open and closed disjoint union decomposition
    \begin{equation}
    \ms{M}_0 \times_{\Spec \mc{O}_F} \ms{M}(n - r, r)^{(d)} [1/(d \Delta)] = \coprod_{L''} \mc{M}^{L''}
    \end{equation}
running over representatives $L''$, one for each ad\`elic isomorphism class of non-degenerate Hermitian $\mc{O}_F$-lattices of signature $(n - r, r)$ satisfying $L'' \subseteq L^{\prime \prime \vee}$ and $|L^{\prime \prime \vee} / L''| = d$. We have used flatness of $\ms{M}(n - r, r)^{(d)}[1/(d \Delta)] \ra \Spec \mc{O}_F[1/(d \Delta)]$ in the open and closed decomposition (to lift to characteristic $0$; cf. \cite[{Proposition 2.12}]{KR14} \cite[{Remark 4.2}]{RSZ18}).
With notation as above, the map
    \begin{equation}
    \begin{tikzcd}[row sep = tiny]
    \mc{M}^L \arrow{r} & \mc{M}^{L'} \\
    (A_0, \iota_0, \lambda_0, A, \iota, \lambda) \arrow[mapsto]{r} & (A_0, \iota_0, \lambda_0, A, \iota, b_L \lambda)
    \end{tikzcd}
    \end{equation}
is an isomorphism for any $L$, after restricting to $\Spec \mc{O}_F[1/(d_L \Delta)]$.

\begin{remark}\label{remark:ab_var:integral_models:n=1}
If $L$ has rank $n = 1$, we can construct $\mc{M}$ without discarding any primes. Then $\mc{M} \ra \Spec \mc{O}_F$ is smooth, by smoothness of $\ms{M}(1,0)^{(d)} \ra \Spec \mc{O}_F$ for any $d \in \Z_{>0}$.
\end{remark}

In the next lemma, $A_0^{\s}$ is the abelian scheme $A_0$ but with $\mc{O}_F$-action $\iota_0 \circ \s$.

\begin{lemma}\label{lemma:ab_var:integral_models:nonempty}
Let $L$ be any non-degenerate Hermitian $\mc{O}_F$-lattice of rank $n$ and signature $(n - r, r)$, with associated moduli stack $\mc{M}$. 
There exists a finite degree field extension $E / F$ and $(A_0, \iota_0, \lambda_0, A, \iota, \lambda) \in \mc{M}(\mc{O}_E[1/(d_L \Delta)])$ such that $A$ is $\mc{O}_F$-linearly isogenous to $A_0^{n - r} \times (A_0^{\s})^r$. 
In particular, $\mc{M}$ is nonempty.
\end{lemma}
\begin{proof}
First consider $\kappa = \C$ (equipped with a morphism $\mc{O}_F \ra \C$). Fix the trivializations of roots of unity $\Z / d \Z \xra{\sim} \pmb{\mu}_{d}(\C)$ sending $1 \mapsto e^{- 2 \pi i / d}$.

Choose $\sqrt{\Delta}$ to be the square-root whose image under $F \ra \C$ has positive imaginary part. We pass between Hermitian and alternating forms using the generator $\sqrt{\Delta}$ of the different ideal (as in \crefext{I:ssec:Hermitian_conventions:Hermitian_alternating_symmetric}). Express $L$ as a triple $(L, \iota, \lambda)$ where $\iota \colon \mc{O}_F \ra \End_{\Z}(L)$ is an action and $\lambda$ is a $\mc{O}_F$-compatible alternating pairing on $L$.

Take $(A_0, \iota_0, \lambda_0)$ to be the complex elliptic curve $\C / \mc{O}_F$. If $L_0 \coloneqq \mc{O}_F$ is the rank one Hermitian $\mc{O}_F$-lattice with Hermitian pairing $(x,y) = x^{\s} y$, we have $H_1(A_0,\Z) \cong L_0$ as Hermitian lattices.

Take any orthogonal decomposition $L_F = W \oplus W^{\perp}$ where $W$ is positive definite of rank $n - r$ and $W^{\perp}$ is negative definite of rank $r$. Define the $\C = \mc{O}_F \otimes_{\Z} \R$-action on $L \otimes_{\Z} \R$ to agree with $\iota$ on $W \otimes_F \C$ and to agree with $\iota \circ \s$ on $W^{\perp} \otimes_F \C$. This complex structure gives a tuple $(A, \iota, \lambda)$, where $A \coloneqq (L \otimes_{\Z} \R) / L$ is an abelian variety with $\mc{O}_F$-action $\iota$ and action compatible quasi-polarization $\lambda$. We have $H_1(A, \Z) \cong L$ as Hermitian lattices. By the usual comparison of $H_1(-,\Z)$ with $p$-adic Tate modules \cite[{\S 24 Theorem 1}]{Mumford85}, we conclude $(A_0, \iota_0, \lambda_0, A, \iota, \lambda) \in \mc{M}(\C)$.

We claim that $A$ is $\mc{O}_F$-linearly isogenous to $A_0^{n - r} \times (A^{\s}_0)^r$. Indeed, any $\mc{O}_F$-linear inclusion $\mc{O}_F^{n-r} \hookrightarrow L \cap W$ and any $\s$-linear inclusion $\mc{O}_F^{r} \hookrightarrow L \cap W^{\perp}$ will define an $\mc{O}_F$-linear isogeny $A_0^{n - r} \times (A^{\s}_0)^r \ra A$. 

Since $A_0$ is defined over some number field $\overline{\Q}$, it follows that $A$ and any isogeny $A_0^{n - r} \times (A_0^{\s})^r \ra A$ may also be defined over $\overline{\Q}$ (here using characteristic zero, so the kernel of the isogeny is \'etale). Descend these objects to some number field $E$.

Over a number field, it is a classical fact that any elliptic curve with $\mc{O}_F$-action has everywhere potentially good reduction \cite{Deuring41}. After extending $E$ if necessary, we thus obtain $(A_0, \iota_0, \lambda_0, A, \iota, \lambda) \in \mc{M}(\mc{O}_E[1/(d_L \Delta)])$ (the $\mc{O}_F$-actions extend by the N\'eron mapping property, and the polarizations extend to polarizations as in the proof of \cite[{Theorem 1.9}]{FC90}). The N\'eron mapping property extends the isogeny $A_0^{n - r} \times (A_0^{\s})^r \ra A$ over $\Spec \mc{O}_E[1/(d_L \Delta)]$.
\end{proof}

The preceding lemma will provide a base point for non-Archimedean uniformization (Section \ref{ssec:non-Arch_uniformization:framing}).
For arbitrary $L$ and specializing to any geometric point of characteristic $p > 0$, the abelian variety $A$ of Lemma \ref{lemma:ab_var:integral_models:nonempty} is supersingular (resp. ordinary) if $p$ is nonsplit (resp. split) by classical results of Deuring \cite{Deuring41} on endomorphism rings of elliptic curves.

\begin{notation}
Given a commutative ring $R$ with an automorphism $\s \colon R \ra R$ (e.g. $R = \mc{O}_F$), given a presheaf of modules $\mc{F}$ on a scheme $S$ over $\Spec R$, and given an action $\iota \colon R \ra \End(\mc{F})$ (with $\mc{F}$ viewed as a presheaf of abelian groups), we say the $R$ action via $\iota$ is \emph{$R$-linear} (resp. \emph{$\s$-linear}) if $\iota(a) = a$ (resp. $\iota(a) = a^{\s}$) for all $a \in R$. Here we view $\mc{O}_S$ as an $R$-algebra via the structure morphism $S \ra \Spec R$.
\end{notation}

Given any $(A, \iota, \lambda) \in \ms{M}(n - r, r)[1/\Delta]$ for a base scheme $S$, there is a canonical eigenspace decomposition
    \begin{equation}\label{equation:ab_var:integral_models:Lie_decomp_unr}
    \Lie A = (\Lie A)^+ \oplus (\Lie A)^-
    \end{equation}
characterized by $(\Lie A)^+$ (resp. $(\Lie A)^-$) being rank $n - r$ (resp. rank $r$) and the $\mc{O}_F$ action via $\iota$ on $(\Lie A)^+$ (resp. $(\Lie A)^-$) being $\mc{O}_F$-linear (resp. $\s$-linear). 

\begin{definition}\label{definition:ab_var:integral_models:tautological_bundle}
By the \emph{tautological bundle} on $\ms{M}(n-r,r)[1/\Delta]$, we mean the rank $r$ locally free sheaf $\ms{E}$ (for the fppf topology) whose dual is given by the assignment $\ms{E}^{\vee} \coloneqq (\Lie A)^-$ for $(A, \iota, \lambda) \in \ms{M}(n - r, r)[1/\Delta](S)$ for $\mc{O}_F$-schemes $S$.
\end{definition}

        \subsection{Exotic smoothness}
        \label{ssec:ab_var:exotic_smooth}
            Our main results at ramified primes eventually restrict to even $n$ and residue characteristic $\neq 2$. For this reason, we require $n$ even and $2 \nmid \Delta$ in Section \ref{ssec:ab_var:exotic_smooth}.

\begin{notation}[Exotic smooth setup, even rank]
A non-degenerate Hermitian $\mc{O}_F$-lattice $L$ of even rank $n$ satisfies the \emph{even rank exotic smooth} setup if $2 \nmid \Delta$, the signature of $L$ is $(n - 1, 1)$, and $L$ is self-dual (for the trace pairing).
\end{notation}

In the exotic smooth setup, we recall how $\mc{M}$ can be extended to a certain smooth integral model over $\Spec \mc{O}_F$.
We consider arbitrary signature $(n - r, r)$. For Hermitian abelian schemes $(A, \iota, \lambda)$, we consider
    \begin{enumerate}[label = {}]
        \item (Polarization condition $\circ$) The quasi-polarization $|\Delta| \cdot \lambda$ is a polarization and we have $\ker(|\Delta| \cdot \lambda) = A[\sqrt{\Delta}]$.
    \end{enumerate}
We write $\ms{M}(n - r, r)^{\mrm{Kot},\circ}$ for the substack of $\ms{M}(n - r, r)$ consisting of Hermitian abelian schemes $(A, \iota, \lambda)$ where $\lambda$ satisfies polarization condition $\circ$ from above. Here ``$\mrm{Kot}$'' indicates that we have ``only'' imposed the Kottwitz signature condition. For $d = |\Delta|^n$, the map
    \begin{equation}
    \begin{tikzcd}[row sep = tiny]
    \ms{M}(n - r, r)^{\mrm{Kot}, \circ} \arrow{r} & \ms{M}(n - r, r)^{(d)} \\
    (A, \iota, \lambda) \arrow[mapsto]{r} & (A, \iota, |\Delta| \cdot \lambda)
    \end{tikzcd}
    \end{equation}
is representable by closed immersions of schemes. In particular, $\ms{M}(n - r, r)^{\mrm{Kot}, \circ}$ is also a separated Deligne--Mumford stack which is finite type over $\Spec \mc{O}_F$. The restriction $\ms{M}(n - r, r)^{\mrm{Kot}, \circ}[1/\Delta] \ra \ms{M}(n - r, r)^{(d)}[1/\Delta]$ is an open immersion.
If $\kappa$ is any algebraically closed field of characteristic $0$, an object $(A, \iota, \lambda) \in \ms{M}(n - r, r)$ lies in $\ms{M}(n - r, r)^{\mrm{Kot}, \circ}$ if and only if the Hermitian $\mc{O}_F \otimes_{\Z} \Z_p$-lattice $T_p(A)$ is self-dual (for the trace pairing) for all $p$.
In particular, we have
    \begin{equation}
    \ms{M}_0 \times_{\Spec \mc{O}_F} \ms{M}(n - r, r)^{\mrm{Kot}, \circ}[1/\Delta] = \mc{M}^{L}
    \end{equation}
where $L$ is a representative for the unique ad\`elic isomorphism class of self-dual signature $(n - r, r)$ non-degenerate Hermitian $\mc{O}_F$-lattices (if it exists). Such $L$ exists if and only if $n \equiv 2^r \pmod{4}$ due to the global product formula for local invariants of Hermitian spaces.

Now we restrict to signature $(n - 1, 1)$ and $n \geq 2$. Let $\ms{M}(n - 1, 1)^{\circ}$ be the flat part of $\ms{M}(n - 1, 1)^{\mrm{Kot}, \circ}$, i.e. the scheme-theoretic image of the generic fiber. Equivalently, this is the largest closed substack which is flat over $\Spec \mc{O}_F$.

\begin{example}\label{example:integral_models:Serre_tensor_global}
Suppose $\mc{O}_F^{\times} = \{ \pm 1 \}$ (i.e. further exclude $F = \Q[\sqrt{-3}]$). If $\ms{M}_{\text{ell}} \ra \Spec \mc{O}_F$ denotes the moduli stack of elliptic curves base-changed to $\Spec \mc{O}_F$, the \emph{Serre tensor construction} $E \mapsto E \otimes_{\Z} \mc{O}_F$ defines an open and closed immersion $i_{\mrm{Serre}} \colon \ms{M}_{\text{ell}} \ra \ms{M}(1,1)^{\circ}$ \crefext{IV:ssec:arithmetic_Siegel-Weil:Serre_tensor}. If we replace $\mc{O}_F$ by (representatives of) fractional ideal classes for $\mc{O}_F$, we obtain an isomorphism $\coprod_{\mrm{Cl}(\mc{O}_F)} \ms{M}_{\text{ell}} \ra \ms{M}(1,1)^{\circ}$, where $\mrm{Cl}(\mc{O}_F)$ is the class group. This is \cite[{Proposition 14.4}]{KR14}. The local analogue (e.g. \crefext{I:lemma:moduli_pDiv:Serre_tensor}) plays an important role in this work.
In \crefext{IV:ssec:arithmetic_Siegel-Weil:Serre_tensor}, we revisit this description of $\ms{M}(1,1)^{\circ}$ to restate our main theorem in the simplest case.
\end{example}

Rapoport--Smithling--Zhang have given a moduli description \cite[{\S 6}]{RSZ21} for $\ms{M}(n - 1, 1)^{\circ}$. They define\footnote{Strictly speaking, Rapoport--Smithling--Zhang normalize their polarization differently (i.e. our $\lambda$ is their $|\Delta|^{-1} \lambda$). Their convention is more common elsewhere in the literature, and is of course equivalent to our formulation. We prefer our normalization, which seems more natural for our main results on the comparison with Eisenstein series Fourier coefficients. A related remark is \cite[{Footnote 9}]{LL22II}. \label{footnote:pi_inverse_modular}} a closed substack $\ms{M}(n - 1, 1)^{\mrm{RSZ}} \subseteq \ms{M}(n - 1, 1)^{\mrm{Kot}, \circ}$ with the same generic fiber, with $\ms{M}(n - 1, 1)^{\mrm{RSZ}} \ra \Spec \mc{O}_F$ smooth of relative dimension $n - 1$. This moduli description is as follows: given any scheme $S$ over $\Spec \mc{O}_F$, the groupoid $\ms{M}(n - 1, 1)^{\mrm{RSZ}}(S) \subseteq \ms{M}(n - 1, 1)^{\mrm{Kot}, \circ}(S)$ is the full subcategory consisting of tuples $(A, \iota, \lambda)$ such that the action $\iota \colon \mc{O}_F \ra \End(A)$ satisfies:
    \begin{enumerate}[(1)]
        \item (Pappas wedge condition) For all $a \in \mc{O}_F$, the action of $\iota(a)$ on $\Lie A$ satisfies
            \begin{equation*}
            \bigwedge^{2} (\iota(a) - a) = 0 \quad \text{ and } \quad \bigwedge^{n}(\iota(a) - a^{\s}) = 0.
            \end{equation*}
        \item (PRRSZ spin condition) For every geometric point $\overline{s}$ of $S$, the action of $(\iota(a) - a)$ on $\Lie A_{\overline{s}}$ is nonzero for some $a \in \mc{O}_F$.
\end{enumerate}
The signature condition implies that the equation involving $\bigwedge^{n}$ in the wedge condition is automatic, and that the wedge condition is empty if $n = 2$.
The wedge and spin conditions are automatic (given the signature condition) over $\Spec \mc{O}_F[1/\Delta]$, i.e. $\ms{M}(n-1,1)^{\mrm{Kot}, \circ}[1/\Delta] = \ms{M}(n-1,1)^{\mrm{RSZ}}[1/\Delta]$. For closedness of the spin condition, we refer to the closedness assertion in \cite[{Theorem 5.4}]{RSZ21}. The acronym PRRSZ stands for Pappas, Rapoport, Richarz, Smithling, and Zhang.
We have $\ms{M}(n-1,1)^{\circ} = \ms{M}(n-1,1)^{\mrm{RSZ}}$ by agreement in the generic fiber, flatness, and closedness.

We define the \emph{exotic smooth moduli stack} 
    \begin{equation}\label{equation:ab_var:exotic_smooth:M_circ}
    \mc{M}^{\circ} \coloneqq \ms{M}_0 \times_{\Spec \mc{O}_F} \ms{M}(n - 1, 1)^{\circ}
    \end{equation}
associated to any self-dual lattice $L$ of signature $(n - 1, 1)$.
The structure morphism $\mc{M}^{\circ} \ra \Spec \mc{O}_F$ is smooth, by the discussion above.

\begin{notation}\label{notation:exotic_smooth:notation_override}
From here on, we always write $\mc{M}$ instead of $\mc{M}^{\circ}$ if $L$ satisfies the even rank exotic smooth setup (we are overriding previous notation, see Warning \ref{warning:exotic_smooth_notation}; i.e. $\mc{M}[1/\Delta]$ recovers the moduli stack in Definition \ref{definition:stack_ass_to_lattice}). Recall that we have set $d_L \coloneqq 1$ for $L$ satisfying the even rank exotic smooth setup.
\end{notation}

\begin{remark}\label{remark:ab_var:exotic_smooth:nonempty}
Suppose $L$ satisfies the even rank exotic smooth setup, and form the associated moduli stack $\mc{M} \ra \Spec \mc{O}_F$. Then Lemma \ref{lemma:ab_var:integral_models:nonempty} holds for every geometric point $\Spec \kappa \ra \Spec \mc{O}_F$ by the same proof verbatim (replacing $d_L \Delta$ in loc. cit. with the number $1$).
\end{remark}

We have the following analogue of \eqref{equation:ab_var:integral_models:Lie_decomp_unr}: set
    \begin{equation}\label{equation:ab_var:integral_models:Lie_decomp_ram}
    (\Lie A)^+ \coloneqq \bigcap_{a \in \mc{O}_F} \ker (\iota(a) - a)|_{\Lie A}
    \end{equation}
for objects $(A, \iota, \lambda) \in \ms{M}(n - 1, 1)^{\circ}(S)$ over $\mc{O}_F$-schemes $S$. 

\begin{lemma}\label{lemma:unique_Kramer_hyperplane}
For objects $(A, \iota, \lambda) \in \ms{M}(n - 1, 1)^{\circ}(S)$ over $\mc{O}_F$-schemes $S$, the subsheaf $(\Lie A)^+ \subseteq \Lie A$ is a local direct summand of rank $n - 1$ whose formation commutes with arbitrary base change. The $\mc{O}_F$ action via $\iota$ on $(\Lie A)^+$ (resp. the line bundle $(\Lie A)/(\Lie A)^+$) is $\mc{O}_F$-linear (resp. $\s$-linear).
\end{lemma}
\begin{proof}
This lemma (and its proof) is a global analogue of \cite[{Lemma 2.36}]{LL22II} (the latter is an analogous statement on a Rapoport--Zink space).

Fix $a \in \mc{O}_F$ such that $\{1, a\}$ forms a $\Z$-basis of $\mc{O}_F$.
We have exact sequences
    \begin{equation}\label{equation:lemma:unique_Kramer_hyperplane:1}
    \begin{tikzcd}[row sep = small]
    0 \arrow{r} & (\Lie A)^+ \arrow{r} & \Lie A \arrow{r}{\iota(a) - a} \arrow{r} &  \operatorname{im}(\iota(a) - a) \arrow{r} & 0
    \\
    0 \arrow{r} & \operatorname{im}(\iota(a) - a) \arrow{r} & \Lie A \arrow{r} & \coker(\iota(a) - a) \arrow{r} & 0
    \end{tikzcd}
    \end{equation}
of quasi-coherent sheaves on $S$. The wedge and spin conditions imply that $(\Lie A)^+$ has rank $n - 1$ if $S = \Spec k$ for a field $k$. If $S$ is an arbitrary reduced scheme, the rank constancy of $\coker(\iota(a) - a)$ on geometric points implies that $\coker(\iota(a) - a)$ is finite locally free of rank $n - 1$ (e.g. by \cite[\href{https://stacks.math.columbia.edu/tag/0FWG}{Lemma 0FWG}]{stacks-project}). Hence, when $S$ is reduced, every sheaf appearing in \eqref{equation:lemma:unique_Kramer_hyperplane:1} is finite locally free, with $(\Lie A)^+$, $\operatorname{im}(\iota(a) - a)$, and $\coker(\iota(a) - a)$ having ranks $n-1$, $1$, and $n-1$ respectively. Thus the exact sequences of \eqref{equation:lemma:unique_Kramer_hyperplane:1} remain exact after pullback along any morphism of schemes $S' \ra S$ (where $S$ is reduced but $S'$ is not necessarily reduced). For arbitrary $S$ (not necessarily reduced), the morphism
$S \ra \ms{M}(n-1, 1)^{\circ}$ corresponding to $(A, \iota, \lambda)$ factors through a regular (hence reduced) locally Noetherian scheme fppf locally on $S$ (since the moduli stack $\ms{M}(n-1,1)^{\circ}$ is smooth over $\Spec \mc{O}_F$). These considerations show that $(\Lie A)^+ \subseteq \Lie A$ is a local direct summand of rank $n - 1$ whose formation commutes with arbitrary base change.

It is clear that the $\mc{O}_F$ action via $\iota$ on $(\Lie A)^+$ is $\mc{O}_F$-linear. To show that the action on $(\Lie A)/(\Lie A)^+$ is $\s$-linear, it is enough to check the case where $S$ is an integral scheme (argue fppf locally as above). When $S$ is integral, the $\s$-linearity follows from the $(n-1,1)$ signature condition on $\Lie A$.
\end{proof}

\begin{definition}\label{definition:ab_var:exotic_smooth:tautological_bundle}
By the \emph{tautological bundle} on $\ms{M}(n - 1, 1)^{\circ}$, we mean the rank one locally free sheaf $\ms{E}$ whose dual is $\ms{E}^{\vee} \coloneqq (\Lie A)/(\Lie A)^+$ for $(A, \iota, \lambda) \in \ms{M}(n - 1, 1)^{\circ}(S)$ for $\mc{O}_F$-schemes $S$.
\end{definition}

The restriction of $\ms{E}$ to $\ms{M}(n - 1, 1)^{\circ}[1/\Delta]$ coincides with the restriction of the tautological bundle defined in Definition \ref{definition:ab_var:integral_models:tautological_bundle}.

\begin{remark}\label{remark:ab_var:integral_models:other_models}
We mention how $\ms{M}(n-1,1)^{\mrm{RSZ}}$ relates to other moduli stacks. We caution that the terms ``Pappas model'' and ``Kr\"amer model'' in the literature may refer to variants, e.g. using principal polarizations. 

Let $\ms{M}(n - 1, 1)^{\mrm{Pap}} \subseteq \ms{M}(n-1,1)^{\mrm{Kot}, \circ}$ be the closed substack (``Pappas model'', named for the work \cite{Pappas00}) where we impose the wedge condition (for both $n$ even and odd) but not the spin condition. 

Let $\ms{M}(n - 1, 1)^{\mrm{Kr\ddot a}}$ be the stack in groupoids over $\Spec \mc{O}_F$ (``Kr\"amer model'', named for the work \cite{Kramer03})
consisting of tuples $(A, \iota, \lambda, \mc{F})$ for $(A, \iota, \lambda) \in \ms{M}(n - 1, 1)^{\mrm{Kot}, \circ}(S)$ and $\mc{F} \subseteq \Lie A$ a $\iota$-stable local direct summand of rank $n - 1$, such that the $\mc{O}_F$ action via $\iota$ on $\mc{F}$ (resp. $(\Lie A)/\mc{F}$) is $\mc{O}_F$-linear (resp. $\s$-linear). 

We have a diagram
    \begin{equation}
    \begin{tikzcd}
    \ms{M}(n - 1, 1)^{\mrm{RSZ}} \times_{\ms{M}(n - 1, 1)^{\mrm{Kot}, \circ}} \ms{M}(n - 1, 1)^{\mrm{Kr\ddot a}} \arrow[hook]{rr} \arrow{d} \arrow[dr, phantom, "\lrcorner", very near start] & & \ms{M}(n - 1, 1)^{\mrm{Kr\ddot a}} \arrow{d} 
    \\
    \ms{M}(n-1, 1)^{\mrm{RSZ}} \arrow[hook]{r} & \ms{M}(n - 1, 1)^{\mrm{Pap}} \arrow[hook]{r} & \ms{M}(n - 1, 1)^{\mrm{Kot}, \circ}
    \end{tikzcd}
    \end{equation}
where the horizontal arrows are closed immersions, the vertical arrows are forgetful, and the outer square is $2$-Cartesian. The left vertical arrow is an isomorphism (by Lemma \ref{lemma:unique_Kramer_hyperplane}, i.e. $\mc{F} = (\Lie A)^+$)
and the inclusion $\ms{M}(n-1, 1)^{\mrm{RSZ}} \hookrightarrow \ms{M}(n-1,1)^{\mrm{Pap}}$ is also an open immersion. All arrows are isomorphisms after base-change to $\Spec \mc{O}_F[1/\Delta]$.
\end{remark}
        
        \subsection{Special cycles} 
        \label{ssec:ab_var:special_cycles}
            Fix a non-degenerate Hermitian $\mc{O}_F$-lattice $L$ of rank $n$, with associated moduli stack $\mc{M}$.
The following definition of special cycles is due to Kudla--Rapoport \cite[{Definition 2.8}]{KR14} (there in a principally polarized situation). The notation $(\underline{x}, \underline{x})$ was explained in Section \ref{appendix:quasi-compactness_lemma:terminology}.

\begin{definition}[Kudla--Rapoport special cycles]\label{definition:global_special_cycles}
Given an integer $m \geq 0$, let $T \in \mrm{Herm}_m(\Q)$ be an $m \times m$ Hermitian matrix (with coefficients in $F$). The \emph{Kudla--Rapoport (KR) special cycle} $\mc{Z}(T)$ is the stack in groupoids over $\Spec \mc{O}_F$ defined as follows: for schemes $S$ over $\Spec \mc{O}_F$, we take $\mc{Z}(T)(S)$ to be the groupoid
    \begin{equation}
    \mc{Z}(T)(S) \coloneqq 
    \left \{ (A_0, \iota_0, \lambda_0, A, \iota, \lambda, \underline{x}) : \begin{array}{c} (A_0, \iota_0, \lambda_0, A, \iota, \lambda) \in \mc{M}(S) \\ \underline{x} = [x_1,\ldots,x_m] \in \Hom_{\mc{O}_F}(A_0,A)^m \\ (\underline{x}, \underline{x}) = T \end{array} \right \}.
    \end{equation}
We sometimes refer to elements $x \in \Hom_{\mc{O}_F}(A_0, A)$ as \emph{special homomorphisms}.
\end{definition}

\begin{example}
Suppose $2 \nmid \Delta$, and consider $L$ which is self-dual of signature $(1,1)$. Let $j \in \Z_{>0}$ be any positive integer. If $\mc{O}_F^{\times} = \{ \pm 1 \}$, consider the inclusion
    \begin{equation}
    \ms{M}_0 \times_{\Spec \mc{O}_F} \ms{M}_{\text{ell}} \xra{1 \times i_{\mrm{Serre}}} \ms{M}_0 \times_{\Spec \mc{O}_F} \ms{M}(1, 1)^{\circ} = \mc{M}
    \end{equation}
with $i_{\mrm{Serre}}$ as in Example \ref{example:integral_models:Serre_tensor_global}. Then $\mc{Z}(j) \ra \mc{M}$ pulls back to the $j$-th Hecke correspondence over the left-hand side, parameterizing triples $(E_0, E, w)$ where $E_0$ and $E$ are elliptic curves, $E_0$ has $\mc{O}_F$ action of signature $(1,0)$, and $w \colon E \ra E_0$ is an isogeny of degree $j$.
This is \cite[{Proposition 14.5}]{KR14}. 
We revisit this example in \crefext{IV:ssec:arithmetic_Siegel-Weil:Serre_tensor}, where we restate our main theorem in the simplest case via this description.
\end{example}

In the situation of Definition \ref{definition:global_special_cycles}, recall $\End_{\mc{O}_F}(A_0) = \mc{O}_F$ (if the right-hand side is abuse of notation for global sections of the constant sheaf $\mc{O}_F$ on $S$). 
If the Hermitian pairing on $L$ is $\mc{O}_F$-valued, we thus have $\mc{Z}(T) = \emptyset$ unless $T$ has coefficients in $\mc{O}_F$. If $L$ is self-dual and $2 \nmid \Delta$, we have $\mc{Z}(T) = \emptyset$ unless $\sqrt{\Delta} \cdot T$ has coefficients in $\mc{O}_F$.
Positivity of the Rosati involution also implies that the special cycle $\mc{Z}(T)$ is empty unless $T$ is positive semi-definite of rank $\leq n$.

By Lemma \ref{lemma:quasi-compactness_lemma}, the forgetful map $\mc{Z}(T) \ra \mc{M}$ is representable by schemes, finite, and unramified (and of finite presentation). Hence $\mc{Z}(T)$ is a separated Deligne--Mumford stack of finite type over $\Spec \mc{O}_F$. We will soon verify that $\mc{Z}(T)$ is smooth (after base change) over an explicit nonempty open subscheme of $\Spec \mc{O}_F$ (Lemma \ref{lemma:special_cycles_generically_smooth}).

We refer to $\mc{Z}(T) \ra \mc{M}$ as being a cycle over $\mc{M}$, although it is not literally a cycle (where the precise version of \emph{cycle} means a formal linear combination of integral closed substacks). We also do not wish to take pushforward, which may lose information. However, since $\mc{Z}(T) \ra \mc{M}$ is finite and unramified, this morphism is \'etale locally on the target a disjoint union of closed immersions \cite[\href{https://stacks.math.columbia.edu/tag/04HJ}{Lemma 04HJ}]{stacks-project}. For a more explicit version with level structure, see Lemma \ref{lemma:special_cycles_level_structure} below.

We record a few miscellaneous facts which will be used later. 
If $T_i \in \mrm{Herm}_{m_i}(\Q)$ for $i = 1, \ldots, j$ with $m \coloneqq m_1 + \cdots + m_j$ and all $m_i > 0$, then there is an identification
    \begin{equation}\label{equation:special_cycles_diagonal_entry_intersection}
    \mc{Z}(T_1) \times_{\mc{M}} \cdots \times_{\mc{M}} \mc{Z}(T_j) \cong \coprod_{\substack{T \in \mrm{Herm}_m(\Q) \\ \text{satisfying \eqref{equation:special_cycles_diagonal_entry_intersection_2}}}} \mc{Z}(T)
    \end{equation}
where the disjoint union runs over $T$ of the form
    \begin{equation}\label{equation:special_cycles_diagonal_entry_intersection_2}
    T = \begin{pmatrix} T_1 & & * \\ & \ddots & \\ * & & T_j \end{pmatrix}
    \end{equation}
(i.e. whose block diagonal entries are given by the $T_i$).

We write $\mc{Z}(T)_{\ms{H}} \subseteq \mc{Z}(T)$ for the largest closed substack flat over $\Spec \mc{O}_F$, and call $\mc{Z}(T)_{\ms{H}}$ a \emph{horizontal special cycle} or the \emph{flat part}. 
There is a decomposition of $\mc{Z}(T)$ as a scheme-theoretic union of closed substacks\footnote{By the \emph{scheme-theoretic union} of finitely many closed substacks $\mc{Z}_i$ of a Deligne--Mumford stack $\mc{Z}$, we mean the closed substack whose ideal sheaf is given by intersecting the ideal sheaves of $\mc{Z}_i$ on the small \'etale site of $\mc{Z}$.}
    \begin{equation}\label{equation:setup:ab_var:horizontal_vertical_decomp}
    \mc{Z}(T) = \mc{Z}(T)_{\ms{H}} \cup \bigcup_{p} \mc{Z}(T)_{\ms{V}, p}
    \end{equation}
where $\mc{Z}(T)_{\ms{V}, p} \coloneqq \mc{Z}(T) \times_{\Spec \Z} \Spec \Z/p^{e_p}\Z$ for a choice of $e_p \gg 0$ (notation $e_p$ is temporary). This follows from quasi-compactness of $\mc{Z}(T)$ (which also ensures that we may take $e_p = 0$ for all but finitely many $p$). We think of \eqref{equation:setup:ab_var:horizontal_vertical_decomp} as decomposition into a ``horizontal part'' and ``vertical parts''. A similar horizontal/vertical decomposition for local special cycles on Rapoport--Zink spaces is \cite[{\S 2.9}]{LZ22unitary} (inert case).

While the horizontal part $\mc{Z}(T)_{\ms{H}}$ is defined canonically, the vertical parts $\mc{Z}(T)_{\ms{V}, p}$ depend on $e_p$ as above. We will mostly work with the ``derived vertical special cycle classes'' from Section \ref{ssec:arith_cycle_classes:vertical}, which do not depend on such a choice of $e_p$.

Given any $T \in \mrm{Herm}_m(\Q)$ and $\gamma \in M_{n,n}(\mc{O}_F)$, there is a commutative diagram
    \begin{equation}\label{equation:ab_var:special_cycles:gamma_action}
    \begin{tikzcd}[row sep = tiny]
    (A_0, \iota_0, \lambda_0, A, \iota, \lambda, \underline{x}) \arrow[mapsto]{rr} & & (A_0, \iota_0, \lambda_0, A, \iota, \lambda, \underline{x} \cdot \gamma)
    \\
    \mc{Z}(T) \arrow{ddr} \arrow{rr} & & \mc{Z}({}^t \overline{\gamma} T \gamma) \arrow{ddl} 
    \\
    & &
    \\
    & \mc{M} &
    \end{tikzcd}
    \end{equation}
induced by $\gamma$. Below, we set $\mc{O}_{F,(p)} \coloneqq \mc{O}_F \otimes_{\Z} \Z_{(p)}$.

\begin{lemma}\label{lemma:ab_var:special_cycles:gamma_action}
Fix any $T \in \mrm{Herm}_m(\Q)$ and $\gamma \in M_{n,n}(\mc{O}_F)$. Consider the induced map
    \begin{equation}
    \mc{Z}(T) \ra \mc{Z}({}^t \overline{\gamma} T \gamma).
    \end{equation}
    \begin{enumerate}[(1)]
        \item This map is a finite morphism of algebraic stacks. If moreover $\gamma \in \GL_m(F)$ (resp. $\gamma \in \GL_m(\mc{O}_F)$) then the map is a closed immersion (resp. isomorphism).
        \item Assume $\gamma \in \GL_m(F)$, and let $N \in \Z$ be the product of primes $p$ such that $\gamma \not \in \GL_m(\mc{O}_{F,(p)})$. Then the restriction $\mc{Z}(T)[1/N] \ra \mc{Z}({}^t \overline{\gamma} T \gamma)[1/N]$ is an open immersion.
    \end{enumerate}
\end{lemma}
\begin{proof}
(1) The map $\mc{Z}(T) \ra \mc{Z}(T)({}^t \overline{\gamma} T \gamma)$ is finite because the projections to $\mc{M}$ are finite (finiteness for morphisms of algebraic stacks may be checked fppf locally on the target, so we reduce to the case of schemes). If $\gamma \in \GL_m(F)$, then $\mc{Z}(T) \ra \mc{Z}(T)({}^t \overline{\gamma} T \gamma)$ is a monomorphism of algebraic stacks (check via the moduli description), and any proper monomorphism of algebraic stacks is a closed immersion. If $\gamma \in \GL_m(\mc{O}_F)$, there is an inverse $\mc{Z}({}^t \overline{\gamma} T \gamma) \ra \mc{Z}(T)$ sending $\underline{x} \mapsto \underline{x} \cdot \gamma^{-1}$.

(2) Consider the substack $\mc{Z} \subseteq \mc{Z}({}^t \overline{\gamma} T \gamma)$ consisting of tuples $(A_0, \iota_0, \lambda_0, A, \iota, \lambda, \underline{w})$ such that $\underline{w} \cdot \gamma^{-1} \in \Hom_{\mc{O}_F}(A_0, A)^m$ (i.e. that the tuple of quasi-homomorphisms $\underline{w} \cdot \gamma^{-1}$ is a tuple of homomorphisms). The closed immersion $\mc{Z}(T) \ra \mc{Z}({}^t \overline{\gamma} T \gamma)$ maps isomorphically onto $\mc{Z}$. 

A quasi-homomorphism $f \colon B \ra B'$ of abelian schemes (over any base scheme $S$) is a homomorphism if and only the induced quasi-homomorphisms of $p$-divisible groups $f[p^{\infty}] \colon B[p^{\infty}] \ra B'[p^{\infty}]$ are homomorphisms for all primes $p$. This is a closed condition on $S$ for each prime $p$ (e.g. the quasi-homomorphism version of \cite[{Proposition 2.9}]{RZ96}). This is also an open condition for any prime $p$ which is invertible on $S$ (by \'etaleness of the $p$-divisible groups). If $p$ is a prime such that $\gamma \in \GL_m(\mc{O}_{F,(p)})$, then the tuple $\underline{w} \cdot \gamma^{-1}$ induces a tuple of quasi-homomorphisms $A_0[p^{\infty}] \ra A[p^{\infty}]$ consisting of homomorphisms.
\end{proof}
        
        \subsection{Level structure} 
        \label{ssec:ab_var:level_structure}
            Let $L$ be any non-degenerate Hermitian $\mc{O}_F$-lattice of rank $n$, and form the associated moduli stack $\mc{M}$. We discuss level structure for $\mc{M}$.

Set $V = L \otimes_{\mc{O}_F} F$, and form the unitary group $U(V)$ (over $\Spec \Q$). Let $L_0 = \mc{O}_F$ any self-dual Hermitian $\mc{O}_F$-lattice of rank $1$.
We set
    \begin{align*}
    & K_{L_0,p} \coloneqq \mrm{Stab}_{U(V_0)(\Q_p)}(L_0 \otimes \Z_p) & & K_{L,p} \coloneqq \mrm{Stab}_{U(V)(\Q_p)}(L \otimes \Z_p) 
    \\
    & K_{L_0,f} = \prod_{p} K_{L_0,p} && K_{L,f} = \prod_{p} K_{L,p}
    \end{align*}
for all $p$, where $\mrm{Stab}_{U(V)(\Q_p)}(L \otimes \Z_p)$ denotes the stabilizer of $L \otimes \Z_p$ in $U(V)(\Q_p)$, etc.. We say that $K_{L,f} \subseteq U(V)(\A_f)$ is the \emph{ad\`elic stabilizer} of $L$. Set $K'_{L,f} = K_{L_0,f} \times K_{L,f}$. Note that there is no dependence (up to functorial isomorphism) on the choice of $L_0$, or the choice of $L$ within its ad\`elic isomorphism class. We use the usual notation where $K_{L,f}^p$ means to omit the $p$-th factor in the product, etc..

For integers $N \geq 1$, we define the ``principal congruence subgroups''
    \begin{align*}
    & K_{p}(N) \coloneqq \ker(K_{L,p} \ra \GL(L \otimes \Z_p/N \Z_p)) && K_{f}(N) = \prod_{p} K_{L,p}(N)
    \end{align*}
(suppressing $L$ dependence from notation) and similarly define $K_{0,p}(N_0)$ and $K_f(N_0)$ for $N_0 \geq 1$. Given a pair $N' = (N_0, N)$ of integers $N_0, N \geq 1$, we set $K'_f(N') \coloneqq K_{0,f}(N_0) \times K_f(N)$ and $K'_p(N') \coloneqq K_{0,p}(N_0) \times K_p(N)$, etc..
Given $N' = (N_0, N)$, we sometimes abuse notation, e.g. $N' \geq a$ means $N_0, N \geq a$, and the notation $\mc{X}[1/N']$ for an algebraic stack $\mc{X}$ will mean inverting all primes dividing $N'$. If $N_0 = N$, we write $K'_f(N) \coloneqq K'_f(N')$.

Let $K'_f = K_{0,f} \times K_f \subseteq K'_{L,f}$ be any open compact subgroup which admits product factorizations $K_{0,f} = \prod_{p} K_{0,p}$ and $K_f = \prod_p K_p$. We set $K'_{p} \coloneqq K_{0,p} \times K_p$, etc.. Let $N_{K'_f}$ be the product of primes $p$ for which $K_{0,p} \neq K_{L_0,p}$ or $K_{p} \neq K_{L,p}$. We say that $K'_f$ is \emph{standard} at $p$ if $p \nmid N_{K'_f}$.

\begin{notation}
For $K'_f \subseteq K'_{L,f}$ as above, we reserve the term \emph{small} or \emph{small level} to mean that $K'_f \subseteq K'_f(N')$ for some $N' \geq 3$.
\end{notation}

Consider $\a = (A_0, \iota_0, \lambda_0, A, \iota, \lambda) \in \mc{M}(S)$ for some scheme $S$. Suppose $p$ is a prime which is invertible on $S$. For any integer $e \geq 0$, we consider the fppf sheaf
    \begin{equation}\label{equation:ab_var:level_structure:level_p^e}
    \underline{\mrm{Level}}(p^e) \subseteq \underline{\mrm{Isom}}(A_0[p^e], L_0 \otimes_{\Z} \Z / p^e \Z) \times \underline{\mrm{Isom}}(A[p^e], (L_0 \otimes_{\mc{O}_F} L) \otimes_{\Z} \Z / p^e \Z)
    \end{equation}
on $S$, where $\underline{\mrm{Level}}(p^e)$ is the open and closed subfunctor corresponding to isomorphisms $A_0[p^e] \ra L_0 / p^e L_0$ and $A[p^e] \ra L / p^e L$ which lift to $\mc{O}_F$-linear isomorphisms
    \begin{align*}
    T_p(A_{0,\overline{s}}) \xra{\sim} L_0 && ~ & \text{unitary up to scalar}
    \\
    \Hom_{\mc{O}_F}(T_p(A_{0,\overline{s}}),T_p(A_{\overline{s}})) \xra{\sim} L && ~ & \text{unitary}
    \end{align*}
over every geometric point $\overline{s}$ of $S$. Since $L_0$ is rank one, the ``unitary up to scalar'' condition is automatic.

Consider any open compact subgroup $K'_p \subseteq K'_{L,p}$. For each $e \geq 0$, we write $K'_p ~\mrm{mod}~ p^e$ (temporary notation) for the image of $K'_p$ in $\GL_1(L_0 \otimes_{\Z} \Z / p^e \Z) \times \GL_n(L \otimes_{\Z} \Z / p^e \Z)$. There is a canonical action of $K'_p ~\mrm{mod}~ p^e$ on $\underline{\mrm{Level}}(p^e)$.

\begin{definition}[Level structure]\label{definition:ab_var:level_structure}
Let $K'_f \subseteq K'_{L,p}$ be any factorizable open compact subgroup, as above. Consider an object $\a = (A_0, \iota_0, \lambda_0, A, \iota, \lambda) \in \mc{M}(S)$ for some scheme $S \ra \Spec \mc{O}_F[1/N_{K'_f}]$.

If $p$ is prime which is invertible on $S$, the \emph{sheaf of level $K'_p$ structures for $\a$} is the quotient (coequalizer)\footnote{Note that this quotient $\underline{\mrm{Level}}_{K'_p}$ is (isomorphic to) a sheaf (not just a pro-sheaf), so it is sensible to refer to $\underline{\mrm{Level}}_{K'_p}$ an fppf sheaf on $S$ (rather than a pro-object).}
    \begin{equation}
    \underline{\mrm{Level}}_{K'_p} \coloneqq (K'_p ~\mrm{mod}~ p^e)_{e \geq 0} \backslash (\underline{\mrm{Level}}(p^e))_{e \geq 0}
    \end{equation}
in the category of pro-objects for the category of fppf sheaves on $S$. If $p$ is not invertible on $S$, we let $\underline{\mrm{Level}}_{K'_p}$ be the constant sheaf valued in a singleton set.

The \emph{sheaf of level $K'_f$ structures for $\a$} is the product
    \begin{equation}
    \underline{\mrm{Level}}_{K'_f} \coloneqq \prod_p \underline{\mrm{Level}}_{K'_p}
    \end{equation}
over all primes $p$. This is an fppf sheaf on $S$, and is locally constant in the \'etale topology.

A \emph{$K'_p$ level structure} (rep. \emph{$K'_f$ level structure}) for $\a$ is a global section of $\underline{\mrm{Level}}_{K'_p}$ (resp. $\underline{\mrm{Level}}_{K'_f}$).
\end{definition}

\begin{remark}\label{remark:ab_var:level_structure_connected}
Let $K'_f$, $\a$, and $S$ be as in \cref{definition:ab_var:level_structure}. Fix a geometric point $\overline{s}$ of $S$, with $\mrm{char}(\overline{s}) = p \geq 0$. Assume moreover that $S$ is connected. In this case, giving a $K'_f$-level structure for $\a$ is (canonically) same as giving a pair $(\tilde{\eta}_0, \tilde{\eta})$ where $\tilde{\eta}_0$ (resp. $\tilde{\eta}$) is a $\pi_{1, \et}(S, \overline{s})$-stable $K_{0,f}^p$-orbit (resp. $K_f^p$-orbit) of isomorphisms
    \begin{align*}
    \eta_0 \colon T^p(A_{0, \overline{s}}) \xra{\sim} L_0 \otimes_{\Z} \hat{\Z}^p && ~ & \text{ unitary up to scalar}
    \\
    \eta \colon \Hom_{\mc{O}_F \otimes_{\Z} \hat{\Z}_p}(T^p(A_{0, \overline{s}}, A_{\overline{s}}) \xra{\sim} L \otimes_{\Z} \hat{\Z}^p && ~ & \text{ unitary}.
    \end{align*}
\end{remark}

In the notation of \cref{definition:ab_var:level_structure}, note that $K'_p = K'_p(p^e)$ implies $\underline{\mrm{Level}}_{K'_p} = \underline{\mrm{Level}}(p^e)$. 
In \cref{remark:ab_var:level_structure_connected}, note that the ``unitary up to scalar'' condition is automatic because $L_0$ has rank $1$. Even when $S$ is not connected, we abuse notation and write $(\tilde{\eta}_0, \tilde{\eta})$ for level $K'_f$ structure in the sense of \cref{ssec:ab_var:level_structure}

Given an open compact $K'_f$ as in Definition \ref{definition:ab_var:level_structure}, we now define a stack in groupoids $\mc{M}_{K'_f}$ over $\Spec \mc{O}_F[1/(d_L N_{K'_f})]$ with
    \begin{equation}
     \mc{M}_{K'_f}(S) \coloneqq \{ (\a, \tilde{\eta}_0, \tilde{\eta}) : \a \in \mc{M}(s) \text{ and } (\tilde{\eta}_0, \tilde{\eta}) \text{ a level $K'_f$ structure for $\a$} \}
    \end{equation}
for schemes $S$ over $\Spec \mc{O}_F[1/(d_L N_{K'_f})]$.
Given $T \in \mrm{Herm}_m(\Q)$, we write $\mc{Z}(T)_{K'_f} \coloneqq \mc{Z}(T) \times_{\mc{M}} \mc{M}_{K'_f}$ (``level $K'_f$ special cycle'').

Write $\mc{A}_{n,d,N}$ for the moduli stack over $\Spec \mc{O}_F[1/N]$ of (relative) $n$-dimensional abelian schemes $A$ with degree $d$ polarization and a chosen isomorphism $A[N] \ra L \otimes \Z / N \Z$ of group schemes (not necessarily compatible with symplectic pairings). We similarly form $\mc{A}_{1,1,N_0}$ using the lattice $L_0$ (and pick a basis of $L_0$ for convenience). Recall that $\mc{A}_{n,d,N}$ is representable by a separated Deligne--Mumford stack of finite type over $\Spec \mc{O}_F[1/N]$, and that $\mc{A}_{n,d,N}$ is a scheme quasi-projective over $\Spec \mc{O}_F[1/N]$ if $N \geq 3$ (see \cite[{\S 7.2 Theorem 7.9}]{MFK94}).

Let $b_L, d_L \in \Z_{>0}$ be associated to $L$, as discussed before Definition \ref{definition:stack_ass_to_lattice}. If $K'_f(N')$ is the principal congruence subgroup of some level $N' = (N_0, N)$, consider the forgetful morphism 
    \begin{equation}
    \begin{tikzcd}[row sep = tiny]
    \mc{M}_{K'_f(N')} \arrow{r} & \mc{A}_{1,1,N_0} \times \mc{A}_{n,d_L,N}
    \end{tikzcd}
    \end{equation}
which forgets the $\mc{O}_F$-actions and sends $\lambda \mapsto b_L \lambda$. For level structure, see the description in \cref{equation:ab_var:level_structure:level_p^e}. The induced map
    \begin{equation}\label{equation:ab_var:level_structure:map_to_Ag}
    \begin{tikzcd}[row sep = tiny]
    \mc{M}_{K'_f(N')} \ra \mc{M}[1/N'] \times_{(\mc{A}_{1,1} \times \mc{A}_{n,d_L})} (\mc{A}_{1,1,N_0} \times \mc{A}_{n,d_L,N})[1/d_L]
    \end{tikzcd}
    \end{equation}
is representable by schemes and is an open and closed immersion. Hence $\mc{M}_{K'_f(N')} \ra \mc{A}_{1,1,N_0} \times \mc{A}_{n,d_L,N}[1/d_L]$ is finite (and representable by schemes).

\begin{lemma}\hfill
\begin{enumerate}[(1)]
    \item For any open compact subgroup $K'_f$ as in Definition \ref{definition:ab_var:level_structure}, the stack $\mc{M}_{K'_f}$ is a separated Deligne--Mumford stack of finite type over $\Spec \mc{O}_F$. If $K'_f$ is small, then $\mc{M}_{K'_f}$ is a quasi-projective scheme over $\Spec \mc{O}_F$.
    \item For any inclusion $K'_f \subseteq K''_f$, the forgetful morphism $\mc{M}_{K'_f} \ra \mc{M}_{K''_f}[1/N_{K'_f}]$ (i.e. expand a $K'_f$-orbit to a $K''_f$ orbit) is finite \'etale of degree $|K''_f/K'_f|$. If $K'_f \subseteq K''_f$ is a normal subgroup, then $\mc{M}_{K'_f} \ra \mc{M}_{K''_f}[1/N_{K'_f}]$ is a torsor for the finite discrete group $K''_f / K'_f$.
\end{enumerate}
\end{lemma}
\begin{proof}
The second sentence in part (2) is clear from construction (and makes sense before we know these stacks are algebraic). When $K'_f = K'_f(N')$ for some $N'$, the claims in part (1) follow from \eqref{equation:ab_var:level_structure:map_to_Ag}.

For general $K'_f$, select $N' = (N_0, N)$ such that $\mc{M}_{K'_f(N')}$ is a scheme and $K'_f(N') \subseteq K'_f$. Then $\mc{M}_{K'_f(N')} \ra \mc{M}_{K'_f}[1/N']$ is a torsor for the finite discrete group $K'_f/K'_f(N')$ (in particular, finite \'etale), and hence admits the stack quotient presentation $\mc{M}_{K'_f}[1/N'] \cong [\mc{M}_{K'_f(N')}/ (K'_f/K'_f(N'))]$, which shows that $\mc{M}_{K'_f}[1/N']$ is Deligne--Mumford. Picking another $M' = (M_0, M)$ such that $\gcd(N_0 N, M_0 M)$ is only divisible by primes dividing $N_{K'_f}$, we find that $\mc{M}_{K'_f}[1/M']$ is Deligne--Mumford as well. These two charts show that $\mc{M}_{K'_f}$ is Deligne--Mumford, as well as separated and finite type over $\Spec \mc{O}_F$.

If $K'_f \subseteq K''_f$, then for any scheme $S$ with a morphism $S \ra \mc{M}_{K''_f}$, the $2$-fiber product $\mc{M}_{K'_f} \times_{\mc{M}_{K''_f}} S$ is fibered in setoids, hence equivalent to a sheaf (of sets). But since $\mc{M}_{K'_f(N')} \ra \mc{M}_{K'_f}[1/N']$ is a $K'_f/K'_f(N')$-torsor and affine morphisms satisfy fpqc descent \cite[\href{https://stacks.math.columbia.edu/tag/0244}{Section 0244}]{stacks-project}, we conclude that $\mc{M}_{K'_f}[1/N'] \times_{\mc{M}_{K''_f}} S$ is represented by a scheme. As above, we may pick some other $M'$ to patch and show that the morphism $\mc{M}_{K'_f} \ra \mc{M}_{K''_f}[1/N_{K'_f}]$ is representable by schemes. Since $\mc{M}_{K'_f(N')} \ra \mc{M}_{K'_f}[1/N']$ and $\mc{M}_{K'_f(N')} \ra \mc{M}_{K''_f}[1/N']$ are both finite \'etale surjections, we conclude that $\mc{M}_{K'_f} \ra \mc{M}_{K''_f}[1/N_{K'_f}]$ is finite \'etale by varying $N'$ again (using standard facts like \cite[\href{https://stacks.math.columbia.edu/tag/02K6}{Lemma 02K6}, \href{https://stacks.math.columbia.edu/tag/01KV}{Lemma 01KV}, \href{https://stacks.math.columbia.edu/tag/0AH6}{Lemma 0AH6}, \href{https://stacks.math.columbia.edu/tag/02LS}{Lemma 02LS}]{stacks-project}). The remaining claims follow from this.
\end{proof}

\begin{lemma}\label{lemma:special_cycles_level_structure}
Fix any prime $p$ and a matrix $T \in \mrm{Herm}_m(\Q)$ with $m \geq 0$. The morphism $\mc{Z}(T)_{K'_f(p^e)} \ra \mc{M}_{K'_f(p^e)}$ is a disjoint union of closed immersions for all $e \gg 0$.
\end{lemma}
\begin{proof}

For $e \in \Z_{\geq 0}$, we define a stack (used only in this proof) $\mc{M}(p^e)$ over $\Spec \mc{O}_F[1/(d_L p)]$ as follows. For schemes $S$ over $\Spec \mc{O}_F[1/(d_L p)]$, we take $\mc{M}(p^e)(S)$ to be the groupoid
    \begin{equation}
    \mc{M}(p^e)(S) \coloneqq \left \{ (A_0, \iota_0, \lambda_0, A, \iota, \lambda, \underline{x}) : \begin{array}{c} (A_0, \iota_0, \lambda_0, A, \iota, \lambda) \in \mc{M}(S) \\ \underline{x} = [x_1,\ldots,x_m] \in \Hom_{\mc{O}_F}(A_0[p^e],A[p^e])^m \end{array} \right \}.
    \end{equation}
We have a commutative diagram
    \begin{equation}\label{equation:special_cycles_level_structure_1}
    \begin{tikzcd}
    & \mc{M}(p^e) \arrow{d} \\
    \mc{Z}(T)[1/p] \arrow{ur} \arrow{r} & \mc{M} \rlap{~.}
    \end{tikzcd}
    \end{equation}
The forgetful morphism $\mc{M}(p^e) \ra \mc{M}[1/p]$ is representable by schemes and a finite \'etale surjection. Thus, $\mc{M}(p^e)$ is representable by a separated Deligne--Mumford stack of finite type over $\Spec \mc{O}_F$.

We claim that $\mc{Z}(T) \ra \mc{M}(p^e)$ is a closed immersion for $e \gg 0$. 
This may be checked fppf locally on the target. Suppose $S \ra \mc{M}[1/p]$ is an fppf cover by a Noetherian scheme $S$ (possible since $\mc{M}$ is locally Noetherian and quasi-compact). It is enough to check that $\mc{Z}(T) \times_{\mc{M}} S \ra \mc{M}(p^e) \times_{\mc{M}} S$ is a closed immersion of schemes.
Since the morphism $\mc{Z}(T) \ra \mc{M}$ (resp. $\mc{M}(p^e) \ra \mc{M}[1/p]$) is finite and unramified (resp. finite), we conclude that $\mc{Z}(T)[1/p] \ra \mc{M}(p^e)$ is also finite and unramified.

To prove the claim, it remains only to check that the morphism of schemes $\mc{Z}(T) \times_{\mc{M}} S \ra \mc{M}(p^e) \times_{\mc{M}} S$ is universally injective for $e \gg 0$ (for morphisms of schemes, being a closed immersion is the same as being proper, unramified, and universally injective \cite[\href{https://stacks.math.columbia.edu/tag/04XV}{Lemma 04XV}]{stacks-project}).

We first show that universal injectivity holds fiber-wise over every point $s \in S$ for $e$ sufficiently large (with $e$ possibly depending on $s$). For any point $s$ on $S$ with residue field $k(s)$, we know that $\mc{Z}(T)_{k(s)} \ra \mc{M}(p^e)_{k(s)}$ is universally injective for $e \gg 0$ (possibly depending on $s$) because the map $\Hom(A_1, A_2) \ra \Hom(T_{p}(A_1), T_{p}(A_2))$ is injective for abelian varieties $A_1$, $A_2$ over any field of characteristic $\neq p$ (apply this over a geometric point mapping to $s$ and use finiteness of $\mc{Z}(T)$).

Being universally injective may be checked fiber-wise over $S$, so we need to show that there is a value of $e$ which works for all points $s \in S$ simultaneously. We can select $e \gg 0$ so that $\mc{Z}(T) \times_{\mc{M}} S \ra \mc{M}(p^e) \times_{\mc{M}} S$ is universally injective (hence a closed immersion) over the generic point of each irreducible component of $S$. For such $e$, a limiting argument (``spreading out'') implies that $\mc{Z}(T) \times_{\mc{M}} S \ra \mc{M}(p^e) \times_{\mc{M}} S$ is a closed immersion over an open dense subset of $S$. Applying Noetherian induction on $S$ proves the claim.

To finish the proof of the lemma, we observe that $\mc{M}(p^e) \times_{\mc{M}} \mc{M}_{K'_f(p^e)} \ra \mc{M}_{K'_f(p^e)}$ is a finite disjoint union of isomorphisms, corresponding to the constant sheaf valued in $\Hom_{\mc{O}_F}(L_0 \otimes_{\Z} \Z/p^e \Z, L \otimes_{\Z} \Z/p^e \Z)^m$. Hence $\mc{Z}(T)_{K'_f(p^e)} \ra \mc{M}_{K'_f(p^e)}$ is a disjoint union of closed immersions.
\end{proof}
        
        \subsection{Generic smoothness} 
        \label{ssec:ab_var:generic_smoothness}
            We explain a generic smoothness result for special cycles (Lemma \ref{lemma:special_cycles_generically_smooth}).
The other lemmas are auxiliary. The proof proceeds by reducing to $p$-divisible groups over a base where $p$ is locally nilpotent, and then checking formal smoothness using Serre--Tate and Grothendieck--Messing deformation theory. 

We first consider $p$-divisible groups (see Appendix \ref{appendix:pDiv_prelim:terminology} for a review of terminology). For primes $p$, set $\mc{O}_{F_p} \coloneqq \mc{O}_{F} \otimes_{\Z} \Z_p$. Suppose $p \nmid \Delta$ and consider schemes $S$ over $\Spf \mc{O}_{F_p}$, i.e. $S$ is a scheme over $\Spec \mc{O}_F$ with $p$ locally nilpotent on $S$. We consider tuples $(Y, \iota, \lambda)$ over $S$ where
    \begin{align}
    & Y & & \parbox[t]{0.5 \textwidth}{is a $p$-divisible group over $S$ of height $2 n$ and dimension $n$} 
    \notag
    \\
    & \iota \colon \mc{O}_{F_p} \ra \End(Y) & & \parbox[t]{0.5\textwidth}{is an action satisfying the $(n-r,r)$ Kottwitz signature condition, i.e. for all $a \in \mc{O}_{F_p}$, the characteristic polynomial of $\iota(a)$ acting on $\Lie Y$ is $(x - a)^{n - r}(x - a^{\s})^{r} \in \mc{O}_S[x]$} 
    \label{equation:p_div_generic_smoothness}    
    \\
    & \lambda \colon Y \xra{\sim} Y^{\vee} & & \parbox[t]{0.5\textwidth}{is a principal polarization whose Rosati involution $\dagger$ on $\End(Y)$ satisfies $\iota(a)^{\dagger} = \iota(a^{\sigma})$ for all $a \in \mc{O}_{F_p}$.}
    \notag
    \end{align}
In the signature condition described above, we view $\mc{O}_S$ as an $\mc{O}_{F_p}$-algebra via the structure map $S \ra \Spf \mc{O}_{F_p}$.

Parts of the next formal smoothness result (Lemma \ref{lemma:formally_smooth_deformations_p_div}) may exist in some form in the literature, see e.g. discussion about formal smoothness for ``unramified Rapoport--Zink data'' in \cite[{3.82}]{RZ96} and the reference to \cite[{\S 5}]{Kottwitz92} given there.

Following \cite[\href{https://stacks.math.columbia.edu/tag/04EW}{Section 04EW}]{stacks-project}, we use the term \emph{thickening} to refer to a closed immersion which is a homeomorphism on underlying topological spaces, and the term \emph{first order thickening} for a thickening defined by a square zero ideal.

Let $S$ be a scheme over $\Spf \mc{O}_{F_p}$, and suppose $(Y, \iota, \lambda)$ is a tuple over $S$ as in \eqref{equation:p_div_generic_smoothness}. There is an associated deformation functor $\mrm{Def}_{(Y,\iota,\lambda)}$ (possibly non-standard usage, and it will not appear after Lemma \ref{lemma:formally_smooth_deformations_p_div}) which sends a thickening $S \ra S'$ to the set of (isomorphism classes of) lifts of $(Y, \iota, \lambda)$ to $S'$. Write $S[\e]$ and $S[\e,\e']$ as shorthand for $S \times_{\Spec \Z} \Spec \Z[\e]/(\e^2)$ and $S \times_{\Spec \Z} \Spec \Z[\e, \e']/(\e^2, \e\e', \e^{\prime 2})$, respectively. In the proof of Lemma \ref{lemma:formally_smooth_deformations_p_div}, we will see that the canonical map 
    \begin{equation}\label{equation:tangent_space_condition_E}
    \mrm{Def}_{(Y,\iota,\lambda)}(S[\e, \e']) \ra \mrm{Def}_{(Y,\iota,\lambda)}(S[\e]) \times_{\mrm{Def}_{(Y,\iota,\lambda)}(S)} \mrm{Def}_{(Y,\iota,\lambda)}(S[\e'])
    \end{equation}
is an isomorphism. More generally, if $M$ is a finite rank free $\mc{O}_S$-module and $\mc{O}_S \oplus M$ denotes the quasi-coherent $\mc{O}_S$-algebra with $M$ an ideal of square zero, we will see that the functor $M \mapsto \mrm{Def}_{(Y,\iota,\lambda)}(\underline{\Spec}_S(\mc{O}_S \oplus M))$ preserves fiber products over the base $M = 0$ (note that this holds when $\mrm{Def}_{(Y,\iota,\lambda)}$ is replaced by any scheme, and this is essentially the method of proof). Here $\underline{\Spec}_S$ denotes relative $\Spec$.

For any scheme $S$ over $\Spf \mc{O}_{F_p}$, the above considerations imply that the set $\mrm{Def}_{(Y,\iota,\lambda)}(S[\e])$ has the natural structure of a $\Gamma(S, \mc{O}_S)$-module in the standard way (as a ``tangent space'') as in \cite[{Proposition 3.6}]{SGA3II} or \cite[\href{https://stacks.math.columbia.edu/tag/06I2}{Section 06I2}]{stacks-project}.

\begin{lemma}\label{lemma:formally_smooth_deformations_p_div}
Let $p$ be a prime which is unramified in $\mc{O}_F$. The deformation problem for triples as in \eqref{equation:p_div_generic_smoothness} is formally smooth of relative dimension $(n - r) r$ in the following sense. Let $S$ be any scheme over $\Spf \mc{O}_{F_p}$, and let $(Y,\iota,\lambda)$ be a triple over $S$ as in \eqref{equation:p_div_generic_smoothness}. 
    \begin{enumerate}[(1)]
        \item The triple $(Y, \iota, \lambda)$ lifts along any first order thickening of affine schemes $S \ra S'$, i.e. the map $\mrm{Def}_{(Y,\iota,\lambda)}(S') \ra \mrm{Def}_{(Y,\iota,\lambda)}(S)$ is surjective. 
        \item When $S = \Spec \kappa$ for a field $\kappa$, the $\kappa$ vector space $\mrm{Def}_{(Y,\iota,\lambda)}(\kappa[\e])$ has dimension $(n - r) r$.
    \end{enumerate}
If $(n - r) r = 0$, then $(Y,\iota,\lambda)$ lifts uniquely along any first order thickening of schemes $S \ra S'$.
\end{lemma}
\begin{proof}
We study this lifting problem for $p$-divisible groups in terms of Grothendieck--Messing deformation theory.
Let $S \ra S'$ be a first order thickening of schemes (not necessarily affine).  View $S \hookrightarrow S'$ as a PD thickening, with trivial PD structure on the square zero ideal of the thickening.

Write $\mbb{D}(Y)$ for the covariant Dieudonn\'e crystal of $Y$, and write $\mbb{D}(Y)(S)$ and $\mbb{D}(Y)(S')$ for the evaluation of this crystal on the PD thickenings $\mrm{id} \colon S \ra S$ and $S \hookrightarrow S'$ respectively. We have a short exact sequence of $\mc{O}_S$-modules given by the Hodge filtration
    \begin{equation}\label{equation:deformation_Hodge_filtration}
    0 \ra \Omega_{Y^{\vee}} \ra \mbb{D}(Y)(S) \ra \Lie_Y \ra 0
    \end{equation}
with $\Omega_{Y^{\vee}} = (\Lie_{Y^{\vee}})^{\vee}$ and each $\mc{O}_S$-module above being finite locally free. 

We may decompose the Hodge filtration into eigenspaces with respect to the action $\iota \colon \mc{O}_{F_p} \ra \End(Y)$ (and the structure morphism $S \ra \Spec \mc{O}_{F_p}$). We use superscripts $(-)^+$ and $(-)^-$ to denote these eigenspaces, where $\mc{O}_{F_p}$ acts linearly (resp. $\s$-linearly) on $(-)^+$ (resp. $(-)^-$) via $\iota$.
Then we have short exact sequences 
    \begin{align}
    & 0 \ra \Omega_{Y^{\vee}}^+ \ra \mbb{D}(Y)(S)^+ \ra \Lie_{Y}^+ \ra 0 \label{equation:0_eigenspace_deformation_Hodge_filtration}
    \\
    & 0 \ra \Omega_{Y^{\vee}}^- \ra \mbb{D}(Y)(S)^-\ra \Lie_{Y}^- \ra 0
    \label{equation:1_eigenspace_deformation_Hodge_filtration}
    \end{align}
where each $\mc{O}_S$-module above is finite locally free and, for example, we have $\mbb{D}(Y) = \mbb{D}(Y)^+ \oplus \mbb{D}(Y)^-$. From left to right, the modules in \eqref{equation:0_eigenspace_deformation_Hodge_filtration} have ranks $r$, $n$, and $n - r$, and the modules in \eqref{equation:1_eigenspace_deformation_Hodge_filtration} have ranks $n - r$, $n$, and $r$ respectively.

Using the polarization $\lambda$, we may identify \eqref{equation:0_eigenspace_deformation_Hodge_filtration} with the dual of \eqref{equation:1_eigenspace_deformation_Hodge_filtration}. There is a choice of sign in this identification, which plays essentially no role in this proof.

We have $\mbb{D}(Y)(S')|_S \cong \mbb{D}(Y)(S)$ canonically (as $\mbb{D}(Y)$ is a crystal), and Grothendieck--Messing theory implies that lifting $(Y,\iota,\lambda)$ to $S'$ is the same as lifting the Hodge filtration \eqref{equation:deformation_Hodge_filtration} compatibly with the action $\iota$ and the polarization $\lambda$. Compatibility with the $\iota$ action means that we should lift the eigenspace decomposition in \eqref{equation:0_eigenspace_deformation_Hodge_filtration} and \eqref{equation:1_eigenspace_deformation_Hodge_filtration}, and compatibility with the polarization $\lambda$ means that the resulting exact sequences should again be dual to each other (as determined by $\lambda$). It is equivalent to lift either one of the exact sequences of \eqref{equation:0_eigenspace_deformation_Hodge_filtration} and \eqref{equation:1_eigenspace_deformation_Hodge_filtration} (one determines the other) to a filtration of $\mbb{D}(Y)(S')^+$ or $\mbb{D}(Y)(S')^-$ respectively (with no additional restrictions).

Consider the lifting problem for, say, the $+$ eigenspace of the Hodge filtration as in \eqref{equation:0_eigenspace_deformation_Hodge_filtration}. Zariski locally on $S'$, this lifting problem may be identified with the problem of lifting an $S$ point to an $S'$ point on the Grassmannian parametrizing rank $r$ subbundles of the rank $n$ trivial bundle. This Grassmannian is smooth of relative dimension $(n - r) r$, which proves the claims in the lemma statement.
\end{proof}

The next three lemmas are used to prove Lemma \ref{lemma:smooth_Jacobson_formal_completions}. This latter lemma is in turn used in the proof of generic smoothness in Lemma \ref{lemma:special_cycles_generically_smooth}, to reduce to bases where $p$ is locally nilpotent for some unramified prime $p$. This will allow us to reduce to formal smoothness for deformations of $p$-divisible groups (with certain additional structure) as proved in Lemma \ref{lemma:formally_smooth_deformations_p_div}.

\begin{lemma}\label{lemma:flatness_smoothness_Spf_schemes}
Let $A$ be an adic Noetherian ring, and let $X$ be a locally Noetherian scheme over $\Spec A$. If $X_{\Spf A} \ra \Spf A$ is flat, then $X \ra \Spec A$ is flat at every point of $X$ which lies over $\Spf A$. If $X \ra \Spec A$ is locally of finite type, then the same holds with ``flat'' replaced by ``smooth''.
\end{lemma}
\begin{proof}

Here, flatness (resp. smoothness) of $X_{\Spf A} \ra \Spf A$ is equivalent to the requirement that, for every scheme $T$ with a map $T \ra \Spf A$, the base changed map $X_T \ra T$ is flat (resp. smooth).

We first check the flatness assertion. Passing to an affine open of $X$, we may reduce to the case where $X = \Spec B$ for a Noetherian ring $B$. Let $I \subseteq A$ be an ideal of definition. Then $X_{\Spf A}$ is described by a completed tensor product, and we have $X_{\Spf A} = \hat{B}$ where $\hat{B}$ is the $I$-adic completion of $B$. Since $B$ is a Noetherian ring, the canonical map $B \ra \hat{B}$ is flat. Since $X_{\Spf A} \ra \Spf A$ is flat, we know that $A \ra \hat{B}$ is a flat ring map. We conclude that $B$ is flat over $A$ at every prime in the image of $\Spec \hat{B} \ra \Spec B$. These are precisely the points of $X$ lying over $\Spf A$.

Next, assume $X_{\Spf A} \ra \Spf A$ is smooth. By Noetherianity of $A$, the map $f \colon X \ra \Spec A$ is locally of finite presentation. We have just shown that $X \ra \Spec A$ is flat at every point $x \in X$ which lies over $\Spf A$. Thus, for such $x \in X$, the map $f \colon X \ra \Spec A$ is smooth at $x$ if and only if $X_{f(x)} \ra \Spec k(f(x))$ is smooth at $x$, where $k(f(x))$ denotes the residue field of $f(x)$. But since $\Spec k(f(x)) \ra \Spec A$ factors through $\Spf A \ra \Spec A$, we conclude that $X_{f(x)} \ra \Spec k(f(x))$ is indeed smooth.
\end{proof}

\begin{lemma}\label{lemma:smooth_closed_points_Jacobson}
Let $f \colon X \ra Y$ be a locally of finite type (resp. locally of finite presentation) morphism of schemes, and assume that $Y$ is a Jacobson scheme. Then $f$ is smooth (resp. flat) if and only if $f$ is smooth (resp. flat) at every point of $X$ which lies over a closed point of $Y$.
\end{lemma}
\begin{proof}
Since $f$ is locally of finite type (resp. locally of finite presentation), we know that $X$ is a Jacobson scheme (i.e. closed points are dense in every closed subset). Since $f$ is smooth (resp. flat) on an open subset of $X$, it is enough to check that $f$ is smooth (resp. flat) at every closed point of $X$. As $f$ is locally of finite type and $Y$ is Jacobson, we know that $f$ maps closed points to closed points \cite[\href{https://stacks.math.columbia.edu/tag/01TB}{Lemma 01TB}]{stacks-project} which gives the lemma claim.
\end{proof}

\begin{lemma}\label{lemma:smooth_Jacobson_formal_completions}
Let $\mc{X}$ be an algebraic stack, let $Y$ be a Jacobson locally Noetherian scheme, and let $f \colon \mc{X} \ra Y$ be a morphism which is locally of finite type. For points $y \in Y$, write $\widehat{\mc{O}}_{Y,y}$ for the completion of the local ring at $y$. Then $\mc{X} \ra Y$ is smooth (resp. flat) if and only if $\mc{X}_{\Spf \widehat{\mc{O}}_{Y,y}} \ra \Spf \widehat{\mc{O}}_{Y,y}$ is smooth (resp. flat) for every closed point $y \in Y$.
\end{lemma}
\begin{proof}
Select any scheme $U$ with a surjective smooth morphism $U \ra \mc{X}$. Then $U \ra Y$ is a locally of finite type morphism of Jacobson locally Noetherian schemes, and $\mc{X} \ra Y$ is smooth (resp. flat) if and only if $U \ra Y$ is smooth (resp. flat). By Lemma \ref{lemma:smooth_closed_points_Jacobson}, we may check smoothness (resp. flatness) of $U \ra Y$ at points of $U$ lying over closed points of $Y$. If $x \in U$ and $y = f(x)$, then $U \ra Y$ is smooth (resp. flat) at $x$ if and only if $U_{\Spec \widehat{\mc{O}}_{Y,y}} \ra \Spec \widehat{\mc{O}}_{Y,y}$ is smooth at $x$ (first checking flatness, then checking smoothness in the fiber over the closed point). 
For any $y \in Y$, Lemma \ref{lemma:flatness_smoothness_Spf_schemes} implies that $U \ra Y$ is smooth (resp. flat) at all points $x \in U$ lying over $y$ if and only if $U_{\Spf \widehat{\mc{O}}_{Y,y}} \ra \Spf \widehat{\mc{O}}_{Y,y}$ is smooth (resp. flat). By Lemma \ref{lemma:smooth_closed_points_Jacobson}, we then see that $U \ra Y$ is smooth (resp. flat) if and only if $U_{\Spf \widehat{\mc{O}}_{Y,y}} \ra \Spf \widehat{\mc{O}}_{Y,y}$ is smooth (resp. flat) for every closed point $y \in Y$. This is equivalent to the condition that $\mc{X}_{\Spf \widehat{\mc{O}}_{Y,y}} \ra \Spf \widehat{\mc{O}}_{Y,y}$ is smooth (resp. flat) for all closed points $y \in Y$, since $U \ra \mc{X}$ is a smooth surjection.
\end{proof}

\begin{lemma}\label{lemma:special_cycles_generically_smooth}
Let $L$ be any non-degenerate Hermitian $\mc{O}_F$-lattice, with associated moduli stack $\mc{M}$. Fix $T \in \mrm{Herm}_m(F)$.

Then there exists $N \in \Z$ such that $\mc{Z}(T)[1/(N d_L \Delta)]$ is either empty\footnote{Following the Stacks project \cite[\href{https://stacks.math.columbia.edu/tag/0055}{Definition 0055}]{stacks-project}, our convention is that $\dim \emptyset = - \infty$.} or smooth of relative dimension $(n - r - \rank(T))r$ over $\Spec \mc{O}_F[1/(N d_L \Delta)]$.

We may take $N$ such that for $p \nmid N d_L \Delta$, there exists $g \in \GL_m(\mc{O}_{F_p})$ with
    \begin{equation}\label{equation:special_cycles_generically_smooth}
    {}^t \overline{g} T g = \begin{pmatrix} \mrm{Id}_{\rank(T)} & 0 \\ 0 & 0 \end{pmatrix}
    \end{equation}
where ${}^t \overline{g}$ denotes the conjugate transpose of $g$.
\end{lemma}
\begin{proof}
Fix a prime $p \nmid d_L \Delta$ such that there exists $g \in \GL_m(\mc{O}_{F_p})$ satisfying \eqref{equation:special_cycles_generically_smooth}. By Lemma \ref{lemma:smooth_Jacobson_formal_completions}, it is enough to check that the base change $\mc{Z}(T)_{\Spf \mc{O}_{F_p}}$ is either empty or smooth of relative dimension $(n - r - \rank(T))r$ over $\Spf \mc{O}_{F_p}$.

The morphism $\mc{Z}(T)_{\Spf \mc{O}_{F_p}} \ra \Spf \mc{O}_{F_p}$ is representable by algebraic stacks and locally of finite presentation. Thus $\mc{Z}(T)_{\Spf \mc{O}_{F_p}} \ra \Spf \mc{O}_{F_p}$ is smooth if and only if it is formally smooth \cite[\href{https://stacks.math.columbia.edu/tag/0DP0}{Lemma 0DP0}]{stacks-project}. 

Let $S \ra S'$ be a first order thickening of affine schemes, and assume $S'$ is equipped with a morphism to $\Spf \mc{O}_{F_p}$. To check formal smoothness of $\mc{Z}(T)_{\Spf \mc{O}_{F_p}} \ra \Spf \mc{O}_{F_p}$, we need to show that every object $(A_0, \iota_0, \lambda_0, A, \iota, \lambda, \underline{x}) \in \mc{Z}(T)(S)$ admits a lift to $S'$.

Form $(X_0, \iota_0, \lambda_0)$, where $X_0 = A_0[p^{\infty}]$ is the $p$-divisible group of $A_0$ with induced action $\iota_0 \colon \mc{O}_{F_p} \ra \End(X_0)$ and principal polarization $\lambda_0 \colon X_0 \ra X^{\vee}_0$. Similarly associate $(X, \iota, \lambda)$ to $(A, \iota, \lambda)$, where $X = A[p^{\infty}]$ is the $p$-divisible group of $A$. Note that the polarization $\lambda \colon X \ra X^{\vee}$ is principal because $p \nmid d_L \Delta$. Write also $\underline{x} = [x_1,\ldots,x_m]$ for the corresponding $m$-tuple of morphisms $x_i \colon X_0 \ra X$. By Serre--Tate, lifting $(A_0, \iota_0, \lambda_0, A, \iota, \lambda, \underline{x})$ from $S$ to $S'$ is the same as lifting $(X_0, \iota_0, \lambda_0, X, \iota, \lambda, \underline{x})$ from $S$ to $S'$. 

Using an element $g \in \GL_m(\mc{O}_{F_p})$ satisfying \eqref{equation:special_cycles_generically_smooth} as a ``change of basis'' for $X_0^m$, we obtain an $\mc{O}_F$-linear ``orthogonal splitting'' $X \cong X_0^{\rank(T)} \times Y$. That is, $Y$ is a $p$-divisible group with action $\iota_Y \colon \mc{O}_{F_p} \ra \End(Y)$, and a principal polarization $\lambda_Y \colon Y \ra Y^{\vee}$ whose Rosati involution $\dagger$ satisfies $\iota_Y(a)^{\dagger} = \iota_Y(a^{\s})$ for all $a \in \mc{O}_F$. Under the described identification $X \cong X_0^{\rank(T)} \times Y$, the polarization $\lambda$ on $X$ is given by $(\lambda_0)^{\rank(T)} \times \lambda_Y$. The map $\underline{x} \colon X_0^m \ra X$ may be identified with the projection $X_0^m \ra X_0^{\rank(T)}$ onto the first $\rank(T)$ factors, followed by the canonical inclusion $X_0^{\rank(T)} \ra X_0^{\rank(T)} \times Y$.

Note that the actions of $\mc{O}_{F_p}$ on $X_0$, $X$, and $Y$ have signatures $(1,0)$, $(n-r,r)$, and $(n - r - \rank(T),r)$ respectively (in the sense of \eqref{equation:p_div_generic_smoothness}). These considerations also show that $\rank(T) \leq n - r$ if $\mc{Z}(T)_{\Spf \mc{O}_{F_p}}$ is nonempty.

The triple $(X_0, \iota_0, \lambda_0)$ admits a unique lift to $S'$ as in Lemma \ref{lemma:formally_smooth_deformations_p_div}. The projection map $\underline{x} \colon X_0^m \ra X_0^{\rank(T)}$ clearly lifts to $S'$ as well. So it remains only to lift $(Y, \iota_Y, \lambda_Y)$ from $S$ to $S'$. Such a lift exists by formal smoothness of the corresponding deformation problem described in Lemma \ref{lemma:formally_smooth_deformations_p_div}. We apply the same lemma to compute tangent spaces (e.g. after passing to an \'etale cover by a scheme), which shows that the relative dimension is $(n -r - \rank(T))r$.
\end{proof}

\begin{remark}\label{remark:smoothness_moduli_stack}
Taking $T = \emptyset$ (or $T = 0$) in Lemma \ref{lemma:special_cycles_generically_smooth} and varying over non-degenerate Hermitian $\mc{O}_F$-lattices $L$ satisfying $L \subseteq L^{\vee}$ and $|L^{\vee} / L| = d$, we see that $\ms{M}(n-r,r)^{(d)} \ra \Spec \mc{O}_F[1/(d \Delta)]$ is smooth of relative dimension $(n - r)r$ for every $d \in \Z_{\geq 0}$. If $\mc{M}$ is associated with any non-degenerate Hermitian $\mc{O}_F$-lattice $L$ (not necessarily with $L \subseteq L^{\vee}$), this then implies that $\mc{M} \ra \Spec \mc{O}_F[1/(d_L \Delta)]$ is smooth of relative dimension $(n-r)r$.
\end{remark}

%


    \section{Arithmetic special cycle classes} 
    \label{sec:arith_cycle_classes}

        \subsection{Hermitian vector bundles}
        \label{ssec:arith_cycle_classes:Hermitian_bundles}
            Given a smooth algebraic stack $\mc{X}$ over $\Spec \C$, a \emph{Hermitian vector bundle} $\widehat{\mc{E}}$ on $\mc{X}$ is the following functorial assignment: for every morphism $S \ra \mc{X}$ with $S$ a scheme smooth over $\Spec \C$, the assignment gives a vector bundle on $S$ with a smooth Hermitian metric on the analytification. We sometimes write $\widehat{\mc{E}} = (\mc{E}, \norm{-})$ where $\mc{E}$ is the underlying vector bundle on $\mc{X}$ and $\norm{-}$ is the norm associated with the smooth Hermitian metric (on the analytification), defined functorially.

Let $(R,\Sigma,c_{\infty})$ be an arithmetic ring in the sense of Gillet--Soul\'e \cite[{\S 3.1}]{GS90}, i.e. $R$ is an excellent regular Noetherian integral domain (e.g. Dedekind domains with fraction field of characteristic $0$ or fields), $\Sigma$ is a finite nonempty set of injective homomorphisms $\t \colon R \ra \C$, and $c_{\infty} \colon \C^{\Sigma} \ra \C^{\Sigma}$ is a conjugate-linear involution of $\C \otimes_{\Z} R$-algebras. Write $K$ for the fraction field of $R$. 

Suppose $\mc{X}$ is an algebraic stack which is flat and finite type over $\Spec R$. Write $\mc{X}_{\t} \coloneqq \mc{X} \times_{\Spec K, \t} \Spec \C$. Assume that the generic fiber $\mc{X}_K$ is smooth over $\Spec K$. A \emph{Hermitian vector bundle} on $\mc{X}$ is a vector bundle $\mc{E}$ on $\mc{X}$ equipped with a smooth Hermitian metric on $\mc{E}|_{\mc{X}_{\t}}$ for each $\t \in \Sigma$, such that this collection of metrics is conjugation invariant (meaning $c_{\infty}$-invariant). We write $\widehat{\Pic}(\mc{X})$ for the group of (isomorphism classes of) Hermitian line bundles, with group structure given by the tensor product. We use the notation $\widehat{\mc{L}} = (\mc{L}, \norm{-})$ for a Hermitian line bundle with underlying line bundle $\mc{L}$ and $\norm{-}$ its norm (meaning a collection of norms $\{\norm{-}_{\t}\}_{\t \in \Sigma}$). We also write $\norm{-}_{\infty} = \prod_{\t \in \Sigma} \norm{-}_{\t}$.

Next, we discuss stacky degrees of Hermitian line bundles. We fix the arithmetic ring $(R, \Sigma, c_{\infty})$ associated with $R = \mc{O}_K[1/N]$ for a number field $K$ and an integer $N \in \Z$, i.e. $\Sigma$ is the set of all embeddings $\t \colon K \ra \C$, and $c_{\infty}$ is induced by complex conjugation. For the rest of Section \ref{ssec:arith_cycle_classes:Hermitian_bundles}, we assume that 
    \begin{equation}\label{equation:stacky_degrees_stack_hypotheses}
    \parbox[t]{.8 \textwidth}{$\mc{X}$ is a reduced $1$-dimensional Deligne--Mumford stack which is proper and flat over $\Spec R$.} \end{equation}
Here, dimension is used in an absolute sense rather than a relative one (e.g. we could have $\mc{X} = \Spec R$).

Let $\widehat{\mc{L}} = (\mc{L}, \norm{-})$ be a Hermitian line bundle on $\mc{X}$. For each complex place $v$ of $\mc{O}_K$, we set $\norm{-}_v \coloneqq \norm{-}_{\t_1} \norm{-}_{\t_2}$ where $\t_1, \t_2 \colon K \ra \C$ are the two embeddings corresponding to $v$.
The Arakelov \emph{arithmetic degree} $\widehat{\deg}(\widehat{\mc{L}})$ of $\widehat{\mc{L}}$ on $\mc{X}$ is valued in $\R_N \coloneqq \R/(\sum_{p \mid N} \Q \cdot \log p)$, and may be described as follows. If $\mc{X} = \Spec \mc{O}_E[1/N]$ for a number field $E$, we have the standard definition
    \begin{equation}
    \widehat{\deg}(\widehat{\mc{L}}) \coloneqq - \sum_{v \nmid N} \log \norm{s}_v \quad \text{with} \quad \norm{s}_v \coloneqq q_v^{-\mrm{ord}_v(s)} \text{ if $v < \infty$}
    \end{equation}
where the sum runs over all places $v$ of $E$ (Archimedean included) with $v \nmid N$, the quantity $q_v$ is the cardinality of the residue field at $v$, and $s$ is any rational section of $\mc{L}$.

If $\mc{X}$ is integral (equivalently, reduced and irreducible by quasi-separatedness), select any number field $E$ with a finite surjection $\Spec \mc{O}_E[1/N] \ra \mc{X}$ and set
    \begin{equation} \label{equation:stacky_degree_Hermitian_bundle}
    \widehat{\deg}(\widehat{\mc{L}}) \coloneqq \frac{1}{\deg(E/\mc{X}_K)} \widehat{\deg}(\widehat{\mc{L}}|_{\Spec \mc{O}_{E}[1/N]}) 
    \end{equation}
where $\deg(E/\mc{X}_K)$ denotes the degree of the finite \'etale morphism $\Spec E \ra \mc{X}_K$. This definition does not depend on the choice of $E$ or the morphism $\Spec \mc{O}_E[1/N] \ra \mc{X}$, since any two such choices may be covered by finite surjections from a third such choice $\Spec \mc{O}_{E'}[1/N]$ (and these finite surjections can be made compatible with the maps to $\mc{X}$).

\begin{remark}\label{remark:existence_finite_covers_stack}
Existence of $E$ as in the preceding paragraph follows from some general theory. Indeed, a general fact about Noetherian Deligne--Mumford stacks with separated diagonal \cite[{Th\'eor\`eme 16.6}]{LMB00} implies that there exists a scheme $Z$ with a morphism $Z \ra \mc{X}$ which is finite, surjective, and generically \'etale (in the sense that $Z_{\mc{U}} \ra \mc{U}$ is \'etale for a dense open substack $\mc{U} \subseteq \mc{X}$). Selecting an irreducible component of $Z$ which surjects onto $\mc{X}$, we may assume that $Z$ is also integral. Thus $Z$ is a $1$-dimensional integral scheme which is proper and flat over $\Spec R$. Such a map $Z \ra \Spec R$ must be quasi-finite, hence also finite. If $E$ denotes the fraction field of $Z$, the normalization of $Z$ must be $\tilde{Z} = \Spec \mc{O}_{E}[1/N]$.
\end{remark}

One can check that the definition of stacky arithmetic degree in \eqref{equation:stacky_degree_Hermitian_bundle} recovers the definition of \cite[{(4.4)}]{KRY04} and \cite[{\S 2.1}]{KRY06} which counts geometric points weighted by orders of automorphism groups.\footnote{In loc. cit., the functorial definition of Hermitian line bundle should also include the additional assumption of complex conjugation invariance as above.}

In the general case when $\mc{X}$ is not necessarily irreducible, we take
    \begin{align*}
    \widehat{\deg}(\widehat{\mc{L}}) \coloneqq \sum_{\xi} \widehat{\deg}(\widehat{\mc{L}}|_{\mc{X}_{\xi}})
    \end{align*}
where the sum runs over generic points $\xi$ of irreducible components $\mc{X}_{\xi}$ of $\mc{X}$.

The preceding discussion showed that $\mc{X}$ admits a finite surjection from a scheme which is finite over $\Spec R$, hence $\mc{X} \ra \Spec R$ is also quasi-finite (in the sense of \cite[\href{https://stacks.math.columbia.edu/tag/0G2M}{Definition 0G2M}]{stacks-project}).

Suppose $\mc{X}'$ and $\mc{X}$ are Deligne--Mumford stacks which both satisfy the hypotheses from \eqref{equation:stacky_degrees_stack_hypotheses}, and consider a morphism $f \colon \mc{X}' \ra \mc{X}$ over $\Spec R$. Let $\widehat{\mc{L}}$ be a Hermitian line bundle on $\mc{X}$. First consider the case where $\mc{X}$ is irreducible. We say that the morphism $\mc{X}'_K \ra \mc{X}_K$ has \emph{degree} $\deg(\mc{X}'_K/\mc{X}_K) \coloneqq \deg(\mc{X}'_K/K) / \deg(\mc{X}_K/K)$ (with stacky degrees of $0$-cycles over fields as in \eqref{equation:zero_cycle_degree_over_field}). We have $\widehat{\deg}(f^* \widehat{\mc{L}}) = \deg(\mc{X}'_K / \mc{X}_K) \widehat{\deg}(\widehat{\mc{L}})$.
Next, consider the general case where $\mc{X}$ is not necessarily irreducible. We say that $\mc{X}'_K \ra \mc{X}_K$ has \emph{constant degree $d$} if $\mc{X}'_K \times_{\mc{X}_K} \mc{X}_{\xi,K} \ra \mc{X}_{\xi,K}$ has degree $d$ for every irreducible component $\mc{X}_{\xi}$ of $\mc{X}$. In this case, we have
    \begin{equation} \label{equation:stacky_degrees_pullback_equation}
    \widehat{\deg}(f^* \widehat{\mc{L}}) = d \cdot  \widehat{\deg}(\widehat{\mc{L}}).
    \end{equation}

        \subsection{Arithmetic Chow rings}
        \label{ssec:arith_cycle_classes:arithmetic_Chow}
            We fix definitions for arithmetic Chow rings with rational coefficients.

Let $(R,\Sigma, c_{\infty})$ be an arithmetic ring (Section \ref{ssec:arith_cycle_classes:Hermitian_bundles}) and write $K$ for the fraction field of $R$.
Suppose $X$ is a scheme which is separated, flat, and finite type over $\Spec R$ with smooth and quasi-projective generic fiber $X_K$. There are associated Gillet--Soul\'e \emph{arithmetic Chow groups} $\arithCh^m(X)$ in codimensions $m \geq 0$. If $X$ is moreover regular, these groups form an \emph{arithmetic Chow ring} $\arithCh^*(X)_{\Q}$ (with $\Q$-coefficients) \cite[{Theorem 4.2.3}]{GS90}.

Let $L$ be any non-degenerate Hermitian $\mc{O}_F$-lattice of rank $n$, with associated moduli stack $\mc{M}$. Consider the arithmetic ring $(R, \Sigma, c_{\infty})$ associated with $R = \Spec \mc{O}_F[1/d_L]$. We define \emph{arithmetic Chow groups} for $\mc{M}$ by limiting over level structure: for any nonzero integer $N$, set
    \begin{equation}
    \arithCh^*(\mc{M}[1/N])_{\Q} \coloneqq \lim_{K'_f} \arithCh^*(\mc{M}_{K'_f}[1/N])_{\Q}
    \end{equation}
where $K'_f$ varies over all small levels as in Section \ref{ssec:ab_var:level_structure} (so that each $\mc{M}_{K'_f}$ is a scheme). Similar limiting procedures appeared in \cite{BBK07} and \cite[{\S 4.4}]{BH21}; see also \cite{Gillet09} for more on arithmetic Chow rings of Deligne--Mumford stacks.

Since $\mc{M} \ra \Spec \mc{O}_F[1/d_L]$ is smooth, we know that $\mc{M}$ is regular. Hence we obtain an \emph{arithmetic Chow ring} $\arithCh^*(\mc{M}[1/N])_{\Q}$ via the intersection product for each $\arithCh^*(\mc{M}_{K'_f}[1/N])_{\Q}$.

Suppose $\mc{Z} \ra \mc{M}$ is a finite morphism of algebraic stacks with $\mc{Z} \ra \Spec \mc{O}_F[1/d_L]$ proper and $\mc{Z}$ equidimensional of dimension $d$. Then we define the \emph{height} of $\mc{Z}$ with respect to any Hermitian line bundle $\widehat{\mc{L}}$ on $\mc{M}$ as follows: if $\mc{Z}_{K'_f} \coloneqq \mc{Z} \times_{\mc{M}} \mc{M}_{K'_f}$, the quantity
    \begin{equation}
    \widehat{\deg}(\widehat{\mc{L}}^d|_{\mc{Z}}) \coloneqq \frac{1}{[K'_f(1) : K'_f]} \widehat{\deg}(\widehat{\mc{L}}^d|_{\mc{Z}_{K'_f}}) \in \R_{d_L N_{K'_f}} = \R/(\sum\nolimits_{p \mid d_L N_{K'_f}} \Q \cdot \log p)
    \end{equation}
does not depend on the choice of small level $K'_f$, where $\widehat{\deg}(\widehat{\mc{L}}^d|_{\mc{Z}_{K'_f}})$ is the arithmetic height as in \cite[{Proposition 2.3.1, Remarks(ii)}]{BGS94} (see also \cite{Zhang95}) calculated by replacing $\mc{Z}_{K'_f}$ with its pushforward cycle on $\mc{M}_{K'_f}$. Varying $K'_f$, we obtain the \emph{height} $\widehat{\deg}(\widehat{\mc{L}}^d|_{\mc{Z}}) \in \R_{d_L}$. We will be primarily interested in the case where $d = 1$ with $\mc{Z}$ reduced and flat over $\Spec \mc{O}_F[1/d_L]$. In this case, $\widehat{\deg}(\widehat{\mc{L}}|_{\mc{Z}})$ is the (stacky) arithmetic degree of $\widehat{\mc{L}}$ restricted to $\mc{Z}$, as discussed in Section \ref{ssec:arith_cycle_classes:Hermitian_bundles}.

        \subsection{Hodge bundles}
        \label{ssec:arith_cycle_classes:Hodge_bundles}
            Given an abelian scheme $\pi \colon A \ra S$ of relative dimension $n$, we consider the \emph{Hodge bundles} $\Omega_A \coloneqq \pi_* \Omega^1_{A/S}$ and $\omega_A \coloneqq \pi_* \bigwedge^n \Omega_A$. Here $\Omega_A$ and $\omega_A$ are locally free of ranks $n$ and $1$ respectively. If $e \colon S \ra A$ denotes the identity section, there are canonical isomorphisms $\Omega_A \cong e^* \Omega_{A/S} \cong (\Lie A)^{\vee}$. These Hodge bundles are (contravariantly) functorial in $A$, and their formation commutes with arbitrary base change \cite[{Proposition 2.5.2}]{BBM82}.

When $S$ is a scheme which is smooth over $\Spec \C$, the analytification of $\omega_A$ may be equipped with a Hermitian metric (\emph{Faltings} or \emph{Hodge} metric), normalized as follows. The fiber of $\omega_A$ over any $s \in S(\C)$ is canonically identified with $H^0(A_s, \omega_{A_s})$. Viewing $H^0(A_s, \omega_{A_s})$ as the $1$-dimensional $\C$-vector space of holomorphic $n$-forms on $A_s(\C)$, we take the metric on $\omega_A$ at $s$ to be
    \begin{equation}\label{equation:arith_cycle_classes:Hodge_bundles:Faltings_metric}
    (\a, \b) = \left ( \frac{i}{2 \pi} \right ) ^n \int_{A_s(\C)} \b \wedge \bar{\a}
    \end{equation}
for $\a, \b \in H^0(A_s, \omega_{A_s})$. We call the resulting Hermitian line bundle $\widehat{\omega}_A \coloneqq (\omega_A, \norm{-})$ a \emph{metrized Hodge bundle}. This metric on $\omega_A$ is functorial for isomorphisms between abelian schemes $A$ over $S$.

Next, suppose $S = \Spec \C$ and suppose $\lambda \colon A \ra A^{\vee}$ is a quasi-polarization. There is an associated Hermitian metric on $\Omega_A^{\vee} \cong \Lie(A)$ which we normalize as follows: 
if $\lambda$ is a polarization and $\lambda(a) = (t_a^* \mc{L}) \otimes \mc{L}^{-1}$ for an ample line bundle $\mc{L}$ on $A$ (where $t_a$ is translation by $a$), then the Chern class $c_1(\mc{L}) \in H^2(A,\Z)$ defines a $\Z$-valued alternating form $\psi$ on $H_1(A,\Z)$ (sign conventions etc. as in \cite[{\S I}]{Mumford85}, except our Hermitian pairings are conjugate linear in the first variable) and a positive definite Hermitian pairing
\begin{equation}\label{equation:arith_cycle_classes:Hodge_bundles:Hermitian_normalization}
    (x,y) \coloneqq \pi \sqrt{|\Delta|} \big ( \psi(ix, y) - i \psi(x,y) \big )
    \end{equation}
on $\Lie A$. 
If instead $\lambda$ is a quasi-polarization, the associated Hermitian pairing is $m^{-1}$ times the Hermitian pairing of $m \lambda$ for any $m \in \Z_{>0}$ such that $m \lambda$ is a polarization. If $\dim A = 1$ and if $\lambda$ is the unique principal polarization of $A$, the induced dual metric on $\Omega_A \cong (\Lie A)^{\vee}$ is $\sqrt{|\Delta|}^{-1}$ times the Faltings metric (cf. the proof of \cite[{Lemma 5.1.4}]{BHKRY20II}).

Over an arbitrary smooth scheme $S \ra \Spec \C$, any quasi-polarization $\lambda \colon A \ra A^{\vee}$ defines a smooth Hermitian metric on $\Lie A$, given fiberwise by the construction above. We also call the resulting Hermitian vector bundle $\widehat{\Omega}_A^{\vee}$ a \emph{metrized Hodge bundle}.

Next, let $(R, \Sigma, c_{\infty})$ be an arithmetic ring, and suppose $\mc{X}$ is an algebraic stack which is flat and finite type over $\Spec R$, and whose generic fiber is smooth. Suppose we are given a relative abelian scheme $\mc{A}$ over $\mc{X}$ (equivalently, a functorial assignment of abelian schemes $A \ra S$ to objects $x \in \mc{X}(S)$, for $R$-schemes $S$). Formation of the metrized Hodge bundle $\widehat{\omega}_{A}$ is functorial, hence defines a \emph{metrized Hodge bundle} $\widehat{\omega}$ on $\mc{X}$. If $\mc{A}$ is equipped with a quasi-polarization, then there is similarly a \emph{metrized Hodge bundle} $\widehat{\Omega}^{\vee}$ on $\mc{X}$.

Let $L$ be a non-degenerate Hermitian $\mc{O}_F$-lattice of rank $n$, with associated moduli stack $\mc{M}$. If $L$ is self-dual and signature $(n - 1, 1)$, we write $\ms{E}$ for the pullback of the tautological bundle (Definition \ref{definition:ab_var:exotic_smooth:tautological_bundle}) along $\mc{M} \ra \ms{M}(n - 1, 1)^{\circ}$. Otherwise, we write $\ms{E}$ for the pullback of the tautological bundle (Definition \ref{definition:ab_var:integral_models:tautological_bundle}) along $\mc{M} \ra \ms{M}(n - r, r)[1/\Delta]$.
We similarly write $\Omega$ for the pullback to $\mc{M}$ of the Hodge bundle from $\ms{M}(n - r, r)[1/\Delta]$ in the non exotic smooth case (resp.  $\ms{M}(n - 1, 1)^{\circ}$ in the exotic smooth case) to $\mc{M}$ (where $\Omega$ is the rank $n$ Hodge bundle).
After base change to $\Spec \mc{O}_{F}[1/\Delta]$, the (dual) Hodge bundle $\Omega^{\vee}$ on $\mc{M}$ admits a (canonical) decomposition $\Omega^{\vee}[1/\Delta] = \Omega^{\vee}[1/\Delta]^+ \oplus \Omega^{\vee}[1/\Delta]^-$ (with $\Omega^{\vee}[1/\Delta]^- = \ms{E}^{\vee}[1/\Delta]$) where the $\mc{O}_F$-action (via $\iota$) on $\Omega^{\vee}[1/\Delta]^+$ (resp. $\Omega^{\vee}[1/\Delta]^-$) is $\mc{O}_F$-linear (resp. $\s$-linear). 

Equip $\Omega^{\vee}[1/\Delta]$ with the Hermitian metric which is 
    \begin{equation}\label{equation:arith_cycle_classes:Hodge_bundles:4_pi_gamma_normalization}
    (x,y) \coloneqq 4 \pi^2 e^{\gamma} \sqrt{|\Delta|} \big ( \psi(ix, y) - i \psi(x,y) \big )
    \end{equation}
in complex fibers (i.e. multiply the Hermitian metric from \eqref{equation:arith_cycle_classes:Hodge_bundles:Hermitian_normalization} by $4 \pi e^{\gamma}$) where $\gamma$ is the Euler--Mascheroni constant. We remark that the normalization constant $4 \pi e^{\gamma}$ has appeared previously in similar contexts, e.g. \cite[{(0.4)}]{KRY04} \cite[{\S 7}]{KRY06} \cite[{\S 7.2}]{BHKRY20}. We refer to loc. cit. for possible conceptual explanations of this constant.

Then $\ms{E}^{\vee}[1/\Delta] \subseteq \Omega^{\vee}[1/\Delta]$ (via decomposition in previous paragraph) inherits a Hermitian metric by restriction. This makes $\ms{E}^{\vee}$ into a Hermitian vector bundle $\widehat{\ms{E}}^{\vee}$.

Write $\Omega_0$ for the Hodge bundle on $\ms{M}_0$. Equip the dual $\Omega_0^{\vee}$ with the metric described fiberwise by \eqref{equation:arith_cycle_classes:Hodge_bundles:Hermitian_normalization}, giving a Hermitian line bundle $\widehat{\Omega}_0^{\vee}$. By the \emph{metrized dual tautological bundle} on $\mc{M}$, we mean
    \begin{equation}
    \widehat{\mc{E}}^{\vee} \coloneqq \widehat{\Omega}_0^{\vee} \otimes \widehat{\ms{E}}^{\vee}
    \end{equation}
where we have suppressed pullbacks from notation. The metric on $\widehat{\mc{E}}^{\vee}$ is the tensor product of the metrics described above.
This definition of $\widehat{\mc{E}}^{\vee}$ is similar to \cite[{\S 2.4, \S 7.2}]{BHKRY20}, though in a different setup (we are considering not-necessarily principal polarizations). Taking a dual gives the \emph{metrized tautological bundle} $\widehat{\mc{E}}$. For more discussion on these metric normalization choices, see Section \ref{ssec:Arch_uniformization:framing}. We write $\mc{E}^{\vee}$ for underlying line bundle of $\widehat{\mc{E}}^{\vee}$.

For readers interested in Faltings height, we also consider the metrized Hodge (determinant) bundle $\widehat{\omega}$ on $\mc{M}$, which is pulled back from the Hodge determinant bundle $\omega$ on $\ms{M}(n - r, r)$ and with metric normalized as in \eqref{equation:arith_cycle_classes:Hodge_bundles:Faltings_metric}.

Suppose $A_0 \ra \Spec \mc{O}_E$ is any (relative) elliptic curve with $\mc{O}_F$-action, where $E$ is a number field. If $\widehat{\omega}_{A_0}$ denotes the associated metrized Hodge bundle (normalized as in \eqref{equation:arith_cycle_classes:Hodge_bundles:Faltings_metric}), we recall that the \emph{Faltings height} of $A_0$ is 
    \begin{equation}\label{equation:arith_cycle_classes:Hodge_bundles:CM_Faltings_height}
    h_{\mrm{Fal}}^{\mrm{CM}} \coloneqq \frac{1}{[E : \Q]} \widehat{\deg}(\widehat{\omega}_{A_0}) = \frac{1}{2} \frac{L'(1,\eta)}{L(1, \eta)} + \frac{1}{2} \frac{\Gamma'(1)}{\Gamma(1)} + \frac{1}{4} \log |\Delta| - \frac{1}{2} \log(2 \pi)
    \end{equation}
where $\eta$ is the quadratic character associated to $F / \Q$, and $\Gamma$ is the usual gamma function. This comes from the classical Chowla--Selberg formula (the statement above is as in \cite[{Proposition 10.10}]{KRY04}). It will be convenient to define the height constants
    \begin{equation}\label{equation:arith_cycle_classes:Hodge_bundles:taut_height_constants}
    h_{\mrm{tau}}^{\mrm{CM}} \coloneqq - h_{\mrm{Fal}}^{\mrm{CM}} + \frac{1}{4} \log|\Delta| - \frac{1}{2} \log(4 \pi e^{\gamma}) \quad \quad h_{\widehat{\mc{E}}^{\vee}}^{\mrm{CM}} \coloneqq - h^{\mrm{CM}}_{\mrm{Fal}} - \frac{1}{4} \log|\Delta| + h_{\mrm{tau}}^{\mrm{CM}}.
    \end{equation}
These will re-appear in \cref{ssec:non-Arch_uniformization:horizontal}.

        \subsection{Arithmetic special cycle classes}
        \label{ssec:arith_cycle_classes:arithmetic_cycle_classes}
            Let $L$ be a non-degenerate Hermitian $\mc{O}_F$-lattice of rank $n$, with associated moduli stack $\mc{M}$. For the rest of Section \ref{sec:arith_cycle_classes}, we assume $L$ has signature $(n - 1, 1)$. Consider an $m \times m$ Hermitian matrix $T \in \mrm{Herm}_m(\Q)$, assume $m \leq n$, and form the associated special cycle $\mc{Z}(T) \ra \mc{M}$. One expects to be able to construct an associated \emph{arithmetic special cycle class} $[\widehat{\mc{Z}}(T)] \in \arithCh^{m}(\mc{M})_{\Q}$.

For arbitrary singular $T$, there is no proposed definition of $[\widehat{\mc{Z}}(T)]$ in the literature. In general, $\mc{Z}(T)_{\ms{H}}$ has larger-than-expected dimension. The stack $\mc{Z}(T)$ could also have components with larger-than-expected dimension in positive characteristic (occurs already for nonsingular $T$). Available methods in the literature for treating the non-Archimedean theory ($K$-theoretic and derived algebro-geometric) do not incorporate the Archimedean place in general, as needed for arithmetic intersection theory (see introduction).

The analogue of the ``linear invariance'' approach of \cite[{\S 6.4}]{KRY04} (there for Shimura curves) is to first define $[\widehat{\mc{Z}}(T^{\flat})]$ for nonsingular $T^{\flat}$, to consider ${}^t \overline{\gamma} T \gamma = \mrm{diag}(0, T^{\flat})$ for some $\gamma \in \GL_m(\mc{O}_F)$ with $T^{\flat}$ nonsingular, and to define $[\widehat{\mc{Z}}(T)]$ by intersecting $[\widehat{\mc{Z}}(T^{\flat})]$ with a power of some metrized tautological bundle (possibly with additional Archimedean adjustment). This is not literally possible in the unitary setting, where $\mc{O}_F$ may have class number $\neq 1$ (in particular, $\gamma$ as above may not exist). One also needs to verify independence of the choice of $\gamma$.

For arbitrary $T \in \mrm{Herm}_m(\Q)$, we propose to define $[\widehat{\mc{Z}}(T)]$ as a sum
    \begin{equation}
    [\widehat{\mc{Z}}(T)] \coloneqq [\widehat{\mc{Z}}(T)_{\ms{H}}] + \sum_{\substack{p \text{ prime} \\ p \nmid d_L}} [{}^{\mathbb{L}}\mc{Z}(T)_{\ms{V},p}] \in \arithCh^m(\mc{M})_{\Q}.
    \end{equation}
We construct $[\widehat{\mc{Z}}(T)_{\ms{H}}]$ using the horizontal part $\mc{Z}(T)_{\ms{H}}$ and an appropriate Green current $g_{T,y}$ \eqref{equation:arith_cycle_classes:horizontal} with an additional parameter $y \in \mrm{Herm}_m(\R)_{>0}$. The element $[{}^{\mbb{L}} \mc{Z}(T)_{\ms{V},p}]$ arises from a class ${}^{\mbb{L}} \mc{Z}(T)_{\ms{V},p} \in \mrm{gr}^m_{\mc{M}} K_0'(\mc{Z}(T)_{\F_p})_{\Q}$ corresponding to the ``vertical part'' of $\mc{Z}(T)$ at $p$ \eqref{equation:arith_cycle_classes:vertical}. The classes ${}^{\mbb{L}} \mc{Z}(T)_{\ms{V},p}$ will be zero for all but finitely many primes $p$. We define ${}^{\mbb{L}} \mc{Z}(T)_{\ms{V},p}$ using a ``$p$-local'' variant of the linear invariance strategy above.

We will show that $[\widehat{\mc{Z}}(T)]$ satisfies the ``linear invariance'' property
    \begin{equation}
    [\widehat{\mc{Z}}(T)] = [\widehat{\mc{Z}}({}^t \overline{\gamma} T \gamma)]
    \end{equation}
for all $\gamma \in \GL_m(\mc{O}_F)$, where $[\widehat{\mc{Z}}(T)]$ is formed with respect to $y \in \mrm{Herm}_m(\R)_{>0}$ and $[\widehat{\mc{Z}}(T)]$ is formed with respect to $\gamma^{-1} y {}^t \overline{\gamma}^{-1}$. 

In fact, we prove refined statements. We show
    \begin{equation}
    [\widehat{\mc{Z}}(T)_{\ms{H}}] = [\widehat{\mc{Z}}({}^t \overline{\gamma} T \gamma)_{\ms{H}}]
    \end{equation}
for the Green currents $g_{T,y}$ defined in Section \ref{ssec:Arch_uniformization:Archimedean} (where the current $g_{{}^t \overline{\gamma} T \gamma, \gamma^{-1} y {}^t \overline{\gamma}^{-1}}$ is used on the right-hand side). Moreover, we show $g_{T,y} = g_{{}^t \overline{\gamma} T \gamma, \gamma^{-1} y {}^t \overline{\gamma}^{-1}}$ (Section \ref{ssec:Arch_uniformization:Archimedean}); this property is also satisfied for the Garcia--Sankaran currents in \cite[{(4.38)}]{GS19} (which we do not use for our arithmetic Siegel--Weil results).

For any $\gamma \in \GL_m(\mc{O}_{F,(p)}) \cap M_{m,m}(\mc{O}_F)$, we show that the pullback
    \begin{equation}\label{equation:arith_cycle_classes:arithmetic_cycle_classes:vertical_invariance}
    \mrm{gr}^m_{\mc{M}} K'_0(\mc{Z}(T)_{\F_p})_{\Q} \leftarrow \mrm{gr}^m_{\mc{M}} K'_0(\mc{Z}({}^t \overline{\gamma} T \gamma)_{\F_p})_{\Q}
    \end{equation}
along $\mc{Z}(T)_{\F_p} \ra \mc{Z}({}^t \overline{\gamma} T \gamma)_{\F_p}$ (defined in \eqref{equation:ab_var:special_cycles:gamma_action}) sends ${}^{\mbb{L}} \mc{Z}({}^t \overline{\gamma} T \gamma)_{\ms{V},p}$ to ${}^{\mbb{L}} \mc{Z}(T)_{\ms{V},p}$ (see \eqref{equation:arith_cycle_classes:vertical:linear_invariance}).

        \subsection{Horizontal arithmetic special cycle classes}
        \label{ssec:arith_cycle_classes:horizontal}
            Consider any $m \times m$ Hermitian matrix $T \in \mrm{Herm}_m(\Q)$. The horizontal arithmetic special cycle class $[\widehat{\mc{Z}}(T)_{\ms{H}}]$ should involve some extra Archimedean data, e.g. from a Green current $g_{T,y}$ (which we allow to depend on a parameter $y \in \mrm{Herm}_m(\R)_{>0}$, as is typical in the literature).

Given an equidimensional complex manifold $X$, recall that a \emph{$(p,q)$-current} on $X$ is a continuous linear map $\Omega^{\dim X - p, \dim X - q}_c(X) \ra \C$ on compactly supported smooth forms of degree $(\dim X - p, \dim X - q)$, where $\Omega^{\dim X - p, \dim X - q}_c(X)$ has the usual colimit topology. A $(p,p)$-current is \emph{real} if it is induced by a continuous real-valued linear map on real $(p,p)$-forms. Given a top degree current $g$ on $X$ (i.e. a distribution), we say that $g$ is \emph{integrable} or that $\int_X g$ \emph{converges} (possibly non-standard usage) if $g$ extends (necessarily uniquely) to a continuous map $C^{\infty}_{b_1}(X) \ra \C$, where
    \begin{equation}
    C^{\infty}_{b_1}(X) \coloneqq \{f \in C^{\infty}(X) : |f(x)| \leq 1 \text{ for all $x \in X$}\}
    \end{equation} 
with topology given by $\sup$-norms ranging over all compact subsets $K \subseteq X$. In this case, we write $\int_X g$ for the value of $g$ on $1 \in C^{\infty}_{b_1}(X)$.
Suppose $\a$ is a (measurable) locally $L^1$ form of top degree on $X$. If $\a$ is integrable, then the associated distribution $[\a]$ on $X$ is integrable, and we have $\int_X [\a] = \int_X \a$.
We use the orientation and sign conventions of \cite{GS90}.

Returning to the moduli stack $\mc{M} \ra \Spec \mc{O}_F$ from above, choose any embedding $F \ra \C$ and form the base changes $\mc{M}_{\C} \coloneqq \mc{M} \times_{\Spec \mc{O}_F} \Spec \C$ and $\mc{Z}(T)_{\C} \coloneqq \mc{Z}(T) \times_{\Spec \mc{O}_F} \Spec \C$, etc.. By a \emph{$(p,q)$-current} on $\mc{M}_{\C}$, we mean a system of currents $g = (g_{K'_f})_{K'_f} = (\Omega^{n - 1 - p, n - 1 - q}_c(\mc{M}_{K'_f,\C}) \ra \C)_{K'_f}$ compatible with pullback of currents as we vary $K'_f$ among all small levels.
We say a $(p,q)$-current on $\mc{M}_{\C}$ is \emph{real} if the associated current at each small level $K'_f$ is real. If $g$ is a current of top degree on $\mc{M}$, its \emph{integral} is defined as
    \begin{equation}\label{equation:arith_cycle_classes:horizontal:integral}
    \int_{\mc{M}_{\C}} g \coloneq \frac{1}{[K'_{L,f} : K'_f]} \int_{\mc{M}_{K'_f, \C}} g_{K'_f}
    \end{equation}
for any sufficiently small level $K'_f$ (conditional on convergence). This definition does not depend on the choice of small level.

 Suppose $g_{T,y}$ is any real $(m - 1, m - 1)$ current on $\mc{M}_{\C}$ which satisfying a \emph{modified current equation}, i.e. such that
    \begin{equation}\label{equation:arith_cycle_classes:horizontal:Green_current}
    -\frac{1}{2 \pi i} \partial \overline{\partial} g_{T,y} + \delta_{\mc{Z}(T)_{\C}} \wedge [c_1(\widehat{\mc{E}}^{\vee}_{\C})^{m - \rank(T)}]
    \end{equation}
is (represented by) a smooth $(m, m)$-form (for all small levels $K'_f$), where $c_1(\widehat{\mc{E}}^{\vee}_{\C})$ is the Chern form of the Hermitian line bundle $\widehat{\mc{E}}^{\vee}_{\C}$. We call $g_{T,y}$ a \emph{Green current}. We typically write $g_{T,y}$ instead of $g_{T,y,K'_f}$ to lighten notation, for understood level $K'_f$.

For each small level $K'_f$, pick a representative $(\mc{Z}_0, g_0)$ for the self-intersection arithmetic cycle class $\widehat{c}_1(\widehat{\mc{E}}^{\vee})^{m - \rank(T)} \in \arithCh^{m - \rank(T)}(\mc{M}_{K'_f})_{\Q}$. We can assume that $\mc{Z}_0$ intersects $\mc{Z}(T)_{K'_f}$ properly in the generic fiber (moving lemma) and that $g_0$ is a Green form of logarithmic type for $\mc{Z}_0$ (in the sense of \cite[{\S 4}]{GS90}).

The intersection pairing for Chow groups with supports (as in \cite[{\S 4}]{GS90}) gives a class $\mc{Z}(T)_{K'_f, \ms{H}} \cdot \mc{Z}_0 \in \mrm{Ch}^m_{Z(T)_{K'_f, \ms{H}} \cap \mc{Z}_0}(\mc{M}_{K'_f})_{\Q}$. We set
    \begin{equation}
    [\widehat{\mc{Z}}(T)_{K'_f, \ms{H}}] \coloneqq [(\mc{Z}(T)_{K'_f, \ms{H}} \cdot \mc{Z}_0, g_{T,y} + g_0 \wedge \delta_{\mc{Z}(T)_{K'_f,\C}})]\in \arithCh^m(\mc{M}_{K'_f}[1/N])_{\Q}
    \end{equation}
(where we have suppressed the $1/N$ notation from the left). As in \cite[{(5.158)}]{GS19}, a short computation (using well-definedness of arithmetic intersection products) shows that this class does not depend on the choice of $(\mc{Z}_0,g_0)$. One can verify that $g_{T,y} + g_0 \wedge \delta_{\mc{Z}(T)_{K'_f,\C}}$ satisfies a Green current equation for $(\mc{Z}(T)_{K'_f} \cap \mc{Z}_0)_{\C}$ by combining the Green current equation for $g_0$ with the modified current equation of $g_{T,y}$ (see also \cite[{\S 5.4}]{GS19}).

These classes $[\widehat{\mc{Z}}(T)_{K'_f, \ms{H}}]$ thus form a compatible system as $K'_f$ varies, and hence give an element
    \begin{equation}\label{equation:arith_cycle_classes:horizontal}
    [\widehat{\mc{Z}}(T)_{\ms{H}}] \coloneqq ([\widehat{\mc{Z}}(T)_{K'_f, \ms{H}}])_{K'_f} \in \arithCh^m(\mc{M}[1/N])_{\Q}.
    \end{equation}
This construction of $[\widehat{\mc{Z}}(T)_{\ms{H}}]$ is essentially that of \cite[{\S 5.4}]{GS19}.
If $g_{T,y} = g_{{}^t \overline{\gamma} T \gamma, \gamma^{-1} y {}^t \overline{\gamma}^{-1}}$, note that we automatically have the ``linear invariance'' equality
    \begin{equation}\label{equation:arith_cycle_classes:horizontal:linear_invariance}
    [\widehat{\mc{Z}}(T)_{\ms{H}}] = [\widehat{\mc{Z}}({}^t \overline{\gamma} T \gamma)_{\ms{H}}].
    \end{equation}

Currents $g_{T,y}$ satisfying \eqref{equation:arith_cycle_classes:horizontal:Green_current} were studied by Garcia--Sankaran \cite[{Theorem 1.1 and \S 4}]{GS19}. 
We choose to use the star-product approach of Kudla \cite{Kudla97b} to define currents $g_{T,y}$ for our main results (for $\rank T \geq n - 1$ or $\det T \neq 0$ with $T$ not positive definite), and postpone the explicit description of $g_{T,y}$ to Section \ref{ssec:Arch_uniformization:Archimedean} \eqref{equation:Arch_uniformization:Archimedean:current_definition}.
Our definition of $g_{T,y}$ is that of \cite[{Theorem 4.20}]{Liu11} in the nonsingular cases. When $T \in \mrm{Herm}_n(\Q)$ is singular with $\rank(T) = n - 1$, our definition is new (still based on star products). These Green currents satisfy $g_{T,y} = g_{{}^t \overline{\gamma} T \gamma, \gamma^{-1} y {}^t \overline{\gamma}^{-1}}$ (Section \ref{ssec:Arch_uniformization:Archimedean}), so linear invariance \eqref{equation:arith_cycle_classes:horizontal:linear_invariance} is satisfied.

        \subsection{Vertical special cycle classes}
        \label{ssec:arith_cycle_classes:vertical}
            We define vertical special cycle classes via $K_0$ groups. 
We remind the reader that $d_L \in \Z$ is an integer associated to the lattice $L$ as discussed before Definition \ref{definition:stack_ass_to_lattice}. Recall that $d_L = 1$ if $L$ is self-dual of signature $(n - 1, 1)$ with $2 \nmid \Delta$.

For our notation and definitions regarding $K_0$-groups for Deligne--Mumford stacks, we refer to Appendix \ref{appendix:K0}. Note that the stacky $K_0$ groups we use are different from those used in \cite{HM22}. It is also possible to avoid stacky $K_0$ groups entirely by working with compatible systems of classes in towers of level structure.

Fix any prime $p \nmid d_L$ and set
    \begin{equation}
    \mc{M}_{(p)} \coloneqq \mc{M} \times_{\Spec \Z} \Spec \Z_{(p)} \quad \quad \mc{Z}(T)_{(p)} \coloneqq \mc{Z}(T) \times_{\Spec \Z} \Spec \Z_{(p)}
    \end{equation}
for any $T \in \mrm{Herm}_m(\Q)$. Since $\mc{M}_{(p)}$ admits a finite \'etale cover by a scheme (add away-from-$p$ level structure), we may consider filtrations for $K_0'$ groups as in Definition \ref{definition:stack_dimension_filtration}.

Let $T \in \mrm{Herm}_m(\Q)$ be any $m \times m$ Hermitian matrix (with entries in $F$). We first describe a ``$p$-local'' derived special cycle class ${}^{\mbb{L}} \mc{Z}(T)_{(p)} \in F^m_{\mc{M}_{(p)}} K'_0(\mc{Z}(T)_{(p)})_{\Q}$ (with $F^m_{\mc{M}_{(p)}} \coloneqq F_{n - m}$ denoting the $m$-th step of the codimension filtration) before extracting a ``vertical'' piece.

For any $t \in \Q$, we define ${}^{\mbb{L}} \mc{Z}(t)_{(p)} \in F^1_{\mc{M}_{(p)}} K'_0(\mc{Z}(t)_{(p)})_{\Q}$ to be the element
    \begin{align}
    {}^{\mbb{L}} \mc{Z}(t)_{(p)} \coloneqq
    \begin{cases}
    [\mc{O}_{\mc{Z}(t)_{(p)}}] & \text{if $t \neq 0$} \\
    [\mc{O}_{\mc{M}_{(p)}}] - [\mc{E}] & \text{if $t = 0$}. \notag
    \end{cases}
    \end{align}
Write $t_1,\ldots,t_m$ for the diagonal entries of $T$. Using the intersection pairing of Lemma \ref{lemma:K0_relative} as well as compatibility with dimension filtrations from Lemma \ref{lemma:K0_relative_with_filtrations} via Lemma \ref{lemma:special_cycles_level_structure}, we form the intersection ${}^{\mbb{L}} \mc{Z}(t_1)_{(p)} \cdots {}^{\mbb{L}} \mc{Z}(t_m)_{(p)}$ and define ${}^{\mbb{L}} \mc{Z}(T)_{(p)}$ by the restriction
    \begin{equation}\label{equation:setup:arith_cycle_classes:vertical:derived_definition} 
    \begin{tikzcd}[row sep = tiny]    
    F^m_{\mc{M}} K_0'(\mc{Z}(t_1)_{(p)} \times_{\mc{M}_{(p)}} \cdots \times_{\mc{M}_{(p)}} \mc{Z}(t_m)_{(p)})_{\Q} \arrow{r} & F^m_{\mc{M}_{(p)}} K'_0(\mc{Z}(T)_{(p)})_{\Q} \\
    {}^{\mbb{L}} \mc{Z}(t_1)_{(p)} \cdots {}^{\mbb{L}} \mc{Z}(t_m)_{(p)} \arrow[mapsto]{r} & {}^{\mbb{L}} \mc{Z}(T)_{(p)}.
    \end{tikzcd}
    \end{equation}
This displayed restriction map comes from the disjoint union decomposition in \eqref{equation:special_cycles_diagonal_entry_intersection}.
We call ${}^{\mbb{L}} \mc{Z}(T)_{(p)}$ the \emph{$p$-local derived special cycle class}\footnote{This is the construction of \cite[{Definition 5.1.3}]{HM22} (there for orthogonal Shimura varieties). This construction also underlies the intersection numbers considered in \cite{KR14} for non-degenerate $T$. We differ slightly from those references by localizing at $p$, since we will only be interested in the ``vertical'' part of ${}^{\mbb{L}} \mc{Z}(T)_{(p)}$. The ``horizontal part'' is accounted for by Section \ref{ssec:arith_cycle_classes:horizontal}.
} associated with $T$. 

The following lemma is a ``$p$-local'' version of linear invariance, and is proved using a variant on ideas from \cite{Howard19,HM22}. The map in \eqref{equation:arith_cycle_classes:vertical:p-local_linear_invariance} was defined in \eqref{equation:ab_var:special_cycles:gamma_action}.

\begin{lemma}\label{lemma:arith_cycle_classes:vertical:p-local_linear_invariance}
Given any $T \in \mrm{Herm}_m(\Q)$ and any $\gamma \in \GL_m(\mc{O}_{F,(p)}) \cap M_{m,m}(\mc{O}_F)$, the pullback along
    \begin{equation}\label{equation:arith_cycle_classes:vertical:p-local_linear_invariance}
    \mc{Z}(T)_{(p)} \ra \mc{Z}({}^t \overline{\gamma} T \gamma)_{(p)}
    \end{equation}
sends ${}^{\mbb{L}} \mc{Z}({}^t \overline{\gamma} T \gamma)_{(p)}$ to ${}^{\mbb{L}} \mc{Z}(T)_{(p)}$.
\end{lemma}
\begin{proof}
By Lemma \ref{lemma:ab_var:special_cycles:gamma_action}, we know that \eqref{equation:arith_cycle_classes:vertical:p-local_linear_invariance} is an open and closed immersion.

The ring $\mc{O}_{F,(p)}$ is a Euclidean domain, with Euclidean function $\phi(a) \coloneqq \sum_{\mf{p}_i} v_{\mf{p}_i}(a) \cdot f_i$ for nonzero $a \in \mc{O}_{F,(p)}$ (summing over primes $\mf{p}_i$ in $\mc{O}_F$ lying over $p$, with residue cardinality $p^{f_i}$). Row reducing via the Euclidean algorithm shows that $\GL_m(\mc{O}_{F,(p)})$ is generated by elementary matrices.

Any $\gamma \in \GL_m(\mc{O}_{F,(p)})$ may thus be expressed as $\gamma = \gamma_1 \gamma_2^{-1}$ where each $\gamma_1$ and $\gamma_2$ are products of elementary matrices lying in $\GL_m(\mc{O}_{F,(p)}) \cap M_{m,m}(\mc{O}_F)$. If moreover $\gamma \in \GL_m(\mc{O}_{F,(p)}) \cap M_{m,m}(\mc{O}_F)$, the commutative diagram
\begin{equation}
    \begin{tikzcd}
    & \mc{Z}({}^t \overline{\gamma}_1 T \gamma_1)_{(p)} & \\
    \mc{Z}(T)_{(p)} \arrow{rr} \arrow{ur} & & \arrow{ul} \mc{Z}({}^t \overline{\gamma} T \gamma)_{(p)}
    \end{tikzcd}
    \end{equation}
shows that it is enough to prove the lemma when $\gamma \in \GL_m(\mc{O}_{F,(p)}) \cap M_{m,m}(\mc{O}_F)$ is an elementary matrix.

If $\gamma$ is a permutation matrix, the lemma is clear. Next, consider $a \in \mc{O}_{F,(p)}^{\times} \cap \mc{O}_F$. For any $t \in \Q$, note that $\mc{Z}(t)_{(p)} \ra \mc{Z}(\overline{a} t a)_{(p)}$ is an open and closed immersion (by Lemma \ref{lemma:ab_var:special_cycles:gamma_action} again). This fact implies the present lemma for the case where $\gamma = \mrm{diag}(a, 1, \ldots, 1)$.

It remains to check the case where $\gamma$ is an elementary unipotent matrix. This case follows as in the analogous result \cite[{Proposition 5.4.1}]{HM22} (there for $\mrm{GSpin}$).\footnote{Strictly speaking, our setup for stacky $K'_0$ groups is slightly different from that of \cite{HM22}, see Appendix \ref{appendix:K0}. This makes no difference in the proof of the cited result. Alternatively, one can replace $\mc{M}_{(p)}$ by a finite \'etale cover by a scheme to reduce to the case of schemes, where our setup agrees with \cite[{\S A.2}]{HM22}.} 
The latter is proved using methods from \cite{Howard19} (the analogous local linear invariance result on Rapoport--Zink spaces). We are also using global analogues of \cite[{Lemma 2.36, Lemma 2.37, Lemma 2.41}]{LL22II} (there about a tautological bundle on an exotic smooth Rapoport--Zink space) which may be proved similarly, e.g. our Lemma \ref{lemma:unique_Kramer_hyperplane} replaces \cite[{Lemma 2.36}]{LL22II} in the global setup.
Alternatively, linear invariance for $\gamma \in \GL_m(\mc{O}_F)$ should also follow from the derived algebro-geometric methods in \cite{Madapusi22}.
\end{proof}

Next, we define a \emph{derived vertical special cycle class}
    \begin{equation}
    {}^{\mbb{L}} \mc{Z}(T)_{\ms{V},p} \in \mrm{gr}^m_{\mc{M}} K'_0(\mc{Z}(T)_{\F_p})_{\Q}
    \end{equation}
at $p$, where $\mc{Z}(T)_{\F_p} \coloneqq \mc{Z}(T) \times_{\Spec \Z} \Spec \F_p$.

First consider the case where $\det T \neq 0$. Using Lemmas \ref{lemma:special_cycles_generically_smooth} and \ref{lemma:K'0_component_decomp} as well as \eqref{equation:setup:ab_var:horizontal_vertical_decomp}, we decompose
    \begin{equation}\label{equation:nonsingular_horizontal_vertical_decomp}
    \mrm{gr}^{m}_{\mc{M}_{(p)}} K'_0(\mc{Z}(T)_{(p)})_{\Q} = \mrm{gr}^{m}_{\mc{M}_{(p)}} K'_0(\mc{Z}(T)_{{(p)},\ms{H}})_{\Q} \oplus \mrm{gr}^{m}_{\mc{M}} K'_0(\mc{Z}(T)_{\F_p})_{\Q}
    \end{equation}
into a ``horizontal'' part and a ``vertical'' part. This uses nonsingularity of $T$ (via Lemma \ref{lemma:special_cycles_generically_smooth}), so that $\mc{Z}(T)_{{(p)},\ms{H}} \cap \mc{Z}(T)_{\F_p}$ is of dimension $< n - m$.
We are also using the d\'evissage pushforward identification $K'_0(\mc{Z}(T)_{\F_{p}}) \xra{\sim} K'_0(\mc{Z}(T)_{\ms{V},p})$, with $\mc{Z}(T)_{\ms{V},p}$ as in \eqref{equation:setup:ab_var:horizontal_vertical_decomp}. The above decomposition of $\mrm{gr}^{m}_{\mc{M}_{(p)}} K'_0(\mc{Z}(T)_{(p)})_{\Q}$ is independent of the choice of $e_{p}$ in \eqref{equation:setup:ab_var:horizontal_vertical_decomp}.
We define ${}^{\mbb{L}} \mc{Z}(T)_{\ms{V},p}$ to be given by the projection
    \begin{equation}
    \begin{tikzcd}[row sep = tiny]
    \mrm{gr}^{m}_{\mc{M}} K'_0(\mc{Z}(T)_{(p)})_{\Q} \arrow{r} & \mrm{gr}^{m}_{\mc{M}} K'_0(\mc{Z}(T)_{\F_p})_{\Q} \\
    {}^{\mbb{L}} \mc{Z}(T)_{(p)} \arrow[mapsto]{r} & {}^{\mbb{L}} \mc{Z}(T)_{\ms{V},p}.
    \end{tikzcd}
    \end{equation}

If $T = \mrm{diag}(0,T^{\flat})$ (with $\det T^{\flat} \neq 0$), we set
    \begin{equation}
    {}^{\mbb{L}} \mc{Z}(T)_{\ms{V},p} \coloneqq (\mc{E}^{\vee})^{m - \rank(T)} \cdot {}^{\mbb{L}} \mc{Z}(T^{\flat})_{\ms{V},p} \in \mrm{gr}^m_{\mc{M}} K'_0(\mc{Z}(T)_{\F_p})_{\Q}
    \end{equation}
where $\mc{E}^{\vee}$ stands for the class $[\mc{O}_{\mc{M}_{(p)}}] - [\mc{E}] \in F^1_{\mc{M}_{(p)}}(\mc{M}_{(p)})_{\Q}$. Given arbitrary $T$ (not necessarily block diagonal), select any $\gamma \in \GL_m(\mc{O}_{F,(p)}) \cap M_{m,m}(\mc{O}_F)$ such that
    \begin{equation}
    {}^t \overline{\gamma} T \gamma = \mrm{diag}(0,T^{\prime\flat})
    \end{equation}
where $\det T^{\prime \flat} \neq 0$. Set $T' \coloneqq \mrm{diag}(0,T^{\prime\flat})$. We define ${}^{\mbb{L}} \mc{Z}(T)_{\ms{V},p}$ to be the pullback class
    \begin{equation}\label{equation:arith_cycle_classes:vertical}
    \begin{tikzcd}[row sep = tiny]
    \mrm{gr}^m_{\mc{M}}(\mc{Z}(T)_{\F_p})_{\Q} & \arrow{l} \mrm{gr}^m_{\mc{M}}(\mc{Z}(T')_{\F_p})_{\Q} \\
    {}^{\mbb{L}} \mc{Z}(T)_{\ms{V},p} & \arrow[mapsto]{l} {}^{\mbb{L}} \mc{Z}(T')_{\ms{V},p}
    \end{tikzcd}
    \end{equation}
along the open and closed immersion $\mc{Z}(T)_{\F_p} \ra \mc{Z}(T')_{\F_p}$ induced by $\gamma$.

By Lemma \ref{lemma:arith_cycle_classes:vertical:p-local_linear_invariance} (applied to $T^{\prime \flat}$, in the notation above), the preceding definition of ${}^{\mbb{L}} \mc{Z}(T)_{\ms{V},p}$ does not depend on the choice of $\gamma$. Moreover, the class $\mc{Z}(T)_{\ms{V},p}$ is \emph{linearly invariant} in the following sense: given any $\gamma \in \GL_m(\mc{O}_{F,(p)}) \cap M_{m,m}(\mc{O}_F)$ (no additional assumptions on ${}^t \overline{\gamma} T \gamma$), the pullback along
    \begin{equation}\label{equation:arith_cycle_classes:vertical:linear_invariance}
    \mc{Z}(T)_{\F_p} \ra \mc{Z}({}^t \overline{\gamma} T \gamma)_{\F_p}
    \end{equation}
sends ${}^{\mbb{L}} \mc{Z}({}^t \overline{\gamma} T \gamma)_{\ms{V},p}$ to ${}^{\mbb{L}} \mc{Z}(T)_{\ms{V},p}$. This follows from the construction of ${}^{\mbb{L}} \mc{Z}(T)_{\ms{V},p}$.

We see that ${}^{\mbb{L}}\mc{Z}(T)_{\ms{V},p} = 0$ for all but finitely many primes $p$. Taking pushforward, we obtain associated cycle classes $[{}^{\mbb{L}} \mc{Z}(T)_{\ms{V},p}] \in \arithCh^m(\mc{M})_{\Q}$.
The preceding constructions may be repeated (essentially verbatim) with $K'_f$ level structure away from $p$. The resulting classes ${}^{\mbb{L}} \mc{Z}(T)_{(p),K'_f}$ and ${}^{\mbb{L}} \mc{Z}(T)_{\ms{V},p,K'_f}$ will be compatible with pullback for varying level.

        \subsection{Degrees of arithmetic special cycles}
        \label{ssec:arith_cycle_classes:degrees}

The moduli stack $\mc{M} \ra \Spec \mc{O}_F[1/d_L]$ considered in Section \ref{ssec:arith_cycle_classes:arithmetic_cycle_classes} may not be proper. For a robust arithmetic degree theory via arithmetic Chow groups for arbitrary $T \in \mrm{Herm}_m(\Q)$, one might instead consider arithmetic special cycle classes on a suitably compactified moduli space.

If the special cycle $\mc{Z}(T)$ is already proper over $\Spec \mc{O}_F[1/d_L]$, we can directly define the \emph{arithmetic degree without boundary contributions} which should result from a compactification: set

    \begin{align}\label{equation:arith_cycle_classes:degrees:if_proper}
        \widehat{\deg}([\widehat{\mc{Z}}(T)] \cdot \widehat{c}_1(\widehat{\mc{E}}^{\vee})^{n-m}) & \coloneqq \left ( \int_{\mc{M}_{\C}} g_{T,y} \wedge c_1(\widehat{\mc{E}}^{\vee}_{\C})^{n-m} \right ) 
        \\
        & \hphantom{\coloneqq} + \widehat{\deg}((\widehat{\mc{E}}^{\vee})^{n - \rank(T)}|_{\mc{Z}(T)_{\ms{H}}}) \notag
        \\
        & \hphantom{\coloneqq} + \sum_{\substack{p \text{ prime} \\ p \nmid d_L}} \deg_{\F_p}({}^{\mathbb{L}}\mc{Z}(T)_{\ms{V},p} \cdot (\mc{E}^{\vee})^{n-m}) \log p \notag
    \end{align}
conditional on convergence of the integral. Since compactification of $\mc{M}$ plays no other role in this work, we take this approach. As in Section \ref{ssec:arith_cycle_classes:horizontal}, the notation $\mc{M}_{\C}$ mean $\mc{M} \times_{\Spec \mc{O}_F} \Spec \C$ for a choice of $F \ra \C$ (the choice does not matter).

The quantity in \eqref{equation:arith_cycle_classes:degrees:if_proper} is an element of $\R_{d_L} = \R / (\sum_{p \mid d_L} \Q \cdot \log p)$. Here $\widehat{\deg}((\widehat{\mc{E}}^{\vee})^{n-m}|_{\mc{Z}(T)_{\ms{H}}})$ is the height of $\mc{Z}(T)_{\ms{H}}$ with respect to the metrized tautological bundle $\widehat{\mc{E}}^{\vee}$ (Sections \ref{ssec:arith_cycle_classes:arithmetic_Chow} and \ref{ssec:arith_cycle_classes:Hodge_bundles}).
The symbol ${}^{\mbb{L}} \mc{Z}(T)_{\ms{V},p} \cdot (\mc{E}^{\vee})^{n-m}$ is shorthand for the intersection product
    \begin{equation}
    {}^{\mbb{L}} \mc{Z}(T)_{\ms{V},p} \cdot ([\mc{O}_{\mc{M}}] - [\mc{E}])^{n-m} \in \mrm{gr}^n_{\mc{M}} K'_0 (\mc{Z}(T)_{\F_p})_{\Q} = \mrm{gr}_0 K'_0(\mc{Z}(T)_{\F_p})_{\Q},
    \end{equation}
defined in Lemma \ref{lemma:K0_relative}. With $\mc{Z}(T)_{\F_p}$ viewed as a proper scheme over $\F_p$, the notation $\deg_{\F_p}$ refers to the degree map $\deg \colon \mrm{gr}_0 K'_0(\mc{Z}(T)_{\F_p})_{\Q} \ra \Q$ as defined in \eqref{equation:euler_char_degree}.

Certainly $\mc{Z}(T) \ra \Spec \mc{O}_F[1/d_L]$ is proper if $\mc{Z}(T)$ is empty (e.g. if $T$ is not positive semidefinite). In this case, the right-hand side of \eqref{equation:arith_cycle_classes:degrees:if_proper} consists only of the integral $\int_{\mc{M}_{\C}} g_{T,y} \wedge c_1(\widehat{\mc{E}}^{\vee}_{\C})^{n-m}$.

We show below that $\mc{Z}(T) \ra \Spec \mc{O}_F[1/d_L]$ is also proper if $\rank(T) \geq n - 1$, so \eqref{equation:arith_cycle_classes:degrees:if_proper} applies in this case as well.

\begin{lemma}\label{lemma:corank_1_isogeny_decomposition}
Fix a Hermitian matrix $T \in \mrm{Herm}_m(\Q)$ with $\rank(T) \geq n-1$ and $m \geq 0$. Let $\kappa$ be a field, and consider $(A_0, \iota_0, \lambda_0, A, \iota, \lambda, \underline{x}) \in \mc{Z}(T)(\kappa)$. There exists an $\mc{O}_F$-linear isogeny $(A_0)^{n-1} \times A^- \ra A$ where $A^-$ is an elliptic curve with $\mc{O}_F$-action. After replacing $\kappa$ by a finite extension, we may take $A^- = A_0^{\s}$ where $A_0^{\s} = A_0$ but with $\mc{O}_F$-action $\iota_0 \circ \s$.
\end{lemma}
\begin{proof}
Write $\underline{x} = [x_1,\ldots,x_m]$, where the $x_i$ are $\mc{O}_F$-linear homomorphisms $x_i \colon A_0 \ra A$. Since $T$ has rank $\geq n-1$, we may assume (rearranging the elements $x_i$ if necessary) that $\underline{x}^{\flat} = [x_1, \ldots, x_{n-1}]$ has nonsingular Gram matrix $(\underline{x}^{\flat}, \underline{x}^{\flat})$. Then the map
    \begin{equation}
    f \colon A \xra{\sqrt{\Delta} \circ (x_1^{\dagger} \times \cdots \times x_{n-1}^{\dagger})} A_0^{n-1}
    \end{equation}
is a homomorphism and a surjection of fppf sheaves. 
We form the ``isogeny complement'' in the standard way, i.e. we let $A^-$ be the reduced connected component of $\ker f$. If $j \colon A^- \ra A$ is the natural inclusion, then the map $(A_0)^{n-1} \times A^- \xra{x_1 \times \cdots \times x_{n-1} \times j} A$ is an $\mc{O}_F$-linear isogeny.

Note that $A^-$ is $\mc{O}_F$-linearly isogenous to an elliptic curve with $\mc{O}_F$ action of signature $(0,1)$: if $\mrm{char}(k) = p > 0$ with $p$ nonsplit in $\mc{O}_F$, then $A^-$ is supersingular, so apply Skolem--Noether to $\End(A^-) \otimes \Q$; otherwise, $A^-$ automatically has signature $(0,1)$ because $A$ has signature $(n - 1, 1)$.

If $\kappa$ is algebraically closed, any two elliptic curves over $\kappa$ with $\mc{O}_F$-action of the same signature are $\mc{O}_F$-linearly isogenous. This is classical: lift to characteristic zero to reduce to $\kappa = \C$ (the moduli stack $\ms{M}_0 \ra \Spec \mc{O}_F$ is \'etale; 
more classically, see Deuring \cite{Deuring41}); recall that elliptic curves over $\C$ with $\mc{O}_F$-action are defined and isogenous over $\overline{\Q}$. By a standard limiting argument, we conclude that $A^-$ and $A_0^{\s}$ are $\mc{O}_F$-linearly isogenous over a finite extension of the (not necessarily algebraically closed) original field $\kappa$.
\end{proof}

\begin{remark}\label{remark:arith_cycle_clases:degrees:split_full_rank_empty}
If $p$ is a prime which splits in $\mc{O}_F$ and if $\rank(T) \geq n$, then $\mc{Z}(T)_{\F_p} = \emptyset$. This is a standard argument (e.g. \cite[{Lemma 2.21}]{KR14}): if $\kappa$ is a field of characteristic $p$ and $(A_0, \iota_0, \lambda_0, A, \iota, \lambda, \underline{x}) \in \mc{Z}(T)(\kappa)$, arguing as in Lemma \ref{lemma:corank_1_isogeny_decomposition} shows that $A$ is $\mc{O}_F$-linearly isogenous to $A_0^n$. This contradicts Lemma \ref{lemma:corank_1_isogeny_decomposition}, because there is no nonzero $\mc{O}_F$-linear map $A_0 \ra A_0^{\s}$ as $A_0$ and $A_0^{\s}$ have $\mc{O}_F$-action of opposite signature (e.g. there are no nonzero maps of the underlying ordinary $p$-divisible groups).
\end{remark}

We say a characteristic $p > 0$ geometric point $(A_0,\iota_0, \lambda_0, A, \iota, \lambda)$ of $\mc{M}$ lies in the \emph{supersingular locus} if $A_0$ and $A$ are supersingular abelian varieties (i.e. the associated $p$-divisible groups are supersingular). The following corollary also holds in arbitrary signature $(n - r, r)$ (i.e. all but the last sentence of Lemma \ref{lemma:corank_1_isogeny_decomposition} is valid for arbitrary signature $(n - r, r)$).

\begin{corollary}\label{corollary:special_cycle_supersingular_support}
Let $p$ be a prime which is nonsplit in $\mc{O}_F$. Fix $T \in \mrm{Herm}_m(\Q)$ with $\rank(T) \geq n-1$ and $m \geq 0$. The morphism $\mc{Z}(T)_{\overline{\F}_p} \ra \mc{M}_{\overline{\F}_p}$ factors (set-theoretically) through the supersingular locus on $\mc{M}_{\overline{\F}_p}$.
\end{corollary}
\begin{proof}
Follows from Lemma \ref{lemma:corank_1_isogeny_decomposition} and Deuring's classical results on endomorphisms of elliptic curves in positive characteristic \cite{Deuring41} (i.e. over a field of characteristic $p > 0$, the $p$-divisible group of an elliptic curve with $\mc{O}_F$-action is supersingular (resp. ordinary) if $p$ is nonsplit (resp. split) in $\mc{O}_F$). Here we used the notation $\mc{Z}(T)_{\overline{\F}_p} \coloneqq \mc{Z}(T) \times_{\Spec \Z} \Spec \overline{\F}_p$ and similarly for $\mc{M}$.
\end{proof}

\begin{lemma}\label{lemma:horizontal_global_special_cycles_proper_quasi-finite}
Fix $T \in \mrm{Herm}_m(\Q)$ with $\rank(T) \geq n-1$ and $m \geq 0$. Then the horizontal special cycle $\mc{Z}(T)_{\ms{H}}$ is proper and quasi-finite over $\Spec \mc{O}_F[1/d_L]$.
\end{lemma}
\begin{proof}
By Lemma \ref{lemma:special_cycles_generically_smooth}, we know the generic fiber $\mc{Z}(T)_{\ms{H}, F} \ra \Spec F$ is smooth of relative dimension $0$. Hence each generic point of $\mc{Z}(T)_{\ms{H}}$ is the image of a map $\Spec E \ra \mc{Z}(T)$ for some number field $E$, corresponding to an object $(A_0, \iota_0, \lambda_0, A, \iota, \lambda, \underline{x}) \in \mc{Z}(T)(E)$. By Lemma \ref{lemma:corank_1_isogeny_decomposition}, we know that $A$ is isogenous to a product of elliptic curves with complex multiplication by $\mc{O}_F$. It is a classical result of Deuring that such elliptic curves have everywhere potentially good reduction, so $A_0$ and $A$ have everywhere potentially good reduction \cite{Deuring41,ST68}. Enlarging $E$ if necessary, we can thus extend $\Spec E \ra \mc{Z}(T)$ to a morphism $\Spec \mc{O}_E[1/d_L] \ra \mc{Z}(T)$ (the N\'eron mapping property ensures that the datum $(\iota_0, \lambda_0, \iota, \lambda, \underline{x})$ extends as well; the polarizations must extend to polarizations as in the proof of \cite[{Theorem 1.9}]{FC90}).

Hence each irreducible component of $\mc{Z}(T)_{\ms{H}}$ may be covered by a morphism $\Spec \mc{O}_E[1/d_L] \ra \mc{Z}(T)$ for some number field $E$. Since $\mc{Z}(T)$ is quasi-compact and separated, this implies that $\mc{Z}(T) \ra \Spec \mc{O}_F[1/d_L]$ is proper and quasi-finite.
\end{proof}

\begin{lemma}\label{lemma:global_special_cycles_proper}
For $m \geq 0$, suppose $T \in \mrm{Herm}_m(\Q)$  has $\rank(T) \geq n-1$. Then the structure map $\mc{Z}(T) \ra \Spec \mc{O}_F[1/d_L]$ is proper.
\end{lemma}
\begin{proof}
We already know that the horizontal part $\mc{Z}(T)_{\ms{H}}$ is proper over $\Spec \mc{O}_F[1/d_L]$, so it suffices to check that every irreducible component of $\mc{Z}(T)$ in characteristic $p \nmid d_L$ is proper over $\Spec \F_p$. It is enough to check that $\mc{Z}(T)_{\overline{\F}_p} \ra \Spec \overline{\F}_p$ is proper by fpqc descent (e.g. use away-from-$p$ level structure to replace $\mc{Z}(T)_{\overline{\F}_p}$ with a finite cover by a scheme, then use fpqc descent for morphisms of schemes). It is enough to check properness of the map $\mc{Z}(T)_{\overline{\F}_p,\mrm{red}} \ra \Spec \overline{\F}_p$ on reduced substacks (e.g. by local Noetherianity of these algebraic stacks, or because $\mc{Z}(T)_{\overline{\F}_p} \ra \Spec \overline{\F}_p$ is locally of finite type).

The supersingular locus on $\mc{M}_{\overline{\F}_p}$ is proper (follows from the proof of \cite[{Theorem 1.1a}]{Oort74}; and finiteness of the forgetful map $\ms{M}^{(d)} \ra \mc{A}_{n,d}$ in Section \ref{ssec:ab_var:integral_models}). Properness of $\mc{Z}(T)_{\overline{\F}_p} \ra \Spec \overline{\F}_p$ now follows from Corollary \ref{corollary:special_cycle_supersingular_support}.
\end{proof}

\clearpage


\part{Local change of heights}
\label{part:local_change_heights}
    Fix an imaginary quadratic field $F / \Q$. Write $\Delta$ for the discriminant and $\s$ for the nontrivial involution. We allow $2 \mid \Delta$ in Part \ref{part:local_change_heights} unless otherwise specified. We set $F_p \coloneqq F \otimes_{\Q} \Q_p$ and $\mc{O}_{F_p} \coloneqq \mc{O}_F \otimes_{\Z} \Z_p$.

Throughout Part \ref{part:local_change_heights}, we write $E$ for a number field, with ring of integers $\mc{O}_E$. Given a prime $p$, we set $\mc{O}_{E, (p)} \coloneqq \mc{O}_E \otimes_{\Z} \Z_{(p)}$. We use $\breve{E}$ to denote a finite degree field extension of $\breve{\Q}_p$, with ring of integers $\mc{O}_{\breve{E}}$. We write $\mf{d}_p \subseteq \mc{O}_{F_p}$ for the different ideal of $\mc{O}_{F_p/\Q_p}$. We abuse notation and also mean $\mf{d}_p \coloneqq \sqrt{\Delta}$, which is a generator of the different ideal.

By a \emph{place} $\breve{w}$ of $E \otimes_{\Q} \breve{\Q}_p$, we mean a prime ideal of $E \otimes_{\Q} \breve{\Q}_p$ (equivalently, an element of $\Hom_{\breve{\Q}_p}(E \otimes_{\Q} \breve{\Q}_p, \C_p)$ up to automorphisms of $\C_p$). We write $\breve{E}_{\breve{w}}$ for the residue field of $\breve{w}$, with ring of integers $\mc{O}_{\breve{E}_{\breve{w}}}$. We use the shorthand $\breve{w} \mid p$ to indicate a place of $E \otimes_{\Q} \breve{\Q}_p$, and may use subscripts (e.g. $X'_{\breve{w}}$ and $\phi_{\breve{w}}$ in Section \ref{ssec:qcan_heights:minimal_isogenies}) to indicate base-change from $\Spec \mc{O}_{E,(p)}$ to $\Spf \mc{O}_{\breve{E}_{\breve{w}}}$.

Whenever an $\mc{O}_F$-action or $F$-action is mentioned (e.g. on a sheaf of modules on $\Spec \mc{O}_E$), we assume that $\mc{O}_E$ (resp. $\mc{O}_{\breve{E}}$) is equipped with morphism $\mc{O}_F \ra \mc{O}_E$ (resp. $\mc{O}_F \ra \mc{O}_{\breve{E}}$).

We write $\mf{X}_s$ for a level $s \geq 0$ quasi-canonical lifting of signature $(1,0)$ over $\Spec \mc{O}_{\breve{E}}$ with its $\mc{O}_{F_p}$-action $\iota_{\mf{X}_s}$, as explained in \crefext{I:sec:can_and_qcan}. The framing $\rho_{\mf{X}_s}$ of loc. cit. is unimportant in Part \ref{part:local_change_heights} (and will be omitted). As before, the notation $\mf{X}_s^{\s}$ means $\mf{X}_s$ but with $\mc{O}_{F_p}$-action pre-composed by $\s$.

Given a group scheme $G$ over a base $S$, we typically write $e \colon S \ra G$ for the identity section. We abuse notation and use ``$e$'' simultaneously for different group schemes.

    \section{Faltings and ``tautological'' heights}
    \label{sec:Faltings_and_taut}
    
        \subsection{Heights}
        \label{ssec:Faltings_and_taut:heights}

Suppose $A \ra \Spec \mc{O}_E$ is a semi-abelian N\'eron model of an abelian variety over $E$.
The \emph{Faltings height} of $A$ (or its generic fiber $A_E$) is
    \begin{equation}
    h_{\mrm{Fal}}(A_E) \coloneqq h_{\mrm{Fal}}(A) \coloneqq \frac{1}{[E : \Q]} \widehat{\deg}(\widehat{\omega}_{A})
    \end{equation}
where $\widehat{\omega}_{A} = (\omega_{A}, \norm{-}) = (e^* \bigwedge^n \Omega^1_{A/\mc{O}_E}, \norm{-})$ is the Hermitian line bundle with norm $\norm{-}$ normalized as in \eqref{equation:arith_cycle_classes:Hodge_bundles:Faltings_metric}. The usual arithmetic degree $\widehat{\deg}$ was recalled in Section \ref{ssec:arith_cycle_classes:Hermitian_bundles}. Any abelian variety over a number field has everywhere potentially semi-abelian reduction, and the \emph{Faltings height} of any abelian variety $B$ over $\Spec E$ is defined so that $h_{\mrm{Fal}}(B)$ is remains constant under finite field extensions $E \ra E'$. (We only consider stable Faltings height, as defined above.)

We also consider certain ``tautological heights'' to describe the arithmetic intersections appearing in Section \ref{ssec:arith_cycle_classes:degrees}. The terminology we introduce for this (e.g. ``Kr\"amer datum'') is likely nonstandard.

\begin{definition}\label{definition:Falting_and_taut:heights:Kramer_datum}
\hfill
\begin{enumerate}[(1)]
    \item Given a scheme $S$ over $\Spec \mc{O}_F$, a \emph{Kr\"amer datum} (of signature $(n - 1, 1)$) is a tuple $(A, \iota, \mc{F})$ where $A \ra S$ is an abelian scheme, where $\iota \colon \mc{O}_F \ra \End(A)$ an action of signature $(n - 1, 1)$, and where $\mc{F} \subseteq \Lie A$ is a $\iota$-stable local direct summand of rank $n - 1$ such that the $\mc{O}_F$ action via $\iota$ on $\mc{F}$ (resp. $(\Lie A) / \mc{F}$) is $\mc{O}_F$-linear (resp. $\s$-linear). We say that $\mc{F}$ is the associated \emph{Kr\"amer hyperplane}.
    
    \item Given a formal scheme $S$ over $\Spf \mc{O}_{F_p}$, a \emph{local Kr\"amer datum} (of signature $(n - 1,1)$) is a tuple $(X, \iota, \mc{F})$ where $X$ is a $p$-divisible group over $S$ of height $2n$ and dimension $n$, where $\iota \colon \mc{O}_{F_p} \ra \End(X)$ is an action of signature $(n - 1, 1)$, and where $\mc{F} \subseteq \Lie X$ is a $\iota$-stable local direct summand of rank $n - 1$ such that the $\mc{O}_F$ action via $\iota$ on $\mc{F}$ (resp. $(\Lie X) / \mc{F}$) is $\mc{O}_F$-linear (resp. $\s$-linear). We say that $\mc{F}$ is the associated \emph{Kr\"amer hyperplane}.
    
    \item A \emph{quasi-polarized} Kr\"amer datum (resp. \emph{quasi-polarized} local Kr\"amer datum) is a tuple $(A, \iota, \lambda, \mc{F})$ (resp. $(X, \iota, \lambda, \mc{F})$) where $(A, \iota, \lambda)$ is a Hermitian abelian scheme (\cref{definition:ab_var:integral_models:Hermitian_abelian_scheme}) (resp. $(X, \iota, \lambda)$ is a Hermitian $p$-divisible group (\crefext{I:definition:moduli_pDiv:RZ:Hermitian_p-divisible}, but we allow $p = 2$ even if $F / \Q_p$ is ramified)) and $(A, \iota, \mc{F})$ is a Kr\"amer datum (resp. $(X, \iota, \mc{F})$ is a Kr\"amer datum).
\end{enumerate}
\end{definition}

The name ``Kr\"amer datum'' refers to the Kr\"amer model mentioned in Remark \ref{remark:ab_var:integral_models:other_models}.
For an understood Kr\"amer datum $(A, \iota, \lambda, \mc{F})$, we will use the shorthand $\ms{E}^{\vee} \coloneqq (\Lie A)/\mc{F}$ (cf. the ``tautological bundles'' of Definition \ref{definition:ab_var:integral_models:tautological_bundle} and Definition \ref{definition:ab_var:exotic_smooth:tautological_bundle}). We use the same notation $\ms{E}^{\vee} \coloneqq (\Lie X) / \mc{F}$ given an understood local Kr\"amer datum $(X, \iota, \mc{F})$. In both cases, the sheaf $\ms{E}^{\vee}$ is locally free of rank $1$, and we call it the associated \emph{Kr\"amer hyperplane quotient}.

\begin{definition}
A \emph{morphism} (resp. \emph{isogeny}) of Kr\"amer data $(A_1, \iota_1, \mc{F}_1) \ra (A_2, \iota_2, \mc{F}_2)$ is an $\mc{O}_F$-homomorphism (resp. isogeny) $A_1 \ra A_2$ such that $\mrm{im}(\mc{F}_1) \subseteq \mc{F}_2$, where $\mrm{im}(\mc{F}_1)$ is the image of $\mc{F}_1$ under $\Lie A_1 \ra \Lie A_2$. A \emph{morphism} (resp. \emph{isogeny}) of local Kr\"amer data is defined in the same way.
\end{definition}

\begin{lemma}\label{lemma:Faltings_and_taut:height:canonical_Kramer_data}
Let $S$ be a scheme over $\Spec \mc{O}_F$. Assume either that $S$ is a scheme over $\Spec \mc{O}_F[1/\Delta]$ or that $S = \Spec R$ where $R$ is a Dedekind domain with fraction field of characteristic $0$.
\begin{enumerate}[(1)]
    \item Suppose $A \ra S$ is an abelian scheme with an action $\iota \colon \mc{O}_F \ra \End(X)$ of signature $(n - 1, 1)$. Then the pair $(A, \iota)$ extends uniquely to a Kr\"amer datum $(A, \iota, \mc{F})$ over $S$.
    \item Given pairs $(A_1, \iota_1)$ and $(A_2, \iota_2)$ as above, any $\mc{O}_F$-linear homomorphism (resp. isogeny) $A_1 \ra A_2$ induces a morphism (resp. isogeny) of Kr\"amer data.
    \item If $S$ is a scheme over $\Spec \mc{O}_F[1 / \Delta]$, the exact sequence
    \begin{equation}
    0 \ra \mc{F} \ra \Lie A \ra \ms{E}^{\vee} \ra 0
    \end{equation}
has a unique $\mc{O}_F$-linear splitting.
\end{enumerate}
\end{lemma}
\begin{proof}
If $S$ is a scheme over $\Spec \mc{O}_F[1/\Delta]$, the claims hold because there is a unique decomposition $\Lie A = (\Lie A)^+ \oplus (\Lie A)^-$ characterized by $\iota$ acting $\mc{O}_F$-linearly on the rank $n - 1$ subbundle $(\Lie A)^+$ (resp. $\s$-linearly on the rank $1$ bundle $(\Lie A)^-$).

Suppose instead that $S = \Spec R$ is a Dedekind domain with fraction field $K$ of characteristic $0$. By localizing, it is enough to verify the lemma when $R$ is a discrete valuation ring. Then part (1) amounts to the following fact: given a finite free $R$-module $M$ and any $K$-subspace $W \subseteq M \otimes K$, there is a unique summand $M' \subseteq M$ such that $W = M' \otimes K$ (namely $M' = M \cap W$; note that $M' \subseteq M$ is a saturated sublattice). We are applying this when $M = \Lie A$ and $W = (\Lie A \otimes K)^+$, in the notation above (and taking $\mc{F} = M'$). The signature $(n - 1, 1)$ condition forces the $\mc{O}_F$-action on $(\Lie A) / \mc{F}$ to be $\s$-linear. 
These considerations also verify the claim in part (2) (since it holds in the generic fiber).
\end{proof}

\begin{lemma}\label{lemma:Faltings_and_taut:height:canonical_local_Kramer_data}
Let $S$ be a formal scheme over $\Spf \mc{O}_{F_p}$. Assume either that $p$ is unramified in $\mc{O}_F$ or that $S = \Spf R$ for an adic ring which is a Dedekind domain with fraction field of characteristic $0$. Then the following conclusions hold.
    \begin{enumerate}[(1)]
    \item Suppose $X$ is a $p$-divisible group of height $2n$ over $S$ with an action $\iota \colon \mc{O}_{F_p} \ra \End(X)$ of signature $(n - 1, 1)$. Then the pair $(X, \iota)$ extends uniquely to a local Kr\"amer datum $(X, \iota, \mc{F})$ over $S$.
    \item Given pairs $(X_1, \iota_1)$ and $(X_2, \iota_2)$ as above, any $\mc{O}_{F_p}$-linear homomorphism (resp. isogeny) $X_1 \ra X_2$ induces a morphism (resp. isogeny) of local Kr\"amer data.
    \item If $p$ is unramified, the exact sequence
        \begin{equation}
        0 \ra \mc{F} \ra \Lie A \ra \ms{E}^{\vee} \ra 0
        \end{equation}
    has a unique $\mc{O}_{F_p}$-linear splitting.
    \end{enumerate}
\end{lemma}
\begin{proof}
This may be proved in the same way as Lemma \ref{lemma:Faltings_and_taut:height:canonical_Kramer_data}.
If $S = \Spf R$ for $R$ a Dedekind domain with fraction field of characteristic $0$, note that $R$ must be a complete discrete valuation ring.
\end{proof}

In the situations of Lemma \ref{lemma:Faltings_and_taut:height:canonical_Kramer_data} and \ref{lemma:Faltings_and_taut:height:canonical_Kramer_data}, we also use the alternative terminology \emph{dual tautological bundle} for the Kr\"amer hyperplane quotient $\ms{E}^{\vee}$.

If $(A, \iota, \lambda)$ is a Hermitian abelian scheme of signature $(n - 1, 1)$ over $\Spec \mc{O}_E$ with associated quasi-polarized Kr\"amer datum $(A, \iota, \lambda, \mc{F})$, we thus obtain a Hermitian line bundle $\widehat{\ms{E}^{\vee}} = (\ms{E}^{\vee}, \norm{-})$ on $\Spec \mc{O}_E$ as follows: the metric $\norm{-}$ is given by restricting the metric on $\Lie A$ induced by $\lambda$ (which we take to be normalized as in \eqref{equation:arith_cycle_classes:Hodge_bundles:4_pi_gamma_normalization}) along the $\mc{O}_F$-linear splitting $\ms{E}^{\vee}[1/\Delta] \hookrightarrow (\Lie A)[1/\Delta]$ (where $(-)[1/\Delta]$ means restriction to $\Spec \mc{O}_E[1/\Delta]$). We say $\widehat{\ms{E}}^{\vee}$ is the associated \emph{metrized dual tautological bundle}. We also make the same construction over $\Spec \mc{O}_{E,(p)}$ and $\Spec E$.

\begin{definition}
Let $(A, \iota, \lambda)$ be a Hermitian abelian scheme of signature $(n - 1, 1)$ over $\Spec \mc{O}_E$. The associated \emph{tautological height} is
    \begin{equation}
    h_{\mrm{tau}}(A_E) \coloneqq h_{\mrm{tau}}(A) \coloneqq \frac{1}{[E : \Q]} \widehat{\deg}(\widehat{\ms{E}}^{\vee}).
    \end{equation}
\end{definition}

The tautological height depends on the auxiliary data in the definition (not just $A_E$ or $A$), which we have suppressed from notation. If $(A, \iota, \lambda)$ is a Hermitian abelian scheme of signature $(n - 1, 1)$ over $\Spec E$ such that $A$ has everywhere potentially good reduction, we define the \emph{tautological height} $h_{\mrm{tau}}(A)$ so that it is invariant under finite degree field extension $E \ra E'$.

\begin{remark}\label{remark:Faltings_and_taut:heights:punctured_heights}
If we instead work over $\mc{O}_E[1/N]$ for some integer $N \geq 1$, we may define \emph{Faltings height} and \emph{tautological height} as above, but where $\widehat{\deg}$ now takes values in $\R_N \coloneqq \R / \sum_{p \mid N} \Q \cdot \log p$ (as explained in Section \ref{ssec:arith_cycle_classes:Hermitian_bundles}).
\end{remark}

        \subsection{Change along global isogenies}
        \label{ssec:Faltings_and_taut:isogeny_change_global}
            Let $A_1 \ra \Spec \mc{O}_E$ and $A_2 \ra \Spec \mc{O}_E$ be semi-abelian N\'eron models of abelian varieties over $E$. We have
    \begin{equation}
    [E : \Q] (h_{\mrm{Fal}}(A_2) - h_{\mrm{Fal}}(A_1)) = \widehat{\deg}(\widehat{\omega}_{A_2}) - \widehat{\deg}(\widehat{\omega}_{A_1}) = -\widehat{\deg}(\underline{\Hom}(\widehat{\omega}_{A_2}, \widehat{\omega}_{A_1})).
    \end{equation}
Any isogeny $\phi \colon A_1 \ra A_2$ defines a section $\phi$ of the Hermitian line bundle $\underline{\Hom}(\widehat{\omega}_{A_2}, \widehat{\omega}_{A_1})$, which gives
    \begin{align}\label{equation:Faltings_and_taut:isogeny_change_global:Faltings}
    h_{\mrm{Fal}}(A_2) - h_{\mrm{Fal}}(A_1) & = \frac{1}{[E : \Q]} \left ( \log \norm{\phi}_{\infty} + \sum_{v < \infty} \log \norm{\phi}_v \right )
    \\
    & = \frac{1}{2} \log(\deg \phi) - \frac{1}{[E : \Q]} \log|e^* \Omega^1_{\ker \phi / \mc{O}_E}|
    \end{align}
(sum is over places $v$ of $E$) as in \cite[{Lemma 5}]{Faltings86}, where $|e^* \Omega^1_{\ker \phi / \mc{O}_E}|$ denotes the cardinality of the finite length $\mc{O}_E$-module $e^* \Omega^1_{\ker \phi / \mc{O}_E}$. Note $|e^*\Omega^1_{\ker \phi / \mc{O}_E}| = |\coker(\phi^* \colon \omega_{A_2} \ra \omega_{A_1})| = |\coker(\phi_* \colon \Lie A_1 \ra \Lie A_2)|$. Also note
    \begin{equation}\label{equation:Faltings_and_taut:isogeny_change_global:Faltings:log_p}
    h_{\mrm{Fal}}(A_2) - h_{\mrm{Fal}}(A_1) = \sum_p a_p \log p = \sum_{p \mid \deg \phi} a_p \log p
    \end{equation}
for some $a_p \in \Q$ independent of $\phi$.

Given Hermitian abelian schemes $(A_1, \iota_1, \lambda_1)$ and $(A_2, \iota_2, \lambda_2)$ of signature $(n - 1, 1)$ over $\Spec \mc{O}_E$ with associated Hermitian line bundles $\widehat{\ms{E}}^{\vee}_1$ and $\widehat{\ms{E}}^{\vee}_2$, we similarly have
    \begin{equation}
    h_{\mrm{tau}}(A_2) - h_{\mrm{tau}}(A_1) = \frac{1}{[E : \Q]} \widehat{\deg}(\underline{\Hom}(\widehat{\ms{E}}_{1}^{\vee}, \widehat{\ms{E}}_{2}^{\vee})).
    \end{equation}
Any $\mc{O}_F$-linear isogeny $\phi \colon A_1 \ra A_2$ defines a section $\phi$ of the Hermitian line bundle $\underline{\Hom}(\widehat{\ms{E}}_{1}^{\vee}, \widehat{\ms{E}}_{2}^{\vee})$, and we have
    \begin{align}\label{equation:Faltings_and_taut:isogeny_change_global:taut}
    h_{\mrm{tau}}(A_2)- h_{\mrm{tau}}(A_1) & = \frac{1}{[E : \Q]} \left ( - \log \norm{\phi}_{\infty} - \sum_{v < \infty} \log \norm{\phi}_v \right ) 
    \\
    & = \frac{1}{[E : \Q]} \left ( - \log \norm{\phi}_{\infty} +  \log |\coker(\phi_* \colon \ms{E}_1^{\vee} \ra \ms{E}_2^{\vee})| \right ).
    \end{align}

        \subsection{Change along local isogenies: Faltings}
        \label{ssec:Faltings_and_taut:isogeny_change_local_Faltings}

Given an isogeny $\phi \colon A_1 \ra A_2$ of abelian schemes over $\Spec \mc{O}_{E,(p)}$, we define the \emph{semi-global change of Faltings height}
    \begin{equation}
    \delta_{\mrm{Fal},(p)}(\phi) \coloneqq - \frac{1}{2}\log|\deg \phi|_p - \frac{1}{[E : \Q]} \log|e^* \Omega^1_{\ker \phi / \mc{O}_{E,(p)}}|_p
    \end{equation}
where $|-|_p$ is the usual $p$-adic norm. We have $\delta_{\mrm{Fal},(p)}(\phi) \in \Q \cdot \log p$. The formula for change of Faltings height \eqref{equation:Faltings_and_taut:isogeny_change_global:Faltings} shows that $\delta_{\mrm{Fal},(p)}(\phi) = a_p \log p$, in the notation of \eqref{equation:Faltings_and_taut:isogeny_change_global:Faltings:log_p}.
In particular, $\delta_{\mrm{Fal},(p)}(\phi)$ does not depend on the choice of isogeny $\phi$ (and depends only on $A_1$ and $A_2$).
If $A_1$ and $A_2$ have everywhere potentially good reduction, we have
    \begin{equation}\label{equation:Faltings_and_taut:isogeny_change_local_Faltings:decomp_to_semi-global}
    h_{\mrm{Fal}}(A_{2,E}) - h_{\mrm{Fal}}(A_{1,E}) = \sum_{\ell} \delta_{\mrm{Fal},(\ell)}(\phi) = \sum_{\ell \mid \deg \phi} \delta_{\mrm{Fal},(\ell)}(\phi)
    \end{equation}
where $\phi$ also denotes the induced isogeny on N\'eron models over $\Spec \mc{O}_{E,(\ell)}$ for each prime $\ell$ (after enlarging $E$ if necessary).

Given any isogeny $\phi \colon X_1 \ra X_2$ of $p$-divisible groups over $\Spf \mc{O}_{\breve{E}}$, we have $(\Lie X_i)^{\vee} \cong e^* \Omega^1_{X_i[p^N]/\Spec \mc{O}_{\breve{E}}}$ (canonically) for $N \gg 0$ by \cite[{Corollary II.3.3.17}]{Messing72} (passing to the limit over $\mc{O}_{\breve{E}} / p^k \mc{O}_{\breve{E}}$ as $k \ra \infty$), so there is a canonical exact sequence
    \begin{equation}\label{equation:Faltings_and_taut:isogeny_change_local:coker_and_differentials}
    0 \ra (\Lie X_2)^{\vee} \xra{\phi^*} (\Lie X_1)^{\vee} \ra e^* \Omega^1_{\ker \phi / \Spec \mc{O}_{\breve{E}}} \ra 0
    \end{equation}
of finite free $\mc{O}_{\breve{E}}$-modules (note that $\Lie X_1 \ra \Lie X_2$ is injective, e.g. by Lemma \ref{lemma:isogeny_criterion}). If $X_1$ and $X_2$ are moreover height $2n$ and dimension $n$, we define the \emph{local change of Faltings height}
    \begin{equation}
    \breve{\delta}_{\mrm{Fal}}(\phi) \coloneqq \frac{1}{2} \log(\deg \phi) - \frac{1}{[\breve{E} : \breve{\Q}_p]} \mrm{length}_{\mc{O}_{\breve{E}}}(e^* \Omega^1_{\ker \phi / \Spec \mc{O}_{\breve{E}}}) \cdot \log p.
    \end{equation}
We have $\breve{\delta}_{\mrm{Fal}}(\phi) = \Q \cdot \log p$, as well as
    \begin{align}
    & \breve{\delta}_{\mrm{Fal}}(\phi' \circ \phi) = \breve{\delta}_{\mrm{Fal}}(\phi') + \breve{\delta}_{\mrm{Fal}}(\phi) && \breve{\delta}_{\mrm{Fal}}([N]) = 0
    \end{align}
where $\phi' \colon X_2 \ra X_3$ is any isogeny of $p$-divisible groups and $[N] \colon X_1 \ra X_1$ is the multiplication-by-$N$ isogeny (follows from \eqref{equation:Faltings_and_taut:isogeny_change_local:coker_and_differentials}). 
Unlike $\delta_{\mrm{Fal},(p)}(-)$ from above, the quantity $\breve{\delta}_{\mrm{Fal}}(\phi)$ may depend on the isogeny $\phi$.

Given isogenous abelian schemes over $\Spec \mc{O}_{E,(p)}$ and an isogeny $\phi \colon X_1 \ra X_2$ of the associated $p$-divisible groups, set
    \begin{equation}\label{equation:Faltings_and_taut:isogeny_change_local_Faltings:pDiv_semi-global_def}
    \delta_{\mrm{Fal},(p)}(\phi) \coloneqq \frac{1}{[E : \Q]} \sum_{\breve{w} \mid p} [\breve{E}_{\breve{w}} : \breve{\Q}_p] \breve{\delta}_{\mrm{Fal}}(\phi_{\breve{w}})
    \end{equation}
where $\phi_{\breve{w}}$ denotes the base-change of $\phi$ to $\Spf \mc{O}_{\breve{E}_{\breve{w}}}$.

\begin{lemma}\label{lemma:Faltings_and_taut:isogeny_change_local_Faltings:global-to-local}
Let $A_1, A_2$ be isogenous abelian schemes over $\Spec \mc{O}_{E,(p)}$. Let $X_i$ be the associated $p$-divisible groups. Given any isogenies $\tilde{\phi} \colon A_1 \ra A_2$ and $\phi \colon X_1 \ra X_2$, we have
    \begin{equation}\label{equation:Faltings_and_taut:isogeny_change_local:global-to-local}
    \delta_{\mrm{Fal},(p)}(\tilde{\phi}) = \delta_{\mrm{Fal},(p)}(\phi).
    \end{equation}
\end{lemma}
\begin{proof}
The lemma is clear if $\phi$ is the isogeny associated with $\tilde{\phi}$. If $\phi' \colon X_1 \ra X_2$ is another isogeny, we have $[p^N] \circ \phi = \phi' \circ \phi''$ for some isogeny $\phi'' \colon X_1 \ra X_1$ (Lemma \ref{lemma:isogeny_criterion}). By additivity of $\breve{\delta}_{\mrm{Fal}}$ and since $\breve{\delta}_{\mrm{Fal}}([p^N]) = 0$, it is enough to show $\delta_{\mrm{Fal},(p)}(\phi) = 0$ if $X_1 = X_2$. For this purpose, we may also assume $A_1 = A_2$.

Write $A \coloneqq A_1$ and $X \coloneq X_1$ to lighten notation. As usual, $A_E$ and $X_E$ denote the respective generic fibers (over $\Spec E$). We write $\mrm{Isog}(A)$ and $\mrm{Isog}(X)$ for the set of self-isogenies of $A$ and $X$.

We have canonical identifications
    \begin{equation}
    \End(X) = \End(X_E) = \End(T_p(X_E)).
    \end{equation}
The first equality holds by a theorem of Tate \cite[{Theorem 4}]{Tate67} (base-change along $\Spec E \ra \Spec \mc{O}_{E,(p)}$ is fully faithful on $p$-divisible groups) and the second equality holds because $X_E$ is an \'etale $p$-divisible group. Here, the notation $\End(T_p(X_E))$ means endomorphisms of $T_p(X_E)$ as a Galois module.

Equip the finite $\Z_p$-module $\End(T_p(X_E))$ with the $p$-adic topology, and give $\mrm{Isog}(X)$ the subspace topology.
We have $\delta_{\mrm{Fal},(p)}(\phi \circ \phi') = \delta_{\mrm{Fal},(p)}(\phi)$ for any $\phi \in \mrm{Isog}(X)$ and $\phi' \in \End(X)$ with $\phi' \equiv 1 \pmod{p}$, since any such $\phi'$ is an automorphism of $X$. The map $\mrm{Isog}(X) \ra \R$ given by $\phi \mapsto \delta_{\mrm{Fal},(p)}(\phi)$ is thus locally constant.

We also have canonical identifications
    \begin{equation}
    \End(A) \otimes_{\Z} \Z_p = \End(A_E) \otimes_{\Z} \Z_p = \End(T_p(A_E)) = \End(T_p(X_E)).
    \end{equation}
The first equality holds by the N\'eron mapping property, and the second equality holds by Faltings's theorem \cite[{\S 5 Corollary 1}]{Faltings86}. Hence $\mrm{Isog}(A) = \mrm{Isog}(X) \cap \End(A)$ is a dense subset of $\mrm{Isog}(X)$. Since $\delta_{\mrm{Fal},(p)}(\phi) = 0$ for any $\phi \in \mrm{Isog}(A)$, this proves the lemma.
\end{proof}

\begin{corollary}
In the situation of Lemma \ref{lemma:Faltings_and_taut:isogeny_change_local_Faltings:global-to-local}, the quantity $\delta_{\mrm{Fal},(p)}(\phi)$ does not depend on the choice of isogeny $\phi \colon X_1 \ra X_2$.
\end{corollary}
\begin{proof}
In the notation of the lemma, this follows immediately from $\tilde{\phi}$ independence of $\delta_{\mrm{Fal},(p)}(\tilde{\phi})$ (discussed above).
\end{proof}

We will use Lemma \ref{lemma:Faltings_and_taut:isogeny_change_local_Faltings:global-to-local} to compute Faltings heights without producing isogenies on abelian varieties, only isogenies on underlying $p$-divisible groups over $\Spec \mc{O}_{E,(p)}$. The analogous lemma for tautological height (Lemma \ref{lemma:Faltings_and_taut:isogeny_change_local_taut:global-to-local}) serves a similar purpose.

        \subsection{Change along local isogenies: tautological}
        \label{ssec:Faltings_and_taut:isogeny_change_local_taut}
            To locally decompose the change of tautological height along an isogeny, we impose an additional condition.

\begin{definition}\label{definition:Faltings_and_taut:isogeny_change_local_taut:special}
\hfill
\begin{enumerate}[(1)]
    \item A Hermitian abelian scheme $(A, \iota, \lambda)$ of signature $(n - 1, 1)$ over $E$ is \emph{special} if $A$ is $\mc{O}_F$-linearly isogenous to a product of elliptic curves, each with $\mc{O}_F$-action. A Hermitian abelian scheme of signature $(n - 1, 1)$ over $\Spec \mc{O}_E$ or $\Spec \mc{O}_{E,(p)}$ is \emph{special} if its generic fiber is special.
    \item A Hermitian $p$-divisible group $(X, \iota, \lambda)$ of signature $(n - 1, 1)$ over $\Spf \mc{O}_{\breve{E}}$ is \emph{special} if $X$ is $\mc{O}_{F_p}$-linearly isogenous to $\mf{X}_0^{n - 1} \times \mf{X}_0^{\s}$.
\end{enumerate}
\end{definition}

We only use the term ``special'' this way in Part \ref{part:local_change_heights} (but we have global special cycles in mind, cf. Lemma \ref{lemma:corank_1_isogeny_decomposition}). The norm $\norm{-}_{\infty}$ below is as in \eqref{equation:Faltings_and_taut:isogeny_change_global:taut}.

\begin{lemma}\label{lemma:Faltings_and_taut:isogeny_change_local_taut:infty_norm}
Let $(A_1, \iota_1, \lambda_1)$ and $(A_2, \iota_2, \lambda_2)$ be special Hermitian abelian schemes of signature $(n - 1, 1)$ over $\Spec E$. For any $\mc{O}_F$-linear isogeny $\phi \colon A_1 \ra A_2$, we have $\norm{\phi}_{\infty}^2 \in \Q_{>0}$.
\end{lemma}
\begin{proof}
Given such $\phi$, form a diagram
    \begin{equation}
    B_1 \times B_1^{\perp} \xra{\phi_1} A_1 \xra{\phi} A_2 \xra{\phi_2} B_2 \times B_2^{\perp}
    \end{equation}
where each $\phi_i$ is an $\mc{O}_F$-linear isogeny, each $B_i$ is a product of $(n
- 1)$ elliptic curves each with $\mc{O}_F$-action of signature $(1,0)$, and each $B_i^{\perp}$ is an elliptic curve with $\mc{O}_F$-action of signature $(0,1)$. Signature incompatibility implies that $\lambda_1$ pulls back to a diagonal quasi-polarization $\lambda_{B_1} \times \lambda_{B_1^{\perp}}$ on $B_1 \times B_1^{\perp}$ (e.g. $\Hom_{\mc{O}_F}(B_1^{\s}, B_1^{\perp \vee}) = \Hom_{\mc{O}_F}(B_1^{\perp \s}, B_1^{\vee}) = 0$). Similarly, $\lambda_2$ pulls back along the quasi-isogeny $\phi_2^{-1}$ to a diagonal quasi-polarization $\lambda_{B_2} \times \lambda_{B_2^{\perp}}$.

With these quasi-polarizations, we have $\norm{\phi_2 \circ \phi \circ \phi_1}_{\infty} = \norm{\phi_2}_{\infty} \norm{\phi}_{\infty} \norm{\phi_1}_{\infty} = \norm{\phi}_{\infty}$ since $\norm{\phi_1}_{\infty} = \norm{\phi_2}_{\infty} = 1$ (because $\phi_1$ and $\phi_2$ preserve quasi-polarizations, by construction).
On the other hand, if $\phi' \colon B_1^{\perp} \ra B_2^{\perp}$ is the induced isogeny (signature incompatibility again implies $\Hom_{\mc{O}_F}(B_1, B_2^{\perp}) = \Hom_{\mc{O}_F}(B_1^{\perp}, B_2) = 0$), we must have $\norm{\phi_2 \circ \phi \circ \phi_1}_{\infty} = \norm{\phi'}_{\infty}$ (the latter norm is taken with respect to $\lambda_{B_1^{\perp}}$ and $\lambda_{B_2^{\perp}}$). For each embedding $\t \colon E \ra \C$, the quantity $\norm{\phi'}^2_{\t}$ must be the element of $\Q_{>0}$ satisfying
    \begin{equation}
    \phi^{\prime *} \lambda_{B_2^{\perp}} = \norm{\phi'}_{\t}^2 \lambda_{B_1^{\perp}},
    \end{equation}
(quasi-polarizations on elliptic curves are unique up to $\Q_{>0}$ scalar), so we have $\norm{\phi'}_{\infty}^2 = \prod_{\t \colon E \ra \C} \norm{\phi'}_{\t}^2 \in \Q_{>0}$.
\end{proof}

For the rest of Section \ref{ssec:Faltings_and_taut:isogeny_change_local_taut}, we let $(A_i, \iota_i, \lambda_i)$ for $i = 1, 2$ be special Hermitian abelian schemes of signature $(n - 1, 1)$ over $\Spec \mc{O}_{E,(p)}$, with associated Kr\"amer hyperplanes $\mc{F}_i$ and dual tautological bundles $\ms{E}^{\vee}_i$. We also let $(X_i, \iota_i, \lambda_i)$ for $i = 1, 2, 3$ be special Hermitian $p$-divisible groups of signature $(n - 1, 1)$ over $\Spf \mc{O}_{\breve{E}}$, and reuse the notation $\mc{F}_i$ and $\ms{E}^{\vee}_i$ for the respective Kr\"amer hyperplanes and dual tautological bundles.

Given $(X_1, \iota_1, \lambda_1)$ and an $\mc{O}_{F_p}$-linear isogeny $Y_1 \times Y_1^{\perp} \ra X_1$ with $Y_1$ being a product of $n - 1$ canonical liftings of signature $(1,0)$ and $Y_1^{\perp}$ being a canonical lifting of signature $(0,1)$, there is an induced decomposition
    \begin{equation}\label{equation:Faltings_and_taug:isogeny_change_local_taut:T_p-decomp}
    T_p(X_1)^0 = T_p(Y_1)^0 \oplus T_p(Y_1^{\perp})^0
    \end{equation}
on rational Tate modules (of the generic fibers). Equip $Y_1 \times Y_1^{\perp}$ with the pullback of $\lambda_1$. This gives a product quasi-polarization $\lambda_{Y_1} \times \lambda_{Y_1^{\perp}}$ on $Y_1 \times Y_1^{\perp}$ (by signature incompatibility as in the abelian scheme case, i.e. $\Hom_{\mc{O}_{F_p}}(Y_1^{\s}, Y_1^{\perp \vee}) = \Hom_{\mc{O}_{F_p}}(Y_1^{\perp \s}, Y_1^{\vee}) = 0$).
Hence the decomposition in \eqref{equation:Faltings_and_taug:isogeny_change_local_taut:T_p-decomp} is orthogonal for the Hermitian pairing on $T_p(X_1)^0$.

Consider $(X_2, \iota_2, \lambda_2)$ with $\mc{O}_{F_p}$-linear isogeny $Y_2 \times Y_2^{\perp} \ra X_2$ as above, and suppose $\phi \colon X_1 \ra X_2$ is an $\mc{O}_{F_p}$-linear isogeny. Then the induced map $\phi_* \colon T_p(X_1)^0 \ra T_p(X_2)^0$ sends $T_p(Y_1)^0$ to $T_p(Y_2)^0$ and similarly for $T_p(Y_i^{\perp})^0$ (again by signature incompatibility, i.e. $\Hom_{\mc{O}_{F_p}}(Y_1, Y_2^{\perp}) = \Hom_{\mc{O}_{F_p}}(Y_1^{\perp}, Y_2) = 0$). In particular, the decomposition in \eqref{equation:Faltings_and_taug:isogeny_change_local_taut:T_p-decomp} does not depend on the choice of $Y_1 \times Y_1^{\perp} \ra X_1$.

Any $\mc{O}_{F_p}$-linear isogeny $\phi \colon X_1 \ra X_2$ thus gives a nonzero element $\phi \in \Hom_{F_p}(T_p(Y_1^{\perp})^0, T_p(Y_2^{\perp})^0)$. We then set
    \begin{equation}\label{equation:Faltings_and_taut:isogeny_change_local_taut:infinity_p_norm}
    \norm{\phi}_{\infty,p} \coloneqq \norm{\phi}
    \end{equation}
where $\norm{-}$ on the right means the norm for the (one-dimensional and non-degenerate) $F_p$-Hermitian space $\Hom_{F_p}(T_p(Y_1^{\perp})^0, T_p(Y_2^{\perp})^0)$.

We may now proceed as in the Faltings height case.
Given an $\mc{O}_F$-linear isogeny $\phi \colon A_1 \ra A_2$, we define the \emph{semi-global change of tautological height}
    \begin{equation}
    \delta_{\mrm{tau},(p)}(\phi) \coloneqq \frac{1}{[E : \Q]} \left (\log | \norm{\phi}_{\infty} |_p +  \log|\coker(\phi_* \colon \ms{E}_1^{\vee} \ra \ms{E}_2^{\vee})| \right )
    \end{equation}
where $|-|_p$ is the usual $p$-adic norm (well-defined by Lemma \ref{lemma:Faltings_and_taut:isogeny_change_local_taut:infty_norm}).
We have $\delta_{\mrm{tau},(p)}(\phi) \in \Q \cdot \log p$. Since $A_1$ and $A_2$ have everywhere potentially good reduction (implied by the special hypothesis: elliptic curves with $\mc{O}_F$-action over number fields have everywhere potentially good reduction)
the formula for change of tautological height \eqref{equation:Faltings_and_taut:isogeny_change_global:taut} implies
    \begin{equation}\label{equation:Faltings_and_taut:isogeny_change_local_taut:decomp_to_semi-global}
    h_{\mrm{tau}}(A_{2,E}) - h_{\mrm{tau}}(A_{1,E}) = \sum_{\ell} \delta_{\mrm{tau},(\ell)}(\phi) = \sum_{\ell \mid \deg \phi} \delta_{\mrm{tau},(\ell)}(\phi)
    \end{equation}
where $\phi$ also denotes the induced isogeny on N\'eron models over $\Spec \mc{O}_{E,(\ell)}$ for each prime $\ell$ (after enlarging $E$ if necessary). In particular, $\delta_{\mrm{tau},(p)}(\phi)$ does not depend on the choice of isogeny $\phi$.

Given any $\mc{O}_{F_p}$-linear isogeny $\phi \colon X_1 \ra X_2$, we define the \emph{local change of tautological height}
    \begin{equation}
    \breve{\delta}_{\mrm{tau}}(\phi) \coloneqq \log \norm{\phi}_{\infty, p} + \frac{1}{[\breve{E} : \breve{\Q}_p]} \mrm{length}_{\mc{O}_{\breve{E}}}(\coker(\phi_* \colon \ms{E}_1^{\vee} \ra \ms{E}_2^{\vee})) \cdot \log p.
    \end{equation}
We have $\breve{\delta}_{\mrm{tau}}(\phi) \in \Q \cdot \log p$, as well as
    \begin{align}\label{equation:isogeny_change_local:taut:linear_relations}
    & \breve{\delta}_{\mrm{tau}}(\phi' \circ \phi) = \breve{\delta}_{\mrm{tau}}(\phi') + \breve{\delta}_{\mrm{tau}}(\phi)
    && \breve{\delta}_{\mrm{tau}}([N]) = 0
    \end{align}
where $\phi' \colon X_2 \ra X_2$ is any $\mc{O}_{F_p}$-linear isogeny and $[N] \colon X_1 \ra X_1$ is the multiplication-by-$N$ isogeny. 
For use in later calculations, we note the identity
    \begin{align}\label{equation:Faltings_and_taut:isogeny_change_local_taut:length_additivity}
    & \operatorname{length}_{\mc{O}_{\breve{E}}}(\coker(\phi_* \colon \Lie(X_1) \ra \Lie(X_2)))
    \\
    & = \operatorname{length}_{\mc{O}_{\breve{E}}}(\coker(\phi_* \colon \mc{F}_1 \ra \mc{F}_2)) + \operatorname{length}_{\mc{O}_{\breve{E}}}(\coker(\phi_* \colon \ms{E}^{\vee}_1 \ra \ms{E}^{\vee}_2)) \notag
    \end{align}
(by the snake lemma).

\begin{lemma}\label{lemma:isogeny_change_local_taut:nonsplit_independence}
If $F_p / \Q_p$ is nonsplit, we have $\breve{\delta}_{\mrm{tau}}(\phi) = \breve{\delta}_{\mrm{tau}}(\phi')$ for any two $\mc{O}_F$-linear isogenies $\phi, \phi' \colon X_1 \ra X_2$.
\end{lemma}
\begin{proof}
Set $X = \mf{X}_0^{n - 1} \times \mf{X}_0^{\s}$, and equip $X$ with any $\mc{O}_F$-action-compatible quasi-polarization. Select any $\mc{O}_F$-linear isogeny $\phi'' \colon X \ra X_1$. Using the additivity property $\breve{\delta}_{\mrm{tau}}(\phi \circ \phi'') = \breve{\delta}_{\mrm{tau}}(\phi) + \breve{\delta}_{\mrm{tau}}(\phi'')$ and similarly for $\phi'$, this reduces us to the case where $X = X_1$.

As in the proof of Lemma \ref{lemma:Faltings_and_taut:isogeny_change_local_Faltings:global-to-local}, there exists an isogeny $\phi'' \colon X \ra X$ such that $[p^N] \circ \phi = \phi' \circ \phi''$ for some $N \geq 0$, so the additivity properties of $\breve{\delta}_{\mrm{tau}}$ reduce us to showing $\breve{\delta}_{\mrm{tau}}(\phi) = 0$ when $(X_1, \iota_1, \lambda_1) = (X_2, \iota_2, \lambda_2)$.

Since $\Hom_{\mc{O}_F}(\mf{X}_0, \mf{X}_0^{\s}) = \Hom_{\mc{O}_F}(\mf{X}_0^{\s}, \mf{X}_0) = 0$, we must have $\phi = f \times f^{\perp}$ where $f \colon \mf{X}_0^{n - 1} \ra \mf{X}_0$ and $f^{\perp} \colon \mf{X}_0^{\s} \ra \mf{X}_0^{\s}$. We find $\breve{\delta}_{\mrm{tau}}(\phi) = \breve{\delta}_{\mrm{tau}}(f^{\perp}) = 0$ since $f^{\perp} \colon \mf{X}_0^{\s} \ra \mf{X}_0^{\s}$ is an automorphism times $[p^N]$ for some $N \geq 0$.
\end{proof}

\begin{remark}\label{remark:Faltings_and_taut:isogeny_change_local_taut:split_complication}
If $F_p / \Q_p$ is split, then Lemma \ref{lemma:isogeny_change_local_taut:nonsplit_independence} fails (consider multiplication by $(1,p)$ and $(p,1)$ in $\mc{O}_{F_p} \cong \Z_p \times \Z_p$). This is the reason for Lemma \ref{lemma:Faltings_and_taut:isogeny_change_local_taut:global-to-local} below, which allows us to uniformly treat all cases of $F_p / \Q_p$.
\end{remark}

Continuing to allow $F_p/\Q_p$ inert/ramified/split, now suppose that $(X_i, \iota_i, \lambda_i)$ is the Hermitian $p$-divisible group associated with $(A_i, \iota_i, \lambda_i)$, for $i = 1, 2$. Since each $(A_i, \iota_i, \lambda_i)$ is special, there automatically exists an $\mc{O}_F$-linear isogeny $A_1 \ra A_2$ after possibly replacing $E$ by a finite extension (by the theory of complex multiplication for elliptic curves). 
Given any $\mc{O}_{F_p}$-linear isogeny $\phi \colon X_1 \ra X_2$, set
    \begin{equation}\label{equation:Faltings_and_taut:isogeny_change_local_taut:pDiv_semi-global_def}
    \delta_{\mrm{tau},(p)}(\phi) \coloneqq \frac{1}{[E : \Q]} \sum_{\breve{w} \mid p} [\breve{E}_{\breve{w}} : \breve{\Q}_p] \breve{\delta}_{\mrm{tau}}(\phi_{\breve{w}})
    \end{equation}
where $\phi_{\breve{w}}$ denotes the base-change of $\phi$ to $\Spf \mc{O}_{\breve{E}_{\breve{w}}}$.

\begin{lemma}\label{lemma:Faltings_and_taut:isogeny_change_local_taut:global-to-local}
Suppose that $(X_i, \iota_i, \lambda_i)$ is the Hermitian $p$-divisible group associated with $(A_i, \iota_i, \lambda_i)$, for $i = 1, 2$.
For any $\mc{O}_F$-linear isogenies $\tilde{\phi} \colon A_1 \ra A_2$ and $\phi \colon X_1 \ra X_2$, we have
    \begin{equation}
    \delta_{\mrm{tau},(p)}(\tilde{\phi}) = \delta_{\mrm{tau}, (p)}(\phi).
    \end{equation}
\end{lemma}
\begin{proof}
This may be proved exactly as in Lemma \ref{lemma:Faltings_and_taut:isogeny_change_local_Faltings:global-to-local}, now requiring isogenies and endomorphisms to be $\mc{O}_F$-linear.
\end{proof}

\begin{corollary}
In the situation of Lemma \ref{lemma:Faltings_and_taut:isogeny_change_local_taut:global-to-local}, the quantity $\delta_{\mrm{tau},(p)}(\phi)$ does not depend on the choice of isogeny $\phi \colon X_1 \ra X_2$.
\end{corollary}
\begin{proof}
In the notation of the lemma, this follows immediately from $\tilde{\phi}$ independence of $\delta_{\mrm{tau},(p)}(\tilde{\phi})$ (discussed above).
\end{proof}

        \subsection{Serre tensor}
        \label{ssec:Faltings_and_taut:Serre_tensor}
            We compute local changes of Faltings and tautological heights for isogenies involving the Serre tensor $p$-divisible groups $\mf{X}_s \otimes_{\Z_p} \mc{O}_{F_p}$. These results will later be used to compute heights of arithmetic special $1$-cycles.

Given $s \in \Z_{\geq 0}$ and a quasi-canonical lifting $\mf{X}_s$ over $\Spf \mc{O}_{\breve{E}}$, we write $\lambda_{\mf{X}_s}$ for an understood principal polarization of $\mf{X}_s$. Recall that $\lambda_{\mf{X}_s}$ exists and is unique up to $\Z_p^{\times}$ scalar (see \crefext{I:lemma:moduli_pDiv:Serre_tensor:lift_polarization} and its proof). As in \crefext{I:ssec:moduli_pDiv:RZ}, we consider the map $\lambda_{\mrm{tr}} \colon \mc{O}_{F_p} \ra \mc{O}_{F_p}^{*}$ determined by the $\Z_p$-bilinear pairing $\mrm{tr}_{F_p / \Q_p}(x^{\s} y)$ on $\mc{O}_{F_p}$, where $\mc{O}_{F_p}^{*} \coloneqq \Hom_{\Z_p}(\mc{O}_{F_p}, \Z_p)$.

We equip $\mf{X}_s \otimes_{\Z_p} \mc{O}_{F_p}$ with its Serre tensor $\mc{O}_{F_p}$-action $\iota$ and the polarization $-\iota(\mf{d}_p^2)^{-1} \circ (\lambda_{\mf{X}_s} \otimes \lambda_{\mrm{tr}})$. We equip $\mf{X}_0 \times \mf{X}_0^{\s}$ with its diagonal $\mc{O}_{F_p}$ action $\iota_{\mf{X}_0} \times \iota_{\mf{X}_0}^{\s}$ (of signature $(1,1)$) and the diagonal quasi-polarization $-\iota(\mf{d}_p^2)^{-1} \circ (\lambda_{\mf{X}_0} \times \lambda_{\mf{X}_0})$.

\begin{lemma}\label{lemma:Faltings_and_taut:Serre_tensor:product}
For the $\mc{O}_{F_p}$-linear isogeny
    \begin{equation}
    \begin{tikzcd}[row sep = tiny]
    \mf{X}_0 \otimes_{\Z_p} \mc{O}_{F_p} \arrow{r}{\phi} & \mf{X}_0 \times \mf{X}_0^{\s}
    \\
    x \otimes a \arrow[mapsto]{r} & (\iota_{\mf{X}_0}(a) x, \iota_{\mf{X}_0}(a^{\s}) x)
    \end{tikzcd}
    \end{equation}
we have $\breve{\delta}_{\mrm{Fal}}(\phi) = 0$. Assuming $p \neq 2$ if $F_p / \Q_p$ is ramified, we also have $\breve{\delta}_{\mrm{tau}}(\phi) = 0$. 
\end{lemma}
\begin{proof}
We already know $\deg \phi = |\Delta|_p^{-1}$ (see \crefext{I:equation:moduli_pDiv:RZ:isogeny_from_tensor:polarizations} and surrounding discussion). 

Pick any $\mc{O}_{F_p}$-linear isomorphism $\Lie \mf{X}_0 \cong \mc{O}_{\breve{E}}$. Then the map $\phi_* \colon \Lie(\mf{X}_0 \otimes_{\Z_p} \mc{O}_{F_p}) \ra \Lie(\mf{X}_0 \times \mf{X}_0^{\s})$ may be identified with the map of $\mc{O}_{\breve{E}}$-modules $f \colon \mc{O}_{\breve{E}} \otimes_{\Z_p} \mc{O}_{F_p} \ra \mc{O}_{\breve{E}} \oplus \mc{O}_{\breve{E}}$ given by $f(x \otimes a) = (a x, a^{\s} x)$. Thus $\phi_*$ is given by the matrix in \crefext{I:equation:moduli_pDiv:RZ:isogeny_from_tensor:matrix} (the same matrix describing $\phi$ after identifying $\mf{X}_0 \otimes_{\Z_p} \mc{O}_{F_p} \cong \mf{X}_0^2$ using a $\Z_p$-basis of $\mc{O}_{F_p}$). That matrix has determinant which generates the different ideal $\mf{d}_p$, hence 
    \begin{equation}
    \operatorname{length}_{\mc{O}_{\breve{E}}}(\coker(\phi_* \Lie(\mf{X}_0 \otimes_{\Z_p} \mc{O}_{F_p}) \ra \Lie(\mf{X}_0 \times \mf{X}_0^{\s}))) = \frac{1}{2} [\breve{E} : \breve{\Q}_p] v_p(\Delta). 
    \end{equation}
This gives $2 \breve{\delta}_{\mrm{Fal}}(\phi) = \log \deg \phi - v_p(\Delta) \log p = 0$.

We also know that $\phi^*(\lambda_{\mf{X}_0} \times \lambda_{\mf{X}_0}) = \lambda_{\mf{X}_0} \otimes \lambda_{\mrm{tr}}$ (see discussion surrounding \crefext{I:equation:moduli_pDiv:RZ:isogeny_from_tensor:polarizations} again). Thus $\norm{\phi}_{\infty, p} = 1$, in the notation of \eqref{equation:Faltings_and_taut:isogeny_change_local_taut:infinity_p_norm}.

Let $\mc{F}_1 \subseteq \Lie (\mf{X}_0 \otimes_{\Z_p} \mc{O}_{F_p})$ and $\mc{F}_2 \subseteq \Lie(\mf{X}_0 \times \mf{X}_0^{\s})$ be the (unique) associated Kr\"amer hyperplanes, with associated Kr\"amer hyperplane quotients $\ms{E}^{\vee}_1 $ and $\ms{E}^{\vee}_2$. 
If $F_p / \Q_p$ is unramified, then $\phi$ is an isomorphism, hence $\coker(\phi_* \colon \ms{E}^{\vee}_1 \ra \ms{E}^{\vee}_2) = 0$. 
If $F_p / \Q_p$ is ramified, assume $p \neq 2$ and select a uniformizer $\varpi \in \mc{O}_{F_p}$ satisfying $\varpi^{\s} = -\varpi$. 
Then $(\varpi \otimes 1 + 1 \otimes \varpi) \in \mc{O}_{\breve{E}} \otimes_{\Z_p} \mc{O}_{F_p}$ is a generator of $\mc{F}_1$. 
We thus find $\coker(\phi_* \colon \mc{F}_1 \ra \mc{F}_2) \cong \mc{O}_{\breve{E}} / \varpi \mc{O}_{\breve{E}}$.

By \eqref{equation:Faltings_and_taut:isogeny_change_local_taut:length_additivity}, the previous computations imply $\coker(\phi_* \colon \ms{E}^{\vee}_1 \ra \ms{E}^{\vee}_2) = 0$, and hence
    \begin{equation}
    \breve{\delta}_{\mrm{tau}}(\phi) = \log \norm{\phi}_{\infty, p} + \frac{1}{[\breve{E} : \breve{\Q}_p]} \mrm{length}_{\mc{O}_{\breve{E}}}(\coker(\phi_* \colon \ms{E}^{\vee}_1 \ra \ms{E}^{\vee}_2)) \cdot \log p = 0.
    \end{equation}
\end{proof}

For $s \in \Z_{\geq 0}$, and given any isogeny $\psi_s \colon \mf{X}_0 \ra \mf{X}_s$ of minimal degree
(i.e. degree $p^s$), the discussion surrounding \crefext{I:equation:can_and_qcan:qcan:local_change_taut}
implies
    \begin{equation}
    \breve{\delta}_{\mrm{tau}}(\psi_s)
    =
    \delta_{\mrm{tau}}(s) \log p
    =
    -\frac{1}{2} \left ( s - \frac{(1 - p^{-s})(1 - \eta^*(p))}{(1 - p^{-1})(p - \eta^*(p))} \right ) \log p
    \end{equation}
where the constant $\delta_{\mrm{tau}}(s)$ was defined in
\crefext{I:equation:can_and_qcan:qcan:local_change_taut},
and $\eta^*(p) \coloneqq -1, 0, 1$ for $F_p / \Q_p$ in the inert, ramified, and split cases, respectively.
As in loc. cit., we also define $\delta_{\mrm{Fal}}(s) \coloneqq -2 \delta_{\mrm{tau}}(s)$.

\begin{lemma}\label{lemma:Faltings_and_taut:Serre_tensor:qcan_isog}
Let $\psi_s \colon \mf{X}_0 \ra \mf{X}_s$ be any isogeny of degree $p^s$. For the $\mc{O}_{F_p}$-linear isogeny
    \begin{equation}
    \begin{tikzcd}[row sep = tiny]
    \mf{X}_0 \otimes_{\Z_p} \mc{O}_{F_p} \arrow{r}{\phi} & \mf{X}_s \otimes_{\Z_p} \mc{O}_{F_p}
    \\
    x \otimes a \arrow[mapsto]{r} & \psi_s(x) \otimes a
    \end{tikzcd}
    \end{equation}
we have 
    \begin{align}
    & \breve{\delta}_{\mrm{Fal}}(\phi) = -2 \breve{\delta}_{\mrm{tau}}(\phi) = -2 \delta_{\mrm{tau}}(s) \cdot \log p.
    \end{align}
\end{lemma}
\begin{proof}
Recall that $\psi_s$ is unique up to pre-composition by elements of $\mc{O}_{F_p}^{\times}$ \crefext{I:equation:can_and_qcan:qcan:psi_s}. Write $\mc{F}_1$ and $\mc{F}_2$ (resp. $\ms{E}^{\vee}_1$ and $\ms{E}^{\vee}_2$) for the associated Kr\"amer hyperplanes (resp. dual tautological bundles) of $\mf{X}_0 \otimes_{\Z_p} \mc{O}_{F_p}$ and $\mf{X}_s \otimes_{\Z_p} \mc{O}_{F_p}$ respectively.

We have $\deg \phi = (\deg \psi_s)^2 = p^{2s}$. Since quasi-polarizations on $\mf{X}_0$ are unique up to $\Q_p^{\times}$ scalar (follows from Drinfeld rigidity and the corresponding statement for $\mbf{X}_0$ in \crefext{I:ssec:moduli_pDiv:RZ}), we have $\psi_s^* \lambda_{\mf{X}_s} = b \lambda_{\mf{X}_0}$ for some $b \in p^s \Z_p^{\times}$. Hence we have $\phi^*(\lambda_{\mf{X}_s} \otimes \lambda_{\mrm{tr}}) = b (\lambda_{\mf{X}_0} \otimes \lambda_{\mrm{tr}})$, so $\norm{\phi}_{\infty,p} = p^{-s/2}$.

Pick any identifications $\Lie \mf{X}_0 \cong \Lie \mf{X}_s \cong \mc{O}_{\breve{E}}$ of $\mc{O}_{\breve{E}}$-modules. With these identifications, the map $\psi_{s,*} \colon \Lie \mf{X}_0 \ra \Lie \mf{X}_s$ is multiplication by some $c \in \mc{O}_{\breve{E}}$ satisfying $[\breve{E} : \breve{\Q}_p] v_p(c) = \operatorname{length}_{\mc{O}_{\breve{E}}}(\coker (\psi_{s,*} \colon \Lie \mf{X}_0 \ra \Lie \mf{X}_s))$.

We also obtain identifications $\Lie (\mf{X}_0 \otimes_{\Z_p} \mc{O}_{F_p}) \cong \Lie (\mf{X}_s \otimes_{\Z_p} \mc{O}_{F_p})$ of $\mc{O}_{\breve{E}} \otimes_{\Z_p} \mc{O}_{F_p}$-modules, with induced identifications $\mc{F}_1 \cong \mc{F}_2$ and $\ms{E}^{\vee}_1 \cong \ms{E}_2^{\vee}$. Then $\phi_* \colon \Lie (\mf{X}_0 \otimes_{\Z_p} \mc{O}_{F_p}) \ra \Lie (\mf{X}_s \otimes_{\Z_p} \mc{O}_{F_p})$ is identified with multiplication by $c$, and hence $\phi_* \colon \ms{E}^{\vee}_1 \ra \ms{E}_2^{\vee}$ must also be multiplication by $c$.
Hence
    \begin{align}
    \operatorname{length}_{\mc{O}_{\breve{E}}}(\coker (\phi_* \colon \Lie(\mf{X}_0 \otimes_{\Z_p} \mc{O}_{F_p}) \ra \Lie(\mf{X}_s \otimes_{\Z_p} \mc{O}_{F_p}))) & = 2 [\breve{E} : \breve{\Q}_p] v_p(c)
    \\
    \operatorname{length}_{\mc{O}_{\breve{E}}}(\coker (\phi_* \colon \ms{E}^{\vee}_1 \ra \ms{E}^{\vee}_2)) & = [\breve{E} : \breve{\Q}_p] v_p(c).
    \end{align}
The lemma now follows from the formula for $\mrm{length}_{\mc{O}_{\breve{E}}}(\coker (\psi_{s,*} \colon \Lie \mf{X}_0 \ra \Lie \mf{X}_s))$ in \crefext{I:equation:can_and_qcan:qcan:Nakkajima--Taguchi}.
\end{proof}

    \section{Heights and quasi-canonical liftings}
    \label{sec:qcan_heights}

        \subsection{A descent lemma}
        \label{ssec:qcan_heights:descent_lemma}
            To compute Faltings and tautological heights, we will produce isogenies of $p$-divisible groups over $\Spec \mc{O}_{E,(p)}$ from isogenies over $\Spf \mc{O}_{\breve{E}_{\breve{w}}}$ for any choice of $\breve{w} \mid p$. We now explain this descent procedure, in a more general setup.

\begin{lemma}\label{lemma:qcan_heights:descent_lemma}
Let $S' \ra S$ be a morphism of schemes whose scheme-theoretic image is all of $S$. Suppose $X$ is a $p$-divisible group over $S$ which satisfies $\End^{0}(X) = \End^{0}(X_{S'})$. Let $Y$ and $Z$ be $p$-divisible groups over $S$ which are isogenous to $X$. The base-change maps
    \begin{align*}
    & \! \Hom^0(Y,Z) \ra \Hom^0(Y_{S'}, Z_{S'}) & & \! \Hom(Y,Z) \ra \Hom(Y_{S'}, Z_{S'}) \\
    & \mrm{Isog}^0(Y,Z) \ra \mrm{Isog}^0(Y_{S'}, Z_{S'}) & & \mrm{Isog}(Y,Z) \ra \mrm{Isog}(Y_{S'}, Z_{S'})
    \end{align*}
are bijections. 
\end{lemma}
\begin{proof}
Choose isogenies $\phi_Y \colon X \ra Y$ and $\phi_Z \colon X \ra Z$. There is a commutative diagram
    \[
    \begin{tikzcd}
    \a \arrow[mapsto]{d} & \End^0(X) \arrow{r}{\sim} \arrow{d}[vertLabel]{\sim} & \End^0(X_{S'}) \arrow{d}[vertLabelSwap]{\sim} & \a' \arrow[mapsto]{d} \\
    \phi_Z \circ \a \circ \phi_Y^{-1} & \Hom^0(Y,Z) \arrow{r} & \Hom^0(Y_{S'}, Z_{S'}) & \phi_{Z,S'} \circ \a' \circ \phi_{Y,S'}^{-1}
    \end{tikzcd}
    \]
where horizontal arrows are base-change. The vertical arrows are isomorphisms, and the upper horizontal arrow is an isomorphism by hypothesis. Hence the bottom arrow is an isomorphism. Suppose $\b \in \Hom^0(Y,Z)$ is any quasi-homomorphism. The functor $T \mapsto \{ \phi \in \Hom(T,S) : \phi^* \b \text{ is a homomorphism} \}$ is represented by a closed subscheme of $T$, see \cite[{Proposition 2.9}]{RZ96}. If $\b|_{S'}$ is a homomorphism, then $\b$ must also be a homomorphism, since the smallest closed subscheme of $S$ through which $S'$ factors is all of $S$ (by hypothesis). Hence $\Hom(Y,Z) \ra \Hom(Y_{S'},Z_{S'})$ is an isomorphism. The statements about (quasi-)isogenies follow from an essentially identical argument, replacing $\End$ and $\Hom$ with $\mrm{Isog}$, and noting $\mrm{Isog}^0(X) = (\End^0(X))^\times$ (e.g. by Lemma \ref{lemma:qisog_invertible}).
\end{proof}

\begin{remark}
We will be interested in the case where $S = \Spec \mc{O}_{E,(p)}$ and $S' = \Spec \mc{O}_{\breve{E}}$ for some finite extension $\breve{E}$ of $\breve{E}_{\breve{w}}$ for some $\breve{w} \mid p$. In this case, Lemma \ref{lemma:qcan_heights:descent_lemma} admits an alternative proof: a quasi-homomorphism of $p$-divisible groups over $\Spec E$ is a homomorphism if and only if the map on rational Tate modules preserves (integral) Tate modules, and this can be checked after base-change to $\Spec \breve{E}$. Then apply the theorem of Tate \cite[{Theorem 4}]{Tate67} which states that the generic fiber functor for $p$-divisible groups over $\Spec \mc{O}_{E,(p)}$ (similarly, for $\Spec \mc{O}_{\breve{E}}$) is fully faithful.
\end{remark}

\begin{lemma}\label{lemma:qcan_heights:descent_lemma:product_decomp_pdiv}
Let $X$ be a $p$-divisible group over a formal scheme $S$. Suppose there is a decomposition $X = X_1 \times X_2$ as fppf sheaves of abelian groups (on $(\Sch/S)_{fppf}$). Then $X_1$ and $X_2$ are both $p$-divisible groups.
\end{lemma}
\begin{proof}
Write $e_1, e_2 \in \End(X)$ for the projections to $X_1$ and $X_2$ respectively. As being a $p$-divisible group can be checked locally on $(\Sch/S)_{fppf}$, assume $S$ is a usual scheme.

It is clear that the multiplication by $p$ map $[p] \colon X \ra X$ is a surjection if and only if $[p] \colon X_1 \ra X_1$ and $[p] \colon X_2 \ra X_2$ are surjections. We also have $X[p^n] = X_1[p^n] \times X_2[p^n]$ for all $n \geq 1$. Thus the natural map $\varinjlim X[p^n] \ra X$ is an isomorphism if and only if $\varinjlim X_1[p^n] \ra X_1$ and $\varinjlim X_2[p^n] \ra X_2$ are isomorphisms.

Next, note $X_1[p] = \ker(e_2 \colon X[p] \ra X[p])$ and similarly $X_2[p] = \ker(e_1 \colon X[p] \ra X[p])$. Since $X[p]$ is representable by a finite locally free scheme over $S$, we conclude that $X_1[p]$ and $X_2[p]$ are represented by schemes which are finite and finitely presented over $S$. We also have short exact sequences
    \begin{align*}
    & 0 \ra X_1[p] \ra X[p] \ra X_2[p] \ra 0 \\
    & 0 \ra X_2[p] \ra X[p] \ra X_1[p] \ra 0
    \end{align*}
so Lemma \ref{lemma:surj_finite_locally_free_fppf} implies that $X_1[p]$ and $X_2[p]$ are finite locally free over $S$.
\end{proof}

\begin{corollary}\label{corollary:qcan_heights:descent_lemma:hom_lemma_product_decomp}
Let $S' \ra S$ and $X$ be as in Lemma \ref{lemma:qcan_heights:descent_lemma}. Suppose $Y$ and $Z$ are $p$-divisible groups over $S$ isogenous to $X$. 

If $Y_{S'} = Y'_1 \times \cdots \times Y'_r$ and $Z_{S'} = Z'_1 \times \cdots \times Z'_r$ for $p$-divisible groups $Y'_i$ and $Z'_i$ over $S'$, then there are unique decompositions $Y = Y_1 \times \cdots \times Y_r$ and $Z = Z_1 \times \cdots \times Z_r$ such that $Y_i|_{S'} = Y'_i$ and $Z_i|_{S'} = Z'_i$ for all $i$. For any $i$, the base-change maps
    \begin{align*}
    & \! \Hom^0(Y_i, Z_i) \ra \Hom^0(Y_{i, S'}, Z_{i, S'}) & & \! \Hom(Y_i,Z_i) \ra \Hom(Y_{i, S'}, Z_{i, S'}) \\
    & \mrm{Isog}^0(Y_i, Z_i) \ra \mrm{Isog}^0(Y_{i, S'}, Z_{i, S'}) & & \mrm{Isog}(Y_i ,Z_i) \ra \mrm{Isog}(Y_{i, S'}, Z_{i, S'})
    \end{align*}
are bijective.
\end{corollary}
\begin{proof}
The decomposition $Y_{S'} = Y'_1 \times \cdots \times Y'_r$ corresponds to a system of orthogonal idempotents $d'_1, \ldots, d'_r \in \End(Y_{S'})$, i.e. $d_i^{\prime 2} = d'_i$ for all $i$ and $d'_i d'_j = 0$ for all $i \neq j$. Lifting to a decomposition $Y = Y_1 \times \cdots \times Y_r$ is the same as lifting $\{d'_i\}_{i}$ to a system of orthogonal idempotents $\{d_i\}_{i}$ in $\End(Y)$. Such a lift exists and is unique by Lemma \ref{lemma:qcan_heights:descent_lemma}. The same applies for $Z$, and we write $\{e'_i\}_i$ and $\{e_i\}_i$ for the corresponding systems of idempotents. 
Using Lemma \ref{lemma:qcan_heights:descent_lemma}, we have 
    \begin{align*}
    & \Hom^{0}(Y_i,Z_i) = d_i \Hom^{0}(Y,Z) e_i = d'_i \Hom^{0} (Y_{S'}, Z_{S'}) e'_i = \Hom^{0}(Y_{i, S'},Z_{i, S'})
    \\
    & \Hom(Y_i,Z_i) = d_i \Hom(Y,Z) e_i = d'_i \Hom (Y_{S'}, Z_{S'}) e'_i = \Hom(Y_{i, S'}, Z_{i, S'}).
    \end{align*}
The statement about $\mrm{Isog}^{0}$ then follows from Lemma \ref{lemma:qisog_invertible}, and the statement about $\mrm{Isog}$ follows from the relation $\mrm{Isog}(-,-) = \mrm{Isog}^{0}(-,-) \cap \Hom(-,-)$.
\end{proof}
        
        \subsection{Minimal isogenies}
        \label{ssec:qcan_heights:minimal_isogenies}
            Given any abelian scheme $A \ra S$ over some base $S$, we can form the \emph{Serre tensor} abelian scheme $A \otimes_{\Z} \mc{O}_F$ given by $(A \otimes_{\Z} \mc{O}_F)(T) \coloneqq A(T) \otimes_{\Z} \mc{O}_F$ for $S$-schemes $T$. There is a natural action of $\mc{O}_F$ on $A \otimes_{\Z} \mc{O}_F$, as we have discussed for $p$-divisible groups \eqref{equation:appendix:Serre_tensor}. If $\lambda \colon A \ra A^{\vee}$ is a quasi-polarization, then $\lambda \otimes \lambda_{\mrm{tr}} \colon A \otimes_{\Z} \mc{O}_F \ra A^{\vee} \otimes_{\Z} \mc{O}_F^{\vee} \cong (A \otimes_{\Z} \mc{O}_F)^{\vee}$ is a polarization, where $\lambda_{\mrm{tr}} \colon \mc{O}_F \ra \mc{O}_F^{\vee}$ is induced by the trace pairing, as above.

Let $A_0 \ra \Spec \mc{O}_{E,(p)}$ be any (relative) elliptic curve with $\mc{O}_F$-action $\iota_0$ of signature $(1,0)$, and let $\lambda_0$ be the unique principal polarization of $A_0$. For $n \geq 2$, set
    \begin{equation}
    A \coloneqq A_0^{n - 2} \times (A_0 \otimes_{\Z} \mc{O}_F)
    \end{equation}
with $\mc{O}_F$ action $\iota$ which is diagonal on $A_0^{n - 2}$ and the Serre tensor action on $A_0 \otimes_{\Z} \mc{O}_F$, and polarization $\lambda_0^{n - 2} \times (|\Delta|^{-1} (\lambda_0 \otimes \lambda_{\mrm{tr}}))$. Then $(A, \iota, \lambda)$ is a special Hermitian abelian scheme of signature $(n - 1, 1)$. We write $(X, \iota, \lambda)$ for the associated special Hermitian $p$-divisible group of signature $(n - 1, 1)$, with
    \begin{equation}
    X = X_0^{n - 2} \times (X_0 \otimes_{\Z_p} \mc{O}_{F_p})
    \end{equation}
where $X_0$ is the $p$-divisible group of $A_0$. For any $\breve{w} \mid p$, the base-change $X_{0, \breve{w}}$ is a canonical lifting. The preceding notation (e.g. for $A_0$ and $X_0$) will be fixed for all of Section \ref{ssec:qcan_heights:minimal_isogenies}.

In Proposition \ref{proposition:qcan_heights:minimal_isogenies} and Corollary \ref{corollary:qcan_heights:minimal_isogenies:local_decomp} below, we equip $\mf{X}_0^{n - 2} \times (\mf{X}_s \otimes_{\Z_p} \mc{O}_{F_p})$ with the diagonal $\mc{O}_{F_p}$ action (which is the Serre tensor action on $\mf{X}_s \otimes_{\Z_p} \mc{O}_{F_p}$) and a product quasi-polarization, for some quasi-polarization of $\mf{X}_0^{n - 2}$ and the quasi-polarization $-\iota(\mf{d}_p^2)^{-1} \circ (\lambda_{\mf{X}_s} \otimes \lambda_{\mrm{tr}})$ on $(\mf{X}_s \otimes_{\Z_p} \mc{O}_{F_p})$.

\begin{proposition}\label{proposition:qcan_heights:minimal_isogenies}
Let $(A', \iota', \lambda')$ be a special Hermitian abelian scheme of signature $(n - 1, 1)$ over $\Spec \mc{O}_{E,(p)}$, with associated Hermitian $p$-divisible group $(X', \iota', \lambda')$. Replace $E$ with a finite extension if necessary, so that $A$ and $A'$ are $\mc{O}_F$-linearly isogenous.

Suppose there exists a $\mc{O}_{F_p}$-linear quasi-polarization preserving isomorphism
    
    \begin{equation}\label{equation:qcan_heights:minimal_isogenies:local_decomp}
    X'_{\Spf \mc{O}_{\breve{E}}} \cong \mf{X}_0^{n - 2} \times (\mf{X}_s \otimes_{\Z_p} \mc{O}_{F_p})
    \end{equation}
over $\Spf \mc{O}_{\breve{E}}$, where $\breve{E}$ is a finite extension of $\mc{O}_{\breve{E}_{\breve{w}'}}$, for some $\breve{w}' \mid p$ and $s \geq 0$.
Fix an isomorphism $X_{0, \Spf \mc{O}_{\breve{E}}} \cong \mf{X}_0$.

\begin{enumerate}[(1)]
    \item Then there exists an $\mc{O}_{F_p}$-linear quasi-polarization preserving isomorphism
        \begin{equation}\label{equation:qcan_heights:minimal_isogenies:semi-global_decomp}
        X' \cong X_0^{n - 2} \times (X_s \otimes_{\Z_p} \mc{O}_{F_p})
        \end{equation}
    over $\Spec \mc{O}_{E,(p)}$, for some $p$-divisible group $X_s$ of height $2$ and dimension $1$ with fixed identification $X_{s, \Spf \mc{O}_{\breve{E}}} \cong \mf{X}_s$, such that \eqref{equation:qcan_heights:minimal_isogenies:semi-global_decomp} recovers \eqref{equation:qcan_heights:minimal_isogenies:local_decomp} upon base-change to $\Spf \mc{O}_{\breve{E}}$. 
    On the right-hand side of \eqref{equation:qcan_heights:minimal_isogenies:semi-global_decomp}, the polarization is the product of a polarization on $X_0^{n - 2}$ and a quasi-polarization $- (\mf{d}^2_p)^{-1} \cdot (\lambda_s \otimes_{\Z_p} \lambda_{\mrm{tr}})$ on $X_s \otimes_{\Z_p} \mc{O}_{F_p}$, where $\lambda_s$ is a principal polarization on $X_s$.
    
    \item For any $\breve{w} \mid p$, the base-change $X_{s, \breve{w}}$ is a quasi-canonical lifting of level $s$, and hence there is an identification as in \eqref{equation:qcan_heights:minimal_isogenies:semi-global_decomp} for all $\breve{w} \mid p$.

    \item There exists an isogeny $\psi_s \colon X_0 \ra X_s$ of degree $p^s$. The $\mc{O}_{F_p}$-linear product isogeny $\phi \colon X \ra X'$ given by
        \begin{equation}
        \phi \coloneqq \mrm{id}_{X_0^{n - 2}} \times (\psi_s \otimes 1) \colon X_0^{n - 2} \times (X_0 \otimes_{\Z_p} \mc{O}_{F_p}) \ra X_0^{n - 2} \times (X_s \otimes_{\Z_p} \mc{O}_{F_p})
        \end{equation}
    over $\Spec \mc{O}_{E,(p)}$ satisfies
        \begin{align}
        & \breve{\delta}_{\mrm{Fal}}(\phi_{\breve{w}}) = -2 \breve{\delta}_{\mrm{tau}}(\phi_{\breve{w}}) = -2 \delta_{\mrm{tau}}(s) \cdot \log p && \text{for all $\breve{w} \mid p$} \label{equation:qcan_heights:minimal_isogenies:local_change}
        \\
        & \delta_{\mrm{Fal}, (p)}(\phi) = -2 \delta_{\mrm{tau}, (p)}(\phi) = -2 \delta_{\mrm{tau}}(s) \cdot \log p. \label{equation:qcan_heights:minimal_isogenies:semi-global_change}
        \end{align}

\end{enumerate}
\end{proposition}
\begin{proof}
Note that $X$ satisfies the hypotheses of Lemma \ref{lemma:qcan_heights:descent_lemma} with $S = \Spec \mc{O}_{E,{(p)}}$ and $S' = \Spec \mc{O}_{\breve{E}_{\breve{w}}}$ for any $\breve{w} \mid p$, as $\End(X) \cong M_{n,n}(\mc{O}_{F_p})$ over both $\Spec \mc{O}_{E,(p)}$ and $\Spf \mc{O}_{\breve{w}}$ for any $\breve{w}$. The same holds for $S' = \Spec \mc{O}_{\breve{E}}$. Again, we pass between $\Spf \mc{O}_{\breve{E}_{\breve{w}}}$ and $\Spec \mc{O}_{\breve{E}_{\breve{w}}}$ as in Appendix \ref{appendix:pDiv_prelim:Spec_v_Spf}.

The proposition then follows from repeated applications of Lemma \ref{lemma:qcan_heights:descent_lemma} and Corollary \ref{corollary:qcan_heights:descent_lemma:hom_lemma_product_decomp}, as we now explain.

\begin{enumerate}[(1)]
    \item and (2) Corollary \ref{corollary:qcan_heights:descent_lemma:hom_lemma_product_decomp} implies that \eqref{equation:qcan_heights:minimal_isogenies:local_decomp} descends to a $\mc{O}_{F_p}$-linear product decomposition $X' \cong X_0^{\prime n - 2} \times (X_s \otimes_{\Z_p} \mc{O}_{F_p})$ over $\Spec \mc{O}_{E, (p)}$ for some $X'_0$ descending $\mf{X}_0$ (first pick any identification of $p$-divisible groups $\mf{X}_s \otimes_{\Z_p} \mc{O}_{F_p} \cong \mf{X}_s^2$, then descend the $\mc{O}_{F_p}$-action), and the fully-faithfulness in Corollary \ref{corollary:qcan_heights:descent_lemma:hom_lemma_product_decomp} implies $\End(X_s) = \mc{O}_{F_p,s}$ (with $\mc{O}_{F_p,s} = \Z_p + p^s \mc{O}_{F_p}$ as in \crefext{I:ssec:can_and_qcan:qcan}) over $\Spec \mc{O}_{E,(p)}$ and also over $\Spf \mc{O}_{\breve{w}}$ for any $\breve{w} \mid p$. The fully-faithfulness in Corollary \ref{corollary:qcan_heights:descent_lemma:hom_lemma_product_decomp} also implies that the fixed $\mc{O}_{F_p}$-linear isomorphism $X_{0, \breve{w}} \ra \mf{X}_0$ lifts to an isomorphism $X_0 \ra X'_0$. The polarization on $\mf{X}_0^{n - 2} \times (\mf{X}_s \otimes_{\Z_p} \mc{O}_{F_p})$ descends to $X_0^{n - 2} \times (X_s \otimes_{\Z_p} \mc{O}_{F_p})$ by Corollary \ref{corollary:qcan_heights:descent_lemma:hom_lemma_product_decomp} again (applied to $X'$ and $X^{\prime \vee}$; note that the property of being a polarization is represented by a closed subfunctor of $\Spec \mc{O}_{E,(p)}$, hence can be checked in the generic fiber or over $\Spec \breve{E}$)
    \addtocounter{enumi}{1}

    \item If $\psi_s \colon \mf{X}_0 \ra \mf{X}_s$ is any isogeny of degree $p^s$ (exists and is unique up to precomposition by $\mc{O}_{F_p}^{\times}$, as discussed in \crefext{I:ssec:can_and_qcan:qcan}), we apply Corollary \ref{corollary:qcan_heights:descent_lemma:hom_lemma_product_decomp} to descend to an isogeny $\psi_s \colon X_0 \ra X_s$ of degree $p^s$. Equation \eqref{equation:qcan_heights:minimal_isogenies:local_change} now follows from Lemma \ref{lemma:Faltings_and_taut:Serre_tensor:qcan_isog}. Equation \eqref{equation:qcan_heights:minimal_isogenies:semi-global_change} follows from this (by the definitions in \eqref{equation:Faltings_and_taut:isogeny_change_local_Faltings:pDiv_semi-global_def} and \eqref{equation:Faltings_and_taut:isogeny_change_local_taut:pDiv_semi-global_def}). \qedhere
\end{enumerate}
\end{proof}

We will use the following reformulation (tailored to our intended application for global heights via local special cycles). In the corollary statement and proof, $A_0^{\s}$ and $A_0^{n - 1} \times A_0^{\s}$ are equipped with the product quasi-polarizations $- \Delta^{-1} \lambda_0$ and $-\Delta^{-1} (\lambda_0 \times \cdots \times \lambda_0)$ (where $A_0^{\s} = A_0$ but with $\mc{O}_F$-action $\iota \circ \s$, as above).

\begin{corollary}\label{corollary:qcan_heights:minimal_isogenies:local_decomp}
Let $S$ be a reduced scheme which is finite flat over $\Spec \mc{O}_F$. Let $(A', \iota', \lambda', \mc{F}')$ be a quasi-polarized Kr\"amer datum over $S$ (of signature $(n - 1, 1)$) for $n \geq 2$, with associated metrized line bundles $\widehat{\omega}$ and $\widehat{\ms{E}}^{\vee}$ on $S$. Assume that $(A', \iota', \lambda')$ is special at all generic points of $S$.
Let $(X', \iota', \lambda')$ be the associated Hermitian $p$-divisible group. 

Suppose we are given a finite \'etale surjection
    \begin{equation}
    \coprod_j \mc{Z}_j \ra S \times_{\Spec \Z} \Spec \breve{\Z}_p
    \end{equation}
such that each restricted map $\Theta_j \colon \mc{Z}_j \ra S \times_{\Spec \Z} \Spec \breve{\Z}_p$ has constant degree $\deg(j)$ onto its image. Assume that $\Theta_j$ and $\Theta_{j'}$ have disjoint images for $j \neq j'$.

For each irreducible component $\mc{Z} \hookrightarrow \coprod_j \mc{Z}_j$, write $\breve{E}_{\mc{Z}}$ for the residue field of its generic point. Assume there exists an isomorphism of Hermitian $p$-divisible groups
    \begin{equation}
    X'|_{\Spf \mc{O}_{\breve{E}_{\mc{Z}}}} \cong \mf{X}_0^{n - 2} \times (\mf{X}_{s_{\mc{Z}}} \otimes_{\Z_p} \mc{O}_{F_p})
    \end{equation}
for all $\mc{Z}$, where $s_{\mc{Z}} \in \Z_{\geq 0}$ is an integer depending on $\mc{Z}$.

We then have
    \begin{align}\label{equation:qcan_heights:minimal_isogenies:degree_decomp}
    \widehat{\deg}(\widehat{\ms{E}}^{\vee}) - (\deg_{\Z} S) \cdot h_{\mrm{tau}}(A_0^{\s}) & =  \sum_{j \in J} \frac{1}{\deg(j)} \sum_{\mc{Z} \hookrightarrow \mc{Z}_j} (\deg_{\breve{\Z}_p} \mc{Z}) \delta_{\mrm{tau}}(s_{\mc{Z}}) \log p
    \\
    \widehat{\deg}(\widehat{\omega}) - (\deg_{\Z} S) \cdot h_{\mrm{Fal}}(A_0^{n - 1} \times A_0^{\s}) & = \sum_{j \in J} \frac{1}{\deg(j)} \sum_{\mc{Z} \hookrightarrow \mc{Z}_j}  (\deg_{\breve{\Z}_p} \mc{Z}) \delta_{\mrm{Fal}}(s_{\mc{Z}}) \log p \notag
    \end{align}
modulo $\sum_{\ell \neq p} \Q \cdot \log \ell$, where the inner sums run over all irreducible components $\mc{Z} \hookrightarrow \mc{Z}_j$.
\end{corollary}
\begin{proof}
In the corollary statement, the expression ``modulo $\sum_{\ell \neq p} \Q \cdot \log \ell$'' means an equality of elements in the additive quotient $\R / (\sum_{\ell \neq p} \Q \cdot \log \ell)$. The notation $\deg_{\Z} S$ (resp. $\deg_{\breve{\Z}_p} \mc{Z}$) denotes the degree of $S \ra \Spec \Q$ (resp. $\mc{Z} \ra \Spec \breve{\Z}_p$) in the generic fiber.

By additivity, we immediately reduce to the case where $S$ is irreducible. Then $J$ consists of a single element $j$. By normalization, we may assume $S = \Spec \mc{O}_E$ for a number field $E$. We may also enlarge $E$ as necessary so that $(A_0, \iota, \lambda_0)$ also extends to $\Spec \mc{O}_E$, and such that there exists an $\mc{O}_F$-linear isogeny $\phi \colon A_0^{n - 2} \times (A_0 \otimes_{\Z} \mc{O}_F) \ra A$. We also consider the $\mc{O}_F$-linear isogeny $\phi' \colon A_0^{n - 2} \times (A_0 \otimes_{\Z} \mc{O}_F) \ra A_0^{n - 2} \times (A_0 \times A_0^{\s})$ which is the identity on $A_0^{n - 2}$ and given by $(x \otimes a) \mapsto (ax, a^{\s} x)$ for $(A_0 \otimes_{\Z} \mc{O}_F) \ra A_0 \times A_0^{\s}$.

Since $\delta_{\mrm{tau}}(\phi') = 0$ (Lemma \ref{lemma:Faltings_and_taut:Serre_tensor:product}, along with the local decomposition \eqref{equation:Faltings_and_taut:isogeny_change_local_taut:pDiv_semi-global_def}, also Lemma \ref{lemma:Faltings_and_taut:isogeny_change_local_taut:global-to-local}), the decomposition in \eqref{equation:Faltings_and_taut:isogeny_change_local_taut:decomp_to_semi-global} shows
    \begin{align}
    \widehat{\deg}(\widehat{\ms{E}}^{\vee}) - [E : \Q] \cdot h_{\mrm{tau}}(A_0^{n - 1} \times A_0^{\s}) = [E : \Q] \delta_{\mrm{tau},(p)}(\phi) \mod \sum_{\ell \neq p} \Q \cdot \log \ell.
    \end{align}
We have $\sum_{\mc{Z} \hookrightarrow \mc{Z}_j} (\deg_{\breve{\Z}_p} \mc{Z}) = \deg(j) \cdot [E : \Q]$. Applying Proposition \ref{proposition:qcan_heights:minimal_isogenies}(3) (combined with the ``isogeny independence'' result of Lemma \ref{lemma:Faltings_and_taut:isogeny_change_local_taut:global-to-local}) now shows $\delta_{\mrm{tau},(p)}(\phi) = \delta_{\mrm{tau}}(s_{\mc{Z}}) \log p $ for any $\mc{Z} \hookrightarrow \mc{Z}_j$. This also shows that all $s_{\mc{Z}}$ are equal (when $S$ is irreducible): the quantity $\delta_{\mrm{tau}}(s)$ takes distinct values for distinct $s \in \Z_{\geq 0}$ (in the nonsplit cases, note $\delta_{\mrm{tau}}(s)$ has strictly decreasing $p$-adic valuation as $s$ increases, for $s > 2$). We also have $h_{\mrm{tau}}(A_0^{n - 1} \times A_0^{\s}) = h_{\mrm{tau}}(A_0^{\s})$ (straightforward from the definition). This verifies \eqref{equation:qcan_heights:minimal_isogenies:degree_decomp} for $\widehat{\deg}(\widehat{\ms{E}}^{\vee})$ and the tautological height.

Since $\delta_{\mrm{Fal}}(\phi') = 0$ (Lemma \ref{lemma:Faltings_and_taut:Serre_tensor:product}, along with the local decomposition \eqref{equation:Faltings_and_taut:isogeny_change_local_Faltings:pDiv_semi-global_def}, also Lemma \ref{lemma:Faltings_and_taut:isogeny_change_local_Faltings:global-to-local}), \eqref{equation:Faltings_and_taut:isogeny_change_local_taut:decomp_to_semi-global} similarly shows
    \begin{align}
    \widehat{\deg}(\widehat{\omega}) - [E : \Q] \cdot h_{\mrm{Fal}}(A_0^{n - 1} \times A_0^{\s}) = [E : \Q] \delta_{\mrm{Fal},(p)}(\phi) \mod \sum_{\ell \neq p} \Q \cdot \log \ell.
    \end{align}
Applying Proposition \ref{proposition:qcan_heights:minimal_isogenies}(3) (combined with the ``isogeny independence'' result of Lemma \ref{lemma:Faltings_and_taut:isogeny_change_local_Faltings:global-to-local}) verifies \eqref{equation:qcan_heights:minimal_isogenies:degree_decomp} for $\widehat{\deg}(\widehat{\omega})$ and the Faltings height, just as for tautological height above.
\end{proof}

In the situation above, we have 
    \begin{align}\label{equation:qcan_heights:minimal_isogenies:CM_heights}
    h_{\mrm{tau}}(A_0^{\s}) = h_{\mrm{tau}}^{\mrm{CM}} \quad \quad
    h_{\mrm{Fal}}(A_0^{n - 1} \times A_0^{\s}) = n \cdot h_{\mrm{Fal}}^{\mrm{CM}}
    \end{align}
in the notation of \eqref{equation:arith_cycle_classes:Hodge_bundles:taut_height_constants} and \eqref{equation:arith_cycle_classes:Hodge_bundles:CM_Faltings_height}.

\begin{remark}
In Proposition \ref{proposition:qcan_heights:minimal_isogenies}(3), it was important that $\psi_s$ was an isogeny of minimal degree $p^s$. If $\psi_s$ were replaced by an arbitrary isogeny $f \colon \mf{X}_0 \ra \mf{X}_s$, we would not be able to determine $\breve{\delta}_{\mrm{Fal}}(f)$ or $\breve{\delta}_{\mrm{tau}}(f)$ using only $\deg f$ in the case when $F_p / \Q_p$ is split (due to Remark \ref{remark:Faltings_and_taut:isogeny_change_local_taut:split_complication}).
\end{remark}

\clearpage


    \part{Uniformization}
    \label{part:part_III:uniformization}

        In \cref{sec:non-Arch_uniformization,sec:Arch_uniformization}, we use global notation as in \cref{sec:ab_var,sec:arith_cycle_classes}, e.g. $F$ is an imaginary quadratic field extension of $\Q$ with nontrivial involution $a \mapsto a^{\s}$ and discriminant $\Delta$. The notation $\A_f$ (resp. $\A_f^p$) will always denote the finite ad\`ele ring (resp. finite ad\`ele ring away from $p$) for $\Q$.

For all of \cref{sec:non-Arch_uniformization,sec:Arch_uniformization}, let $L_0 \coloneqq \mc{O}_F$ be the rank one Hermitian $\mc{O}_F$-lattice with pairing $(x,y) \coloneqq x^{\s} y$. Let $L$ be any non-degenerate Hermitian $\mc{O}_F$-lattice of rank $n$ and signature $(n - r, r)$, with associated moduli stack $\mc{M}$ (Definition \ref{definition:stack_ass_to_lattice} and surrounding discussion).

We fix some group-theoretic setup (common in the literature, e.g. \cite{RSZ20,BHKRY20}). 
Set
    \begin{gather*}
    V_0 \coloneqq L_0 \otimes_{\mc{O}_F} F \quad \quad V \coloneqq L \otimes_{\mc{O}_F} F
    \\
    G' \coloneqq \{ (g_0, g) \in GU(V_0) \times GU(V) : c(g_0) = c(g) \} \subseteq GU(V_0) \times GU(V)
    \end{gather*}
where $c \colon GU(V_0) \ra \G_m$ and $c \colon GU(V) \ra \G_m$ are the similitude characters. We use the shorthand
    \begin{equation*}
    L_p \coloneqq L \otimes_{\Z} \Z_p \quad \quad V_p \coloneqq V \otimes_{\Q} \Q_p \quad \quad V_{\R} \coloneqq V \otimes_{\Q} \R
    \end{equation*}
and use similar notation for local versions of other Hermitian spaces. Given a tuple $\underline{x} \in V^m$, we write $\underline{x}_p \in V_p^m$ and $\underline{x}_{\infty} \in V_{\R}^m$ and $\underline{x}_f \in (V \otimes_{\Q} \A_f)^m$ and $\underline{x}^p \in (V \otimes_{\Q} \A_f^p)^m$ for the corresponding projections (and similarly for other Hermitian spaces).

There is an isomorphism
    \begin{equation}\label{equation:uniformization:group_iso} \tag{$*$}
    \begin{tikzcd}[row sep = tiny]
    G' \arrow{r} & GU(V_0) \times U(V) \\
    (g_0, g) \arrow[mapsto]{r} & (g_0, g_0^{-1} g).
    \end{tikzcd}
    \end{equation}
To avoid potential confusion: whenever we write $(g_0, g) \in G'$, we mean $g_0 \in GU(V_0)$ and $g \in GU(V)$ with the same similitude factor.

We use factorizable open compact subgroups $K'_f = K_{0,f} \times K_f \subseteq G'(\A_f)$ as in Section \ref{ssec:ab_var:level_structure}, where $K_{0,f} \subseteq GU(V_0)(\A_f)$ and $K_f \subseteq U(V)(\A_f)$ (using also \eqref{equation:uniformization:group_iso}).

Recall the moduli stack with level structure $\mc{M}_{K'_f}$ defined in Section \ref{ssec:ab_var:level_structure}. We do not require $K'_f$ to be a small level for non-Archimedean uniformization (\cref{sec:non-Arch_uniformization}), so $\mc{M}_{K'_f}$ is allowed to be a stack. We will take small level to have manifolds rather than orbifolds in complex uniformization (\cref{sec:Arch_uniformization}).

\begin{notation*}
In \cref{sec:non-Arch_uniformization,sec:Arch_uniformization}, we implicitly fix an open compact subgroup $K'_f \subseteq G'(\A_f)$ as above. We abusively suppress $K'_f$ from notation: we write
    \begin{equation*}
    \mc{M} \quad \quad \mc{Z}(T) \quad \quad {}^{\mbb{L}}\mc{Z}(T) \quad \quad {}^{\mbb{L}}\mc{Z}(T)_{\ms{V},p}
    \end{equation*}
instead of $\mc{M}_{K'_f}$, $\mc{Z}(T)_{K'_f}$, ${}^{\mbb{L}} \mc{Z}(T)_{K'_f}$, ${}^{\mbb{L}} \mc{Z}(T)_{\ms{V},p,K'_f}$, etc..
\end{notation*}

For example, given a Hermitian matrix $T \in \mrm{Herm}_m(\Q)$ (with entries in $F$) and an appropriate scheme $S$, our notation entails
    \begin{equation*}
    \mc{Z}(T)(S) = \{(A_0, \iota_0, \lambda_0, A, \iota, \lambda, \tilde{\eta}_0, \tilde{\eta}, \underline{x}) \text{ over $S$}\}
    \end{equation*}
where $(\tilde{\eta}_0, \tilde{\eta})$ is a $K'_f$ level structure and $\underline{x} \in \Hom_{\mc{O}_F}(A_0, A)^m$ satisfies $(\underline{x}, \underline{x}) = T$.

        \section{Non-Archimedean}
        \label{sec:non-Arch_uniformization}
            Fix a prime $p$. 
If $p$ is not inert, we assume the signature is $(n - r, r) = (n - 1, 1)$. 
We assume that $L \otimes_{\Z} \Z_p$ is self-dual. If $p$ is ramified, we assume $n$ is even, $L$ is self-dual (for the trace pairing), and $p \neq 2$.

In all cases, we assume the implicit level $K'_f = K_{0,f} \times K_f$ at $p$ is
    \begin{equation}
    K_{0,p} = K_{L_0, p} \quad \quad K_{p} = K_{L,p}.
    \end{equation}
Recall that these denote the stabilizers of $L_0$ and $L$, respectively. 

Set $\mc{O}_{F_p} \coloneqq \mc{O}_F \otimes_{\Z} \Z_p$ and $\mc{O}_{F,(p)} \coloneqq \mc{O}_F \otimes_{\Z} \Z_{(p)}$ and $F_p \coloneqq F \otimes_{\Q} \Q_p$. We write $\breve{F}_p$ for the completion of the maximal unramified extension of $F_p$ is $p$ is nonsplit (resp. $\breve{F}_p \coloneqq \breve{\Q}_p$ if $p$ is split, with a choice of morphism $F_p \ra \breve{F}_p$). In all cases, $\mc{O}_{\breve{F}_p}$ (resp. $\overline{k}$) will denote the ring of integers (resp. residue field) of $\breve{F}_p$.

We discuss Rapoport--Zink uniformization \cite[{\S 6}]{RZ96}, as applied to supersingular loci on special cycles by Kudla--Rapoport \cite[{\S 5, \S 6}]{KR14} (there in the inert case, $p \neq 2$). For $p$ inert or ramified, the material in \cref{ssec:non-Arch_uniformization:formal_completion,ssec:non-Arch_uniformization:away_p_special_cycle,ssec:non-Arch_uniformization:framing,ssec:non-Arch_uniformization:framed_stack,ssec:non-Arch_uniformization:quotient,ssec:non-Arch_uniformization:proof} is essentially a repackaging of Rapoport--Zink and Kudla--Rapoport. However, we need modified arguments at split places: the abelian varieties will be ordinary. We give a mostly uniform treatment for inert/ramified/split places. We also allow $p = 2$ if $p$ is inert, except in Section \ref{ssec:non-Arch_uniformization:horizontal}.

Section \ref{ssec:non-Arch_uniformization:horizontal} is the newest part of Section \ref{sec:non-Arch_uniformization}. Here, we explain how to use uniformization to reduce (global) Faltings or ``tautological'' heights to quantities expressed in terms of local special cycles and the ``local change of heights'' from \cref{part:local_change_heights} (with the main input being \cref{corollary:qcan_heights:minimal_isogenies:local_decomp}).

Section \ref{ssec:non-Arch_uniformization:global_to_local} is the next newest part of Section \ref{sec:non-Arch_uniformization}. We use global special cycles and an ``approximation'' argument to prove certain properties of local special cycles. Some of these results are available or implicit in the literature (for $p$ nonsplit, sometimes with $p \neq 2$ hypotheses and signature $(n - 1, 1)$ hypotheses); we indicate this where relevant. Our methods of proof are different, based on the approximation argument mentioned above.

Section \ref{ssec:non-Arch_uniformization:vertical} is the next newest part of Section \ref{sec:non-Arch_uniformization}. We explain how to reduce global ``vertical intersection numbers'' to local ``vertical intersection numbers''.

\Cref{ssec:non-Arch_uniformization:global_to_local,ssec:non-Arch_uniformization:vertical,ssec:non-Arch_uniformization:horizontal} will need detailed information on the construction of Rapoport--Zink uniformization. This is our other reason for giving an exposition of uniformization in Sections \ref{ssec:non-Arch_uniformization:formal_completion}-\ref{ssec:non-Arch_uniformization:proof}, as we need to explain the relevant maps (and fix notation) to give precise statements.

\Cref{ssec:non-Arch_uniformization:formal_completion,ssec:non-Arch_uniformization:away_p_special_cycle,ssec:non-Arch_uniformization:framing,ssec:non-Arch_uniformization:framed_stack,ssec:non-Arch_uniformization:quotient} state the precise Rapoport--Zink uniformization map for special cycles. The proof of uniformization appears in Section \ref{ssec:non-Arch_uniformization:proof} (and allows $p = 2$ inert). We differ slightly from \cite{RZ96} by working directly with formal algebraic stacks (rather than requiring sufficient level structure) in the sense of \cite{Emerton20}. We occasionally need some notions on formal algebraic stacks which are not defined in \cite{Emerton20}; we will define these as needed.

Throughout \cref{sec:non-Arch_uniformization}, we freely use the relevant Rapoport--Zink spaces and their (Kudla--Rapoport) local special cycles as in \crefext{I:sec:moduli_pDiv}.
    
            \subsection{Formal completion}
            \label{ssec:non-Arch_uniformization:formal_completion}
                Throughout Section \ref{sec:non-Arch_uniformization}, the notation $T$ will always mean an $m \times m$ Hermitian matrix with $F$-coefficients, i.e. $T \in \mrm{Herm}_{m}(\Q)$. If $p$ is split in $\mc{O}_F$, we assume $\rank(T) \geq n - 1$. Form the special cycle $\mc{Z}(T) \ra \mc{M}$. 

Suppose $p$ is nonsplit. The \emph{supersingular locus} on $\mc{Z}(T)_{\overline{k}} \coloneqq \mc{Z}(T) \times_{\Spec \mc{O}_F} \Spec \overline{k}$ is the subset $\mc{Z}(T)^{ss} \subseteq |\mc{Z}(T)_{\overline{k}}|$ of the underlying topological space\footnote{By the \emph{underlying topological space} $|\mc{X}|$ of a formal algebraic stack $\mc{X}$, we mean the underlying topological space of its reduced substack $\mc{X}_{\red}$. As $\mc{X}_{\red}$ is an algebraic stack, it has an underlying topological space in the sense of \cite[\href{https://stacks.math.columbia.edu/tag/04XE}{Section 04XE}]{stacks-project}.} consisting of geometric points
    $(A_0, \iota_0, \lambda_0, A, \iota, \lambda, \tilde{\eta}_0, \tilde{\eta}, \underline{x})$
with $A$ supersingular.
The supersingular locus $\mc{Z}(T)^{ss}$ is a closed subset of $|\mc{Z}(T)_{\overline{k}}|$ (by the Katz--Grothendieck theorem on specialization for Newton polygons).
The \emph{formal completion of $\mc{Z}(T)_{\Spec \mc{O}_{\breve{F}_p}} \coloneqq \mc{Z}(T) \times_{\Spec \mc{O}_F} \Spec \mc{O}_{\breve{F}_p}$ along its supersingular locus} is the (strictly full) substack $\breve{\mc{Z}}(T) \subseteq \mc{Z}(T)_{\Spec \mc{O}_{\breve{F}_p}}$ given by 
    \begin{equation}
    \breve{\mc{Z}}(T) \coloneqq \{ \a \in \mc{Z}(T)_{\Spec \mc{O}_{\breve{F}_p}}(S) : \a(|S|) \subseteq \mc{Z}(T)^{ss} \}
    \end{equation}
for schemes $S$ over $\Spec \mc{O}_{\breve{F}_p}$, where the condition $\a(|S|) \subseteq \mc{Z}(T)^{ss}$ means that the associated map on underlying topological spaces $|S| \ra |\mc{Z}(T)_{\Spec \mc{O}_{\breve{F}_p}}|$ factors through $\mc{Z}(T)^{ss}$ (with $\a \in \mc{Z}(T)_{\Spec \mc{O}_{\breve{F}_p}}(S)$ ``viewed'' as a morphism $S \ra \mc{Z}(T)_{\Spec \mc{O}_{\breve{F}_p}}$ by the $2$-Yoneda lemma).

If $p$ is split, we define
    \begin{equation}
    \breve{\mc{Z}}(T) \coloneqq \mc{Z}(T)_{\Spf \mc{O}_{\breve{F}_p}} \coloneqq \mc{Z}(T) \times_{\Spec \mc{O}_F} \Spf \mc{O}_{\breve{F}_p}.
    \end{equation}
This is also the formal completion of $\mc{Z}(T)_{\Spec \mc{O}_{\breve{F}_p}}$ along its special fiber. For any geometric point $(A_0, \iota_0, \lambda_0, A, \iota, \lambda, \tilde{\eta}_0, \tilde{\eta}, \underline{x})$ of $\breve{\mc{Z}}(T)$, the abelian variety $A$ is ordinary (because Lemma \ref{lemma:corank_1_isogeny_decomposition} implies $A$ is isogenous to a product of elliptic curves with $\mc{O}_F$ action).

In all cases, $\breve{\mc{Z}}(T)$ is a locally Noetherian formal algebraic stack in the sense of \cite{Emerton20} (formal completion is discussed in \cite[{Example 5.9}]{Emerton20}). The structure morphism $\breve{\mc{Z}}(T) \ra \Spf \mc{O}_{\breve{F}_p}$ is formally smooth,\footnote{Given a morphism $f$ of categories fibered in groupoids over some base scheme, there is a category of dotted arrows \cite[\href{https://stacks.math.columbia.edu/tag/0H18}{Definition 0H18}]{stacks-project} associated to the infinitesimal lifting problem along each square-zero thickening of affine schemes. We say that $f$ is \emph{formally smooth} (resp. \emph{formally \'etale}) (resp. \emph{formally unramified}) if each such category of dotted arrows is nonempty (resp. a setoid with exactly one isomorphism class) (resp. either empty or a setoid with exactly one isomorphism class).} formally locally of finite type,\footnote{We say a morphism of locally Noetherian formal algebraic stacks is \emph{formally locally of finite type} if it is locally of finite type on underlying reduced substacks.} separated, and quasi-compact. If $K'_f$ is a small level, then $\breve{\mc{Z}}(T)$ is a locally Noetherian formal scheme.

If $\mc{M} = \mc{Z}(T)$ (e.g. $T = \emptyset$ or $T = 0$), we set $\breve{\mc{M}} \coloneqq \breve{\mc{Z}}(T)$. If $p$ is nonsplit, this is the \emph{formal completion of $\mc{M}_{\Spec \mc{O}_{\breve{F}_p}}$ along its supersingular locus $\mc{M}^{ss}$}.

            \subsection{Local special cycles away from \texorpdfstring{$p$}{p}}
            \label{ssec:non-Arch_uniformization:away_p_special_cycle}
                Given an $m$-tuple $\underline{x}^p = [x_1, \ldots, x_m] \in (V \otimes_{\Q} \A_f^p)^m$, we consider an \emph{``away-from-$p$'' local special cycle}
    \begin{equation}
    \mc{Z}'(\underline{x}^p) \coloneqq \{ (g_0, g) : G'(\A_f^p) / K'^p_f : g^{-1} g_0 x_i \in L \otimes_{\Z} \hat{\Z}^p \text{ for all $x_i \in \underline{x}^p$}\}.
    \end{equation}
We often view $\mc{Z}'(\underline{x}^p)$ and $G'(\A_f^p) / K'^p_f$ as constant formal schemes over $\Spf \mc{O}_{\breve{F}_p}$. We also define the \emph{``away-from-$p$'' local special cycle}
    \begin{equation}
    \mc{Z}(\underline{x}^p) \coloneqq \{ g : U(V)(\A_f^p) / K_f^p : g^{-1} x_i \in L \otimes_{\Z} \hat{\Z}^p \text{ for all $x_i \in \underline{x}^p$}\}.
    \end{equation}
The isomorphism $G'(\A_f^p) / K'^p_f \ra GU(V_0)(\A^p_f)/K^p_{0,f} \times U(V)(\A_f^p)/K^p_f$ \eqref{equation:uniformization:group_iso} induces an isomorphism
    \begin{equation}\label{equation:non-Arch_uniformization:away_p_special_cycle:iso}
    \mc{Z}'(\underline{x}^p) \xra{\sim} GU(V_0)(\A_f^p) / K^p_{0,f} \times \mc{Z}(\underline{x}^p).
    \end{equation}

            \subsection{Framing objects}
            \label{ssec:non-Arch_uniformization:framing}
                To define the uniformization map, we fix an object $(\mbf{A}_0, \iota_{\mbf{A}_0}, \lambda_{\mbf{A}_0}, \mbf{A}, \iota_{\mbf{A}}, \lambda_{\mbf{A}}, \tilde{\pmb{\eta}}_0, \tilde{\pmb{\eta}}) \in \mc{M}(\overline{k})$ (``basepoint of the uniformization''). If $p$ is nonsplit (resp. split), we assume $\mbf{A}$ is supersingular (resp. $\mbf{A}$ is $\mc{O}_F$-linearly isogenous to $\mbf{A}_0^{n - r} \times (\mbf{A}_0^{\s})^r$); such data exists by Lemma \ref{lemma:ab_var:integral_models:nonempty}. Let $(\mbf{X}_0, \iota_{\mbf{X}_0}, \lambda_{\mbf{X}_0}, \mbf{X}, \iota_{\mbf{X}}, \lambda_{\mbf{X}})$ be the tuple obtained by passing to $p$-divisible groups (e.g. $\mbf{X}$ is the $p$-divisible group of $\mbf{A}$). We use this as the framing object over $\overline{k}$ to define the Rapoport--Zink space $\mc{N}'$ (as in \crefext{I:definition:moduli_pDiv:RZ:N_prime}). Set $\mc{N} \coloneqq \mc{N}(n - r, r)$ (with the latter as in \crefext{I:definition:moduli_pDiv:RZ:RZ_spaces}).

In the supersingular cases, the abelian variety $\mbf{A}$ is automatically $\mc{O}_F$-linearly isogenous to $\mbf{A}_0^{n - r} \times (\mbf{A}_0^{\s})^r$, since
    \begin{equation}
    \Hom^0_F(\mbf{A}_0^{n - r} \times (\mbf{A}_0^{\s})^r, \mbf{A}) \otimes_{\Q} \Q_p \xra{\sim} \Hom^0_F(\mbf{X}_0^{n - r} \times (\mbf{X}_0^{\s})^r, \mbf{X})
    \end{equation}
by Tate's isogeny theorem (for any supersingular abelian variety over a finite field, some power of Frobenius will be a power of $p$, e.g. by \cite[{Lemma 6.28}]{RZ96}); then use uniqueness of the framing object $(\mbf{X}, \iota_{\mbf{X}}, \lambda_{\mbf{X}})$ up to isogeny (see \crefext{I:ssec:moduli_pDiv:RZ}).

Since $\mc{M}_{\Spec \mc{O}_{\breve{F}_p}} \ra \Spec \mc{O}_{\breve{F}_p}$ is smooth, this framing object $(\mbf{A}_0, \ldots)$ admits a lift $(\mf{A}_0, \iota_{\mf{A}_0}, \lambda_{\mf{A}_0}, \mf{A}, \iota_{\mf{A}}, \lambda_{\mf{A}}, \tilde{\mf{y}}_0, \tilde{\mf{y}}) \in \mc{M}(\Spf \mc{O}_{\breve{F}_p})$, which we also fix. We fix representatives
    \begin{equation}\label{equation:non-Arch_uniformization:framing:Tate_module_lattice}
    \pmb{\eta}_0 \colon T^p(\mbf{A}_0) \xra{\sim} L_0 \otimes_{\Z} \hat{\Z}^p \quad \quad \pmb{\eta} \colon \Hom_{\mc{O}_F \otimes_{\Z} \hat{\Z}^p}(T^p(\mbf{A}_0), T^p(\mbf{A})) \xra{\sim} L \otimes_{\Z} \hat{\Z}^p
    \end{equation}
for the $K_{0}^p$-orbit $\tilde{\pmb{\eta}}_0$ and the $K^p$-orbit $\tilde{\pmb{\eta}}$ (see Section \ref{ssec:ab_var:level_structure}). Recall that $\pmb{\eta}$ preserves Hermitian pairings but $\pmb{\eta}_0$ need not. We also write
    \begin{equation}
    \mf{y}_0 \colon T^p(\mf{A}_0) \xra{\sim} L_0 \otimes_{\Z} \hat{\Z}^p \quad \quad \mf{y} \colon \Hom_{\mc{O}_F \otimes_{\Z} \hat{\Z}^p}(T^p(\mf{A}_0), T^p(\mf{A})) \xra{\sim} L \otimes_{\Z} \hat{\Z}^p
    \end{equation}
for the identifications induced by $\pmb{\eta}_0$ and $\pmb{\eta}$.

We define Hermitian $F$-modules
    \begin{align}
    & \mbf{W} \coloneqq \Hom^0_F(\mbf{A}_0, \mbf{A}) && \mbf{W}^{\perp} \coloneqq \begin{cases}
    \Hom^0_F(\mbf{A}_0^{\s}, \mbf{A}) & \text{if $p$ is split} \\
    0 & \text{if $p$ is nonsplit}
    \end{cases}
    \\
    & \mbf{V}_0 \coloneqq \Hom^0_F(\mbf{A}_0, \mbf{A}_0)
    && \mbf{V} \coloneqq \mbf{W} \oplus \mbf{W}^{\perp}
    \end{align}
where the direct sum defining $\mbf{V}$ is orthogonal. In all cases, the Hermitian pairing is $(x, y) \coloneqq x^{\dagger} y \in F$. All of these Hermitian spaces are positive definite (positivity of the Rosati involution).

The canonical maps
    \begin{equation}\label{equation:non-Arch_uniformization:framing:at_p_Herm}
    \mbf{W} \otimes_{\Q} \Q_p \xra{\sim} \Hom_{F_p}^0(\mbf{X}_0, \mbf{X}) \quad \quad \mbf{W}^{\perp} \otimes_{\Q} \Q_p \xra{\sim} \Hom_{F_p}^0(\mbf{X}_0^{\s}, \mbf{X})
    \end{equation}
are isomorphisms of Hermitian spaces. In the nonsplit (hence supersingular) cases, this follows from Tate's isogeny theorem as above. In the split case, this follows because $\mbf{A}$ is $\mc{O}_F$-linearly isogenous to $\mbf{A}_0^{n - r} \times (\mbf{A}_0^{\s})^r$. In particular, the local invariant is $\varepsilon(\mbf{W}_p) = (-1)^r$ if $p$ is nonsplit (resp. $\varepsilon(\mbf{W}_p) = 1$ if $p$ is split).

If $p$ is nonsplit, the natural map
    \begin{equation}\label{equation:non-Arch_uniformization:framing:W_to_hom}
    \mbf{W} \otimes_{\Q} \A_f^p \ra \Hom_{F \otimes_{\Q} \A_f^p}(T^p(\mbf{A}_0)^0, T^p(\mbf{A})^0)
    \end{equation}
is an isomorphism of Hermitian spaces, by similar reasoning.

If $p$ is split, any $\mc{O}_F$-linear isogeny $\mbf{A}_0^{n - 1} \times \mbf{A}_0^{\s} \ra \mbf{A}$ defines an $F$-linear orthogonal decomposition
    \begin{equation}
    T^p(\mbf{A})^0 = T^p(\mbf{A}_0^{n - 1})^0 \oplus T^p(\mbf{A}_0^{\s})^0.
    \end{equation}
This decomposition is independent of the choice of isogeny because $\Hom^0_F(\mbf{A}_0, \mbf{A}_0^{\s}) = \Hom^0_F(\mbf{A}_0^{\s}, \mbf{A}_0) = 0$ (e.g. because $\End^0(\mbf{A}_0) = F$). Then the natural map
    \begin{equation}\label{equation:non-Arch_uniformization:framing:W_to_hom:split}
    \mbf{W} \otimes_{\Q} \A_f^p \ra \Hom_{F \otimes_{\Q} \A_f^p}(T^p(\mbf{A}_0)^0, T^p(\mbf{A}_0^{n - 1})^0)
    \end{equation}
is an isomorphism of Hermitian spaces.

Given a tuple $\underline{\mbf{x}} \in \mbf{W}^m$, we write
    \begin{align}
    \underline{\mbf{x}}_p & \in \mbf{W}^m_p = \Hom^0_{F_p}(\mbf{X}_0, \mbf{X})
    \\
    \underline{\mbf{x}}^p & \in \mbf{W}^m \otimes_{\Q} \A_f^p \subseteq \Hom_{F \otimes \A_f^p}(T^p(\mbf{A}_0)^0, T^p(\mbf{A})^0)^m = V^m \otimes_{\Q} \A_f^p
    \end{align}
for the respective images of $\underline{\mbf{x}}$ (using $\pmb{\eta}$ for the identification with $V^m \otimes_{\Q} \A_f^p$ in the second line). 
    
            \subsection{Framed stack}
            \label{ssec:non-Arch_uniformization:framed_stack}
                We consider the stack $\breve{\mc{Z}}(T)_{\mrm{framed}}$ over $\Spf \mc{O}_{\breve{F}_p}$, given by
    \begin{equation*}
    \breve{\mc{Z}}(T)_{\mrm{framed}}(S) \coloneqq \left \{ (\a, \phi_0, \phi) : \begin{array}{l} \a = (A_0, \iota_0, \lambda_0, A, \iota, \lambda, \tilde{\eta}_0, \tilde{\eta}, \underline{x}) \in \breve{\mc{Z}}(T)(S) \\ \phi_0 \colon A_0 \ra \mf{A}_{0,S} \text{ and } \phi \colon A \ra \mf{A}_S \quad \text{quasi-isogenies} \\ \quad \text{such that } \phi_0^* \lambda_{\mf{A}_0,S} = b \lambda_0 \text{ and } \phi^* \lambda_{\mf{A}, S} = b \lambda \\ \quad \text{for some $b \in \Q_{>0}$} \end{array} \right \}
    \end{equation*}
for schemes $S$ over $\Spf \mc{O}_{\breve{F}_p}$. The similitude factor $b$ is allowed to vary (and is only required to be locally constant). If $\mc{M} = \mc{Z}(T)$, we set $\breve{\mc{M}}_{\mrm{framed}} \coloneqq \breve{\mc{Z}}(T)_{\mrm{framed}}$. There is a canonical forgetful map
    \begin{equation}\label{equation:non-Arch_uniformization:framed_stack:forgetful}
    \Theta \colon \breve{\mc{Z}}(T)_{\mrm{framed}} \ra \breve{\mc{Z}}(T)
    \end{equation}
sending $(\a, \phi_0, \phi) \mapsto \a$. This will be the uniformization map (Section \ref{ssec:non-Arch_uniformization:proof}). 

There is a canonical isomorphism
    \begin{equation}\label{equation:non-Arch_uniformization:framing:iso_framed_and_local_decomp}
    \breve{\mc{Z}}(T)_{\mrm{framed}} \xra{\sim} \coprod_{\substack{\underline{\mbf{x}} \in \mbf{W}^m \\ (\underline{\mbf{x}}, \underline{\mbf{x}}) = T}} \mc{Z}'(\underline{\mbf{x}}_p) \times \mc{Z}'(\underline{\mbf{x}}^p)
    \end{equation}
which we now describe. Here $\mc{Z}'(\underline{\mbf{x}}_p)$ is the local special cycle at $p$ from \crefext{I:ssec:moduli_pDiv:special_cycles}, and $\mc{Z}'(\underline{\mbf{x}}^p)$ is the away-from-$p$ local special cycle from \cref{ssec:non-Arch_uniformization:away_p_special_cycle}.

Consider $(\a, \phi_0, \phi) \in \breve{\mc{Z}}(T)_{\mrm{framed}}(S)$ as above. 
Passing to $p$-divisible groups gives a datum $(X_0, \iota_0, \lambda_0, X, \iota, \lambda)$ (e.g. $X$ is the $p$-divisible group of $A$), along with a framing quasi-isogeny $\rho \colon X_{\overline{S}} \ra \mbf{X}_{\overline{S}}$ induced by $\phi$ (where $\overline{S} \coloneqq  S_{\overline{k}}$) and similarly a framing $\rho_0$ induced by $\phi_0$.
We also obtain $g_0 \coloneqq \mf{y}_0 \circ \phi_{0,*} \circ \tilde{\eta}_0^{-1} \in G_0(\A_f^p) / K_{0,f}^p$ and $g \coloneqq \mf{y} \circ (\phi_0^{-1,*} \phi_{*}) \circ \tilde{\eta}^{-1} \in U(V)(\A_f^p) / K_f^p$ where $\phi_{0,*} \colon T^p(A_0)^0 \ra T^p(\mf{A}_0)^0$ and
    \begin{equation}
    \phi_0^{-1,*} \phi_* \colon \Hom_{F \otimes_{\Q} \A_f^p}(T^p(A_0)^0, T^p(A)^0) \ra \Hom_{F \otimes_{\Q} \A_f^p}(T^p(\mf{A}_0)^0, T^p(\mf{A})^0)
    \end{equation}
is pre- and post-composition (when $S$ is connected, pick any geometric point; there is no dependence on this choice). In general, $g_0$ and $g$ will be locally constant elements.
For any $\underline{x} \in \Hom_F(A_0, A)^m$ over a connected base $S$, we have $\phi \circ \underline{x} \circ \phi_0^{-1} \in \mbf{W}^m$ (canonically), by Mumford's rigidity lemma for morphisms of abelian schemes \cite[{Corollary 6.2}]{MFK94}. In the not-necessarily connected case, we obtain a locally constant element of $\mbf{W}^m$.

The above constructions give a map
    \begin{equation}\label{equation:non-Arch_uniformization:framing:iso_framed_and_local_decomp:ambient}
    \begin{tikzcd}[row sep = tiny]
    \breve{\mc{Z}}(T)_{\mrm{framed}} \arrow{r} & \mc{N}' \times G'(\A_f^p)/K'^p_f \times \mbf{W}^m
    \\
    (\a, \phi_0, \phi) \arrow[mapsto]{r} & ((X_0, \iota_0, \lambda_0, \rho_0, X, \iota, \lambda, \rho), (g_0, g_0 g), \phi \circ \underline{x} \circ \phi_0^{-1})
    \end{tikzcd}
    \end{equation}
which induces an isomorphism from $\breve{\mc{Z}}(T)_{\mrm{framed}}$ to the open and closed subfunctor
    \begin{equation}\label{equation:non-Arch_uniformization:framing:iso_framed_and_local_decomp:inclusion}
    \begin{tikzcd}
    \displaystyle\coprod_{\substack{\underline{\mbf{x}} \in \mbf{W}^m \\ (\underline{\mbf{x}}, \underline{\mbf{x}}) = T}} \mc{Z}'(\underline{\mbf{x}}_p) \times \mc{Z}'(\underline{\mbf{x}}^p) \arrow[hook]{r} &
    \mc{N}' \times G'(\A_f^p) / K'^p_f \times \mbf{W}^m.
    \end{tikzcd}
    \end{equation}
One can verify that the map in \eqref{equation:non-Arch_uniformization:framing:iso_framed_and_local_decomp} is an isomorphism by decomposing the kernels of framing quasi-isogenies (rescale to obtain an isogeny) into their $p$-power and $\ell$-power torsion subgroups. The isomorphism implies that $\breve{\mc{Z}}(T)_{\mrm{framed}}$ is a locally Noetherian formal scheme.

            \subsection{Quotient}
            \label{ssec:non-Arch_uniformization:quotient}
                Consider the algebraic groups
    \begin{align}
    I_0 \coloneqq GU(\mbf{V}_0)
    \quad \quad 
    I_1 \coloneqq U(\mbf{W}) \times U(\mbf{W}^{\perp})
    \quad \quad 
    I' \coloneqq I_0 \times I_1
    \end{align}
over $\Q$.
Unless specified otherwise, an element $(\gamma_0, \gamma) \in I'$ will mean $\gamma_0 \in I_0$ and $\gamma \in GU(\mbf{W}) \times GU(\mbf{W}^{\perp})$ with $\gamma_0^{-1} \gamma \in I_1$.

Uniformization will involve the stack quotient $[I'(\Q) \backslash \breve{\mc{Z}}(T)_{\mrm{framed}}]$ for an action of $I'(\Q)$ on $\breve{\mc{Z}}(T)_{\mrm{framed}}$, which we now describe.
For $\Q$-algebras $R$, there are canonical identifications
    \begin{align}
    I_0(R) & = \{\g_0 \in \End^0_F(\mbf{A}_0) \otimes_{\Q} R : \gamma_0^{\dagger} \gamma_0 \in R^{\times} \}
    \\
    I_1(R) & = \{\g \in \End^0_F(\mbf{A}) \otimes_{\Q} R : \gamma^{\dagger} \gamma = 1\} \label{equation:non-Arch_uniformization:quotient:I_group_unitary}
    \end{align}
(act on $\mbf{V}_0$ and $\mbf{V}$ by post-composition). View $I'(\Q)$ as a discrete group. Then $(\gamma_0, \g) \in I'(\Q)$ acts on $\breve{\mc{Z}}(T)_{\mrm{framed}}$ as $(\a, \phi_0, \phi) \mapsto (\a, \gamma_0 \circ \phi_0, \gamma \circ \phi)$.
We are abusing notation: the elements $\gamma_0$ and $\gamma$ lift (uniquely, by Mumford's rigidity lemma or Drinfeld rigidity and Serre--Tate) to quasi-endomorphisms of $\mf{A}_{0,S}$ and $\mf{A}_{S}$ respectively.

In terms of the isomorphism in \eqref{equation:non-Arch_uniformization:framing:iso_framed_and_local_decomp}, the action of $I'(\Q)$ on $\breve{\mc{Z}}(T)_{\mrm{framed}}$ admits the following (equivalent) description. By the isomorphism in \eqref{equation:non-Arch_uniformization:framing:at_p_Herm}, the group $I'(\Q_p)$ acts on $\mc{N}'$ (discussed in \crefext{I:ssec:moduli_pDiv:action_RZ}). By \eqref{equation:non-Arch_uniformization:quotient:I_group_unitary}, we have a faithful action of $I(\A_f^p)$ on
    \begin{equation}
    \Hom_{F \otimes \A_f^p}(T^p(\mbf{A}_0)^0, T^p(\mbf{A})^0) = V \otimes_{\Q} \A_f^p
    \end{equation}
by post-composition. This induces a homomorphism $I_1(\A_f^p) \ra U(V)(\A_f^p)$ and hence an action of $I'(\A_f^p)$ on $G'(\A_f^p) / K'^p_f$ (left multiplication).
The group $I'(\Q)$ also acts on $\mbf{W}$ by the projection $I'(\Q) \ra U(\mbf{W})$.

Hence $I'(\Q)$ acts on
    \begin{equation}
    \mc{N}' \times G'(\A_f^p) / K'^p_f \times \mbf{W}^m.
    \end{equation}
Under the inclusion \eqref{equation:non-Arch_uniformization:framing:iso_framed_and_local_decomp:inclusion}, this induces the same action on $\breve{\mc{Z}}(T)_{\mrm{framed}}$ described previously. Both descriptions will be useful for us.

We now form the (fppf) stack quotient
    \begin{equation}\label{equation:non-Arch_uniformization:quotient}
    \coprod_{\substack{\underline{\mbf{x}} \in \mbf{W}^m \\ (\underline{\mbf{x}}, \underline{\mbf{x}}) = T}} \mc{Z}'(\underline{\mbf{x}}_p) \times \mc{Z}'(\underline{\mbf{x}}^p)
    \ra
    \Biggl [ I'(\Q) \backslash \Biggl (\coprod_{\substack{\underline{\mbf{x}} \in \mbf{W}^m \\ (\underline{\mbf{x}}, \underline{\mbf{x}}) = T}} \mc{Z}'(\underline{\mbf{x}}_p) \times \mc{Z}'(\underline{\mbf{x}}^p) \Biggr ) \Biggr ]
    \end{equation}
The left-hand side is a locally Noetherian formal scheme, and the right-hand side is a locally Noetherian formal algebraic stack which is formally locally of finite type over $\Spf \mc{O}_{\breve{F}_p}$. The right-hand side is also $[I'(\Q) \backslash \breve{\mc{Z}}(T)_{\mrm{framed}}]$. The quotient map is representable by schemes, separated, \'etale, and surjective.

Using \crefext{I:equation:moduli_pDiv:special_cycles:unitary_iso} (and \crefext{I:equation:moduli_pDiv:discrete_reduced:split_Z(L)}) and \eqref{equation:non-Arch_uniformization:away_p_special_cycle:iso} (various incarnations of the isomorphism $G' \xra{\sim} GU(V_0) \times U(V)$) yields a canonical isomorphism from the left-hand side of \eqref{equation:non-Arch_uniformization:quotient} to
    \begin{equation}
    GU(V_0) / K_{0,f} \times \coprod_{\substack{ \underline{\mbf{x}} \in \mbf{W}^m \\ (\underline{\mbf{x}}, \underline{\mbf{x}}) = T}} \mc{Z}(\underline{\mbf{x}}_p) \times U(\mbf{W}_p^{\perp}) / K_{1, \mbf{L}^{\perp}_p} \times \mc{Z}(\underline{\mbf{x}^p})
    \end{equation}
where $K_{1, \mbf{L}^{\perp}_p} \subseteq U(\mbf{W}^{\perp}_p)$ is the unique maximal open compact subgroup (since $\mbf{W}^{\perp}$ has rank $0$ or $1$). This is a disjoint union of various local special cycles $\mc{Z}(\underline{\mbf{x}}_p)$, indexed by the (discrete) set
    \begin{equation}
    J_p(T) \coloneqq GU(V_0) / K_{0,f} \times \coprod_{\substack{ \underline{\mbf{x}} \in \mbf{W}^m \\ (\underline{\mbf{x}}, \underline{\mbf{x}}) = T}} U(\mbf{W}_p^{\perp}) / K_{1, \mbf{L}^{\perp}_p} \times \mc{Z}(\underline{\mbf{x}^p}).
    \end{equation}
In particular, every element $j \in J_p(T)$ defines a morphism
    \begin{equation}
    \Theta_j \colon \mc{Z}(\underline{\mbf{x}}_p) \ra [I'(\Q) \backslash \mc{Z}(T)_{\mrm{framed}}]
    \end{equation}
which is \'etale, separated, and representable by schemes. Given $j \in J_p(T)$, we let $\Aut(j) \subseteq I'(\Q)$ be the stabilizer for the action of $I'(\Q)$ on $J_p(T)$.

The right-hand side of \eqref{equation:non-Arch_uniformization:quotient} is then identified with
    \begin{align}
    \label{equation:non-Arch_uniformization:quotient:quotient_unitary_iso}
    \Biggl [ GU(V_0)(\Q) \backslash \Biggl ( GU(V_0)(\A_f) / K_{0,f} \Biggr )\Biggr ] \times \Biggl [ I_1(\Q) \backslash \Biggl ( \coprod_{\substack{\underline{\mbf{x}} \in \mbf{W}^m \\ (\underline{\mbf{x}}, \underline{\mbf{x}}) = T}} \mc{Z}(\underline{\mbf{x}}_p) \times U(\mbf{W}^{\perp}_p) / K_{1, \mbf{L}^{\perp}_p} \times \mc{Z}(\underline{\mbf{x}}^p) \Biggr ) \Biggr ].
    \end{align}

We have
    \begin{equation}\label{equation:non-Arch_uniformization:quotient:0_stacky_degree}
    \deg [ GU(V_0)(\Q) \backslash ( GU(V_0)(\A_f) / K_{0,f} ) ] = [K_{L_0,f} : K_{0,f}] \cdot h_F / |\mc{O}_F^{\times}|
    \end{equation}
where the left-hand side denotes (stacky) groupoid cardinality, where $[K_{L_0, f} : K_{0,f}]$ is the index of $K_{0,f}$ in $K_{L_0, f}$, and $h_F$ is the class number of $\mc{O}_F$. In the case where $\rank(T) \geq n - 1$ (we have already assumed $\rank(T) \geq n - 1$ if $p$ is split), the groupoid
    \begin{equation}\label{equation:non-Arch_uniformization:quotient:unitary_stacky_degree}
    \Biggl [ I_1(\Q) \backslash \Biggl ( \coprod_{\substack{\underline{\mbf{x}} \in \mbf{W}^m \\ (\underline{\mbf{x}}, \underline{\mbf{x}}) = T}} U(\mbf{W}^{\perp}_p) / K_{1, \mbf{L}^{\perp}_p} \times \mc{Z}(\underline{\mbf{x}}^p) \Biggr ) \Biggr ]
    \end{equation}
has finite automorphism groups and finitely many isomorphism classes, and its groupoid degree is essentially a product of special values of local Whittaker functions away from $p$ \crefext{IV:lemma:local_Siegel-Weil:uniformization_degree:main}.

In the case $\rank(T) \geq n - 1$, the map $\Theta_j$ associated with any $j \in J_p(T)$ is thus representable by schemes and finite \'etale of constant degree $\deg \Theta_j = |\Aut(j)|$.
    
            \subsection{Uniformization}
            \label{ssec:non-Arch_uniformization:proof}
                We explain how the uniformization morphism $\Theta \colon \breve{\mc{Z}}(T)_{\mrm{framed}} \ra \breve{\mc{Z}}(T)$ \eqref{equation:non-Arch_uniformization:framed_stack:forgetful} descends to an isomorphism of locally Noetherian formal algebraic stacks
    \begin{equation}\label{equation:non-Arch_uniformization:proof:isomorphism}
    \tilde{\Theta} \colon \Biggl [ I'(\Q) \backslash \Biggl ( \coprod_{\substack{\underline{\mbf{x}} \in \mbf{W}^m \\ (\underline{\mbf{x}}, \underline{\mbf{x}}) = T}} \mc{Z}'(\underline{\mbf{x}}_p) \times \mc{Z}'(\underline{\mbf{x}}^p) \Biggr ) \Biggr] \xra{\sim} \breve{\mc{Z}}(T).
    \end{equation} 
The main point is surjectivity on $\overline{k}$-points via the Hasse principle (Lemma \ref{lemma:non-Arch_uniformization:proof:kbar_surjective}).
    
When $p$ is split, we will allow a change of choice of framing data $(\mbf{A}_0, \iota_{\mbf{A}_0}, \lambda_{\mbf{A}_0}, \mbf{A}, \iota_{\mbf{A}}, \lambda_{\mbf{A}}, \tilde{\pmb{\eta}}_0, \tilde{\pmb{\eta}})$, $\pmb{\eta}_0$, and $\pmb{\eta}$, possibly depending on $T$.

\begin{lemma}\label{lemma:non-Arch_uniformization:proof:monomorphism}
The map $\Theta \colon \breve{\mc{Z}}(T)_{\mrm{framed}} \ra \breve{\mc{Z}}(T)$ factors uniquely through a monomorphism\footnote{By a \emph{monomorphism} of formal algebraic stacks, we mean a morphism which is fully faithful on underlying fibered categories.}
    \begin{equation}\label{equation:non-Arch_uniformization:proof:tilde_theta}
    \tilde{\Theta} \colon [I'(\Q) \backslash \breve{\mc{Z}}(T)_{\mrm{framed}}] \ra \breve{\mc{Z}}(T)
    \end{equation}
of formal algebraic stacks. The map $\tilde{\Theta}$ is formally locally of finite type and formally \'etale.
\end{lemma}
\begin{proof}
Suppose $(\a, \phi_0, \phi)$ and $(\a', \phi'_0, \phi')$ are objects of $\breve{\mc{Z}}(T)_{\mrm{framed}}(S)$, and suppose $f' \colon \a \ra \a'$ is an isomorphism of objects in the groupoid $\breve{\mc{Z}}(T)(S)$ (for some base scheme $S$). We claim there is a unique $\gamma' = (\gamma_0, \gamma) \in I'(\Q)$ such that $f'$ induces an isomorphism $\gamma' \cdot (\a, \phi_0, \phi) \xra{\sim} (\a', \phi'_0, \phi')$ in the setoid $\breve{\mc{Z}}(T)_{\mrm{framed}}(S)$.

The map $f$ is given by a pair of isomorphisms $f_0 \colon A_0 \ra A'_0$ and $f \colon A \ra A'$ (where $\a = (A_0, \ldots)$ and $\a' = (A'_0, \ldots)$, with notation as above). Then we take $\gamma_0 = \phi'_0 \circ f_0 \circ \phi_0^{-1}$ and $\gamma = \phi' \circ f \circ \phi^{-1}$. Hence $\tilde{\Theta}$ is a monomorphism.

The map $\tilde{\Theta}$ is a map between locally Noetherian formal algebraic stacks which are formally locally of finite type over $\Spf \mc{O}_{\breve{F}_p}$, so $\tilde{\Theta}$ is formally locally of finite type.

The property of being formally \'etale may be checked ``formally \'etale locally on the source''.
The quotient map $\breve{\mc{Z}}(T)_{\mrm{framed}} \ra [I'(\Q) \backslash \breve{\mc{Z}}(T)_{\mrm{framed}}]$ is representable by schemes and formally \'etale, so it is enough to check that $\Theta \colon \breve{\mc{Z}}(T)_{\mrm{framed}} \ra \breve{\mc{Z}}(T)$ is formally \'etale. This property amounts to the following rigidity statement for abelian schemes: given any first order thickening of schemes $T \ra T'$ on which $p$ is locally nilpotent, and given abelian schemes $A_1$ and $A_2$ over $T'$, any quasi-homomorphism $A_{1,T} \ra A_{2,T}$ lifts uniquely to a quasi-homomorphism $A_1 \ra A_2$ (e.g. by Drinfeld rigidity and Serre--Tate).
\end{proof}

\begin{lemma}\label{lemma:non-Arch_uniformization:proof:kbar_surjective}
The map $\Theta(\overline{k}) \colon \breve{\mc{Z}}(T)_{\mrm{framed}}(\overline{k}) \ra \breve{\mc{Z}}(T)(\overline{k})$ (on groupoids of $\overline{k}$-points) is surjective (resp. surjective for some choice of framing data) on isomorphism classes if $p$ is nonsplit (resp. split).
\end{lemma}
\begin{proof}
If $\breve{\mc{Z}}(T)$ is empty, there is nothing to show, so assume $\breve{\mc{Z}}(T)$ is nonempty. If $p$ is split, we can change the framing object to assume it extends to $(\mbf{A}_0, \iota_{\mbf{A}_0}, \lambda_{\mbf{A}_0}, \mbf{A}, \iota_{\mbf{A}}, \lambda_{\mbf{A}}, \tilde{\pmb{\eta}}_0, \tilde{\pmb{\eta}}, \underline{\mbf{x}}) \in \breve{\mc{Z}}(T)(\overline{k})$ (i.e. $\underline{\mbf{x}} \in \mbf{W}^m$ with $(\underline{\mbf{x}}, \underline{\mbf{x}}) = T$). This implies that $T$ has rank $n - 1$ if $p$ is split (we already assumed $\rank(T) \geq n - 1$ if $p$ is split, then see Remark \ref{remark:arith_cycle_clases:degrees:split_full_rank_empty}). We still know that $\mbf{A}$ is $\mc{O}_F$-linearly isogenous to $\mbf{A}_0^{n - 1} \times \mbf{A}_0^{\s}$ (Lemma \ref{lemma:corank_1_isogeny_decomposition}).

In all cases, the task is to show that any $(\mbf{A}'_0, \iota_{\mbf{A}'_0}, \lambda_{\mbf{A}'_0}, \mbf{A}', \iota_{\mbf{A}'}, \lambda_{\mbf{A}'}, \tilde{\pmb{\eta}}_0, \tilde{\pmb{\eta}}, \underline{\mbf{x}}') \in \breve{\mc{Z}}(T)(\overline{k})$ admits a framing $(\phi_0, \phi)$, i.e. quasi-isogenies $\phi_0 \colon \mbf{A}'_0 \ra \mbf{A}_0$ and $\phi \colon \mbf{A}' \ra \mbf{A}$ which preserve quasi-polarizations up to the same scalar in $\Q_{>0}$.

Fix any $\mc{O}_F$-linear isogeny $\phi_0 \colon \mbf{A}'_0 \ra \mbf{A}_0$, which exists because $\mbf{A}'_0$ and $\mbf{A}_0$ are elliptic curves with $\mc{O}_F$-action of the same signature (see the proof of Lemma \ref{lemma:corank_1_isogeny_decomposition}). Let $b \in \Q_{>0}$ be such that $\phi_0^* \lambda_{\mbf{A}_0} = b \lambda_{\mbf{A}'_0}$. Set $\mbf{W}' \coloneqq \Hom^0_{F}(\mbf{A}'_0, \mbf{A}')$
with the Hermitian pairing $(x,y) = x^{\dagger} y$. 

\textit{Case $p$ is nonsplit:} There is an isomorphism of $F$ vector spaces
    \begin{equation}
    \begin{tikzcd}[row sep = tiny]
    \Hom^0_F(\mbf{A}', \mbf{A}) \arrow{r} & \Hom_F(\mbf{W}', \mbf{W}) 
    \\
    \phi \arrow[mapsto]{r} & (f \mapsto \phi \circ f \circ \phi_0^{-1}). 
    \end{tikzcd}
    \end{equation}
An element $\phi \in \Hom^0_F(\mbf{A}', \mbf{A})$ satisfies $\phi^{\dagger} \phi = b$ if and only if $\phi$ corresponds to an isomorphism of Hermitian spaces $\mbf{W}' \ra \mbf{W}$. But we have $\mbf{W}' \otimes_{\Q} \A_f^p \cong \mbf{W} \otimes_{\Q} \A_f^p \cong V \otimes_{\Q} \A_f^p$ as Hermitian spaces, we have $\varepsilon(\mbf{W}') = \varepsilon(\mbf{W}) = (-1)^r$, and we have $\mbf{W}'_{\mbf{R}} \cong \mbf{W}_{\R}$ (both are positive definite of rank $n$). So we have $\mbf{W}' \cong \mbf{W}$ as Hermitian spaces, by the Hasse principle for Hermitian spaces (Landherr's theorem).

\textit{Case $p$ is split:} Fix $\mc{O}_F$-linear isogenies $\mbf{B} \times \mbf{B}^{\perp} \ra \mbf{A}$ and $\mbf{B}' \times \mbf{B}^{\prime \perp} \ra \mbf{A}'$, where $\mbf{B} \cong \mbf{A}_0^{n -1}$, $\mbf{B}^{\perp} \cong \mbf{A}_0^{\s}$, $\mbf{B}' \cong \mbf{A}_0^{\prime n - 1}$, and $\mbf{B}^{\prime \perp} \cong \mbf{A}_0^{\prime \s}$.
Equip $\mbf{B} \times \mbf{B}^{\perp}$ and $\mbf{B}' \times \mbf{B}^{\prime \perp}$ with the quasi-polarizations pulled back from $\lambda_{\mbf{A}}$ and $\lambda_{\mbf{A}'}$ on $\mbf{A}$ and $\mbf{A}'$, respectively.

Any $F$-linear quasi-isogeny $\phi \colon \mbf{A}' \ra \mbf{A}$ decomposes as a product of quasi-isogenies $\mbf{B}' \ra \mbf{B}$ and $\mbf{B}^{\prime \perp} \ra \mbf{B}^{\perp}$, since $\Hom_F(\mbf{B}', \mbf{B}^{\perp}) = \Hom_F(\mbf{B}^{\prime \perp}, \mbf{B}) = 0$ (because of the opposite signatures). We write $\phi^{\perp} \colon \mbf{B}^{\prime \perp} \ra \mbf{B}^{\perp}$ for the quasi-isogeny induced by $\phi$. By similar reasoning, the quasi-polarization on $\mbf{B} \times \mbf{B}^{\perp}$ is the product of a quasi-polarization on $\mbf{B}$ and a quasi-polarization on $\mbf{B}^{\perp}$.

There is an isomorphism of $F$ vector spaces
    \begin{equation}
    \begin{tikzcd}[row sep = tiny]    
    \Hom^0_F(\mbf{A}', \mbf{A}) \arrow{r} & \Hom_F(\mbf{W}', \mbf{W}) \times \Hom^0_F(\mbf{B}^{\prime \perp}, \mbf{B}^{\perp})
    \\
    \phi \arrow[mapsto]{r} & (f \mapsto (\phi \circ f \circ \phi_0^{-1})), \phi^{\perp}.
    \end{tikzcd}
    \end{equation}
An element $\phi \in \Hom^0_F(\mbf{A}', \mbf{A})$ satisfies $\phi^{\dagger} \phi = b$ if and only if $\phi$ corresponds to an isomorphism of Hermitian spaces $\mbf{W}' \ra \mbf{W}$ and with $\phi^{\perp \dagger} \phi^{\perp} = b$.

We have $\underline{\mbf{x}}' \in \mbf{W}'$ and $\underline{\mbf{x}} \in \mbf{W}$ with $(\underline{\mbf{x}}', \underline{\mbf{x}}') = (\underline{\mbf{x}}, \underline{\mbf{x}}) = T$. Since $\rank(T) = \rank(\mbf{W}') = \rank(\mbf{W}) = n - 1$, this implies $\mbf{W}' \cong \mbf{W}$ as Hermitian spaces.

For every prime $\ell \neq p$, the natural map
    \begin{equation}
    \begin{tikzcd}
    \Hom^0_F(\mbf{B}^{\prime \perp}, \mbf{B}^{\perp}) \otimes_{\Q} \Q_{\ell} \arrow{r} & \Hom_{F_{\ell}}(T_{\ell}(\mbf{B}^{\prime \perp})^0, T_{\ell}(\mbf{B}^{\perp})^0)
    \end{tikzcd}
    \end{equation}
is an isomorphism of (one-dimensional) Hermitian spaces. If we set $\mbf{U}_{\ell} \coloneqq \Hom_{F_{\ell}}(T_{\ell}(\mbf{A}_0)^0, T_{\ell}(\mbf{B}^{\perp})^0)$ and $\mbf{U}'_{\ell} \coloneqq \Hom_{F_{\ell}}(T_{\ell}(\mbf{A}_0')^0, T_{\ell}(\mbf{B}^{\prime \perp})^0)$, there is an isomorphism of $F_{\ell}$ vector spaces
    \begin{equation}
    \begin{tikzcd}[row sep = tiny]
    \Hom_{F_{\ell}}(T_{\ell}(\mbf{B}^{\prime \perp})^0, T_{\ell}(\mbf{B}^{\perp})^0) \arrow{r} & \Hom_{F_{\ell}}(\mbf{U}_{\ell}, \mbf{U}'_{\ell})
    \\
    \phi^{\perp} \arrow[mapsto]{r} & (f \mapsto (\phi^{\perp} \circ f \circ \phi_0^{-1})).
    \end{tikzcd}
    \end{equation}
An element $\phi^{\perp}$ on the left satisfies $\phi^{\perp \dagger} \phi^{\perp} = b$ if and only if $\phi^{\perp}$ corresponds to an isomorphism of Hermitian spaces $\mbf{U}_{\ell} \ra \mbf{U}'_{\ell}$. We have $V_{\ell} \cong \mbf{W}'_{\ell} \oplus \mbf{U}'_{\ell} \cong \mbf{W}_{\ell} \oplus \mbf{U}_{\ell}$ (orthogonal direct sum) as Hermitian spaces, for all $\ell \neq p$. Hence $\mbf{U}'_\ell \cong \mbf{U}_{\ell}$ for all $\ell \neq p$ (consider the Hermitian space local invariants (in $\{ \pm 1\}$) via
$\varepsilon$, with conventions as in \crefext{I:ssec:Hermitian_conventions:lattices}).

The preceding discussion produces an element $\phi^{\perp}_{\ell} \in \Hom^0_F(\mbf{B}^{\prime \perp}, \mbf{B}^{\perp}) \otimes_{\Q} \Q_{\ell}$ satisfying $\phi^{\perp \dagger}_{\ell} \phi^{\perp}_{\ell} = b$ for all primes $\ell \neq p$. Since $p$ is split in $\mc{O}_F$, such an element exists for $\ell = p$ as well (i.e. $N_{F_p/\Q_p}(F_p^{\times}) = \Q_p^{\times}$). Since $b > 0$, such an element also exists if $\Q_{\ell}$ is replaced by $\R$ (positivity of the Rosati involution). By the Hasse principle for Hermitian spaces (or Hasse norm theorem), we obtain $\phi^{\perp} \in \Hom^0_F(\mbf{B}^{\prime \perp}, \mbf{B}^{\perp})$ satisfying $\phi^{\perp \dagger} \phi^{\perp} = b$.
\end{proof}

For the rest of Section \ref{sec:non-Arch_uniformization}, we fix framing data as in Lemma \ref{lemma:non-Arch_uniformization:proof:kbar_surjective} if $p$ is split, so that $\Theta(\overline{k})$ is surjective.

For the supersingular cases, we use the following lifting result to prove surjectivity of $\tilde{\Theta}$ by bootstrapping from surjectivity on $\overline{k}$ points (as in the proof of \cite[{Theorem 6.30}]{RZ96}). Recall that a $p$-divisible group $X$ over a base scheme $S$ is said to be \emph{isoclinic} if for any geometric point $\overline{s}$ of $S$, the isocrystal of $X_{\overline{s}}$ has constant slope independent of $\overline{s}$.

\begin{proposition}[Isoclinic lifting theorem]\label{proposition:isoclinic_lifting_theorem}
For any integer $h$, there exists an integer $c$ with the following property: Let $R$ be a reduced Noetherian Henselian local ring with residue field $\kappa$, and assume that $R$ is an $\F_p$-algebra. Let $X$ and $Y$ be isoclinic $p$-divisible groups of heights $\leq h$ over $\Spec R$. For any homomorphism $f \colon X_{\kappa} \ra Y_{\kappa}$, the homomorphism $p^c f$ lifts to a unique homomorphism $X \ra Y$. 
\end{proposition}
\begin{proof}
See \cite[{Corollary 3.4}]{OZ02}. For the statement when $R = \psring{\kappa}{t}$ for an algebraically closed field $\kappa$ (which is enough for Lemma \ref{lemma:non-Arch_uniformization:proof:surjection}), see also \cite[{Theorem 2.7.1}]{Katz79} combined with Grothendieck--Messing theory as in \cite[{pg. 295}]{RZ96}.
\end{proof}

\begin{lemma}\label{lemma:non-Arch_uniformization:proof:surjection}
The uniformization map $\Theta$ is a surjection\footnote{We say a morphism of formal algebraic stacks is a \emph{surjection} if it is surjective on underlying topological spaces.} of formal algebraic stacks.
\end{lemma}
\begin{proof} The reduced substack $\breve{\mc{Z}}(T)_{\red} \subseteq \breve{\mc{Z}}(T)$ is Jacobson, Deligne--Mumford, with quasi-compact diagonal, and finite type over $\Spec \overline{k}$. This implies that the closed points of $\breve{\mc{Z}}(T)_{\red}$ are dense in every closed subset (e.g. \cite[\href{https://stacks.math.columbia.edu/tag/06G2}{Lemma 06G2}]{stacks-project}; the finite type points are the same as closed points here), each closed point is the image of a map $\Spec \overline{k} \ra \breve{\mc{Z}}(T)_{\red}$, and every such map has image being a closed point.

\textit{Case $p$ is nonsplit:} We already know that $\Theta$ is surjective on $\overline{k}$ points. It is thus enough to prove the following claim: suppose $\a' \rightsquigarrow \a$ is an immediate specialization of points in $|\breve{\mc{Z}}(T)|$ (in the sense of \cite[\href{https://stacks.math.columbia.edu/tag/02I9}{Definition 02I9}]{stacks-project}, i.e. $\a$ is a point of ``codimension one'' in the closure of $\a'$). If $\a$ is in the image of $\Theta$, we claim that $\a'$ is also in the image of $\Theta$.
(This specialization process eventually terminates with a $\overline{k}$ point.)

Let $\kappa$ an algebraically closed field with a morphism $\Spec \psring{\kappa}{t} \ra \breve{\mc{Z}}(T)$, which sends the closed point to $\a$ and the open point to $\a'$.\footnote{The following procedure produces such a morphism $\Spec \psring{\kappa}{t} \ra \breve{\mc{Z}}(T)$. First, take an \'etale cover of $\breve{\mc{Z}}(T)_{\red}$ by a scheme $U$ and lift $x' \rightsquigarrow x$ to an immediate specialization $y' \rightsquigarrow y$ on $U$. Write $Z$ for the normalization of the integral closed subscheme of $U$ with generic point $y'$. Note the normalization map is finite, and lift $y' \rightsquigarrow y$ to an immediate specialization $z' \rightsquigarrow z$ on $Z$. Completion of the local ring at $z$ on $Z$ is a power series ring over a field.} Enlarging $\kappa$ if necessary, we may lift $\a$ to a point $(\a, \phi_0, \phi) \in \breve{\mc{Z}}(T)_{\mrm{framed}}$. The task is to lift the framing pair $(\phi_0, \phi)$ to $\Spec \psring{k}{t}$, which is then a framing pair for $\a'$. Serre--Tate (and formal GAGA as in \cite[{Th\'eor\`eme 5.4.1}]{EGAIII1}) implies that it is enough to lift the induced quasi-isogenies of $p$-divisible groups to $\Spec \psring{k}{t}$. This is possible by the isoclinic lifting theorem (Proposition \ref{proposition:isoclinic_lifting_theorem}).

\textit{Case $p$ is split:} By \crefext{I:lemma:moduli_pDiv:discrete_reduced:local_special_cycle_red_finite} (finiteness of local special cycles), and since the groupoid $[I'(\Q) \backslash J_p(T)]$ has finite automorphism groups and finitely many isomorphism classes (Section \ref{ssec:non-Arch_uniformization:quotient}; we assumed $\rank(T) \geq n - 1$ for $p$ split), we know there is a surjection from finitely many copies of $\Spec \overline{k}$ to $[I'(\Q) \backslash \breve{\mc{Z}}(T)_{\mrm{framed}}]$. Since $\Theta$ is surjective on $\overline{k}$-points, the previous considerations show that $|\breve{\mc{Z}}(T)_{\red}|$ is a finite discrete topological space, and that $\Theta$ is a surjection.
\end{proof}

\begin{lemma}\label{lemma:non-Arch_uniformization:proof:proper}
The map $\tilde{\Theta}$ is proper on underlying reduced substacks, and the reduced substack $[I'(\Q) \backslash \breve{\mc{Z}}(T)_{\mrm{framed}}]_{\red}$ is proper over $\Spec \overline{k}$.
\end{lemma}
\begin{proof}
Since the reduced substack $\breve{\mc{Z}}(T)_{\red}$ is separated over $\Spec \overline{k}$, it is enough to check that $[I'(\Q) \backslash \breve{\mc{Z}}(T)_{\mrm{framed}}]_{\mrm{red}}$ is proper over $\Spec \overline{k}$, by \cite[\href{https://stacks.math.columbia.edu/tag/0CPT}{Lemma 0CPT}]{stacks-project}.

We already saw that $\tilde{\Theta}$ is a monomorphism, hence separated.
Since $\breve{\mc{Z}}(T)_{\red}$ is separated, we see that $[I'(\Q) \backslash \breve{\mc{Z}}(T)_{\mrm{framed}}]_{\mrm{red}}$ is also separated over $\Spec \overline{k}$.

We use the description of $\breve{\mc{Z}}(T)_{\mrm{framed}}$ in \eqref{equation:non-Arch_uniformization:framing:iso_framed_and_local_decomp}. 
We know that every irreducible component of the reduced subscheme $\mc{N}'_{\red}$ is projective over $\overline{k}$ \crefext{I:ssec:moduli_pDiv:RZ}, hence the same holds for $\mc{Z}'(\underline{\mbf{x}}_p)$ for any $\underline{\mbf{x}} \in \mbf{W}^m$ (and $\mc{Z}'(\underline{\mbf{x}}^p)$ is discrete). 
Hence each irreducible component of $\breve{\mc{Z}}(T)_{\mrm{framed}, \mrm{red}}$ has closed image in $\breve{\mc{Z}}(T)_{\mrm{red}}$.
Since $\Theta$ is surjective, we conclude that finitely many irreducible components of $\breve{\mc{Z}}(T)_{\mrm{framed}, \mrm{red}}$ cover $\breve{\mc{Z}}(T)_{\mrm{red}}$ (by Noetherianity of the latter). Since $\tilde{\Theta}$ is a monomorphism, hence injective on underlying topological spaces, we conclude that those finitely many irreducible components cover $[I'(\Q) \backslash \breve{\mc{Z}}(T)_{\mrm{framed}}]_{\mrm{red}}$ as well. Then $[I'(\Q) \backslash \breve{\mc{Z}}(T)_{\mrm{framed}}]_{\mrm{red}}$ is proper over $\Spec \overline{k}$ by \cite[\href{https://stacks.math.columbia.edu/tag/0CQK}{Lemma 0CQK}]{stacks-project}.
\end{proof}

\begin{proposition}\label{proposition:non-Arch_uniformization:proof}
The map $\tilde{\Theta}$ is an isomorphism.
\end{proposition}
\begin{proof}
We have seen that the morphism $\tilde{\Theta}$ of locally Noetherian formal algebraic stacks is formally \'etale, surjective, and a monomorphism. The underlying map of reduced substacks is proper. These properties imply that $\tilde{\Theta}$ is an isomorphism.
\end{proof}

            \subsection{Global and local}
            \label{ssec:non-Arch_uniformization:global_to_local}
                The next lemma (purely linear-algebraic) helps us use uniformization to deduce properties of local special cycles via ``approximating'' them by global special cycles.

\begin{lemma}\label{lemma:non-Arch_uniformization:local_global_approximation}
Let $L \subseteq \mbf{W}_p$ be any non-degenerate Hermitian $\mc{O}_{F_p}$-lattice (of any rank). There exists an element $g_p \in U(\mbf{W}_p)$ such that $g_p(L)$ admits a basis consisting of elements in $\mbf{W}$.
\end{lemma}
\begin{proof}
Set $W = L \otimes_{\mc{O}_{F_p}} F_p$. It is enough to produce $g_p \in U(\mbf{W}_p)$ such that $g_p(W)$ admits an $F_p$-basis consisting of elements in $\mbf{W}$ (since this implies that every full rank $\mc{O}_{F_p}$-lattice in $g_p(W)$ admits a basis consisting of elements of $\mbf{W}$).

Select any $F_p$-basis $\underline{e} =[e_1, \ldots, e_d]$ for $W$. Since $\mbf{W}$ is dense in $\mbf{W}_p$, we may select $\underline{\tilde{e}} = [\tilde{e}_1, \ldots, \tilde{e}_d]$ such that each $\norm{\tilde{e}_i - e_i}_p \ll 1$ for all $i$ (meaning $\tilde{e}_i - e_i$ lies in a small neighborhood of $0$ for the $p$-adic topology on $\mbf{W}_p$). Set $\tilde{W} \coloneqq \mrm{span}_{F_p}\{\tilde{e}_1, \ldots, \tilde{e}_d\}$. When each $\tilde{e}_i - e_i$ lies in a sufficiently small neighborhood of $0$, there exists a (non-canonical) isomorphism of Hermitian spaces $W \cong \tilde{W}$ (the associated Gram matrices $(\underline{e}, \underline{e})$ and $(\underline{\tilde{e}}, \underline{\tilde{e}})$ can be made arbitrarily $p$-adically close; hence the local invariants $\varepsilon((\underline{e}, \underline{e}))$ and $\varepsilon((\underline{\tilde{e}}, \underline{\tilde{e}}))$ will agree). By Witt's theorem for Hermitian spaces, any isometry $W \ra \tilde{W}$ extends to an isometry $g_p \colon \mbf{W}_p \ra \mbf{W}_p$. This element $g_p \in U(\mbf{W}_p)$ satisfies the conditions in the lemma statement.
\end{proof}

\begin{corollary}\label{corollary:non-Arch_uniformization:vertical:approximate_cycle_RZ}
Consider any tuple $\underline{\mbf{x}}_p \in \mbf{W}_p^m$ which spans a non-degenerate Hermitian $\mc{O}_{F_p}$-lattice, and write $m^{\flat}$ for its rank. Assume $m^{\flat} = n - 1$ if $p$ is split.

For some $T \in \mrm{Herm}_m(\Q)$ (still assuming $\rank T \geq n - 1$ if $p$ is split), and some $j \in J_p(T)$ with associated $\underline{\mbf{w}} \in \mbf{W}$, there exists $g_p \in U(\mbf{W}_p)$ inducing an automorphism $\mc{N} \ra \mc{N}$ which takes $\mc{Z}(\underline{\mbf{x}}_p)$ isomorphically to $\mc{Z}(\underline{\mbf{w}}_p)$. In particular, there is an induced morphism
    \begin{equation}
    \mc{Z}(\underline{\mbf{x}}_p) \xra{\sim} \mc{Z}(\underline{\mbf{w}}_p) \xra{\Theta_j} \breve{\mc{Z}}(T)
    \end{equation}
which is representable by schemes, separated, and \'etale.
If $m^{\flat} \geq n - 1$ (equivalently, $\rank(T) \geq n - 1$), this map is finite \'etale.
\end{corollary}
\begin{proof}
By Lemma \ref{lemma:non-Arch_uniformization:local_global_approximation}, we may pick an element $g_p \in U(\mbf{W}_p)$ so that $\mrm{span}_{\mc{O}_{F_p}}( g_p \cdot \underline{\mbf{x}}_p)$ admits an $\mc{O}_{F_p}$-basis $\underline{\mbf{w}}^{\flat}$ of elements in $\mbf{W}$. Extend $\underline{\mbf{w}}^{\flat}$ to any $m$-tuple $\underline{\mbf{w}} \in \mbf{W}^m$, and set $T \coloneqq (\underline{\mbf{w}}, \underline{\mbf{w}})$. Recall that $U(\mbf{W}_p)$ acts on $\mc{N}$, and that $g_p$ gives an automorphism of $\mc{N}$ sending $\mc{Z}(\underline{\mbf{x}}_p) \mapsto \mc{Z}(\underline{\mbf{w}}_p)$ \crefext{I:ssec:moduli_pDiv:action_RZ}. 
By uniformization (Proposition \ref{proposition:non-Arch_uniformization:proof}), any $j \in J_p(T)$ whose associated tuple is $\underline{\mbf{w}}$ will satisfy the conditions of the lemma.
Replacing $\underline{\mbf{w}}$ with $a \cdot \underline{\mbf{w}}$ for suitable $a \in \Z$ with $p \nmid a$ ensures $\mc{Z}'(\underline{\mbf{w}}^p) \neq \emptyset$. Then such $j \in J_p(T)$ will exist.
In Section \ref{ssec:non-Arch_uniformization:quotient}, we saw that $\Theta_j$ is finite \'etale if $\rank(T) \geq n - 1$.
\end{proof}

If $p \neq 2$ and in signature $(n - 1, 1)$, the quasi-compactness proved in the next lemma is also \cite[{Lemma 2.9.}]{LZ22unitary} (inert), proved via Bruhat--Tits stratification. In the exotic smooth ramified case, quasi-compactness should be implicit in \cite{LL22II}, via Bruhat--Tits stratification as discussed in \cite[{\S 2.3}]{LL22II}. In the case when $\underline{\mbf{x}}_p$ spans a lattice of rank $n$ and signature $(n - 1, 1)$, see \cite[{Lemma 5.1.1}]{LZ22unitary} (inert, $p \neq 2$) and \cite[{Remark 2.26}]{LL22II} (ramified, exotic smooth).

\begin{lemma}\label{lemma:non-Arch_uniformization:vertical:local_cycle_quasi-compact}
Let $\underline{\mbf{x}}_p \in \mbf{W}_p^d$ be any tuple which spans a non-degenerate Hermitian $\mc{O}_{F_p}$-lattice of rank $\geq n - 1$. Then the local special cycle $\mc{Z}(\underline{\mbf{x}}_p)$ is quasi-compact and 
the structure map
$\mc{Z}(\underline{\mbf{x}}_p) \ra \Spf \mc{O}_{\breve{F}_p}$ is adic and proper.
\end{lemma}
\begin{proof}
By Corollary \ref{corollary:non-Arch_uniformization:vertical:approximate_cycle_RZ}, we obtain $T \in \mrm{Herm}_m(\Q)$ and a map $\mc{Z}(\underline{\mbf{x}}_p) \ra \breve{\mc{Z}}(T)$ which is representable by schemes, and finite \'etale. In particular, $\mc{Z}(\underline{\mbf{x}}_p)$ is quasi-compact because the (base-changed) global special cycle $\breve{\mc{Z}}(T)$ is quasi-compact.

If $p$ is nonsplit, then $\mc{Z}(T)_{\overline{k}} \ra \mc{M}_{\overline{k}}$ automatically factors through the supersingular locus (Corollary \ref{corollary:special_cycle_supersingular_support}), so we have $\breve{\mc{Z}}(T) = \mc{Z}(T)_{\Spf \mc{O}_{\breve{F}_p}}$.
This formula holds for $p$ split as well, by definition. 
Hence $\breve{\mc{Z}}(T) \ra \Spf \mc{O}_{\breve{F}_p}$ is adic\footnote{We say a morphism of formal algebraic stacks is \emph{adic} if the morphism is representable by algebraic stacks in the sense of \cite[{\S 3}]{Emerton20}.} and proper (Lemma \ref{lemma:global_special_cycles_proper}), 
so $\mc{Z}(\underline{\mbf{x}}_p) \ra \Spf \mc{O}_{\breve{F}_p}$ is adic and proper.
\end{proof}

We write $\breve{\mc{Z}}(T)_{\ms{H}}$ for the flat\footnote{Flatness for morphisms of locally Noetherian formal algebraic stacks was defined in \cite[{Definition 8.42}]{Emerton20}. We are using a different definition, since the definition of loc. cit. does not recover the usual notion of flatness for morphisms of schemes (in the situation of \cite[{Lemma 8.41(1)}]{Emerton20}, consider $\mc{X} = \mc{Y} = \Spec k$ for a field $k$ and any non-flat morphism of Noetherian affine $k$-schemes $U \ra V$).

We define flatness in the style of \cite[\href{https://stacks.math.columbia.edu/tag/06FL}{Section 06FL}]{stacks-project} (there for algebraic stacks), which recovers usual flatness for morphisms of locally Noetherian formal schemes.
Let $f \colon X \ra Y$ be a morphism of locally Noetherian formal algebraic spaces. Consider commutative diagrams
    \begin{equation}\label{equation:formal_algebraic_space_flat}
    \begin{tikzcd}[ampersand replacement = \&]
    U \arrow{r}{h} \arrow{d}{a} \& V \arrow{d}{b} \\
    X \arrow{r}{f} \& Y
    \end{tikzcd}
    \end{equation}
where $U$ and $V$ are locally Noetherian formal schemes and the vertical arrows are representable by schemes, flat, and locally of finite presentation. We say that $f$ is \emph{flat} if it satisfies either of the following equivalent conditions.
    \begin{enumerate}[(1)]
     \item For any diagram as above such that in addition $U \ra X \times_{Y} V$ is flat, the morphism $h$ is flat.
    \item For some diagram as above with $a \colon U \ra X$ surjective, the morphism $h$ is flat.    
    \end{enumerate}
Next, consider a morphism $f \colon X \ra Y$ of locally Noetherian formal algebraic stacks. Consider diagrams as above, but assume instead that $U$ and $V$ are locally Noetherian formal algebraic spaces, and that the arrows $a$ and $b$ are representable by algebraic spaces, flat, and locally of finite presentation. We say that $f$ is \emph{flat} if either of the equivalent conditions (1) and (2) as above are satisfied. If the morphism $f$ is adic, then this agrees with the notion of flatness for adic morphisms as in \cite[{Definition 3.11}]{Emerton20}.
\label{footnote:flatness_formal_alg_stacks}
} 
part (``horizontal'') of $\breve{\mc{Z}}(T)$, i.e. the largest closed substack which is flat over $\Spf \mc{O}_{\breve{F}_p}$. We use similar notation $\mc{Z}(\underline{\mbf{x}}_p)_{\ms{H}}$ for the flat part of the local special cycle $\mc{Z}(\underline{\mbf{x}}_p)$.

Formation of ``flat part'' is flat local on the source. The quotient map $\Theta$ \eqref{equation:non-Arch_uniformization:framed_stack:forgetful} is representable by schemes and \'etale, hence flat. So the uniformization result (Proposition \ref{proposition:non-Arch_uniformization:proof}) implies that there is an induced uniformization morphism
    \begin{equation}\label{equation:non-Arch_uniformization:global_to_loca:horizontal_uniformization}
    \Theta \colon \coprod_{\substack{\underline{\mbf{x}} \in \mbf{W}^m \\ (\underline{\mbf{x}}, \underline{\mbf{x}}) = T}} \mc{Z}'(\underline{\mbf{x}}_p)_{\ms{H}} \times \mc{Z}'(\underline{\mbf{x}}^p)
    \ra
    \breve{\mc{Z}}(T)_{\ms{H}}
    \end{equation}
where $\mc{Z}'(\underline{\mbf{x}}_p)_{\ms{H}}$ is the flat part of $\mc{Z}'(\underline{\mbf{x}}_p)$. The action of $I'(\Q)$ must preserve the flat part, so generalities on stack quotients imply that $\Theta$ induces an isomorphism
    \begin{equation}
    \tilde{\Theta} \colon \Biggl [ I'(\Q) \backslash \Biggl ( \coprod_{\substack{\underline{\mbf{x}} \in \mbf{W}^m \\ (\underline{\mbf{x}}, \underline{\mbf{x}}) = T}} \mc{Z}'(\underline{\mbf{x}}_p)_{\ms{H}} \times \mc{Z}'(\underline{\mbf{x}}^p) \Biggr ) \Biggr ]
    \xra{\sim}
    \breve{\mc{Z}}(T)_{\ms{H}}
    \end{equation}
of formal algebraic stacks. For each $j \in J_p(T)$, the map $\Theta_j \colon \mc{Z}(\underline{\mbf{x}}_p) \ra \breve{\mc{Z}}(T)$ induces a map $\Theta_j \colon \mc{Z}(\underline{\mbf{x}}_p)_{\ms{H}} \ra \breve{\mc{Z}}(T)_{\ms{H}}$ (reusing the notation $\Theta_j$). Since $\Theta_j$ is flat and since formation of flat part is flat local on the source, the ``horizontal'' $\Theta_j$ arises from the original $\Theta_j$ by base-change along $\breve{\mc{Z}}(T)_{\ms{H}} \ra \breve{\mc{Z}}(T)$.

In the case $p \neq 2$ and for signature $(n - 1, 1)$ and $m^{\flat} = n - 1$, the next lemma is a consequence of \cite[{Theorem 4.2.1}]{LZ22unitary} (decomposition into quasi-canonical lifting cycles via Breuil modules) and is explained in \cite[{Lemma 2.49(1)}]{LL22II} (also via decomposition into quasi-canonical lifting cycles). In the case $p \neq 2$, signature $(n - 1, 1)$, and $m^{\flat} = n$, see again \cite[{Lemma 5.1.1}]{LZ22unitary} (inert, $p \neq 2$) and \cite[{Remark 2.26}]{LL22II} (ramified, exotic smooth).

\begin{lemma}\label{lemma:non-Arch_uniformization:global_to_local:horizontal_description}
Let $\underline{\mbf{x}}_p \in \mbf{W}_p^m$ be a tuple which spans a non-degenerate Hermitian $\mc{O}_{F_p}$-lattice, whose rank we denote $m^{\flat}$. Assume $m^{\flat} = n - 1$  if $F / \Q_p$ is split. Form the horizontal part $\mc{Z}(\underline{x})_{\ms{H}}$ of $\mc{Z}(\underline{x})$.
    \begin{enumerate}
        \item If $\mc{Z}(\underline{x})_{\ms{H}}$ is nonempty, then it is equidimensional of dimension $(n - r)r + 1 - m^{\flat} r$.
        \item If $m^{\flat} = n - 1$ and the signature is $(n - r, r) = (n - 1, 1)$, then the structure morphism $\mc{Z}(\underline{\mbf{x}}_p)_{\ms{H}} \ra \Spf \mc{O}_{\breve{F}_p}$ is a finite adic morphism of Noetherian formal schemes. The associated finite scheme over $\Spec \mc{O}_{\breve{F}_p}$ has reduced generic fiber.
        \item If $m^{\flat} = n$ and the signature is $(n - r, r) = (n - 1, 1)$, then $\mc{Z}(\underline{\mbf{x}}_p)_{\ms{H}} = \emptyset$. 
    \end{enumerate}
\end{lemma}
\begin{proof}
\hfill
\begin{enumerate}[(1)]
    \item By Corollary \ref{corollary:non-Arch_uniformization:vertical:approximate_cycle_RZ}, we can find $T \in \mrm{Herm}_m(\Q)$ (with $\rank(T) \geq n - 1$ if $p$ is split) and a morphism $\mc{Z}(\underline{\mbf{x}}_p) \ra \breve{\mc{Z}}(T)$ which is representable by schemes and \'etale. As formation of flat part is flat local on the source, we obtain a morphism $\mc{Z}(\underline{\mbf{x}})_{\ms{H}} \ra \breve{\mc{Z}}(T)_{\ms{H}}$ which is still representable by schemes and \'etale.
    The claim now follows from the corresponding global result for $\mc{Z}(T)_{\ms{H}}$ (Lemma \ref{lemma:special_cycles_generically_smooth}). Note that we may assume $K'_f$ is a small level (deepen the level away from $p$) to reduce to the case when $\breve{\mc{Z}}(T)_{\ms{H}}$ is a formal scheme.
    \item In this case, the map $\mc{Z}(\underline{\mbf{x}}_p)_{\ms{H}} \ra \breve{\mc{Z}}(T)_{\ms{H}}$ from part (1) is finite \'etale. We know that $\mc{\breve{Z}}(T)_{\ms{H}} \ra \Spf \mc{O}_{\breve{F}_p}$ (with $T$ as in the proof of loc. cit.) is proper and quasi-finite (Lemma \ref{lemma:horizontal_global_special_cycles_proper_quasi-finite}). 
    Since proper and quasi-finite implies finite (for morphisms of schemes) and since $\mc{\breve{Z}}(\underline{\mbf{x}}_p) \ra \Spf \mc{O}_{\breve{F}_p}$ is adic, (already proved in Lemma \ref{lemma:non-Arch_uniformization:vertical:local_cycle_quasi-compact}) i.e. representable by schemes, the claimed finiteness holds. The claim on reducedness in the generic fiber follows because $\mc{Z}(T)_{\ms{H}} \ra \Spec \mc{O}_F$ is \'etale in the generic fiber (Lemma \ref{lemma:special_cycles_generically_smooth}).
    We are passing from finite relative schemes over $\Spf \mc{O}_{\breve{F}_p}$ and $\Spec \mc{O}_{\breve{F}_p}$ as in \crefext{appendix:pDiv_prelim:Spec_v_Spf} (i.e. $\Spf R \mapsto \Spec R$).

    \item If $m^{\flat} = n$, then $\mc{\breve{Z}}(T)_{\ms{H}} = \emptyset$ (by Lemma \ref{lemma:special_cycles_generically_smooth}), so $\mc{Z}(\underline{\mbf{x}}_p)_{\ms{H}} = \emptyset$ by existence of the map $\mc{Z}(\underline{\mbf{x}}_p)_{\ms{H}} \ra \breve{\mc{Z}}(T)_{\ms{H}}$. \qedhere
\end{enumerate}
\end{proof}

In the case of $p \neq 2$, signature $(n - 1, 1)$, and $m^{\flat} = n - 1$, the following lemma is \cite[{\S 2.9}]{LZ22unitary} (there proved differently, using their quasi-compactness result via Bruhat--Tits stratification).

\begin{lemma}[Horizontal and vertical decomposition]\label{lemma:non-Arch_uniformization:global_to_local:scheme_union_decomp}
Let $\underline{\mbf{x}}_p \in \mbf{W}_p^m$ be a tuple which span a non-degenerate Hermitian $\mc{O}_{F_p}$-lattice of rank $m^{\flat}$. Assume $m^{\flat} = n - 1$ if $F / \Q_p$ is split. For $e \gg 0$, we have a scheme-theoretic union decomposition
    \begin{equation}
    \mc{Z}(\underline{\mbf{x}}_p) = \mc{Z}(\underline{\mbf{x}}_p)_{\ms{H}} \cup \mc{Z}(\underline{\mbf{x}}_p)_{\ms{V}}
    \end{equation}
where $\mc{Z}(\underline{\mbf{x}}_p)_{\ms{V}} \coloneqq \mc{Z}(\underline{\mbf{x}}_p)_{\Spf \mc{O}_{\breve{F}_p} / p^e \mc{O}_{\breve{F}_p}}$.
\end{lemma}
\begin{proof}
If $\mc{I}$ denotes the ideal sheaf of $\mc{Z}(\underline{\mbf{x}}_p)_{\ms{H}}$ as a closed subscheme of $\mc{Z}(\underline{\mbf{x}}_p)$, it is enough to show that $p^e$ annihilates $\mc{I}$ for $e \gg 0$. By Corollary \ref{corollary:non-Arch_uniformization:vertical:approximate_cycle_RZ}, we can find $T \in \mrm{Herm}_m(\Q)$ (with $\rank(T) \geq n - 1$ if $p$ is split) and a morphism $f \colon \mc{Z}(\underline{\mbf{x}}_p) \ra \breve{\mc{Z}}(T)$ which is representable by schemes and \'etale. We may assume $K'_f$ is small (deepen the level away from $p$) so that $\breve{\mc{Z}}(T)$ is a formal scheme.

If $\ms{I}$ denotes the ideal sheaf of the flat part $\breve{\mc{Z}}(T)_{\ms{H}} \subseteq \breve{\mc{Z}}(T)$, then $f^* \ms{I} \ra \mc{I}$ is surjective (by flatness of $f$, i.e. formation of flat part is flat local on the source). If $p^e$ annihilates $\ms{I}$, then $p^e$ also annihilates $\mc{I}$. We know that $\ms{I}$ consists (locally) of $p$-power torsion elements in the structure sheaf. Since $\breve{\mc{Z}}(T)$ is quasi-compact, we know that $\ms{I}$ is annihilated by $p^e$ for $e \gg 0$.
\end{proof}

For the rest of Section \ref{sec:non-Arch_uniformization}, we restrict to signature $(n - 1, 1)$ in all cases.

\begin{lemma}\label{lemma:non-Arch_uniformization:split_proper_quasi-finite_global}
\hfill
\begin{enumerate}[(1)]
    \item If $p$ is split, then $\breve{\mc{Z}}(T) \ra \Spf \mc{O}_{\breve{F}_p}$ is proper and quasi-finite and we have ${}^{\mbb{L}} \mc{Z}(T)_{\ms{V},p} = 0$. 
    \item Assume $n = 2$ and $\rank(T) \geq 1$. Then $\breve{\mc{Z}}(T) \ra \Spf \mc{O}_{\breve{F}_p}$ is proper and quasi-finite. If $\rank(T) = 1$, then we have ${}^{\mbb{L}} \mc{Z}(T)_{\ms{V},p} = 0$.
\end{enumerate}
\end{lemma}
\begin{proof}
(1) Recall our running assumption that $\rank(T) \geq n - 1$ if $p$ is split. 
Recall also $\breve{\mc{Z}}(T) \coloneqq \mc{Z}(T)_{\Spf \mc{O}_{\breve{F}_p}}$ in the split case, so the map $\breve{\mc{Z}}(T) \ra \Spf \mc{O}_{\breve{F}_p}$ is representable by algebraic stacks and locally of finite type.
This map is proper on reduced substacks by uniformization (Lemma \ref{lemma:non-Arch_uniformization:proof:proper} and Proposition \ref{proposition:non-Arch_uniformization:proof}), hence it is proper.

It remains to check that $\breve{\mc{Z}}(T) \ra \Spf \mc{O}_{\breve{F}_p}$ is quasi-finite in the sense of \cite[\href{https://stacks.math.columbia.edu/tag/0G2M}{Definition 0G2M}]{stacks-project}. It is enough to check that $\breve{\mc{Z}}(T)_{\mrm{red}} \ra \Spec \overline{k}$ is quasi-finite. This follows from the uniformization isomorphism, since $\breve{\mc{Z}}(T)_{\mrm{red}}$ may be covered by finitely many copies of $\Spec \overline{k}$ (combine uniformization with the analogous result for local special cycles, which is \crefext{I:lemma:moduli_pDiv:discrete_reduced:local_special_cycle_red_finite}; since we assume $\rank(T) \geq n - 1$ when $p$ is split, the groupoid $[I'(\Q) \backslash J_p(T)]$ has finitely many isomorphism classes, as discussed in Section \ref{ssec:non-Arch_uniformization:quotient}).

The derived vertical special cycle class ${}^{\mbb{L}} \mc{Z}(T)_{\ms{V},p} \in \mrm{gr}^m_{\mc{M}} K'_0(\mc{Z}(T)_{\F_p})_{\Q}$ was defined in \cref{ssec:arith_cycle_classes:vertical}. If $m \geq n$ then $\mc{Z}(T)$ is empty. If $m = n - 1$, then $\mrm{gr}^m_{\mc{M}} K'_0(\mc{Z}(T)_{\F_p})_{\Q} = 0$ because $\mc{Z}(T)_{\F_p}$ has dimension $0$ (and $\mc{M}$ has dimension $n$).

(2) This may be proved as in part (1). We may assume $p$ is nonsplit. We have $\breve{\mc{Z}}(T) = \mc{Z}(T)_{\Spf \mc{O}_{\breve{F}_p}}$ (Lemma \ref{corollary:special_cycle_supersingular_support}). Then use quasi-compactness of local special cycles (Lemma \ref{lemma:non-Arch_uniformization:vertical:local_cycle_quasi-compact}), uniformization, and discreteness of $\mc{N}_{\red}$ \crefext{I:ssec:moduli_pDiv:discrete_reduced}. Suppose $\rank(T) = 1$. First consider the case $m = 1$. Then ${}^{\mbb{L}} \mc{Z}(T)_{\ms{V},p} \in \mrm{gr}^1_{\mc{M}} K'_0(\mc{Z}(T)_{\F_p})_{\Q} = 0$ because $\mc{Z}(T)_{\F_p}$ has dimension $0$ (and $\mc{M}$ has dimension $2$). If $m = 2$, then ${}^{\mbb{L}} \mc{Z}(t_i)_{\ms{V},p} = 0$ for any nonzero diagonal entry $t_i$ of $T$ by the preceding argument, so ${}^{\mbb{L}} \mc{Z}(T)_{\ms{V},p} = 0$ by construction (defined in \cref{ssec:arith_cycle_classes:vertical} as the projection of a product against ${}^{\mbb{L}} \mc{Z}(t_i)_{\ms{V},p} = 0$ for some $i$).
\end{proof}

\begin{lemma}\label{lemma:non-Arch_uniformization:global_local_derived_cycles}
Assume $p$ is nonsplit. Assume that $K'_f$ is a small level, so that $\mc{M}$ is a scheme. Fix any $j \in J_p(T)$ and consider the map
    \begin{equation}
    \mc{Z}(\underline{\mbf{x}}_p) \xra{\Theta_j} \breve{\mc{Z}}(T) \ra \mc{Z}(T).
    \end{equation}
The class ${}^{\mbb{L}} \mc{Z}(T) \in K'_0(\mc{Z}(T))_{\Q}$ pulls back to the class ${}^{\mbb{L}} \mc{Z}(\underline{\mbf{x}}_p) \in K'_0(\mc{Z}(\underline{\mbf{x}}_p))_{\Q}$.
\end{lemma}
\begin{proof}
The maps $\Theta_j \colon \mc{Z}(\underline{\mbf{x}}_p) \ra \breve{\mc{Z}}(T)$ and $\breve{\mc{Z}}(T) \ra \mc{Z}(T)$ are flat maps of locally Noetherian formal schemes, so we may take the non-derived pullback. The lemma may be proved using the fact that the commutative diagrams
    \begin{equation}\label{lemma:non-Arch_uniformization:global_local_derived_cycles:Cartesian}
    \begin{tikzcd}
    \breve{\mc{Z}}(T)_{\mrm{framed}} \arrow{r} \arrow{d} \pullbackcorner & \breve{\mc{M}}_{\mrm{framed}} \arrow{d} \\
    \mc{Z}(T) \arrow{r} & \mc{M}
    \end{tikzcd}
    \quad \quad
    \begin{tikzcd}
    \breve{\mc{Z}}(t_i)_{\mrm{framed}} \arrow{r} \arrow{d} \pullbackcorner & \breve{\mc{M}}_{\mrm{framed}} \arrow{d} \\
    \mc{Z}(t_i) \arrow{r} & \mc{M}
    \end{tikzcd}
    \end{equation}
are $2$-Cartesian (where the $t_i$ are the diagonal entries of $T$), and the fact that the tautological bundle $\mc{E}$ on $\mc{M}$ pulls back to the tautological bundle $\mc{E}$ on $\mc{N}$.
\end{proof}

\begin{corollary}\label{corollary:non-Arch_uniformization:vertical:local_cycle_filtration}
Assume $p$ is nonsplit. For any $\underline{\mbf{x}}_p \in \mbf{W}_p^m$, we have ${}^{\mbb{L}} \mc{Z}(\underline{\mbf{x}}_p) \in F^m_{\mc{N}} K'_0(\mc{Z}(\underline{\mbf{x}}_p))_{\Q}$.
\end{corollary}
\begin{proof}
By Corollary \ref{corollary:non-Arch_uniformization:vertical:approximate_cycle_RZ}, we can find $T \in \mrm{Herm}_m(\Q)$ (with $\rank(T) \geq n - 1$ if $p$ is split) and a morphism $\mc{Z}(\underline{\mbf{x}}_p) \ra \breve{\mc{Z}}(T)$ which is representable by schemes and \'etale. We can deepen the level $K'_f$ away from $p$ to assume $\breve{\mc{Z}}(T)$ is a formal scheme. Since ${}^{\mbb{L}} \mc{Z}(T) \in F^m_{\mc{N}}K'_0(\mc{Z}(T))_{\Q}$, the corollary follows from Lemma \ref{lemma:non-Arch_uniformization:global_local_derived_cycles}.
\end{proof}

We previously defined derived vertical (global) special cycles ${}^{\mbb{L}} \mc{Z}(T)_{\ms{V},p} \in \mrm{gr}^m_{\mc{M}} K'_0(\mc{Z}(T)_{\F_p})_{\Q}$ (\cref{ssec:arith_cycle_classes:vertical}). In the next lemma, we write $\mc{Z}(T)_{(p)} \coloneqq \mc{Z}(T) \times_{\Spec \Z} \Spec \Z_{(p)}$. We also write $\mc{Z}(T)_{\ms{V},p} \coloneqq \mc{Z}(T) \times_{\Spec \Z} \Spec \Z / p^e \Z$ and $\mc{Z}(\underline{\mbf{x}}_p)_{\ms{V},p} \coloneqq \mc{Z}(\underline{\mbf{x}}_p) \times_{\Spf \mc{O}_{\breve{F}_p}} \Spec \mc{O}_{\breve{F}_p} / p^e \mc{O}_{\breve{F}_p}$ for an understood integer $e \gg 0$. We also set $\breve{\mc{Z}}(T)_{\ms{V}} \coloneqq \breve{\mc{Z}}(T) \times_{\Spf \mc{O}_{\breve{F}_p}} \Spec \mc{O}_{\breve{F}_p} / p^e \mc{O}_{\breve{F}_p}$.

\begin{lemma}\label{lemma:non-Arch_uniformization:vertical:pullback_K_0_vertical}
Fix any $j \in J_p(T)$. Write $\underline{\mbf{x}} \in \mbf{W}^m$ for the associated $m$-tuple. Fix any $e \gg 0$ such that there are scheme-theoretic union decompositions
    \begin{equation}
    \mc{Z}(T)_{(p)} = \mc{Z}(T)_{(p), \ms{H}} \cup \mc{Z}(T)_{\ms{V},p} \quad \mc{Z}(\underline{\mbf{x}}_p) = \mc{Z}(\underline{\mbf{x}}_p)_{\ms{H}} \cup \mc{Z}(\underline{\mbf{x}}_p)_{\ms{V},p}
    \end{equation}
(``horizontal and vertical'').
Pullback along the map
    \begin{equation}
    \mc{Z}(\underline{\mbf{x}}_p)_{\ms{V}} \xra{\Theta_j} \breve{\mc{Z}}(T)_{\ms{V}} \ra \mc{Z}(T)_{\ms{V},p}
    \end{equation}
sends ${}^{\mbb{L}} \mc{Z}(T)_{\ms{V},p} \in \mrm{gr}^{m}_{\mc{M}} K'_0(\mc{Z}(T)_{\ms{V},p})_{\Q}$ to ${}^{\mbb{L}} \mc{Z}(\underline{\mbf{x}}_p)_{\ms{V}} \in \mrm{gr}^m_{\mc{N}} K'_0(\mc{Z}(\underline{\mbf{x}}_p)_{\ms{V}})_{\Q}$.
\end{lemma}
\begin{proof}
If $p$ is split, the derived vertical special cycles (global and local) are zero (Lemma \ref{lemma:non-Arch_uniformization:split_proper_quasi-finite_global} (global) and \crefext{I:ssec:moduli_pDiv:horizontal_vertical} (local)) and the lemma is trivial. We remind the reader of our running assumption that $\rank(T) \geq n - 1$ if $p$ is split. 

We thus assume that $p$ is nonsplit. By the local and global linear invariance results \crefext{I:ssec:moduli_pDiv:horizontal_vertical} and \cref{equation:arith_cycle_classes:vertical:linear_invariance}, it is enough to check the case where $T = \mrm{diag}(0,T^{\flat})$ where $\det T^{\flat} \neq 0$.

First consider the case where $T$ is nonsingular, i.e. $T = T^{\flat}$. If $K'_f$ is a small level, the lemma follows from Lemma \ref{lemma:non-Arch_uniformization:global_local_derived_cycles}, since the projections $\mrm{gr}^{m^{\flat}}_{\mc{M}} K'_0(\mc{Z}(T^{\flat}))_{\Q} \ra \mrm{gr}^{m^{\flat}}_{\mc{M}} K'_0(\mc{Z}(T^{\flat})_{\ms{V},p})_{\Q}$ and $\mrm{gr}^{m^{\flat}}_{\mc{N}} K'_0(\mc{Z}(\underline{\mbf{x}}_p))_{\Q} \ra \mrm{gr}^{m^{\flat}}_{\mc{N}} K'_0(\mc{Z}(\underline{\mbf{x}}_p)_{\ms{V}})_{\Q}$
are given by (non-derived) pullbacks of coherent sheaves, see \cref{lemma:K'0_component_decomp} (Deligne--Mumford stacks) and \cite[{Lemma B.1}]{Zhang21} (locally Noetherian formal schemes). Note that the codimension graded pieces $\mrm{gr}^m$ are preserved, by \'etale-ness of $\Theta_j$. In general, we may reduce to the case where $K'_f$ is a small level by compatibility of ${}^{\mbb{L}} \mc{Z}(T^{\flat})_{\ms{V},p}$ with (finite \'etale) pullback for varying levels.

Next, consider the case where $T$ is possibly singular with $T = \mrm{diag}(0,T^{\flat})$. If $K'_f$ is a small level, this follows as in the proof of Lemma \ref{lemma:non-Arch_uniformization:global_local_derived_cycles}. 
That is, the class ${}^{\mbb{L}} \mc{Z}(T)_{\ms{V},p} \coloneqq (\mc{E}^{\vee})^{m - \rank(T)} \cdot \mc{Z}(T^{\flat})_{\ms{V},p}$ pulls back to ${}^{\mbb{L}} \mc{Z}(\underline{\mbf{x}}_p)_{\ms{V}} \coloneqq (\mc{E}^{\vee})^{m - \rank(T)} \cdot \mc{Z}(\underline{\mbf{x}}_p^{\flat})_{\ms{V}}$ (use the result for $T^{\flat}$ just proved). 
In general, we may reduce to the case where $K'_f$ is a small level (deepen level away from $p$) by compatibility of ${}^{\mbb{L}} \mc{Z}(T^{\flat})_{\ms{V},p}$ with (finite \'etale) pullback for varying levels.
\end{proof}
                
            \subsection{Local intersection numbers: vertical}
            \label{ssec:non-Arch_uniformization:vertical}
                The main purpose of this section is to reduce ``global vertical intersection numbers'' to ``local vertical intersection numbers'' (see end of this section). We continue to assume signature $(n - 1, 1)$.

Consider $T' \in \mrm{Herm}_m(\Q_p)$ (with $F_p$-coefficients) with $\rank(T') = n - 1$, and either $m = n - 1$ or $m = n$. 
For any tuple $\underline{\mbf{x}}_p \in \mbf{W}_p^m$ with Gram matrix $T'$, we define the \emph{local vertical intersection number}
    \begin{equation}\label{equation:non-Arch_uniformization:horizontal:Int_vertical_local}
    \mrm{Int}_{\ms{V},p}(T') \coloneqq
    \begin{cases}
    2 [\breve{F}_p : \breve{\Q}_p]^{-1} \deg_{\overline{k}}({}^{\mbb{L}} \mc{Z}(\underline{\mbf{x}}_p)_{\ms{V}} \cdot \mc{E}^{\vee}) \log p & \text{if $m = n - 1$}
    \\
    2 [\breve{F}_p : \breve{\Q}_p]^{-1} \deg_{\overline{k}}({}^{\mbb{L}} \mc{Z}(\underline{\mbf{x}}_p)_{\ms{V}}) \log p& \text{if $m = n$}
    \end{cases}
    \end{equation}
Here, $\mc{E}^{\vee}$ stands for the class $[\mc{O}_{\mc{N}}] - [\mc{E}] \in K'_0(\mc{N})$. If no such $\underline{\mbf{x}}_p$ exists, we set $\mrm{Int}_{\ms{V},p}(T') \coloneqq 0$. The definition of $\mrm{Int}_{\ms{V},p}(T')$ does not depend on the choice of $\underline{\mbf{x}}_p$ (by the action of $U(\mbf{W}_p)$ on $\mc{N}(n - 1, 1)$ as in \crefext{I:ssec:moduli_pDiv:action_RZ}). The factor $2 [\breve{F}_p : \breve{\Q}_p]^{-1}$ will account for total degree of $\Spec \mc{O}_{F_p} \ra \Spec \Z_p$ on residue fields (e.g. we need to account for both primes in $\mc{O}_F$ over $p$ in the split case).
By local linear invariance \crefext{I:ssec:moduli_pDiv:horizontal_vertical}, we have
    \begin{equation}
    \mrm{Int}_{\ms{V},p}(T') = \mrm{Int}_{\ms{V},p}({}^t \overline{\gamma} T' \gamma)
    \end{equation}
for any $\gamma \in \GL_m(\mc{O}_{F_p})$.

Consider any $T \in \mrm{Herm}_m(\Q)$ (with $F$-coefficients) with $\rank(T) = n - 1$, and either $m = n - 1$ or $m = n$. Pick any set of representatives $J \subseteq J_p(T)$ for the isomorphism classes of the groupoid $[I'(\Q) \backslash J_p(T)]$. By Lemma \ref{lemma:non-Arch_uniformization:vertical:pullback_K_0_vertical}, we have
    \begin{align}\label{equation:non-Arch_uniformization:vertical:local_global_degrees}
    & \mrm{Int}_{\ms{V},p,\mrm{global}}(T) \coloneqq \deg_{\F_p}({}^{\mbb{L}} \mc{Z}(T)_{\ms{V},p} \cdot (\mc{E}^{\vee})^{n - m}) \log p
    \\
    & = \mrm{Int}_{\ms{V},p}(T) \sum_{j \in J} \frac{1}{|\Aut(j)|}
    \\
    & = \mrm{Int}_{\ms{V},p}(T) \frac{[K_{L_0, f} : K_{0,f}]}{|\mc{O}_F^{\times}| / h_F} \cdot \deg \Biggl [ I_1(\Q) \backslash \Biggl ( \coprod_{\substack{\underline{\mbf{x}} \in \mbf{W}^m \\ (\underline{\mbf{x}}, \underline{\mbf{x}}) = T}} U(\mbf{W}^{\perp}_p) / K_{1, \mbf{L}^{\perp}_p} \times \mc{Z}(\underline{\mbf{x}}^p) \Biggr ) \Biggr ]. \notag
    \end{align}

For later use in \crefext{IV:remark:arithmetic_Siegel-Weil:main_results:nonsingular}, consider $T \in \mrm{Herm}_n(\Q)$ (with $F$-coefficients) with $\det T \neq 0$. We consider the \emph{local intersection number}
    \begin{equation}
    \mrm{Int}_p(T) \coloneqq 2 [\breve{F}_p : \breve{\Q}_p]^{-1} \deg_{\overline{k}} ( {}^{\mbb{L}} \mc{Z}(\underline{\mbf{x}}_p)) \log p
    \end{equation}
where $\underline{\mbf{x}}_p \in \mbf{W}_p^n$ is any $n$-tuple with Gram matrix $T$ (since $\rank \mbf{W}_p = n - 1$ when $p$ is split, set $\mrm{Int}_p(T) \coloneqq 0$ in this case). 
Note ${}^{\mbb{L}} \mc{Z}(\underline{\mbf{x}}_p)_{\ms{V}} = {}^{\mbb{L}} \mc{Z}(\underline{\mbf{x}}_p)$ 
by Lemmas \ref{lemma:non-Arch_uniformization:global_to_local:horizontal_description} and \ref{lemma:non-Arch_uniformization:vertical:local_cycle_quasi-compact} (under the d\'evissage pushforward identification $K'_0(\mc{Z}(\underline{\mbf{x}}_p)_{\overline{k}}) \xra{\sim} K'_0(\mc{Z}(\underline{\mbf{x}}_p))$).
By Lemma \ref{lemma:non-Arch_uniformization:vertical:pullback_K_0_vertical}, we have
    \begin{align}
    \mrm{Int}_{p, \mrm{global}}(T) & \coloneqq \deg_{\F_p} ( {}^{\mbb{L}}\mc{Z}(T)_{\ms{V},p}) \log p
    \\
    & = \mrm{Int}_{p}(T)  \frac{[K_{L_0, f} : K_{0,f}]}{|\mc{O}_F^{\times}| / h_F} \cdot \deg \Biggl [ I_1(\Q) \backslash \Biggl ( \coprod_{\substack{\underline{\mbf{x}} \in \mbf{W}^n \\ (\underline{\mbf{x}}, \underline{\mbf{x}}) = T}} \mc{Z}(\underline{\mbf{x}}^p) \Biggr ) \Biggr ].
    \end{align}
    
            \subsection{Local intersection numbers: horizontal}
            \label{ssec:non-Arch_uniformization:horizontal}
                The main purpose of this section is to reduce ``global horizontal intersection numbers'' to ``local horizontal intersection numbers'' (see end of this section).
We continue to assume signature $(n - 1, 1)$.
In Section \ref{ssec:non-Arch_uniformization:horizontal}, we require $p \neq 2$ if $p$ is inert (because we required this for our discussion of quasi-canonical lifting cycles, \crefext{I:ssec:can_and_qcan:qcan_cycles}).

Consider $T' \in \mrm{Herm}_m(\Q_p)$ (with $F_p$-coefficients) with $\rank(T') = n - 1$, and either $m = n - 1$ or $m = n$. Select any $\underline{\mbf{x}}_p \in \mbf{W}_p^m$ with Gram matrix $T'$, and set $L^{\flat}_p \coloneqq \mrm{span}_{\mc{O}_{F_p}}(\underline{\mbf{x}}_p)$. We define the \emph{local horizontal intersection number}
    \begin{equation}\label{equation:non-Arch_uniformization:horizontal:Int_horizontal_local}
    \mrm{Int}_{\ms{H},p}(T') \coloneqq \sum_{\substack{L_p^{\flat} \subseteq M_p^{\flat} \subseteq M_p^{\flat *} \\ t(M_p^{\flat}) \leq 1}} \mrm{Int}_{\ms{H},p}(M^{\flat}_p)^{\circ}
    \end{equation}
where the sum runs over lattices $M_p^{\flat} \subseteq L^{\flat}_p \otimes_{\mc{O}_{F_p}} F_p$, where
    \begin{equation}\label{equation:non-Arch_uniformization:horizontal:local_primitive}
    \mrm{Int}_{\ms{H},p}(M^{\flat}_p)^{\circ} \coloneqq 2 \cdot \deg \mc{Z}(M^{\flat}_p)^{\circ} \cdot \delta_{\mrm{tau}}(\mrm{val}'(M^{\flat}_p))
    \end{equation}
with $\mrm{val}'(M^{\flat}_p) \coloneqq \lfloor \mrm{val}(M^{\flat}_p) \rfloor$ and with $\delta_{\mrm{tau}}(-)$ the ``local change of tautological height'' as defined in \crefext{I:equation:can_and_qcan:qcan:local_change_taut}.
Here $\mc{Z}(M^{\flat}_p)^{\circ} \subseteq \mc{N}(n - 1, 1)$ is the quasi-canonical lifting cycle associated with $M^{\flat}_p$ 
\crefext{I:ssec:can_and_qcan:qcan_cycles}. The local horizontal intersection number should be compared with the decomposition of horizontal local special cycles into quasi-canonical lifting cycles \crefext{I:ssec:can_and_qcan:qcan_cycles}. The notation $\deg \mc{Z}(M^{\flat}_p)^{\circ}$ 
means the degree of the finite flat adic morphism $\mc{Z}(M^{\flat}_p)^{\circ} \ra \Spf \mc{O}_{\breve{F}_p}$. If no such $\underline{\mbf{x}}_p$ exists, we set $\mrm{Int}_{\ms{H},p}(T') \coloneqq 0$.

This definition of $\mrm{Int}_{\ms{H},p}(T')$ does not depend on the choice of $\underline{\mbf{x}}_p$ (again by the action of $U(\mbf{W}_p)$ on $\mc{N}$ 
\crefext{I:ssec:moduli_pDiv:action_RZ} and Witt's theorem). 
The formula for $\mrm{deg}\mc{Z}(M^{\flat}_p)^{\circ}$ (combine \crefext{I:equation:can_and_qcan:qcan:degree} and \crefext{I:equation:can_and_qcan:qcan_cycles:underlying_formal_scheme}) shows $\mrm{Int}_{\ms{H},p}(M^{\flat}_p)^{\circ}\in \Z$. 
The extra factor of $2$ in \eqref{equation:non-Arch_uniformization:horizontal:local_primitive} will account for the fact that $\Spf (\mc{O}_{F} \otimes_{\Z} \breve{\Z}_p) \ra \Spf \breve{\Z}_p$ has degree $2$.

In the above situation, we also set
    \begin{equation}
    \deg_{\ms{H},p}(T') \coloneqq \deg \mc{Z}(\underline{\mbf{x}}_p)_{\ms{H}}
    \end{equation}
where the right-hand side means the degree of the finite flat adic morphism $\mc{Z}(\underline{\mbf{x}}_p)_{\ms{H}} \ra \Spf \mc{O}_{\breve{F}_p}$. If no such $\underline{\mbf{x}}_p$ exists, we set $\deg_{\ms{H},p}(T') = 0$. Again, the definition of $\deg_{\ms{H},p}(T')$ does not depend on the choice of $\underline{\mbf{x}}_p$.

Suppose $T \in \mrm{Herm}_m(\Q)$ (with $F$-coefficients) with $\rank(T) = n - 1$, and either $m = n - 1$ or $m = n$. Then (in the notation of \cref{ssec:arith_cycle_classes:Hodge_bundles,ssec:arith_cycle_classes:Hermitian_bundles}) we have
    \begin{equation}
    \widehat{\deg}(\widehat{\mc{E}}^{\vee}|_{\mc{Z}(T)_{\ms{H}}}) = \widehat{\deg}(\widehat{\Omega}_0^{\vee}|_{\mc{Z}(T)_{\ms{H}}}) + \widehat{\deg}(\widehat{\ms{E}}^{\vee}|_{\mc{Z}(T)_{\ms{H}}}).
    \end{equation}
We have
    \begin{equation}
    \widehat{\deg}(\widehat{\Omega}_0^{\vee}|_{\mc{Z}(T)_{\ms{H}}}) = \deg_{\Z} \mc{Z}(T)_{\ms{H}} \cdot (- h_{\mrm{Fal}}^{\mrm{CM}} - \frac{1}{4} \log |\Delta|)
    \end{equation}
where $\deg_{\Z} \mc{Z}(T)_{\ms{H}}$ means the degree of $\mc{Z}(T)_{\ms{H}} \times_{\Spec \Z} \Spec \Q \ra \Spec \Q$ (stacky degrees as in \cref{equation:zero_cycle_degree_over_field} and surrounding discussion). 

Pick any set of representatives $J \subseteq J_p(T)$ for the isomorphism classes of the groupoid $[I'(\Q) \backslash J_p(T)]$. Using the finite \'etale maps $\Theta_j \colon \mc{Z}(\underline{\mbf{x}}_p)_{\ms{H}} \ra \breve{\mc{Z}}(T)_{\ms{H}}$ for $j \in J_p(T)$ (Section \ref{ssec:non-Arch_uniformization:quotient} and \eqref{equation:non-Arch_uniformization:global_to_loca:horizontal_uniformization}) which cover $\breve{\mc{Z}}(T)_{\ms{H}}$ as $j$ ranges over $J$, we find
    \begin{align}\label{equation:non-Arch_uniformization:horizontal:local_global_degrees:non-arithmetic}
    & \deg_{\Z} \mc{Z}(T)_{\ms{H}} = \deg_{\ms{H},p}(T) \sum_{j \in J_p(T)} \frac{1}{|\Aut(j)|}
    \\
    & = \deg_{\ms{H},p}(T) \frac{[K_{L_0, f} : K_{0,f}]}{|\mc{O}_F^{\times}| / h_F} \cdot \deg \Biggl [ I_1(\Q) \backslash \Biggl ( \coprod_{\substack{\underline{\mbf{x}} \in \mbf{W}^m \\ (\underline{\mbf{x}}, \underline{\mbf{x}}) = T}} U(\mbf{W}^{\perp}_p) / K_{1, \mbf{L}^{\perp}_p} \times \mc{Z}(\underline{\mbf{x}}^p) \Biggr ) \Biggr ]. \notag
    \end{align}

Combining the following: (1) the finite \'etale maps $\Theta_j \colon \mc{Z}(\underline{\mbf{x}}_p)_{\ms{H}} \ra \breve{\mc{Z}}(T)_{\ms{H}}$ for $j \in J_p(T)$ (Section \ref{ssec:non-Arch_uniformization:quotient} and \eqref{equation:non-Arch_uniformization:global_to_loca:horizontal_uniformization}) which cover $\breve{\mc{Z}}(T)_{\ms{H}}$ as $j$ ranges over $J$
(2) \crefext{I:proposition:can_and_qcan:qcan_cycles} (decomposition of horizontal local special cycles into quasi-canonical liftings) and discussion surrounding \crefext{I:equation:qcan_cycles:pDiv_Serre_tensor}, and (3) \cref{corollary:qcan_heights:minimal_isogenies:local_decomp} (decomposition of global height into local ``change of heights'' for $p$-divisible groups), we find
    \begin{align}\label{equation:non-Arch_uniformization:horizontal:local_global_degrees}
    & \mrm{Int}_{\ms{H},p,\mrm{global}}(T) \coloneqq \widehat{\deg}(\widehat{\ms{E}}^{\vee}|_{\mc{Z}(T)_{\ms{H}}}) - (\deg_{\Z} \mc{Z}(T)_{\ms{H}}) \cdot h_{\mrm{tau}}^{\mrm{CM}} \mod ~ \sum_{\ell \neq p} \Q \cdot \log \ell
    \\
    & = \mrm{Int}_{\ms{H},p}(T) \sum_{j \in J_p(T)} \frac{1}{|\Aut(j)|} \notag
    \\
    & = \mrm{Int}_{\ms{H},p}(T) \frac{[K_{L_0, f} : K_{0,f}]}{|\mc{O}_F^{\times}| / h_F} \cdot \deg \Biggl [ I_1(\Q) \backslash \Biggl ( \coprod_{\substack{\underline{\mbf{x}} \in \mbf{W}^m \\ (\underline{\mbf{x}}, \underline{\mbf{x}}) = T}} U(\mbf{W}^{\perp}_p) / K_{1, \mbf{L}^{\perp}_p} \times \mc{Z}(\underline{\mbf{x}}^p) \Biggr ) \Biggr ] \notag
    \end{align}
with $h_{\mrm{tau}}^{\mrm{CM}}$ and $h_{\widehat{\mc{E}}^{\vee}}^{\mrm{CM}}$ as in \eqref{equation:arith_cycle_classes:Hodge_bundles:taut_height_constants}.
The notation ``$\mrm{mod} ~ \sum_{\ell \neq p} \Q \cdot \log \ell$'' means that equality holds as elements of the (additive) quotient $\R / ( \sum_{\ell \neq p} \Q \cdot \log \ell)$. Note $\mrm{Int}_{\ms{H},p,\mrm{global}}(T) \in \Q \cdot \log p$.
To apply \cref{corollary:qcan_heights:minimal_isogenies:local_decomp}, we first consider the case of small level $K'_f$ so that $\mc{Z}(T)_{\ms{H}}$ is a scheme. This immediately implies the case of general (stacky) level, by compatibility of arithmetic degree with finite \'etale covers (see \cref{ssec:arith_cycle_classes:Hermitian_bundles}). We have
    \begin{equation}\label{equation:non-Arch_uniformization:horizontal:local_global_degrees:decomp}
    \widehat{\deg}(\widehat{\mc{E}}^{\vee}|_{\mc{Z}(T)_{\ms{H}}}) - (\deg_{\Z} \mc{Z}(T)_{\ms{H}}) \cdot h_{\widehat{\mc{E}}^{\vee}}^{\mrm{CM}} = \sum_{\ell} \mrm{Int}_{\ms{H},\ell,\mrm{global}}(T)
    \end{equation}
where the sum ranges over all primes $\ell$, with all but finitely many terms equal to $0$. The preceding expression should be understood modulo $\Q \cdot \log \ell$ for those primes $\ell$ for which $L \otimes_{\Z} \Z_{\ell}$ is not self-dual. If $L$ is not self-dual, we also quotient by $\Q \cdot \log \ell$ for primes $\ell \mid \Delta$. We also quotient by $\Q \cdot \log 2$ unless $2$ is split in $\mc{O}_F$.

We define \emph{total ``intersection numbers''}
    \begin{equation}\label{equation:non-Arch_uniformization:horizontal:total_intersection}
    \mrm{Int}_{p}(T) \coloneqq \mrm{Int}_{\ms{H},p}(T) + \mrm{Int}_{\ms{V},p}(T) \quad \quad \mrm{Int}_{p, \mrm{global}}(T) \coloneqq \mrm{Int}_{\ms{H}, p, \mrm{global}}(T) + \mrm{Int}_{\ms{V}, p, \mrm{global}}(T)
    \end{equation}
 (local and global) at $p$. The local version featured in our main non-Archimedean local theorems \crefext{I:sec:non-Arch_identity}, and the global version will appear in the proof of our main global theorem \crefext{IV:theorem:arithmetic_Siegel-Weil:main_results:main}.

For readers interested in Faltings heights, we record the relation
    \begin{align}\label{equation:non-Arch_uniformization:horizontal:local_global_degrees:Faltings}
    & \widehat{\deg}(\widehat{\omega}|_{\mc{Z}(T)_{\ms{H}}}) - (\deg_{\Z} \mc{Z}(T)_{\ms{H}}) \cdot n \cdot h_{\mrm{Fal}}^{\mrm{CM}} = - 2 \sum_{\ell} \mrm{Int}_{\ms{H},\ell, \mrm{global}}(T)
    \end{align}
where $\widehat{\omega}$ is the metrized Hodge determinant bundle (see discussion in \cref{ssec:ab_var:integral_models}), the sum again runs over all primes $\ell$, and where $\mrm{Int}_{\ms{H},\ell, \mrm{global}}(T)$ is the same quantity defined above. This follows by the same argument as above, using \cref{corollary:qcan_heights:minimal_isogenies:local_decomp}. The remarks following \eqref{equation:non-Arch_uniformization:horizontal:local_global_degrees:decomp} (about quotienting by $\Q \cdot \log \ell$ for some primes $\ell$) apply here verbatim.
                
        \section{Archimedean}
        \label{sec:Arch_uniformization}
            We explain complex uniformization for special cycles and Green currents on $\mc{M}$.
Aside from our discussion on Green currents for singular $T$ in Section \ref{ssec:Arch_uniformization:Archimedean} (when $\rank(T) = n - 1$), most of the material in \cref{sec:Arch_uniformization} should be fairly standard, e.g. \cite[{\S 3}]{KR14} (uniformization of special cycles), \cite[{\S 2}]{BHKRY20} (including discussion of metrized tautological bundle), \cite[{\S 4B}]{Liu11} (Green currents via uniformization), etc.. Strictly speaking, however, the references \cite{KR14,BHKRY20} restrict to principal polarizations. We will need non-principal polarizations (this slightly affects how we normalize the metric on the tautological bundle), so we explain the setup.

With notation as explained at the start of \cref{part:part_III:uniformization}, we also assume $L$ has signature $(n - 1, 1)$. For technical convenience, we assume the implicit level $K'_f$ is small so that $\mc{M}$ is a scheme (except at the very end of Section \ref{ssec:Arch_uniformization:Archimedean}). Fix one of the two embeddings $F \ra \C$, write $\mc{M}_{\C} \coloneqq \mc{M} \times_{\Spec \mc{O}_F} \Spec \C$, and let $\mc{M}_{\C}^{\mrm{an}}$ be the analytification (outside of Section \ref{sec:Arch_uniformization}, we often abuse notation and drop the superscript $\mrm{an}$). This is a complex manifold of dimension $n - 1$. Given any Hermitian matrix $T \in \mrm{Herm}_m(\Q)$ (with $F$-coefficients) with associated special cycle $\mc{Z}(T) \ra \mc{M}$, we use similar notation $\mc{Z}(T)_{\C}$ and $\mc{Z}(T)^{\mrm{an}}_{\C}$. Since $\mc{Z}(T)_{\C} \ra \Spec \C$ is smooth (Lemma \ref{lemma:special_cycles_generically_smooth}), we know that $\mc{Z}(T)_{\C}^{\mrm{an}}$ is also a complex manifold.

We view $V_{\R}$ as a complex vector space via the identification $F \otimes_{\Q} \R \cong \C$ (induced by the choice of $F \ra \C$). We use notation from \crefext{II:sec:Hermitian_domain} on the Hermitian symmetric space $\mc{D}$ and its local special cycles $\mc{D}(\underline{x})$ for tuples $\underline{x} \in V^m_{\R}$, etc..

We set $V_{\C} \coloneqq V \otimes_{\Q} \C$ and write $V_{\C} = V_{\C}^+ \oplus V_{\C}^-$ where the $F$-action on $V^+_{\C}$ (resp. $V^-_{\C}$) is $F$-linear (resp. $\s$-linear) with respect to the chosen map $F \ra \C$. We use similar notation for other $F \otimes_{\Q} \C$-modules.
    
            \subsection{Local special cycles away from \texorpdfstring{$\infty$}{infinity}}
            \label{ssec:Arch_uniformization:away_infty_special_cycle}
                Given an $m$-tuple $\underline{x}_f = [x_1, \ldots, x_m] \in (V \otimes_{\Q} \A_f)^m$, we consider an ``\emph{away-from-$\infty$}'' local special cycle
    \begin{equation}
    \mc{D}'(\underline{x}_f) \coloneqq \{ (g_0, g) \in G'(\A_f) / K'_f : g^{-1} g_0 x_i \in L \otimes_{\Z} \hat{\Z}^p \text{ for all $x_i \in \underline{x}_f$} \}.
    \end{equation}
We view $\mc{D}'(\underline{x}_f)$ as a discrete set. We also define the ``\emph{away-from-$\infty$}'' local special cycle
    \begin{equation}
    \mc{D}(\underline{x}_f) \coloneqq \{ g \in U(V)(\A_f) / K_f : g^{-1} x_i \in L \otimes_{\Z} \hat{\Z}^p \text{ for all $x_i \in \underline{x}_f$} \}.
    \end{equation}
The isomorphism $G'(\A_f) / K'_f \ra GU(V_0)(\A_f) / K_{0,f} \times U(V)(\A_f)/K_f$ \eqref{equation:uniformization:group_iso} induces a bijection
    \begin{equation}\label{equation:Arch_uniformization:away_infty_special_cycle:iso}
    \mc{D}'(\underline{x}_f) \xra{\sim} GU(V_0)(\A_f) / K_{0,f} \times \mc{D}(\underline{x}_f).
    \end{equation}
    
            \subsection{Framing}
            \label{ssec:Arch_uniformization:framing}
                Fix the isomorphism of Hermitian $\mc{O}_F$-lattices $\Hom_{\mc{O}_F}(L_0, L) \ra L$ sending $x \mapsto x(1)$ (with $1 \in L_0$). This is analogous to $\pmb{\eta}$ from \eqref{equation:non-Arch_uniformization:framing:Tate_module_lattice}.

Given $\a = (A_0, \iota_0, \lambda_0, A, \iota, \lambda, \tilde{\eta}_0, \tilde{\eta}) \in \mc{M}_{\C}^{\mrm{an}}$, a \emph{framing pair} $(\phi_0, \phi)$ for $\a$ consists of isomorphisms of $F$ vector spaces (singular homology)
    \begin{equation}
    \phi_0 \colon H_1(A_0, \Q) \ra V_0 \quad \quad \phi \colon H_1(A, \Q) \ra V
    \end{equation}
such that the induced map
    \begin{equation}
    \phi_0^{-1} \phi \colon \Hom_F( H_1(A_0, \Q), H_1(A, \Q)) \ra \Hom_F(V_0, V) = V
    \end{equation}
is an isomorphism of Hermitian spaces. 

The Hodge structures of weight $-1$ on $H_1(A_0, \Q)$ and $H_1(A, \Q)$ induce a Hodge structure of weight $0$ on $V$, with an associated complex line $F^1 V_{\C} \subseteq V_{\C}^+$. After pullback along the projection isomorphism $V_{\R} \ra V_{\C}^+$ of $F \otimes_{\Q} \R$ vector spaces, the line $F^1 V_{\C} \subseteq V_{\R}$ is a negative definite subspace and hence defines a point $z \in \mc{D}$.
There is a canonical isomorphism of $\C$ vector spaces
    \begin{equation}\label{equation:Arch_uniformization:framing:taut_Hom}
    \Hom_{F \otimes_{\Q} \R} (\Lie A_0, F^0 H_1(A, \Q)_{\C}^+) \cong F^1 V_{\C}.
    \end{equation}

We use the fixed choice of $\sqrt{\Delta}$ to pass between Hermitian/alternating/symmetric forms \crefext{I:ssec:Hermitian_conventions:Hermitian_alternating_symmetric}. This makes $H_1(A_0, \Q)$ and $H_1(A,\Q)$ into Hermitian $F$-modules. Using the $\C$-bilinear extension of the symmetric $\Q$-bilinear trace pairing on $H_1(A, \Q)$, we obtain an induced $\C$-linear identification $F^0 H_1(A,\Q)^+_{\C} \cong \Hom_{\C}((\Lie A)^-, \C)$. 

We equip $F^1 V_{\C} \subseteq V_{\R}$ with the Hermitian metric obtained by restricting the metric on $V_{\R}$. Equip $\Lie A_0$ (resp. $\Lie A$) with the Hermitian metric as normalized in \cref{equation:arith_cycle_classes:Hodge_bundles:Hermitian_normalization} (resp. \cref{equation:arith_cycle_classes:Hodge_bundles:4_pi_gamma_normalization}). Then $(\Lie A)^- \subseteq \Lie A$ inherits a Hermitian metric as well. Under the isomorphism
    \begin{equation}\label{equation:Arch_uniformization:framing:taut_Hom_Hermitian}
    \Hom_{\C}(\Lie A_0, \C) \otimes \Hom_{\C}((\Lie A)^-, \C) \cong F^1 V_{\C}
    \end{equation}
induced by \eqref{equation:Arch_uniformization:framing:taut_Hom}, the Hermitian metric on the left is $-(16 \pi^3 e^{\gamma})^{-1}$ times the Hermitian pairing on the right.

To the datum $(\a, \phi_0, \phi)$, there are associated elements $g_0 \in GU(V_0)/K_{0,f}$ and $g \in U(V)(\A_f)/K_f$ given by $g_0 \coloneqq \phi_0 \circ \tilde{\eta}_0^{-1}$ and $g \coloneqq (\phi_0^{-1} \phi) \circ \tilde{\eta}^{-1}$ (strictly speaking, $\phi_0$ and $\phi$ are tensored with $\A_f$ here, with $H_1(A,\Q) \otimes_{\Q} \A_f = T(A)^0$ (rational ad\`elic Tate module) and similarly for $A_0$).
    
            \subsection{Uniformization}
            \label{ssec:Arch_uniformization:uniformization}
                For any Hermitian matrix $T \in \mrm{Herm}_m(\Q)$ (with $F$-coefficients), define the set
    \begin{equation}
    \mc{Z}(T)_{\C, \mrm{framed}}^{\mrm{an}} \coloneqq \left \{  (\a, \underline{x}, \phi_0, \phi) : \begin{array}{l} \a \in \mc{M}_{\C}^{\mrm{an}} \text{ with $(\a, \underline{x}) \in \mc{Z}(T)^{\mrm{an}}_{\C}$} \\ \text{and $(\phi_0, \phi)$ a framing for $\a$} \end{array} \right \}.
    \end{equation}
There is a canonical injection of sets
    \begin{equation}
    \begin{tikzcd}[row sep = tiny]
    \mc{Z}(T)^{\mrm{an}}_{\C, \mrm{framed}} \arrow{r} & \mc{D} \times G'(\A_f) / K'_f \times V^m
    \\
    (\a, \underline{x}, \phi_0, \phi) \arrow[mapsto]{r} & (z, (g_0, g_0 g), \phi \circ \underline{x} \circ \phi_0^{-1})
    \end{tikzcd}
    \end{equation}
where the Hodge structure $z \in \mc{D}$ and the elements $g_0 \in GU(V_0)(\A_f)/K_{0,f}$ and $g \in U(V)(\A_f) / K_f$ are associated to $(\a, \phi_0, \phi)$ as in Section \ref{ssec:Arch_uniformization:framing}, and $\phi \circ \underline{x} \circ \phi_0^{-1} \in \Hom_F(V_0, V)^m = V^m$ (using the isomorphism $\Hom_{\mc{O}_F}(L_0, L) \cong L$ fixed above). This induces a bijection
    \begin{equation}
    \mc{Z}(T)^{\mrm{an}}_{\C, \mrm{framed}} \xra{\sim} \coprod_{\substack{\underline{x} \in V^m \\ (\underline{x}, \underline{x}) = T}} \mc{D}(\underline{x}_{\infty}) \times \mc{D}'(\underline{x}_f).
    \end{equation}

There is a forgetful map $\mc{Z}(T)^{\mrm{an}}_{\C, \mrm{framed}} \ra \mc{Z}(T)^{\mrm{an}}_{\C}$ sending $(\a, \phi_0, \phi) \mapsto \a$. This is surjective, by the Hasse principle (Landherr's theorem) for Hermitian spaces, and factors through an isomorphism of complex manifolds
    \begin{equation}\label{equation:Arch_uniformization:quotient_iso}
    G'(\Q) \backslash \Biggl ( \coprod_{\substack{\underline{x} \in V^m \\ (\underline{x}, \underline{x}) = T}} \mc{D}(\underline{x}_{\infty}) \times \mc{D}'(\underline{x}_f) \Biggr ) \xra{\sim} \mc{Z}(T)^{\mrm{an}}_{\C, \mrm{framed}}
    \end{equation}
where $G'(\Q)$ acts on $\mc{Z}(T)^{\mrm{an}}_{\C, \mrm{framed}}$ as $(\a, \phi_0, \phi) \mapsto (\a, \g_0 \circ \phi_0, \g \circ \phi)$ for $(\g_0, \g) \in G'(\Q)$ with $\g_0$ and $\g$ having the same similitude factor. The case $T = \emptyset$ (or $T = 0$) gives complex uniformization of $\mc{M}^{\mrm{an}}_{\C, \mrm{framed}}$.

The isomorphism $G' \xra{\sim} GU(V_0) \times U(V)$ (see \eqref{equation:uniformization:group_iso}) induces an isomorphism
    \begin{align}
    & G'(\Q) \backslash \Biggl ( \coprod_{\substack{\underline{x} \in V^m \\ (\underline{x}, \underline{x}) = T}} \mc{D}(\underline{x}_{\infty}) \times \mc{D}'(\underline{x}_f) \Biggr )
    \\
    & \xra{\sim}  \Biggl ( GU(V_0)(\Q) \backslash ( GU(V_0)(\A_f) / K_{0,f}) \Biggr ) \times \Biggl ( U(V)(\Q) \backslash \Biggl ( \coprod_{\substack{\underline{x} \in V^m \\ (\underline{x}, \underline{x}) = T}} \mc{D}(\underline{x}_{\infty}) \times \mc{D}(\underline{x}_f) \Biggr ) \Biggr ) \notag
    \end{align}
where $U(V)(\Q)$ acts on $\mc{D}$ via the $U(V)(\R)$ action, and on $U(V)(\A_f) / K_f$ by left multiplication.
    
            \subsection{Local intersection numbers: Archimedean}  
            \label{ssec:Arch_uniformization:Archimedean}
                Fix $T \in \mrm{Herm}_m(\Q)$ and $y \in \mrm{Herm}_m(\R)_{>0}$ (i.e. $y$ is any positive definite complex Hermitian matrix). Throughout Section \ref{ssec:Arch_uniformization:Archimedean}, we require $m \geq n - 1$ if $T$ is positive definite. If $T$ is singular, we also require $m = n$ and $\rank(T) = n - 1$. 

For such $T$ which are nonsingular, we recall Kudla's Green current $g_{T,y}$ for $\mc{Z}(T)_{\C}^{\mrm{an}}$ (i.e. the unitary analogue studied by Liu \cite[{Proof of Theorem 4.20}]{Liu11}), which is defined via uniformization and star products. For the case of singular $T$, we propose a definition of $g_{T,y}$ by a ``linear invariance'' method, which has some subtleties in the case where $T$ is not $\GL_m(\mc{O}_F)$-equivalent to $\mrm{diag}(0,T^{\flat})$ for $\det T^{\flat} \neq 0$ (``not diagonalizable''). Our treatment of this non-diagonalizable case seems to be new, and is based on the local version which we proposed in \crefext{II:ssec:part_II:Hermitian_domain:corank1_modification}.

Allowing $T$ singular or not for the moment, define the set
    \begin{equation}
    J_{\infty}(T) \coloneqq GU(V_0)(\A_f) / K_{0,f} \times \coprod_{\substack{\underline{x} \in V^m \\ (\underline{x}, \underline{x}) = T}} \mc{D}(\underline{x}_f)
    \end{equation}
In one of our companion papers, we will see that the groupoid $[G'(\Q) \backslash J_{\infty}(T)]$ has with finite stabilizers and finitely many isomorphism classes \crefext{IV:lemma:local_Siegel-Weil:uniformization_degree:main}. Given $j \in J_{\infty}(T)$, we let $\Aut(j) \subseteq G'(\Q)$ be the stabilizer for the action of $G'(\Q)$ on $J_{\infty}(T)$.

For any $\gamma \in \GL_m(\mc{O}_F)$, recall that there is an induced isomorphism $\mc{Z}(T) \xra{\sim} \mc{Z}({}^t \overline{\gamma} T \gamma)$ (i.e. send the tuple of special homomorphisms $\underline{x}$ to $\underline{x} \cdot \gamma$). Similarly, there is an induced isomorphism $\mc{Z}(T)^{\mrm{an}}_{\C, \mrm{framed}} \ra \mc{Z}({}^t \overline{\gamma} T \gamma)^{\mrm{an}}_{\C, \mrm{framed}}$. There is corresponding a bijection $J_{\infty}(T) \ra J_{\infty}({}^t \overline{\gamma} T \gamma)$ (which we denote $j \mapsto j \cdot \gamma$) sending $\underline{x} \mapsto \underline{x} \cdot \gamma$ for $\underline{x} \in V^m$ (acting trivially on the remaining data, i.e. view $J_{\infty}(T)$ as a subset of $G'(\A_f) / K'_f \times V^m$; note $\mc{D}(\underline{x}_f \cdot \gamma) = \mc{D}(\underline{x}_f)$). Note $\Aut(j) = \Aut(j \cdot \gamma)$. 

For each $j \in J_{\infty}(T)$, there is a corresponding map
    \begin{equation}
    \Theta_j \colon \mc{D} \ra \mc{M}^{\mrm{an}}_{\C}
    \end{equation}
induced by the uniformization morphism $\mc{D} \times G'(\A_f) / K'_f \ra \mc{M}^{\mrm{an}}_{\C}$ (consider the projection $J_{\infty}(T) \ra G'(\A_f) / K'_f$; by uniformization of $\mc{M}^{\mrm{an}}_{\C}$, every element of $G'(\A_f) / K'_f$ determines a map $\mc{D} \ra \mc{M}^{\mrm{an}}_{\C}$). For any $\gamma \in \GL_m(\mc{O}_F)$, we have $\Theta_j = \Theta_{j \cdot \gamma}$.

If $\widehat{\mc{E}}_{\C}$ denotes the metrized tautological bundle on $\mc{M}^{\mrm{an}}_{\C}$ (\cref{ssec:ab_var:integral_models}) we have $\Theta_j^* \widehat{\mc{E}}_{\C} \cong \widehat{\mc{E}}$, where $\widehat{\mc{E}}$ is the metrized tautological bundle on $\mc{D}$ (as in \crefext{II:ssec:Hermitian_domain:local_cycles}). By our normalizations, the metric on $\Theta_j^* \widehat{\mc{E}}_{\C}$ is $(16 \pi^3 e^{\gamma})^{-1}$ times the metric on $\widehat{\mc{E}}$ (this normalization constant does not change the Chern form $c_1(\widehat{\mc{E}})$). 

Given any $\underline{x} \in V$ with $(\underline{x}, \underline{x}) = T$, recall that we previously defined a current $[\xi(\underline{x},y)]$ on the Hermitian symmetric domain $\mc{D}$ \crefext{II:ssec:part_II:Hermitian_domain:corank1_modification}.

\begin{definition}
For $T$ as above (singular or not), we define the real current
    \begin{equation}\label{equation:Arch_uniformization:Archimedean:current_definition}
    g_{T,y} \coloneqq \sum_{j \in J} \frac{1}{|\mrm{Aut}(j)|} \Theta_{j,*} [\xi(\underline{x}, y)]
    \end{equation}
on $\mc{M}_{\C}^{\mrm{an}}$, where the sum runs over a set $J \subseteq J_{\infty}(T)$ of representatives for the isomorphism classes of $[G'(\Q) \backslash J_{\infty}(T)]$, where $\underline{x} \in V^m$ is the tuple associated with $j \in J_{\infty}(T)$.
\end{definition}   

In the preceding definition, $\Theta_{j,*}$ denotes pushforward of currents along $\Theta_j$ (for singular $T$, see the convergence estimates in \crefext{II:ssec:Hermitian_domain:convergence}).
The current $g_{T,y}$ does not depend on the choice of $J$, by compatibility of $\mc{D}(\underline{x})$ and $[\xi(\underline{x})]$ with the $U(V)(\R)$ action on $\mc{D}$ \crefext{II:ssec:Hermitian_domain:local_cycles}. It is also compatible with pullback of currents for varying (small) levels $K'_f$. When $T$ is nonsingular, this $g_{T,y}$ agrees with the formulation in \cite[{Proof of Theorem 4.20}]{Liu11} (see also \cite[{\S 15.3}]{LZ22unitary}) up to our different normalization of the Green current (\crefext{II:footnote:Green_current_normalization}).

For any $\gamma \in \GL_m(\mc{O}_F)$, we have 
    \begin{equation}\label{equation:Arch_uniformization:Archimedean:global_linear_invariance}
    g_{{}^t \overline{\gamma} T \gamma, \gamma^{-1} y {}^t \overline{\gamma}^{-1}} = g_{T,y}
    \end{equation}
(``global linear invariance''). This follows from the definition of $g_{T,y}$, from local linear invariance of the currents on $\mc{D}$, and the formulas $\Aut(j) = \Aut(j \cdot \gamma)$ and $\Theta_j = \Theta_{j \cdot \gamma}$.

In all cases, we define the Archimedean intersection number
    \begin{equation}\label{equation:Arch_uniformization:Archimedean}
    \mrm{Int}_{\infty, \mrm{global}}(T, y) \coloneqq \int_{\mc{M}^{\mrm{an}}_{\C}} g_{T,y} \wedge c_1(\widehat{\mc{E}}_{\C}^{\vee})^{n - m}.
    \end{equation}
This is a real number, and the integral is convergent by the estimates in \crefext{II:lemma:convergence:integral_convergence,II:lemma:convergence:KM_form_growth_bound}. It does not depend on the choice of embedding $F \ra \C$. By the compatibility of $g_{T,y}$ with varying small levels $K'_f$, we can extend \eqref{equation:Arch_uniformization:Archimedean} to the case of not-necessarily small level by \eqref{equation:arith_cycle_classes:horizontal:integral} (i.e. cover by a small level and divide by the degree of the cover). In the notation of loc. cit., the stack $\mc{M}$ implicitly has level $K'_{L,f}$ (while we are using the notation $\mc{M}$ to mean arbitrary level $K'_f$ in Section \ref{ssec:Arch_uniformization:Archimedean}).

In all cases (including possibly $K'_f$ not necessarily small level), we have
    \begin{align}\label{equation:Arch_uniformization:Archimedean:local_global_degrees}
    & \mrm{Int}_{\infty, \mrm{global}}(T, y) = \mrm{Int}_{\infty}(T,y) \frac{[K_{L_0,f} : K_{0,f}]}{|\mc{O}_F^{\times}| / h_F} \cdot \deg \Biggl [ U(V)(\Q) \backslash \coprod_{\substack{\underline{x} \in V^m \\ (\underline{x}, \underline{x}) = T}} \mc{D}(\underline{x}_f) \Biggr ]
    \end{align}
by construction, where $\deg$ means (stacky) groupoid cardinality, where $h_F$ is the class number of $\mc{O}_F$, and where 
    \begin{equation}\label{equation:Arch_uniformization:Archimedean:Int_infty_local}
    \mrm{Int}_{\infty}(T,y) \coloneqq \int_{\mc{D}} [\xi(\underline{x}, y)] \wedge c_1(\widehat{\mc{E}}^{\vee})^{n - m}
    \end{equation}
for any $\underline{x} \in V^m$ satisfying $(\underline{x}, \underline{x}) = T$. If there is no such $\underline{x}$, we set $\mrm{Int}_{\infty}(T,y) \coloneqq 0$.

    \clearpage


    \appendix

    \part*{Appendices}

        \section{\texorpdfstring{$K_0$}{K0} groups} \label{appendix:K0}
    
            \subsection{\texorpdfstring{$K_0$}{K0} groups for Deligne--Mumford stacks} 
            \label{appendix:K0:K0_stacky}
                Suppose $\mc{X}$ is a Noetherian Deligne--Mumford stack. There are at least two different ways one might define $K_0$ groups for $\mc{X}$. One way is to define a $K$-theory spectrum for $\mc{X}$ using the $K$-theory spectra of schemes in the small \'etale site of $\mc{X}$, as in \cite[{\S 2}]{Gillet09}. This is the approach used in \cite{HM22}. Another way is to simply mimic a definition of $K_0$ for schemes and consider perfect complexes on the small \'etale site of $\mc{X}$. These two approaches will in general result in different $K_0$ groups \cite[{Remark A.2.4}]{HM22}. At least if $\mc{X}$ is regular (and, say, with the additional running hypotheses of \cite[{Appendix A}]{HM22}), there is a map from the latter $K_0$ group to the former $K_0$ group \cite[{(A.7),(A.8)}]{HM22}.

In this paper, we take the latter approach and mimic constructions for schemes to define $K_0(\mc{X})$. Our definitions and notation will be analogous to those for schemes in \cite[\href{https://stacks.math.columbia.edu/tag/0FDE}{Section 0FDE}]{stacks-project}. When defining dimension/codimension filtrations on $K_0'(\mc{X})$ (with notation and hypotheses as below), we will require existence of a finite flat cover by a scheme (enough for our intended application). A similar approach appears in \cite[{Appendix A}]{YZ17} (at least for $K'_0$), but there the stacks are over a base field. We need a slightly more general setup which allows base schemes such as $\Spec R$ for Dedekind domains $R$.

Suppose $\mc{X}$ is a Deligne--Mumford stack. By an \emph{$\mc{O}_{\mc{X}}$-module}\footnote{This is one of the only places where our conventions differ from the Stacks project \cite[\href{https://stacks.math.columbia.edu/tag/06TF}{Chapter 06TF}]{stacks-project}, which mostly works with sheaves on big sites (say, fppf and \'etale) for general algebraic stacks. Restriction from these two big sites to the small \'etale site (for Deligne--Mumford stacks) induces equivalences on categories of quasi-coherent sheaves. But the equivalences are not compatible with pushforward, and are also not compatible with exactness for $\mc{O}_{\mc{X}}$-modules on big sites versus small sites.},
we mean a sheaf of modules on the small \'etale site\footnote{The small \'etale site is as defined in \cite[{Definition 4.10}]{DM69}: the underlying category has objects which are pairs $(U,f)$ for $f \colon U \ra \mc{X}$ an \'etale morphism (definition as in \cite[\href{https://stacks.math.columbia.edu/tag/0CIL}{Definition 0CIL}]{stacks-project}) from a scheme $U$, and morphisms are pairs $(g, \xi) \colon (U,f) \ra (U',f')$ where $g \colon U \ra U'$ is a $1$-morphism and $\xi \colon f \ra f' \circ g$ is a $2$-isomorphism.} of $\mc{X}$. 
Similarly, \emph{quasi-coherent $\mc{O}_{\mc{X}}$-modules} will mean quasi-coherent sheaves of modules on the small \'etale site. When $\mc{X}$ is locally Noetherian, we will also speak of \emph{coherent $\mc{O}_{\mc{X}}$-modules} on the small \'etale site, which are the same as finitely presented quasi-coherent $\mc{O}_{\mc{X}}$-modules in this situation.

Suppose $\mc{X}$ is a locally Noetherian Deligne--Mumford stack. The category $\textit{Coh}(\mc{O}_{\mc{X}})$ of coherent $\mc{O}_{\mc{X}}$-modules forms a weak Serre subcategory of the abelian category $\textit{Mod}(\mc{O}_{\mc{X}})$ of $\mc{O}_{\mc{X}}$-modules (reduce to the case of small \'etale sites of schemes and apply \cite[\href{https://stacks.math.columbia.edu/tag/05VG}{Lemma 05VG}, \href{https://stacks.math.columbia.edu/tag/0GNB}{Lemma 0GNB}]{stacks-project}). We may form derived categories such as
    \begin{equation}
    D(\mc{O}_{\mc{X}}) \quad D^b(\mc{O}_{\mc{X}}) \quad D_{\textit{perf}}(\mc{O}_{\mc{X}}) \quad D^b_{\textit{Coh}}(\mc{O}_{\mc{X}}) \quad D^b(\textit{Coh}(\mc{O}_{\mc{X}}))
    \end{equation}
which denote the derived category of $\mc{O}_{\mc{X}}$-modules, bounded derived category of $\mc{O}_{\mc{X}}$-modules, derived category of perfect objects (definition as in \cite[\href{https://stacks.math.columbia.edu/tag/08G4}{Section 08G4}]{stacks-project}) in $D(\mc{O}_{\mc{X}})$, bounded derived category of $\mc{O}_{\mc{X}}$-modules with coherent cohomology, and the bounded derived category of coherent $\mc{O}_{\mc{X}}$-modules, respectively.

If $\mc{X}$ happens to be a scheme, then $D_{\textit{perf}}(\mc{O}_{\mc{X}})$ and $D^b_{\textit{Coh}}(\mc{O}_{\mc{X}})$ and $D^b(\textit{Coh}(\mc{O}_{\mc{X}}))$ will agree with the usual constructions using the small Zariski site instead of the small \'etale site, via comparison results such as \cite[\href{https://stacks.math.columbia.edu/tag/08HG}{Lemma 08HG}, \href{https://stacks.math.columbia.edu/tag/071Q}{Lemma 071Q}, \href{https://stacks.math.columbia.edu/tag/05VG}{Lemma 05VG}]{stacks-project}.

\begin{definition}
Let $\mc{X}$ be a locally Noetherian Deligne--Mumford stack. We set
    \begin{equation}
    K_0(\mc{X}) \coloneqq K_0(D_{\textit{perf}}(\mc{O}_{\mc{X}})) \quad \quad K'_0(\mc{X}) \coloneqq K_0(\textit{Coh}(\mc{O}_{\mc{X}})).
    \end{equation}
\end{definition}
Above, the left expression means $K_0$ of a triangulated category and the right expression means $K_0$ of an abelian category. If $\mc{X}$ is a locally Noetherian Deligne--Mumford stack, we have canonical identifications
    \begin{equation}
    K_0(\textit{Coh}(\mc{O}_{\mc{X}})) = K_0(D^b(\textit{Coh}(\mc{O}_{\mc{X}}))) = K_0(D^b_{\textit{Coh}}(\mc{O}_{\mc{X}}))
    \end{equation}
as in \cite[\href{https://stacks.math.columbia.edu/tag/0FDF}{Lemma 0FDF}]{stacks-project} (the case of schemes) by general facts about derived categories (see also \cite[\href{https://stacks.math.columbia.edu/tag/0FCS}{Lemma 0FCS}]{stacks-project}).

If $\mc{X}$ is a quasi-compact locally Noetherian Deligne--Mumford stack, there is an inclusion $D_{\textit{perf}}(\mc{O}_{\mc{X}}) \ra D^b_{\textit{Coh}}(\mc{O}_{\mc{X}})$ and a corresponding group homomorphism $K_0(\mc{X}) \ra K_0'(\mc{X})$. If $\mc{X}$ is a regular locally Noetherian Deligne--Mumford stack (not necessarily quasi-compact), there is an inclusion $D^b_{\textit{Coh}}(\mc{O}_{\mc{X}}) \ra D_{\textit{perf}}(\mc{O}_{\mc{X}})$ and a corresponding group homomorphism $K_0'(\mc{X}) \ra K_0(\mc{X})$. If $\mc{X}$ is a locally Noetherian Deligne--Mumford stack which is both quasi-compact and regular, we have $D_{\textit{perf}}(\mc{O}_{\mc{X}}) = D^b_{\textit{Coh}}(\mc{O}_{\mc{X}})$ and a corresponding isomorphism 
    \begin{equation}
    K_0(\mc{X}) \xra{\sim} K_0'(\mc{X}).    
    \end{equation}
These claims follow from the corresponding facts for schemes \cite[\href{https://stacks.math.columbia.edu/tag/0FDC}{Lemma 0FDC}]{stacks-project} and comparison results mentioned previously.

The derived tensor product $\otimes^{\mbb{L}}$ on $D(\mc{O}_{\mc{X}})$ gives $K_0(\mc{X})$ the structure of a commutative ring. Compatibility of $\otimes^{\mbb{L}}$ with the case when $\mc{X}$ is also a scheme follows from the displayed equation in the proof of \cite[\href{https://stacks.math.columbia.edu/tag/08HF}{Lemma 08HF}]{stacks-project} (comparison between the small Zariski and small \'etale sites).

We next describe dimension and codimension filtrations. Our setup for dimension theory is as in \cite[\href{https://stacks.math.columbia.edu/tag/02QK}{Section 02QK}]{stacks-project}. That is, we work over a locally Noetherian and universally catenary base scheme $S$ with a dimension function $\delta \colon |S| \ra \Z$ (which we typically suppress). Typical setups will be $S = \Spec R$ for $R$ a field or Dedekind domain, where $\delta$ is the dimension function sending closed points to $0$. Any Deligne--Mumford stack $\mc{X}$ which is quasi-separated and locally of finite type over $S$ inherits a dimension function $\delta_{\mc{X}} \colon |\mc{X}| \ra \Z$ (work \'etale locally to pass to the case of schemes; the case of algebraic spaces is \cite[\href{https://stacks.math.columbia.edu/tag/0EDS}{Section 0EDS}]{stacks-project}). If $\mc{X}$ is equidimensional of dimension $n$, then $n - \delta_{\mc{X}}$ is also the codimension function (given by dimensions of local rings on \'etale covers by schemes).

For a scheme $X$ which is locally of finite type over $S$, consider the full subcategory $\textit{Coh}_{\leq d}(\mc{O}_X) \subseteq \textit{Coh}(\mc{O}_X)$ consisting of coherent $\mc{O}_X$-modules $\mc{F}$ with $\dim(\mrm{Supp}(\mc{F})) \leq d$. Then there is an increasing \emph{dimension filtration} on $K'_0(X) = K_0(\textit{Coh}(\mc{O}_X))$ given by the image
    \begin{equation}
    F_d K'_0(X) \coloneqq \opn{im} (K_0(\textit{Coh}_{\leq d}(\mc{O}_X)) \ra K_0(\textit{Coh}(\mc{O}_X)))
    \end{equation}
as in \cite[\href{https://stacks.math.columbia.edu/tag/0FEV}{Section 0FEV}]{stacks-project}. We similarly consider the full subcategory $\textit{Coh}^{\geq m}(\mc{O}_X) \subseteq \textit{Coh}(\mc{O}_X)$ of coherent sheaves supported in codimension $\geq m$, and form the decreasing \emph{codimension filtration}
    \begin{equation}
    F^m K'_0(X) \coloneqq \opn{im} (K_0(\textit{Coh}^{\geq m}(\mc{O}_X)) \ra K_0(\textit{Coh}(\mc{O}_X))).
    \end{equation}
When $X$ is equidimensional of dimension $n$, we have $F^m K'_0(X) = F_{n-m} K'_0(X)$.

For the case of Deligne--Mumford stacks, one could consider naive dimension/codimension filtrations on $K_0'(\mc{X})$ by mimicking the definition for schemes. This may not be a well-behaved notation, and we instead take the filtration defined in \cite[{A.2.3}]{YZ17} (with $\Q$-cofficients).

\begin{definition}\label{definition:stack_dimension_filtration}
For $S$ as above, let $\mc{X}$ be a Deligne--Mumford stack which is quasi-separated and locally of finite type over $S$. Suppose there exists a finite flat surjection $\pi \colon U \ra \mc{X}$ from a scheme $U$. Pick such a morphism $\pi$. 

The \emph{dimension filtration} on $K'_0(\mc{X})_{\Q}$ is the increasing filtration given by
    \begin{equation}
    F_d K'_0(\mc{X})_{\Q} \coloneqq \{ \b \in K'_0(\mc{X})_{\Q} : \pi^* \b \in F_d K'_0(U)_{\Q} \} \subseteq K'_0(\mc{X})_{\Q}
    \end{equation}
for $d \in \Z$. If $\mc{X}$ is equidimensional, we also consider the decreasing \emph{codimension filtration} on $K'_0(\mc{X})_{\Q}$ given by
    \begin{equation}
    F^m K'_0(\mc{X})_{\Q} \coloneqq \{ \b \in K'_0(\mc{X})_{\Q} : \pi^* \b \in F^m K'_0(U)_{\Q} \} \subseteq K'_0(\mc{X})_{\Q}
    \end{equation}
for $m \in \Z$.
\end{definition}

For $\mc{X}$ as in the preceding definition, the filtrations just defined give rise to graded pieces $\mrm{gr}_d K'_0(\mc{X})_{\Q} \coloneqq F_d K'_0(\mc{X})_{\Q} / F_{d-1} K'_0(\mc{X})_{\Q}$ and $\mrm{gr}^m K'_0(\mc{X})_{\Q} \coloneqq F^m K'_0(\mc{X})_{\Q} / F^{m+1} K'_0(\mc{X})_{\Q}$. If $\mc{X}$ as above is equidimensional of dimension $n$, we have $F^m K'_0(\mc{X})_{\Q} = F_{n-m} K'_0(\mc{X})_{\Q}$ for all $m \in \Z$.

\begin{lemma}
With notation as in Definition \ref{definition:stack_dimension_filtration}, the filtrations $F_d K'_0(\mc{X})_{\Q}$ and $F^m K'_0(\mc{X})_{\Q}$ do not depend on the choice of finite flat surjection $\pi \colon U \ra \mc{X}$. If $\mc{X}$ is a scheme, these filtrations recovers the usual filtrations.
\end{lemma}
\begin{proof}
Suppose $X$ is a scheme which is locally of finite type over $S$. If $Z_d(X)$ is the group of $d$-cycles on $X$, recall that there is an identification 
    \begin{equation}
    K_0(\textit{Coh}_{\leq d}(\mc{O}_X) / \textit{Coh}_{\leq {d-1}}(\mc{O}_X)) \xra{\sim} Z_d(X)
    \end{equation}
which is compatible with flat pullback of constant relative dimension and finite pushforward \cite[\href{https://stacks.math.columbia.edu/tag/02S9}{Lemma 02S9}, \href{https://stacks.math.columbia.edu/tag/0FDR}{Lemma 0FDR}]{stacks-project} (see also \cite[\href{https://stacks.math.columbia.edu/tag/02MX}{Lemma 02MX}]{stacks-project}). For any finite flat surjection $\pi \colon U \ra X$ of constant degree $a$, the map $\pi_* \pi^* \colon Z_d(X) \ra Z_d(X)$ is multiplication by $a$.
It follows that $\pi_* \pi^* \colon F_d K'_0(X) / F_{d-1} K'_0(X) \ra F_d K'_0(X) / F_{d-1} K'_0(X)$ is multiplication by $a$. This is an isomorphism after tensoring by $\Q$. When $X$ is equidimensional, this gives the corresponding statement for the codimension filtration as well. This verifies the lemma when $\mc{X}$ is a scheme.

Let $\mc{X}$ is a Deligne--Mumford stack as in the lemma statement. Let $\pi \colon U \ra \mc{X}$ and $\pi' \colon U' \ra \mc{X}$ be two finite flat surjections, for schemes $U$ and $U'$. Consider the fiber product $U \times_{\mc{X}} U'$ with its finite flat projections to $U$ and $U'$. We then apply the preceding discussion to see that the filtrations do not depend on the choice of finite flat surjection.

These arguments are essentially the same as in \cite[{A.2.3}]{YZ17} (the arguments of loc. cit. are over a base field, so we have used different references).
\end{proof}

\begin{remark}
As in \cite[{A.2.3}]{YZ17}, it is possible to have $F_d K'_0(\mc{X})_{\Q} \neq 0$ for $d < 0$ in the situation of Definition \ref{definition:stack_dimension_filtration}. 
\end{remark}

We are mainly interested in $K_0$ groups for the purpose of intersection theory, so we next discuss degree theory over a field. Suppose $S = \Spec k$ for a field $k$, and suppose $\mc{X}$ is a Deligne--Mumford stack which is proper over $S$. Again assuming that $\mc{X}$ admits a finite flat surjection from a scheme, there is a graded group homomorphism $\mrm{gr}_* K'_0(\mc{X})_{\Q} \ra \mrm{Ch}_*(\mc{X})_{\Q}$ as defined in \cite[{A.2.6}]{YZ17} (pass to a finite flat surjection to reduce to the case of schemes).

There is a degree map $\deg \colon \mrm{Ch}_0(\mc{X})_{\Q} \ra \Q$ on $0$-cycles which may be described as follows. Suppose $\mc{Z}$ is a quasi-separated finite type Deligne--Mumford stack over $\Spec k$ with separated diagonal, and assume the underlying topological space $|\mc{Z}|$ is a single point. If $V \ra \mc{Z}$ is any finite flat surjection from a scheme $V$ (which exists as in Remark \ref{remark:existence_finite_covers_stack}), one can check that $V$ is finite over $\Spec k$, and we take
    \begin{equation}\label{equation:zero_cycle_degree_over_field}
    \deg(\mc{Z}) \coloneqq \deg_k(\mc{Z}) \coloneqq \frac{\deg_k (V)}{\deg_{\mc{Z}} (V)}
    \end{equation}
where $\deg_k(V)$ (resp. $\deg_{\mc{Z}}(V)$) is the degree of the finite flat morphism $V \ra \Spec k$ (resp. $V \ra \mc{Z}$).
It is straightforward to see that $\deg_k(\mc{Z})$ does not depend on the choice of $V \ra \mc{Z}$ (compare \cite[{Definition 1.15}]{Vistoli89}). This generalizes immediately to the case where $|\mc{Z}|$ is instead a discrete finite set (add the degrees of its components). When $\mc{Z} = \emptyset$, we take $\deg(\mc{Z}) \coloneqq 0$.

There is an induced degree map 
    \begin{equation}\label{equation:stacky_degree_K0}
    \deg \colon \mrm{gr}_0 K'_0(\mc{X})_{\Q} \ra \Q.
    \end{equation}
Consider a class $\b = \sum_{i} b_i [\mc{F}_i] \in F_0 K'_0(\mc{X})_{\Q}$ where each $\mc{F}_i$ is a coherent sheaf on $\mc{X}$ (we do not assume $[\mc{F}_i] \in F_0 K'_0(\mc{X})_{\Q}$ for any given $i$). Select any finite flat surjection $\pi \colon U \ra \mc{X}$ with $U$ a scheme. If $\pi$ has constant degree $a$, we have
    \begin{equation}\label{equation:euler_char_degree}
    \deg(\b) = \frac{1}{a} \deg(\pi^* \b) = \frac{1}{a} \sum_i b_i \cdot \chi(\pi^* \mc{F}_i)
    \end{equation}
where $\chi$ denotes Euler characteristic. We can give a similar description for general finite flat surjections $\pi$ by decomposing $\mc{X}$ into its connected components. On account of \eqref{equation:euler_char_degree}, we may write $\chi(\b) \coloneqq \deg(\b)$ and think of $\chi \colon \mrm{gr}_0 K'_0(\mc{X})_{\Q} \ra \Q$ as a ``stacky Euler characteristic'' (compare usage in \cite[{Definition 11.4}]{KR14}). We caution, however, that we have only defined $\chi$ on $\mrm{gr}_0 K'_0(\mc{X})_{\Q}$ and have not defined $\chi(\mc{F})$ for a general coherent sheaf $\mc{F}$ on $\mc{X}$.

We conclude this subsection with a lemma which we will use to decompose $K'_0(\mc{X})$ in terms of irreducible components of $\mc{X}$. A similar lemma for formal schemes is \cite[{Lemma B.1}]{Zhang21}.

\begin{lemma}\label{lemma:K'0_component_decomp}
Let $\mc{X}$ be a locally Noetherian Deligne--Mumford stack. Let $\pi_1 \colon \mc{Z}_1 \ra \mc{X}$ and $\pi_2 \colon \mc{Z}_2 \ra \mc{X}$ be closed immersions of Deligne--Mumford stacks with corresponding ideal sheaves $\mc{I}_1$ and $\mc{I}_2$. Assume that the diagonals of $\mc{X}$, $\mc{Z}_1$, and $\mc{Z}_2$ are representable by schemes (e.g. if $\mc{X}$ is separated).

Assume that $\mc{X} = \mc{Z}_1 \cup \mc{Z}_2$ scheme-theoretically (meaning $\mc{I}_1 \cap \mc{I}_2 = 0$). There are mutually inverse isomorphisms
    \begin{equation}
    \begin{tikzcd}[row sep = tiny]
    \frac{K'_0(\mc{X})}{K'_0(\mc{Z}_1 \cap \mc{Z}_2)} \arrow{r} & \arrow{l} \frac{K'_0(\mc{Z}_1)}{K'_0(\mc{Z}_1 \cap \mc{Z}_2)} \oplus \frac{K'_0(\mc{Z}_2)}{K'_0(\mc{Z}_1 \cap \mc{Z}_2)}
    \\
    {[\mc{F}]} \arrow[mapsto]{r} & ({[\pi_1^* \mc{F}_1]}, {[\pi_2^* \mc{F}_2]}) \\
    {[\pi_{1,*} \mc{F}_1]} + {[\pi_{2,*} \mc{F}_2]} & \arrow[mapsto]{l} ({[\mc{F}_1]}, {[\mc{F}_2]}) \rlap{~.}
    \end{tikzcd}
    \end{equation}
Here, $\mc{F}$, $\mc{F}_1$, and $\mc{F}_2$ stand for coherent sheaves on $\mc{X}$, $\mc{Z}_1$, and $\mc{Z}_2$ respectively.
\end{lemma}
\begin{proof}
The condition about diagonals is included for technical convenience. Some additional explanation on notation in the lemma statement: the symbol $\mc{Z}_1 \cap \mc{Z}_2$ denotes the closed substack $\mc{Z}_1 \times_{\mc{X}} \mc{Z}_2$ of $\mc{X}$, with associated ideal sheaf $\mc{I}_1 + \mc{I}_2$, and we have also written $K'_0(\mc{X})/K'_0(\mc{Z}_1 \cap \mc{Z}_2) \coloneqq \coker(K'_0(\mc{Z}_1 \cap \mc{Z}_2) \ra K'_0(\mc{X}))$ etc. (the latter map may not be injective).

Consider the short exact sequence
    \begin{equation}
    0 \ra \mc{O}_{\mc{X}}/(\mc{I}_1 \cap \mc{I}_2) \ra \mc{O}_{\mc{X}}/\mc{I}_1 \oplus \mc{O}_{\mc{X}}/\mc{I}_2 \ra \mc{O}_{\mc{X}}/(\mc{I}_1 + \mc{I}_2) \ra 0.
    \end{equation}
Tensoring by any coherent sheaf $\mc{F}$ on $\mc{X}$, we find that $\mrm{Tor}_1^{\mc{O}_{\mc{X}}}(\mc{F}, \mc{O}_{\mc{X}}/\mc{I}_1)$ is an $\mc{O}_{\mc{X}}/(\mc{I}_1 + \mc{I}_2)$-module, and similarly with $\mc{I}_2$ instead of $\mc{I}_1$. This shows that the displayed projection maps $\mc{F} \mapsto \pi_1^* \mc{F}$ and $\mc{F} \mapsto \pi_2^* \mc{F}$ are well-defined (i.e. that they are additive in short exact sequences and hence descend to the given quotients of $K'_0$-groups). Since $\mrm{Tor}_1^{\mc{O}_{\mc{X}}}(\mc{F}, \mc{O}_{\mc{X}}/(\mc{I}_1 + \mc{I}_2))$ is an $\mc{O}_{\mc{X}}/(\mc{I}_1 + \mc{I}_2)$-module, the Tor long exact sequence of the displayed short exact sequence also shows that $[\mc{F}] = [\mc{F} \otimes_{\mc{O}_{\mc{X}}} \mc{O}_{\mc{X}}/\mc{I}_1] + [\mc{F} \otimes_{\mc{O}_{\mc{X}}} \mc{O}_{\mc{X}}/\mc{I}_2]$ in $K'_0(\mc{X}) / K'_0(\mc{Z}_1 \cap \mc{Z}_2)$.
\end{proof}
                
            \subsection{\texorpdfstring{$K_0$}{K0} groups with supports along finite morphisms} \label{appendix:K0:K0_relative_supports}
                Suppose $X$ is a separated regular Noetherian scheme. There is an established intersection theory for $K_0$ groups with supports along closed subsets of $X$, and the intersection pairing is multiplicative with respect to codimension filtrations (after tensoring by $\Q$) \cite{GS87}. However, we will need a slightly more general setup which allows for ``supports along finite morphisms''. This is needed because the special cycles $\mc{Z}(T) \ra \mc{M}$ (Section \ref{ssec:ab_var:special_cycles}) are not literally cycles but are instead finite unramified morphisms. 

Intersection theory with supports along finite morphisms is also discussed in \cite[{Appendix A.4}]{HM22} for a similar purpose. They are not able to show the codimension filtration is multiplicative in general \cite[{Remark A.4.2}]{HM22}, but can show multiplicativity for intersections against classes of codimension $1$ (in the case of supports along finite unramified morphisms) \cite[{Proposition A.4.4}]{HM22}.

We have two main objectives in this section (besides fixing notation). Our first objective is to comment on another situation where the codimension filtration is multiplicative (namely, when the finite supports become disjoint unions of closed immersions after finite flat base change to a regular scheme). The (short) proof reduces to the case of supports along closed immersions. This is relevant for us because of Lemma \ref{lemma:special_cycles_level_structure}, which says that each special cycle $\mc{Z}(T) \ra \mc{M}$ becomes a disjoint union of closed immersions after finite \'etale base change, at least after inverting the prime $p$ in the cited lemma. 
For $\mc{M}$ associated to a Hermitian lattice of signature $(n - r, r)$, intersecting special cycles over $\mc{M}$ involves multiplicativity for classes of codimension $r$ (not covered by \cite[{Proposition A.4.4}]{HM22} when $r > 1$).

Our second objective is to explain intersection theory with supports along finite morphisms of Deligne--Mumford stacks in terms of the $K_0$ groups of \cref{appendix:K0:K0_stacky}. A stacky theory is also considered in \cite[{Appendix A.4}]{HM22}, but the $K_0$ groups we use are slightly different (as discussed at the beginning of \cref{appendix:K0:K0_stacky}). The setup we consider agrees with \cite[{Appendix A.4}]{HM22} for schemes.

\begin{lemma}\label{lemma:K0_relative}
Consider a $2$-commutative diagram of algebraic stacks
    \begin{equation}
    \begin{tikzcd}
    \mc{Z} \times_{\mc{X}} \mc{W} \arrow{r} \arrow{d} \arrow{dr}{h} & \mc{W} \arrow{d}{g} \\
    \mc{Z} \arrow{r}{f} & \mc{X}
    \end{tikzcd}
    \end{equation}
with outer square $2$-Cartesian, where $\mc{X}$ is a separated regular Noetherian Deligne--Mumford stack and the morphisms $f$ and $g$ (and hence $h$) are finite.

There is a bilinear pairing
    \begin{equation}
    \begin{tikzcd}[row sep = tiny]
    K'_0(\mc{Z}) \times K'_0(\mc{W}) \arrow{r} & K'_0(\mc{Z} \times_{\mc{X}} \mc{W}) \\
    (\mc{F}, \mc{G}) \arrow[mapsto]{r} & \sum_{i} (-1)^i \mrm{Tor}_i^{\mc{O}_{\mc{X}}}(f_*\mc{F}, g_*\mc{G})
    \end{tikzcd}
    \end{equation}
where $\mc{F}$ and $\mc{G}$ stand for coherent $\mc{O}_{\mc{Z}}$-modules and coherent $\mc{O}_{\mc{W}}$-modules, respectively. We have a commutative diagram
    \begin{equation}
    \begin{tikzcd}
    K'_0(\mc{Z}) \times K'_0(\mc{W}) \arrow{r} \arrow{d}{f_* \times g_*} & K'_0(\mc{Z} \times_{\mc{X}} \mc{W}) \arrow{d}{h_*} \\
    K'_0(\mc{X}) \times K'_0(\mc{X}) \arrow{r} & K'_0(\mc{X})
    \end{tikzcd}
    \end{equation}
where vertical arrows are pushforward and the lower horizontal arrow is the bilinear pairing from the ring structure on $K'_0(\mc{X}) \cong K_0(\mc{X})$.
\end{lemma}
\begin{proof}
If $\mc{F}$ is a coherent $\mc{O}_{\mc{Z}}$-module and $\mc{G}$ is a coherent $\mc{O}_{\mc{W}}$-module, we may form the object $(f_* \mc{F} \otimes^{\mbb{L}} g_* \mc{G})$ in $D_{\textit{perf}}(\mc{O}_{\mc{X}})$. For each object $U \ra \mc{X}$ in the small \'etale site of $\mc{X}$ (i.e. $U$ is a scheme with an \'etale morphism to $\mc{X}$), the restriction $(f_* \mc{F} \otimes^{\mbb{L}} g_* \mc{G})|_U \in D_{\textit{perf}}(\mc{O}_{U})$ carries natural $\mc{O}_{U}$-linear actions of $(f_* \mc{O}_{\mc{Z}})|_U $ and $(f_* \mc{O}_{\mc{W}})|_U$.
The resulting cohomology sheaves $\mrm{Tor}_i^{\mc{O}_{\mc{X}}}(f_*\mc{F}, g_* \mc{G}) = H^{-i}(f_* \mc{F} \otimes^{\mbb{L}} g_* \mc{G})$ (a priori coherent $\mc{O}_{\mc{X}}$-modules) are thus sheaves of $(f_* \mc{O}_{\mc{Z}}) \otimes_{\mc{O}_{\mc{X}}} (g_* \mc{O}_{\mc{W}})$-algebras. There is a canonical isomorphism $(f_* \mc{O}_{\mc{Z}}) \otimes_{\mc{O}_{\mc{X}}} (g_* \mc{O}_{\mc{W}}) \ra h_* \mc{O}_{\mc{Z} \times_{\mc{X}} \mc{W}}$ of $\mc{O}_{\mc{X}}$-algebras. Since $h$ is affine, we obtain a lift (up to canonical isomorphism) of each $\mrm{Tor}_i^{\mc{O}_{\mc{X}}}(f_* \mc{F}, g_* \mc{G})$ to a coherent sheaf of $\mc{O}_{\mc{Z} \times_{\mc{X}} \mc{W}}$ modules (to pass between quasi-coherent $h_* \mc{O}_{\mc{Z} \times_{\mc{X}} \mc{W}}$-modules and quasi-coherent $\mc{O}_{\mc{Z} \times_{\mc{X}} \mc{W}}$-modules, we may take an \'etale surjection of $\mc{X}$ from a scheme, use the corresponding result for the small \'etale site of schemes which is \cite[\href{https://stacks.math.columbia.edu/tag/08AI}{Lemma 08AI}]{stacks-project}, and reduce to a statement about glueing data on the small \'etale sites of $\mc{Z} \times_{\mc{X}} \mc{W}$ and $\mc{X}$).

The procedure just described descends to $K'_0$ groups and gives the pairing in the lemma statement.
\end{proof}

We think of the map $K'_0(\mc{Z}) \times K'_0(\mc{W}) \ra K'_0(\mc{Z} \times_{\mc{X}} \mc{W})$ from the preceding lemma as an \emph{intersection pairing} ``with supports along finite morphisms''. 

Next, fix a base scheme $S$ with dimension function $\delta$ as in \cref{appendix:K0:K0_stacky}. Suppose $\mc{X}$ is a Deligne--Mumford stack which is quasi-separated and locally of finite type over $S$. We assume that $\mc{X}$ is equidimensional of dimension $n$, and we also assume that $\mc{X}$ admits a finite flat surjection from a scheme in order to define dimension and codimension filtrations as in Definition \ref{definition:stack_dimension_filtration}.

Consider a finite morphism $f \colon \mc{Z} \ra \mc{X}$ from a Deligne--Mumford stack $\mc{Z}$. We define a ``relative codimension'' filtration on $K'_0(\mc{Z})_{\Q}$ by setting
    \begin{equation}
    F^m_{\mc{X}} K'_0(\mc{Z})_{\Q} \coloneqq F_{n-m} K'_0(\mc{Z})_{\Q}.
    \end{equation}
We similarly set $\mrm{gr}^m_{\mc{X}} K'_0(\mc{Z})_{\Q} \coloneqq F^m_{\mc{X}} K'_0(\mc{Z})_{\Q} / F^{m+1}_{\mc{X}} K'_0(\mc{Z})_{\Q}$. The subscript $\mc{X}$ is meant to remind of the dependence on $\mc{X}$.

\begin{lemma}\label{lemma:K0_relative_with_filtrations}
Let $\mc{X}$ be a regular Noetherian Deligne--Mumford stack which is separated and finite type over $S$. Assume that $\mc{X}$ is equidimensional. Let $f \colon \mc{Z} \ra \mc{X}$ and $g \colon \mc{W} \ra \mc{X}$ be finite morphisms from Deligne--Mumford stacks $\mc{Z}$ and $\mc{W}$.

Assume that there exists a finite flat surjection $\pi \colon U \ra \mc{X}$ with $U$ a regular Noetherian scheme, such that $\mc{Z} \times_{\mc{X}} U \ra U$ and $\mc{W} \times_{\mc{X}} U \ra U$ are both disjoint unions of closed immersions. Then the intersection pairing of Lemma \ref{lemma:K0_relative} restricts to a pairing
    \begin{equation}
    F^s_{\mc{X}} K'_0(\mc{Z})_{\Q} \times F^t_{\mc{X}} K'_0(\mc{W})_{\Q} \ra F^{s+t}_{\mc{X}} K'_0(\mc{Z} \times_{\mc{X}} \mc{W})_{\Q}
    \end{equation}
for any $s,t \in \Z$.
\end{lemma}
\begin{proof}
We use the shorthand $\mc{Z}_U \coloneqq \mc{Z} \times_{\mc{X}} U$ and $\mc{W}_U \coloneqq \mc{W} \times_{\mc{X}} U$. If we abuse notation and also write $\pi$ for the natural projections $\mc{Z}_U \ra \mc{Z}$ and $\mc{W}_U \ra \mc{W}$ and $\mc{Z}_U \times_U \mc{W}_U \ra \mc{Z} \times_{\mc{X}} \mc{W}$, we have $(\pi^* \a) \cdot (\pi^* \b) = \pi^* (a \cdot \b)$ for any $\a \in K'_0(\mc{Z})_{\Q}$ and $\b \in K'_0(\mc{W})_{\Q}$.
By definition of the dimension filtration (Definition \ref{definition:stack_dimension_filtration}), it is enough to check that the intersection pairing over $U$ restricts to
    \begin{equation}
    F^s_{U} K'_0(\mc{Z}_U)_{\Q} \times F^t_{U} K'_0(\mc{W}_U)_{\Q} \ra F_{U}^{s+t} K'_0(\mc{Z}_U \times_{U} \mc{W}_U)_{\Q}
    \end{equation}
(i.e. respects filtrations). We have thus reduced to the case where $\mc{X}$ is a scheme and $\mc{Z} \ra \mc{X}$ and $\mc{W} \ra \mc{X}$ are disjoint unions of closed immersions, and we assume these conditions hold for the rest of the proof. Write $\mc{Z} = \coprod_i \mc{Z}_i$ where each $\mc{Z}_i \ra \mc{X}$ is a closed immersion of schemes, and similarly write $\mc{W} = \coprod_j \mc{W}_j$. By a result of Gillet--Soul\'e \cite[{Proposition 5.5}]{GS87}, the pairing $F^s_{\mc{X}} K'_0(\mc{Z}_i)_{\Q} \times F^t_{\mc{X}} K'_0(\mc{W}_j)_{\Q} \ra K'_0(\mc{Z}_i \times_{\mc{X}} \mc{W}_j)_{\Q}$ factors through $F^{s+t}_{\mc{X}} K'_0(\mc{Z}_i \times_{\mc{X}} \mc{W}_j)_{\Q}$. We may decompose $F_{\mc{X}}^s K'_0(\mc{Z})_{\Q} = \bigoplus_i F_{\mc{X}}^s K'_0(\mc{Z}_i)_{\Q}$ and $F^t_{\mc{X}} K'_0(\mc{W})_{\Q} = \bigoplus_j F^t_{\mc{X}} K'_0(\mc{W}_j)_{\Q}$. Commutativity of the diagram
    \begin{equation}
    \begin{tikzcd}
    F^s_{\mc{X}} K'_0(\mc{Z}_i)_{\Q} \times F^t_{\mc{X}} K'_0(\mc{W}_j)_{\Q} \arrow{r} \arrow{d} & F^{s+t}_{\mc{X}} K'_0(\mc{Z}_i \times_{\mc{X}} \mc{W}_j)_{\Q} \arrow{d} \\
    ( \bigoplus_i F^s_{\mc{X}} K'_0(\mc{Z}_i)_{\Q} ) \times ( \bigoplus_j F^t_{\mc{X}} K'_0(\mc{W}_j)_{\Q} ) \arrow{r} & K'_0(\mc{Z} \times_{\mc{X}} \mc{W})_{\Q}
    \end{tikzcd}
    \end{equation}
for each $i,j$ gives the claim.
\end{proof}

        \section{Quasi-compactness of special cycles}
        \label{appendix:quasi-compactness_lemma}
            Besides fixing notation, the purpose of this appendix is to prove a quasi-compactness statement for special cycles (explicit proofs of other properties, e.g. having finite fibers, are more readily available in the literature, e.g. \cite[{Proposition 2.9}]{KR14}).\footnote{It is possible to have a locally of finite type morphism of locally Noetherian schemes which has finite fibers and satisfies the valuative criterion for properness, but which is not quasi-compact.} A similar proof of quasi-compactness (in the context of special divisors on some orthogonal Shimura varieties) is \cite[{Proposition 2.7.2}]{AGHMP17}.

            \subsection{Terminology}
            \label{appendix:quasi-compactness_lemma:terminology}
                Suppose $A$ and $B$ are abelian schemes
over a base scheme $S$. We write $\underline{\Hom}(A,B)$ for the fppf sheaf (on $S$) of homomorphisms of abelian schemes. Then the sheaf of \emph{quasi-homomorphisms} is $\underline{\Hom}^0(A,B) \coloneqq \underline{\Hom}(A,B) \otimes_{\Z} \Q$. We write $\Hom^0(A,B)$ for the space of global sections and call elements $x \in \Hom^0(A,B)$ \emph{quasi-homomorphisms}, sometimes writing $x \colon A \ra B$. If $S$ is quasi-compact, we have $\Hom^0(A,B) = \Hom(A,B) \otimes_{\Z} \Q$. When $A = B$, we often use the notation $\underline{\End}(A)$, $\underline{\End}^0(A)$, and $\End^0(A)$ instead, and often use the term \emph{quasi-endomorphism}. We write $\operatorname{Isog}(A,B)$ for the set of isogenies $A \ra B$. We write $\underline{\operatorname{Isog}}(A,B) \subseteq \underline{\Hom}(A,B)$ for the subsheaf of sets consisting of isogenies, and $\underline{\operatorname{Isog}}^0(A,B) \subseteq \underline{\Hom}^0(A,B)$ for the subsheaf of \emph{quasi-isogenies}, meaning those quasi-homomorphisms which are locally of the form $m f$ for some isogeny $f$ and some nonzero integer $m \in \Z$. We write $\operatorname{Isog}(A,B)$ (resp. $\operatorname{Isog}^0(A,B)$) for the set of \emph{isogenies} (resp. \emph{quasi-isogenies}), consisting of global sections of $\underline{\operatorname{Isog}}(A,B)$ (resp. $\underline{\operatorname{Isog}}^0(A,B)$). We write $\mrm{Isog}(A)$ and $\mrm{Isog}^0(A)$ for the self-isogenies and self quasi-isogenies of $A$. A \emph{quasi-polarization} of $A$ is a quasi-isogeny $A \ra A^{\vee}$ which is locally of the form $m \lambda$ for some polarization $\lambda$ and some positive integer $m \in \Z_{>0}$.

Suppose the abelian schemes $A$ and $B$ are equipped with quasi-polarizations $\lambda_A \colon A \ra A^{\vee}$ and $\lambda_B \colon B \ra B^{\vee}$. Given any $x \in \Hom^0(B,A)$ with dual $x^{\vee} \in \Hom^0(A^{\vee}, B^{\vee})$, we set $x^{\dagger} \coloneqq \lambda_B^{-1} \circ x^{\vee} \circ \lambda_A \in \Hom^0(A,B)$, and call the resulting map $\dagger \colon \Hom^0(B,A) \ra \Hom^0(A,B)$ the \emph{Rosati involution}.
Given $m$-tuples $\underline{x}, \underline{\smash{y}} \in \Hom^0(B,A)^m$ with $\underline{x} = [x_1, \ldots, x_m]$ and $\underline{\smash{y}} = [y_1, \ldots, y_m]$, we write $(\underline{x}, \underline{\smash{y}})$ for the $m \times m$ matrix whose $i,j$-th entry is $x^{\dagger}_i y_j$. We say that $(\underline{x}, \underline{x})$ is the \emph{Gram matrix} of $\underline{x}$.
If $S = \Spec k$ for a field $k$, the $\Q$-bilinear pairing
    \begin{equation}
    \begin{tikzcd}[row sep = tiny]
    \Hom^0(B,A) \times \Hom^0(B,A) \arrow{r} & \Q
    \\
    x,y \arrow[mapsto]{r} & \tr(x^{\dagger} y)
    \end{tikzcd}
    \end{equation}
is symmetric and positive definite (``positivity of the Rosati involution''), where $\operatorname{tr} \colon \End^0(A) \ra \Q$ is the trace for $\End^0(A)$ acting on the $\Q$-vector space $\End^0(A)$ by left multiplication.
            
            \subsection{Proof}
            \label{appendix:quasi-compactness_lemma:proof}
                We continue in the setup of Appendix \ref{appendix:quasi-compactness_lemma:terminology}.

Given any $y \in \End^0(B)$, define a functor $\mc{Z}(y) \colon (\Sch/S)^{\op} \ra \Set$ as
    \begin{equation}
    \mc{Z}(y) \coloneqq \{x \in \underline{\Hom}(B,A) : x^{\dagger} x = y \}.
    \end{equation}
We will check that $\mc{Z}(y)$ is representable by a scheme which is finite, unramified, and of finite presentation over $S$.

\begin{lemma}\label{lemma:cycle_finite_unr_representable_revamped}
The functor $\mc{Z}(y)$ is represented by a scheme over $S$. The structure morphism $\mc{Z}(y) \ra S$ is separated and locally of finite presentation.
\end{lemma}
\begin{proof}
By a standard limit argument (e.g. using \cite[\href{https://stacks.math.columbia.edu/tag/01ZM}{Lemma 01ZM}]{stacks-project}) we may reduce to the case where $S$ is Noetherian, affine, and connected. It is also enough to check the case where $\lambda_A$ and $\lambda_B$ are polarizations, not just quasi-polarizations.

Existence of the product polarization $\lambda_B \times \lambda_A$ on $B \times A$ implies that $B \times A$ admits a relatively ample line bundle over $S$. Thus the Hilbert functor $\mrm{Hilb}_{B \times A}$ is represented by a scheme, each of whose connected components is locally projective over $S$ (in the sense of \cite[\href{https://stacks.math.columbia.edu/tag/01W8}{Definition 01W8}]{stacks-project}), see \cite[{Theorem 5.15}]{Nitsure05} and \cite[\href{https://stacks.math.columbia.edu/tag/0DPF}{Lemma 0DPF}]{stacks-project}. By \cite[\href{https://stacks.math.columbia.edu/tag/0D1B}{Lemma 0D1B}]{stacks-project}, we know there is a locally closed immersion 
    \[
    \mc{Z}(y) \ra \mrm{Hilb}_{B \times A}
    \]
which sends $x \colon B \ra A$ to its graph $(1 \times x) \colon B \ra B \times A$. In particular, $\mc{Z}(y)$ is represented by a scheme which is separated and locally of finite presentation over $S$.
\end{proof}

\begin{lemma}\label{lemma:cycle_quasi-compact_revamped}
The structure morphism $\mc{Z}(y) \ra S$ is quasi-compact.
\end{lemma}
\begin{proof}
Again, we may reduce to the case where $S$ is affine, Noetherian, and connected by a standard limit argument. It is also enough to check the case where $\lambda_A$ and $\lambda_B$ are polarizations, not just quasi-polarizations. 

Consider the graph morphisms
    \begin{align*}
    B & \xra{1 \times \lambda_B} B \times B^{\vee} \\
    A & \xra{1 \times \lambda_A} A \times A^{\vee} \rlap{~.}
    \end{align*}
If $\mc{P}_{B}$ and $\mc{P}_A$ denote the Poincar\'e bundles on $B \times B^{\vee}$ and $A \times A^{\vee}$ respectively, we know that $\mc{L}_B \coloneqq (1 \times \lambda_B)^* \mc{P}_{B}$ and $\mc{L}_A \coloneqq (1 \times \lambda_A)^* \mc{P}_A$ are relatively ample line bundles on $B$ and $A$, respectively, over $S$. If $\pi_{B} \colon B \times A \ra B$ and $\pi_A \colon B \times A \ra A$ are the natural projections, we know $\mc{E} \coloneqq \pi_{B}^* \mc{L}_B \otimes \pi_{A}^* \mc{L}_A$ is a relatively ample line bundle on $B \times A$. Moreover, $\mc{E}$ is isomorphic to the pullback of the Poincar\'e bundle $\mc{P}_{B \times A}$ along the graph of the polarization $\lambda_B \times \lambda_A$ of $B \times A$. Let $m \in \Z_{\geq 1}$ be any integer such that $m \cdot \lambda_B$ and $m^2 \cdot y$ are both honest homomorphisms (rather than quasi-homomorphisms).

As above, write $\mrm{Hilb}_{B \times A}$ for the Hilbert scheme associated with $B \times A$. Given a numerical polynomial $P \colon \Z \ra \Z$, we write $\mrm{Hilb}^P_{B \times A} \subseteq \mrm{Hilb}_{B \times A}$ for the open and closed subscheme corresponding to the Hilbert polynomial $P$ with respect to the line bundle $\mc{E}^{\otimes m^2}$ on $B \times A$.
That is, for an $S$-scheme $T$, we have
    \[
    \mrm{Hilb}_{B \times A}^P(T) \coloneqq \{ Z \in \mrm{Hilb}_{B \times A}(T) : \chi(Z_t, \mc{E}^{\otimes m^2 n}|_{Z_t}) = P(n) \text{ for all $n \in \Z$ and $t \in T$} \}
    \]
(where $Z_t$ is the fiber of $Z \ra T$ over $t \in T$ and $\chi$ denotes Euler characteristic). We know that $\mrm{Hilb}_{B \times A}^P(T)$ is locally projective over $S$ \cite[{Theorem 5.15}]{Nitsure05}, hence quasi-compact over $S$. 

As in the proof of Lemma \ref{lemma:cycle_finite_unr_representable_revamped}, there is a locally closed immersion $\mc{Z}(y) \ra \mrm{Hilb}_{B \times A}$ which sends $x \in \mc{Z}(y)$ to its graph $1 \times x \colon B \ra B \times A$. To show that $\mc{Z}(y)$ is quasi-compact, it is enough to check that $\mc{Z}(y) \ra \mrm{Hilb}_{B \times A}$ factors through $\mrm{Hilb}^P_{B \times A}$ for some fixed numerical polynomial $P$ (possibly depending on $y$).

 Consider the line bundle $\mc{F} \coloneqq \mc{L}_B^{\otimes m^2} \otimes ((1 \times \lambda_B)^* (m^2 y \times 1)^* \mc{P}_B)$ on $B$. For any point $s \in S$, there is a numerical polynomial $P \colon \Z \ra \Z$ such that
    \begin{equation}
    P(n) = \chi(B_s, \mc{F}^{\otimes n} |_{B_s}) \quad \text{for all $n \in \Z$}
    \end{equation}
as in \cite[\href{https://stacks.math.columbia.edu/tag/0BEM}{Lemma 0BEM}]{stacks-project}. The polynomial $P$ does not depend on $s$ because $S$ is connected and the Euler characteristics are locally constant as a function of $s$ (using flatness and properness and the standard facts \cite[\href{https://stacks.math.columbia.edu/tag/0BDJ}{Lemma 0BDJ}]{stacks-project} and \cite[\href{https://stacks.math.columbia.edu/tag/07VJ}{Section 07VJ}]{stacks-project}).

Let $T$ be a scheme over $S$, and suppose $x \in \mc{Z}(y)(T)$. View $x$ as an element of $\mrm{Hilb}_{B \times A}(T)$ as above. We claim that $x \in \mrm{Hilb}^P_{B \times A}(T)$. By taking a base-change to $T$, we may assume $T = S$ without loss of generality (to lighten notation). 
It is enough to check $\mc{F} \cong (1 \times x)^* \mc{E}^{\otimes m^2}$.

First observe $(1 \times x)^* \mc{E}^{\otimes m^2} \cong \mc{L}_B^{\otimes m^2} \otimes x^* \mc{L}_A^{\otimes m^2}$. It is thus enough to verify the identity $x^* \mc{L}_A^{\otimes m^2} \cong (1 \times \lambda_B)^* (m^2 y \times 1)^* \mc{P}_{B}$.
Consider the commutative diagram
    \[
    \begin{tikzcd}
    A \arrow{r}{1 \times \lambda_A} & A \times A^{\vee} \\
    B \arrow{u}{m x} \arrow{r}{1 \times \lambda_B} & B \times B^{\vee} \arrow{d}[swap]{m^2 y \times 1} \arrow{u}{m x \times m x^{\dagger \vee}} \arrow{r}{m x \times 1} & A \times B^{\vee} \arrow{dl}{m x^\dagger \times 1} \arrow{ul}[swap]{1 \times m x^{\dagger \vee}} \\
    & B \times B^{\vee}& \rlap{~.}
    \end{tikzcd}
    \]
There exists an isomorphism $(m x^{\dagger} \times 1)^* \mc{P}_B \cong (1 \times m x^{\dagger \vee})^* \mc{P}_A$ (this characterizes $m x^{\dagger \vee}$ as the dual of $m x^{\dagger}$). Recall also that $m^* \mc{L}_A \cong \mc{L}_A^{\otimes m^2}$ (consider a similar diagram as above, with $A = B$ and $x = y = 1$, and recall that the pullback of $\mc{P}_A$ along $(m \times 1) \colon A \times A^{\vee} \ra A \times A^{\vee}$ is isomorphic to $\mc{P}_A^{\otimes m^2}$ because $m = m^{\vee}$). These facts prove the claimed identity $x^* \mc{L}_A^{\otimes m^2} \cong (1 \times \lambda_B)^* (m^2 y \times 1)^* \mc{P}_{B}$.
\end{proof}

\begin{lemma}\label{lemma:quasi-compactness_lemma}
The functor $\mc{Z}(y)$ is represented by a scheme over $S$, and the structure morphism $\mc{Z}(y) \ra S$ is finite, unramified, and of finite presentation.
\end{lemma}
\begin{proof}
Again, we may reduce to the case where $S$ is Noetherian by a standard limit argument. By Lemmas \ref{lemma:cycle_finite_unr_representable_revamped} and \ref{lemma:cycle_quasi-compact_revamped}, we already know that $\mc{Z}(y)$ is represented by a scheme which is separated and of finite presentation over $S$.

To see that $\mc{Z}(y) \ra S$ is proper, we can use the valuative criterion for discrete valuation rings \cite[\href{https://stacks.math.columbia.edu/tag/0207}{Lemma 0207}]{stacks-project} because $S$ is Noetherian. This valuative criterion holds by the N\'eron mapping property for abelian schemes over discrete valuation rings.

For unramifiedness, it is enough to check that $\mc{Z}(y) \ra S$ is formally unramified (i.e. satisfies the infinitesimal lifting criterion of \cite[\href{https://stacks.math.columbia.edu/tag/02HE}{Lemma 02HE}]{stacks-project}). Formal unramifiedness holds because of rigidity for morphisms of abelian schemes as in \cite[{Corollary 6.2}]{MFK94}.

Since unramified morphisms of schemes are locally quasi-finite, and since proper locally quasi-finite morphisms of schemes are finite, the lemma is proved.
\end{proof}

Recall that if $\mc{X}$ and $\mc{Y}$ are categories fibered in groupoids over the fppf site of some base scheme with $\mc{Y}$ being a Deligne--Mumford stack, and if $f \colon \mc{X} \ra \mc{Y}$ is a morphism which is representable by algebraic spaces, then $\mc{X}$ is also a Deligne--Mumford stack \cite[\href{https://stacks.math.columbia.edu/tag/04T0\#comment-2142}{Comment 2142}]{stacks-project}. This can be used in combination with Lemma \ref{lemma:quasi-compactness_lemma} to verify that various stacks in this work are Deligne--Mumford.

        \section{Miscellany on \texorpdfstring{$p$}{p}-divisible groups}
        \label{appendix:pDiv_prelim}
            We collect some terminology/notation and miscellaneous facts about $p$-divisible groups, which we use freely.
    
            \subsection{Terminology}
            \label{appendix:pDiv_prelim:terminology}
                Suppose $S$ is a formal scheme\footnote{
The formal schemes we use are the ``pr\'esch\'emas formels'' of \cite[{\S 10}]{EGAI}. Given a formal scheme, the notation $(\Sch/S)_{fppf}$ means the site whose objects are morphisms $T \ra S$ for schemes $T$, where coverings are fppf.
} and suppose $\mc{P}$ is a property of morphisms of schemes which is fppf local on the target and stable under arbitrary base-change. A sheaf $X$ on $(\Sch/S)_{fppf}$ is \emph{represented by a relative scheme with property $\mc{P}$ over $S$} if, for every scheme $T$ over $S$, the restriction sheaf $X|_T$ is represented by a scheme with property $\mc{P}$ over $T$.

Fix a prime $p$. A \emph{$p$-divisible group} over a formal scheme $S$ is a sheaf $X$ of abelian groups on $(\Sch/S)_{fppf}$ which satisfies the following conditions.
    \begin{enumerate}[(1)]
    \item ($p$-divisibility) The multiplication by $p$ map $[p] \colon X \ra X$ is a surjection of sheaves.
    \item ($p^{\infty}$-torsion) The natural map $X[p^{\infty}] \coloneqq \varinjlim_n X[p^n] \ra X$ is an isomorphism, where $X[p^n] \subseteq X$ are the $p^n$-torsion subsheaves.
    \item (representable $p$-power-torsion) The sheaves $X[p^n]$ are represented by finite locally free relative schemes over $S$ for all $n \geq 1$.
    \end{enumerate}
If $S$ is an adic (e.g. locally Noetherian) formal scheme and $\ms{I}$ is an ideal sheaf of definition on $S$, giving a $p$-divisible group over $S$ is the same as giving $p$-divisible groups $X_n$ over each scheme $S_n \coloneqq (S, \mc{O}_S/\ms{I}^n)$ with isomorphisms $X_{n+1}|_{S_n} \xra{\sim} X_n$.

For a general formal scheme $S$, we say a $p$-divisible group $X$ over $S$ has \emph{height} $h$ if $X[p]$ is finite locally free relative scheme over $S$ of degree $p^h$. In general, $h$ is understood as a locally constant function on $S$.

If $p$ is locally topologically nilpotent on $S$ (equivalently, $S$ is a formal scheme over $\Spf \Z_p$) and if $X$ is a $p$-divisible group over $S$, there is an associated sheaf $\Lie X$ on $(\Sch/S)_{fppf}$ (constructed as in \cite[{Definition 3.2}]{SGA3II}). By work of Messing \cite[{Theorem 3.3.18}]{Messing72}, it is known that $\Lie X$ is a finite locally free sheaf of modules on $(\Sch/S)_{fppf}$. We refer to the dual $\Omega_X \coloneqq (\Lie X)^{\vee}$ as a \emph{Hodge bundle}. If $r$ is the rank of $\Lie X$, we say that $X$ has \emph{dimension $r$} (in general, $r$ is a locally constant $\Z_{\geq 0}$-valued function). In this case, we write $\omega_X \coloneqq \bigwedge^r \Omega_X$ for the top exterior power and also call $\omega_X$ a \emph{Hodge bundle}.

If $p$ is locally topologically nilpotent on the formal scheme $S$, a \emph{formal $p$-divisible group} $X$ over $S$ is a $p$-divisible group over $S$ such that, fppf (equivalently, Zariski) locally on any $T \in \mrm{Obj}(\Sch/S)_{fppf}$, the pointed fppf sheaf $X$ is isomorphic to $\Spf \mc{O}_{T}[[x_1,\ldots,x_r]]$ for some $r$ (possibly varying). 
See \cite[{Proposition II.4.4}]{Messing72} for equivalent characterizations.

Given $p$-divisible groups $X$ and $Y$ over a general formal scheme $S$, a \emph{quasi-homomorphism} is a global section of the sheaf $\underline{\Hom}(X,Y) \otimes_{\Z} \Q$ on $(\Sch/S)_{fppf}$. We write $\Hom^0(X,Y)$ for the space of quasi-homomorphisms $X \ra Y$, and similarly $\End^0(X) = \Hom^0(X,X)$. Given a quasi-compact scheme $T$ with a map $T \ra S$, we have $\Hom^0(X_T, Y_T) = \Hom(X_T, Y_T) \otimes_{\Z} \Q$. If $X$ and $Y$ are equipped with an action by a ring $R$, then $\Hom^0_R(X,Y)$ will denote the $R$-linear quasi-homomorphisms.

A morphism $f \colon X \ra Y$ of $p$-divisible groups over $S$ is an \emph{isogeny} if $f$ is a surjection of fppf sheaves and $\ker f$ is represented by a finite locally free relative scheme over $S$. If $\ker f$ is finite locally free of rank $p^r$, we say that $f$ has \emph{degree} $p^r$ and \emph{height} $r$. A \emph{quasi-isogeny} $f \colon X \ra Y$ is a quasi-homomorphism which, locally on $(\Sch/S)_{fppf}$, is of the form $f = p^n g$ for $n \in \Z$ and an isogeny $g$. If the $p$-divisible group $X$ has height $h$, such a quasi-isogeny $f = p^n g$ is said to have \emph{degree} $p^{nh} \deg(g)$ and \emph{height} $nh + \mathrm{height}(g)$. We write $\mrm{Isog}(X,Y)$ (resp. $\mrm{Isog}^0(X,Y)$) for the isogenies (resp. quasi-isogenies) $X \ra Y$. We write $\mrm{Isog}(X)$ (resp. $\mrm{Isog}^0(X)$) for self-isogenies (resp. self quasi-isogenies) of $X$.

A $p$ divisible group $X$ over $S$ is \emph{\'etale} if $X[p]$ is an \'etale relative scheme. This implies that each $X[p^n]$ is an \'etale relative scheme. If $R$ is a Noetherian Henselian local ring, we say that a $p$-divisible group $X$ over $\Spec R$ is \emph{connected} if $X[p]$ is connected. This implies that each $X[p^n]$ is connected.

Given any $p$-divisible group $X$ over a general formal scheme $S$, there is a \emph{dual} $p$-divisible group $X^{\vee}$. A \emph{polarization} of $X$ is an isogeny $\lambda \colon X \ra X^{\vee}$ satisfying $\lambda^{\vee} = - \lambda$. The polarization is \emph{principal} if $\lambda$ is an isomorphism.
A \emph{quasi-polarization} is a quasi-isogeny $f \colon X \ra X^{\vee}$ such that $m f$ is a polarization for some $m \in \Q_p^{\times}$.
Suppose $X$ and $Y$ are $p$-divisible groups over $S$ with quasi-polarizations $\lambda_X \colon X \ra X^{\vee}$ and $\lambda_Y \colon Y \ra Y^{\vee}$. Given any $x \in \Hom^0(Y,X)$ with dual $x^{\vee} \in \Hom^0(X^{\vee}, Y^{\vee})$, we set $x^{\dagger} \coloneqq \lambda_Y^{-1} \circ x^{\vee} \circ \lambda_X \in \Hom^0(X,Y)$, and call the resulting map $\dagger \colon \Hom^0(Y,X) \ra \Hom^0(X,Y)$ the \emph{Rosati involution}.

Over an algebraically closed field, we say that a $p$-divisible group is supersingular if all slopes of its isocrystal are equal to $1/2$, and we say that it is ordinary if all slopes of its isocrystal are either $0$ or $1$. A $p$-divisible group over an arbitrary formal scheme is \emph{supersingular} (resp. \emph{ordinary}) if it is supersingular (resp. ordinary) for every geometric fiber.

Over any algebraically closed field, there is a unique \'etale $p$-divisible group of height $r$ (namely the constant sheaf $(\underline{\Q_p / \Z_p})^r$). 
Over any algebraically closed field of characteristic $p$, there is also a unique $p$-divisible group of height $r$ with all slopes of its isocrystal being $1$ (namely $\pmb{\mu}_{p^{\infty}}^r \coloneqq (\varinjlim_{e} \pmb{\mu}_{p^e})^r \cong (\underline{\Q_p / \Z_p}^{\vee})^r$, given by $p$-th power roots of unity). 
Since the connected \'etale exact sequence of any $p$-divisible group over a perfect field is (canonically) split, we conclude that $\pmb{\mu}_{p^{\infty}}^{n - r} \times (\underline{\Q_p/\Z_p})^r$ is the unique ordinary $p$-divisible group of height $n$ and dimension $n - r$ over any algebraically closed field.

By \emph{Drinfeld rigidity} we mean the following phenomenon: if $S_0 \ra S$ is a finite order thickening of schemes over $\Spf \Z_p$, and $X, Y$ are $p$-divisible groups over $S$, any quasi-homomorphism of $X \ra Y$ over $S_0$ lifts uniquely to a quasi-homomorphism over $S$ \cite[{Theorem 2.2.3}]{Andre03} (alternative proof: Grothendieck--Messing theory).

If $A$ is a relative abelian scheme over a general formal scheme $S$, there is an associated $p$-divisible group $A[p^{\infty}] \coloneqq \varinjlim_n A[p^n]$, where $A[p^n]$ is the $p^n$-torsion subfunctor of $A$. If $p$ is locally topologically nilpotent on $S$, there is a canonical identification $\Lie A \cong \Lie A[p^{\infty}]$.

Given a $p$-divisible group $X$ over a formal scheme $S$ and given a finite free $\Z_p$-module $M$ of some rank $d \geq 0$, there is the \emph{Serre tensor construction} $p$-divisible group $X \otimes_{\Z_p} M$ given by the functor
    \begin{equation}\label{equation:appendix:Serre_tensor}
    (X \otimes_{\Z_p} M)(T) \coloneqq X(T) \otimes_{\Z_p} M
    \end{equation}
for schemes $T$ over $S$. Any choice of $\Z_p$-basis for $M$ gives an isomorphism $X \otimes_{\Z_p} M \cong X^{d}$ as $p$-divisible groups. This construction is functorial in $M$: in particular, any $\Z_p$-algebra $R$ acting on $M$ also acts on $X \otimes_{\Z_p} M$. The resulting $R$-action on $X \otimes_{\Z_p} M$ is the \emph{Serre tensor $R$-action}.
There is a canonical identification $(X \otimes_{\Z_p} M)^{\vee} \cong X^{\vee} \otimes_{\Z_p} M^{\vee}$ where $M^{\vee} \coloneqq \Hom_{\Z_p}(M, \Z_p)$. More generally, see \cite[{\S 7}]{Conrad04}.

            \subsection{Isogeny criterion}
            \label{appendix:pDiv_prelim:isogeny_criterion}
                We explain a criterion for a morphism of $p$-divisible groups to be an isogeny (Lemma \ref{lemma:isogeny_criterion}). This should be well-known.\footnote{The only reference I know is the sketch in \cite[{Lemme 9}]{Fargues05}.
We spell out the argument for completeness.}

\begin{lemma}\label{lemma:surj_finite_locally_free_fppf}
Let $S$ be a scheme, and let $H$, $G$, and $Q$ be commutative group schemes over $S$ which are locally of finite presentation. Suppose
    \[
    0 \ra H \ra G \xra{f} Q \ra 0
    \]
is an exact sequence of fppf sheaves of abelian groups. If $G \ra S$ is finite locally free and $Q \ra S$ is separated, then
    \begin{enumerate}[(1)]
        \item The map $f \colon G \ra Q$ is finite locally free.
        \item The group schemes $Q$ and $H$ are finite locally free over $S$.
    \end{enumerate}
\end{lemma}
\begin{proof}
Since $f$ is a surjection of fppf sheaves, it is a surjection on underlying topological spaces. We also know that $f$ is locally of finite presentation because both $G$ and $Q$ are locally of finite presentation over $S$ \cite[\href{https://stacks.math.columbia.edu/tag/00F4}{Lemma 00F4}]{stacks-project}. Since $G \ra S$ is flat, the fibral flatness criterion \cite[{11.3.11}]{EGAIV3} implies that flatness of $f$ may be checked fiberwise over $S$, i.e. it is enough to check flatness of the base-change $G_{k(s)} \ra Q_{k(s)}$ for each $s \in S$. The exact sequence
    \[
    0 \ra H_{k(s)} \ra G_{k(s)} \ra Q_{k(s)} \ra 0
    \]

Since $H = \ker(f)$ and $f$ is an fppf morphism, we know $H \ra S$ is fppf as well. Since $Q \ra S$ is separated, the identity section $S \ra Q$ is a closed immersion, hence $H = \ker(f)$ is a closed subscheme of $G$. Since $G \ra S$ is finite, we conclude that $H \ra S$ is also finite, hence finite locally free.

We have already seen that $Q \ra S$ is flat, proper, and locally of finite presentation. To check that $Q \ra S$ is finite, it is enough to check that it has finite fibers, which follows because $G \ra Q$ is surjective and $G \ra S$ is finite.
\end{proof}

\begin{lemma}\label{lemma:isogeny_criterion}
Let $S$ be a formal scheme. Let $f \colon X \ra Y$ be a homomorphism of $p$-divisible groups over $S$. Then $f$ is an isogeny if and only if, locally on $(\Sch/S)_{fppf}$, there exists a homomorphism $g \colon Y \ra X$ such that
    \[
    g \circ f = [p^N] \quad \quad f \circ g = [p^N]
    \]
for some integer $N \geq 0$, where $[p^N]$ denotes multiplication by $p^N$.

Moreover, given an isogeny $f$, such $g,N$ will exist globally on $S$ if $S$ is quasi-compact or has finitely many connected components. If $f$ is an isogeny of constant degree $p^n$, we may take $N = n$.
\end{lemma}
\begin{proof}

If $f \colon X \ra Y$ is an isogeny, then $Y$ is the fppf sheaf quotient of $X$ by $\ker (f)$. If $S$ is a quasi-compact formal scheme or if $S$ has finitely many connected components, we have $\ker f \subseteq X[p^N]$ for $N$ large, so $g \circ f = [p^N]$ for some homomorphism $g \colon Y \ra X$. We also have $f \circ g \circ f = [p^N] \circ f$. Since $f$ is an epimorphism of fppf sheaves, we conclude that $f \circ g = [p^N]$.

Conversely, suppose that locally on $(\Sch/S)_{fppf}$ there exists a homomorphism $g \colon Y \ra X$ and an integer $N \geq 0$ as in the lemma statement. Since the property of being an isogeny may be checked locally on $(\Sch/S)_{fppf}$, we may assume that $S$ is a scheme and that $g,N$ exist globally on $S$. Since $f \circ g = [p^N]$, we see that $f$ is a surjection of fppf sheaves. It remains only to check that $\ker f$ is representable by a finite locally free group scheme over $S$.

We know that $\ker(f) \subseteq X[p^N]$ and $\ker(g) \subseteq Y[p^N]$. We have $\ker(f) = \ker(X[p^N] \ra Y[p^N])$ and $\ker(g) = \ker(Y[p^N] \ra X[p^n])$. Since $X[p^N]$ and $Y[p^N]$ represented by finite locally free group schemes over $S$, we see that $\ker(f)$ and $\ker(g)$ are represented by schemes which are finite and locally of finite presentation over $S$.

We have short exact sequences
    \begin{align*}
    & 0 \ra \ker(f) \ra X[p^N] \xra{f} \ker(g) \ra 0
    \\
    & 0 \ra \ker(g) \ra Y[p^N] \xra{g} \ker(f) \ra 0
    \end{align*}
of fppf sheaves of abelian groups. By Lemma \ref{lemma:surj_finite_locally_free_fppf}, we conclude that $\ker(f)$ and $\ker(g)$ are finite locally free group schemes over $S$.
\end{proof}

\begin{lemma}\label{lemma:qisog_invertible}
Let $S$ be a formal scheme. Let $X$ and $Y$ be $p$-divisible groups over $S$. Then $f \in \Hom^0(X,Y)$ is a quasi-isogeny if and only if it is invertible, meaning there exists $g \in \Hom^0(Y,X)$ (necessarily unique) with $f \circ g = \mrm{id}_Y$ and $g \circ f = \mrm{id}_X$.
\end{lemma}
\begin{proof}
Invertibility and the property of being a quasi-isogeny can both be checked locally on $(\Sch/S)_{fppf}$, so the lemma follows from Lemma \ref{lemma:isogeny_criterion}.
\end{proof}

            \subsection{\texorpdfstring{$p$}{p}-divisible groups over \texorpdfstring{$\Spec A$}{Spec A} and \texorpdfstring{$\Spf A$}{Spf A}}
            \label{appendix:pDiv_prelim:Spec_v_Spf}
                The following facts are implicitly used (including in our companion paper \cite{corank1_ASW_I.pdf}).

\begin{lemma}\label{lemma:p_div_completion_equiv}
Let $A$ be an adic Noetherian ring. There are equivalences of categories
    \begin{align*}
    \{ \text{finite schemes over $\Spec A$} \} & \ra \{ \text{finite relative schemes over $\Spf A$} \}
    \\
    \{ \text{finite locally free schemes over $\Spec A$} \} & \ra \{ \text{finite locally free relative schemes over $\Spf A$} \}
    \\
    \{ \text{$p$-divisible groups over $\Spec A$} \} & \ra \{ \text{$p$-divisible groups over $\Spf A$} \}
    \end{align*}
given by base change, i.e. restriction of fppf sheaves along the inclusion $(\Sch/\Spf A)_{fppf} \ra (\Sch / \Spec A)_{fppf}$.
\end{lemma}
\begin{proof}
For the statements about finite relative schemes, the quasi-inverse functor is given by $\Spf R \mapsto \Spec R$ for finite $A$-algebras $R$ (topologized so that $R$ is an adic ring and the map $A \ra R$ is adic). This also gives the quasi-inverse functor for finite locally free relative schemes (check using the local criterion for flatness). For the statement about $p$-divisible groups (which follows from the other statements), see \cite[{4.15, Lemma II.4.16}]{Messing72} or \cite[{Lemma 2.4.4}]{deJong95}.
\end{proof}

\begin{lemma}\label{lemma:p_div_completion_isogeny}
Let $A$ be an adic Noetherian ring, and let $\phi \colon X \ra Y$ be a homomorphism of $p$-divisible groups over $\Spec A$. Then $\phi$ is an isogeny if and only if $\phi_{\Spf A} \colon X_{\Spf A} \ra Y_{\Spf A}$ is an isogeny.
\end{lemma}
\begin{proof}
Follows from Lemma \ref{lemma:p_div_completion_equiv} and the isogeny criterion from Lemma \ref{lemma:isogeny_criterion}.
\end{proof}

For adic Noetherian rings $A$, we may thus pass between $p$-divisible groups over $\Spec A$ and $\Spf A$ without loss of information, and similarly for finite locally free relative schemes. We abuse notation in this way: for example, if $A$ is a domain, the \emph{generic fiber} of a $p$-divisible group over $\Spf A$ will refer to its generic fiber as a $p$-divisible group over $\Spec A$.

To avoid potential confusion, we remark on three situations where $p$-divisible groups may have different properties when considered over $\Spec A$ versus over $\Spf A$.

\begin{remark}
Let $A$ be an adic Noetherian ring, and suppose $p$ is topologically nilpotent in $A$. Let $X$ be a $p$-divisible group over $\Spec A$. By work of Messing, \cite[{\S II}]{Messing72}, the sheaf $\Lie (X_{\Spf A})$ (in the sense of \cite{SGA3II}) is locally free of finite rank on $(\Sch/\Spf A)_{fppf}$. However, $\Lie X$ (viewed as a sheaf on $(\Sch/\Spec A)_{fppf}$) is \emph{not} necessarily locally free.

For example, consider $A = \Z_p$ and $X = \pmb{\mu}_{p^{\infty}} \coloneqq \varinjlim \pmb{\mu}_{p^n}$, where $\pmb{\mu}_{p^n}$ is the group scheme of $p^n$-th roots of unity. Then the $p$-divisible group $X$ over $\Spec \Z_p$ is \'etale in the generic fiber, but connected of dimension $1$ in the special fiber. We find that $\Lie X|_{\Spec \Q_p} = 0$ but $\Lie X|_{\Spec \F_p}$ is free of rank $1$, so $\Lie X$ cannot be a locally free sheaf of modules on $(\Sch/\Spec A)_{fppf}$.

Thus, when writing $\Lie X$ in this situation, we always mean (by abuse of notation) to view $X$ as a $p$-divisible group over $\Spf A$, so that $\Lie X$ will be a finite locally free sheaf on $(\Sch/\Spf A)_{fppf}$. Similarly, if we say $X$ has dimension $r$, we mean that the finite locally free sheaf $\Lie X$ on $(\Sch/\Spf A)_{fppf}$ has rank $r$.
\end{remark}

\begin{remark}\label{remark:formal_non_torsion_sections}
Let $A$ be an adic Noetherian ring, and let $X$ be a $p$-divisible group over $\Spec A$. In general, there are sections of $X_{\Spf A} \ra \Spf A$ which do not arise as sections of $X \ra \Spec A$. Indeed, sections of $X \ra \Spec A$ correspond precisely to torsion sections of $X_{\Spf A} \ra \Spf A$ (use quasi-compactness of $\Spec A$). 
But $X_{\Spf A} \ra \Spf A$ may have many non-torsion sections, e.g. when $A = \Z_p$ and $X_{\Spf A}$ is a formal $p$-divisible group, hence $X_{\Spf A} \cong \Spf \Z_p[[X_1,\ldots,X_r]]$ as pointed fppf sheaves on $(\Sch/\Spf \Z_p)_{fppf}$. There will be uncountably many non-torsion sections in this situation. This makes a difference in \crefext{I:ssec:more_moduli_pDiv_split:lifting_theory}, for example, where some statements are correct over $\Spf R$ (which is the written version) but incorrect over $\Spec R$.
\end{remark}

\begin{remark}\label{remark:no_lift_quasi_hom}
Let $A$ be an adic Noetherian ring. By our conventions, it is \emph{not} true that any quasi-homomorphism of $p$-divisible groups over $\Spf A$ necessarily lifts to a quasi-homomorphism of $p$-divisible groups over $\Spec A$. See \crefext{I:example:no_lift_quasi_hom}. On the other hand, homomorphisms and isogenies will lift (uniquely) by the preceding lemmas.
\end{remark}

    \clearpage


    \phantomsection
    \addcontentsline{toc}{part}{References}
    \renewcommand{\addcontentsline}[3]{}
    \printbibliography
    
\end{document}